\documentclass[12pt,a4paper,fleqn]{article}
\usepackage{tikz}
\usepackage{caption}
\usepackage{lscape}

\usepackage{fullpage}
\usepackage{amsmath}
\usepackage{amsthm}
\usepackage{amssymb}
\usepackage{amscd}
\usepackage{plain}
\usepackage{graphicx}
\usepackage{t1enc,epsfig}
\usepackage[mathscr]{eucal}
\usepackage{epstopdf}
\usepackage[german,english]{babel}
\usepackage{stmaryrd}

\usepackage{xcolor}

\usepackage{wrapfig}
\usepackage{float}

\renewcommand{\thefootnote}{}
\newtheorem{theorem}{Theorem}[section]

\newtheorem{corollary}[theorem]{Corollary}
\newtheorem{definition}[theorem]{Definition}
\newtheorem{example}[theorem]{Example}
\newtheorem{lemma}[theorem]{Lemma}
\newtheorem{proposition}[theorem]{Proposition}
\newtheorem{remark}[theorem]{Remark}
\newtheorem{conjecture}[theorem]{Conjecture}
\numberwithin{equation}{section}

\newtheorem*{thmI}{Theorem I}
\newtheorem*{thmII}{Theorem II}
\newtheorem*{thmIII}{Theorem III}
\newenvironment{proofI}{{\bf Proof of Theorem I. }}{\hfill$\rule{1ex}{1ex}$\par\medskip}
\newenvironment{proofII}{{\bf Proof of Theorem II. }}{\hfill$\rule{1ex}{1ex}$\par\medskip}
\newenvironment{proofIII}{{\bf Proof of Theorem III. }}{\hfill$\rule{1ex}{1ex}$\par\medskip}

\newtheorem*{lema}{Lemma 4.5.1}
\newtheorem*{lemb}{Lemma 4.5.2}
\newtheorem*{lemc}{Lemma 4.5.3}

\newtheorem*{lemd}{Lemma 4.6}
\newtheorem*{leme}{Lemma 4.7}
\newtheorem*{lemf}{Lemma 4.8}
\newtheorem*{lemg}{Lemma 4.9}
\newtheorem*{lemh}{Lemma 4.10}

\newtheorem*{thma}{Theorem 1}
\newtheorem*{thmb}{Theorem 2}
\newtheorem*{thmc}{Theorem 3}
\newtheorem*{thmd}{Theorem 4}
\newtheorem*{thme}{Theorem 5}
\newtheorem*{thmf}{Theorem 6}
\newtheorem*{thmg}{Theorem 7}
\newtheorem*{thmh}{Theorem 8}
\newtheorem*{thmi}{Theorem 9}
\newtheorem*{thmj}{Theorem 10}
\newtheorem*{thmk}{Theorem 11}
\newtheorem*{thml}{Theorem 12}
\newtheorem*{thmm}{Theorem 13}

\newtheorem*{ra}{Remark 3.1}
\newtheorem*{rb}{Remark 4.1}
\newtheorem*{rc}{Remark 4.2}
\newtheorem*{rdd}{Remark 4.3}
\newtheorem*{ree}{Remark 5.1}
\newtheorem*{rff}{Remark 7.1}
\newtheorem*{rgg}{Remark 7.2}

\renewenvironment{proof}{{\bf Proof. }}{\hfill$\rule{1ex}{1ex}$\par\vskip 5 truemm}
\begin{document}

\newcommand{\bthm}{\begin{theorem}}
\newcommand{\ethm}{\end{theorem}}
\newcommand{\bd}{\begin{definition}}
\newcommand{\ed}{\end{definition}}
\newcommand{\bs}{\begin{proposition}}
\newcommand{\es}{\end{proposition}}
\newcommand{\bp}{\begin{proof}}
\newcommand{\ep}{\end{proof}}
\newcommand{\be}{\begin{equation}}
\newcommand{\ee}{\end{equation}}
\newcommand{\ul}{\underline}
\newcommand{\br}{\begin{remark}}
\newcommand{\er}{\end{remark}}
\newcommand{\bex}{\begin{example}}
\newcommand{\eex}{\end{example}}
\newcommand{\bc}{\begin{corollary}}
\newcommand{\ec}{\end{corollary}}
\newcommand{\bl}{\begin{lemma}}
\newcommand{\el}{\end{lemma}}
\newcommand{\bj}{\begin{conjecture}}
\newcommand{\ej}{\end{conjecture}}

%%%%%%%%%%%%%%%%%%%%%%%%%%%%%%%%%%%%%%%%%%%%%%%%%%%%%%%%%
\newcommand{\bthmI}{\begin{thmI}}
\newcommand{\ethmI}{\end{thmI}}
\newcommand{\bthmII}{\begin{thmII}}
\newcommand{\ethmII}{\end{thmII}}
\newcommand{\bthmIII}{\begin{thmIII}}
\newcommand{\ethmIII}{\end{thmIII}}

\newcommand{\bthma}{\begin{thma}}
\newcommand{\ethma}{\end{thma}}
\newcommand{\bthmb}{\begin{thmb}}
\newcommand{\ethmb}{\end{thmb}}
\newcommand{\bthmc}{\begin{thmc}}
\newcommand{\ethmc}{\end{thmc}}
\newcommand{\bthmd}{\begin{thmd}}
\newcommand{\ethmd}{\end{thmd}}
\newcommand{\bthme}{\begin{thme}}
\newcommand{\ethme}{\end{thme}}
\newcommand{\bthmf}{\begin{thmf}}
\newcommand{\bthms}{\begin{thms}}
\newcommand{\ethmf}{\end{thmf}}
\newcommand{\ethms}{\end{thms}}
\newcommand{\bthmg}{\begin{thmg}}
\newcommand{\ethmg}{\end{thmg}}
\newcommand{\bthmh}{\begin{thmh}}
\newcommand{\ethmh}{\end{thmh}}
\newcommand{\bthmi}{\begin{thmi}}
\newcommand{\ethmi}{\end{thmi}}
\newcommand{\bthmj}{\begin{thmj}}
\newcommand{\ethmj}{\end{thmj}}
\newcommand{\bthmk}{\begin{thmk}}
\newcommand{\ethmk}{\end{thmk}}
\newcommand{\bthml}{\begin{thml}}
\newcommand{\ethml}{\end{thml}}
\newcommand{\bthmm}{\begin{thmm}}
\newcommand{\ethmm}{\end{thmm}}

\newcommand{\blema}{\begin{lema}}
\newcommand{\elema}{\end{lema}}
\newcommand{\blemb}{\begin{lemb}}
\newcommand{\elemb}{\end{lemb}}
\newcommand{\blemc}{\begin{lemc}}
\newcommand{\elemc}{\end{lemc}}
\newcommand{\blemd}{\begin{lemd}}
\newcommand{\elemd}{\end{lemd}}
\newcommand{\bleme}{\begin{leme}}
\newcommand{\eleme}{\end{leme}}
\newcommand{\blemf}{\begin{lemf}}
\newcommand{\elemf}{\end{lemf}}
\newcommand{\blemg}{\begin{lemg}}
\newcommand{\elemg}{\end{lemg}}
\newcommand{\blemh}{\begin{lemh}}
\newcommand{\elemh}{\end{lemh}}

\newcommand{\bra}{\begin{ra}}
\newcommand{\era}{\end{ra}}
\newcommand{\brb}{\begin{rb}}
\newcommand{\erb}{\end{rb}}
\newcommand{\brc}{\begin{rc}}
\newcommand{\erc}{\end{rc}}
\newcommand{\brdd}{\begin{rdd}}
\newcommand{\erdd}{\end{rdd}}
\newcommand{\bree}{\begin{ree}}
\newcommand{\eree}{\end{ree}}
\newcommand{\brff}{\begin{rff}}
\newcommand{\erff}{\end{rff}}
\newcommand{\brgg}{\begin{rgg}}
\newcommand{\ergg}{\end{rgg}}

%%%%%%%%%%%%%%%%%%%%%%%%%%%%%%%%%%%%%%%%%%%%%%%%%%%%%%%

\newcommand{\bprfI}{\begin{proofI}}
\newcommand{\eprfI}{\end{proofI}}

\newcommand{\bprfII}{\begin{proofII}}
\newcommand{\eprfII}{\end{proofII}}

\newcommand{\bprfIII}{\begin{proofIII}}
\newcommand{\eprfIII}{\end{proofIII}}

%%%%%%%%%%%%%%%%%%%%%%%%%%%%%%%%%%%%%%%%%%%%%%%%%%%%%%%%%%%%

%\newcommand{\s2p}{\begin{proof2}}
%\newcommand{\f2p}{\end{proof2}}

\def\diy{\displaystyle}

\def\Gam{{\Gamma}}
\def\g {\gamma}
\def\d {\delta}
\def\s {\phi}
\def\ph {\phi}
\def\f{\vphi}
\def\ta {\tau}
\def\ve {\varepsilon}
\def\rr {\varrho}
\def\Om {\Omega}
\def\pl {\partial}

\def\bD{{\mathbf D}}

\def\tO {B}
\def\tg {B}
\def\bs {\ov\phi}

\def\bbA{{\mathbb A}} \def\bbD{{\mathbb D}} 
\def\bbE{{\mathbb E}} \def\bbF{{\mathbb F}} \def\bbH{{\mathbb H}}
\def\bbL{{\mathbb L}} \def\bbM{{\mathbb M}} \def\bbN{{\mathbb N}}
\def\bbO{{\mathbb O}} \def\bbP{{\mathbb P}} \def\bbQ{{\mathbb Q}}
\def\bbR{{\mathbb R}} \def\bbS{{\mathbb S}} \def\bbT{{\mathbb T}}
\def\bbV{{\mathbb V}} \def\bbW{{\mathbb W}} \def\bbZ{{\mathbb Z}}

\def\cA{{\mathcal A}} \def\cB{{\mathcal B}} \def\cC{{\mathcal C}}
\def\cD{{\mathcal D}} \def\cE{{\mathcal E}} \def\cF{{\mathcal F}}
\def\cG{{\mathcal G}} \def\cH{{\mathcal H}} \def\cJ{{\mathcal J}}
\def\cP{{\mathcal P}} \def\cT{{\mathcal T}} \def\cU{{\mathcal U}}
\def\cV{{\mathcal V}} \def\cW{{\mathcal W}}

\def\cA{\mathscr A} \def\sB{\mathscr TB} \def\sC{\mathscr C}
\def\sD{\mathscr D} \def\sE{\mathscr E}
\def\sF{\mathscr F} \def\sG{\mathscr G} \def\sI{\mathscr I}
\def\sL{\mathscr L} \def\sM{\mathscr M}
\def\sN{\mathscr N} \def\sO{\mathscr O}
\def\sP{\mathscr P} \def\sR{\mathscr R}
\def\sS{\mathscr S} \def\sS{\mathscr T}
\def\sU{\mathscr U} \def\sV{\mathscr V} \def\sW{\mathscr W}
\def\sX{\mathscr X} \def\sY{\mathscr Y} \def\sZ{\mathscr Z}

\def\tA{{\tt A}} \def\tB{{\tt B}} \def\tC{{\tt C}}
\def\tD{{\tt D}} \def\tE{{\tt E}} 
\def\tQ{{\tt Q}} \def\tR{{\tt R}} \def\tS{{\tt S}}
\def\tT{{\tt T}}

\def\fA{{\mathfrak A}} \def\fB{{\mathfrak B}} \def\fC{{\mathfrak C}}
\def\fD{{\mathfrak D}} \def\fE{{\mathfrak E}} \def\fF{{\mathfrak F}}
\def\fW{{\mathfrak W}} \def\fX{{\mathfrak X}} \def\fY{{\mathfrak Y}}
\def\fZ{{\mathfrak Z}}

\def\bfe{{\mathbf e}} \def\bh{{\mathbf h}}
\def\bi{{\mathbf i}} \def\bj{{\mathbf j}}
\def\bn{{\mathbf n}} \def\bu{{\mathbf u}} \def\bv{{\mathbf v}}
\def\bx{{\mathbf x}} \def\by{{\mathbf y}} \def\bz{{\mathbf z}}
\def\B1{{\mathbf 1}} \def\co{\complement}

\def\bA{{\bf A}} \def\bB{{\bf B}} \def\bC{{\bf C}}
\def\bD{{\bf D}} \def\bE{{\bf E}} \def\bS{{\bf S}}

\def\bAH{{\bf{AH}}} \def\bBH{{\bf{BH}}} \def\bCH{{\bf{CH}}}
\def\bDH{{\bf{DH}}} \def\bEH{{\bf{EH}}} \def\bSH{{\bf{SH}}}

\def\BZ{{\mathbf Z}}

\def\bmu{{\mbox{\boldmath${\mu}$}}}
\def\bnu{{\mbox{\boldmath${\nu}$}}}
\def\bPhi{{\mbox{\boldmath${\Phi}$}}}

\def\fC{{\mathfrak C}}

\def\rA{{\rm A}}  \def\rB{{\rm B}}  \def\rC{{\rm C}}
\def\rD{{\rm  D}} \def\rE{{\rm  E}} \def\rF{{\rm  F}} 
 \def\rG{{\rm  G}} 
\def\rH{{\rm  H}}  \def\rM{{\rm M}} \def\rR{{\rm  R}} 
 \def\rS{{\rm  S}} \def\rT{{\rm T}} \def\rV{{\rm  V}} 
\def\rd{{\rm d}} \def\re{{\rm e}} \def\rn{{\rm n}}
\def\rs{{\rm s}} \def\rt{{\rm t}} \def\ru{{\rm u}}

\def\RD{{\rR\rD}}

\def\ov{\overline} \def\ovp{{\overline p}}

\def\es {{\varnothing}}

\def\cc {\circ}

\def\wt{\widetilde} \def\wh{\wideha}

\def\wtD{{\wt D}}
\def\wtm{{\wt m}} \def\wtn{{\wt n}} \def\wtk{{\wt k}}

\def\be{\begin{equation}}
\def\ee{\end{equation}}

\def\beq{\begin{equation}}
\def\eeq{\end{equation}}

\def\beal{\begin{array}{l}} \def\beac{\begin{array}{c}} \def\bear{\begin{array}{r}}
\def\beacl{\begin{array}{cl}} \def\beall{\begin{array}{ll}}
\def\bealllll{\begin{array}{lllll}}
\def\beacr{\begin{array}{cr}}
\def\ena{\end{array}}

\def\bma{\begin{matrix}}
\def\ema{\end{matrix}}

\def\diy{\displaystyle}

\def\cF{\mathcal F} \def\cG{\mathcal G} \def\cI{\mathcal I}
\def\cL{\mathcal L} \def\cO{\mathcal O}
\def\cP{\mathcal P} 

\def\Gam{{\Gamma}} \def\gam{{\gamma}}
\def\Del{{\Delta}} \def\del{{\delta}}  \def\odel{{\overline\delta}}

\def\vphi{{\varphi}}
\def\eps{{\epsilon}} \def\veps{{\varepsilon}}
\def\vrho{{\varrho}} \def\vpi{{\varpi}}
\def\Lam{{\Lambda}} \def\lam{{\lambda}} 
\def\Om{{\Omega}} \def\om{{\omega}} 

\def\pa {\partial}
\def\comp{\complement}

\def\bZ{{\mathbf Z}}

\def\rSp {{\rm Supp}}

\def\D{D}

\def\tO {B}
\def\tg {B}
\def\bs {\overline \phi}

\def\c  {{\mathchar"0\hexnumber@\msafam7B}}
\def\es {{\varnothing}}

\def\cc {\circ}

\def\io{\iota}
\def\rInt{{\rm{Int}}}
\def\rExt{{\rm{Ext}}}

\def\cl{\centerline}

%%%%%%%%%%%%%%%%%%%%%%%%%%%%  201068247927  
%%%%%%%%%%%%%%%%%%%%%%%%%%%%%%%%%%%%%
%% 86602540378443864676372317075294   M251

\def\Sgie{\raise-1.5ex\hbox{${{\displaystyle \sum}^{\sharp} \atop {\scriptstyle \{\g_m\}^{\ex} \in V}}$}}

%options for floating pictures

\makeatletter
 \def\fps@figure{htbp}
\makeatother

%definitions for drawing with tikz
\def\rr{0.86602540378443864676372317075294} %\sqrt{3/4}

\def\hexagongrid{
%\clip[yscale=sqrt(3/4), xslant=0.5] (0, \n) -- (\n, 0) -- (\n, -\n) -- (0, -\n) -- (-\n, 0) -- (-\n, \n) -- cycle;
\draw [yscale=sqrt(3/4), xslant=0.5] (-\n, -\n) grid (\n, \n);
\draw [yscale=sqrt(3/4), xslant=-0.5] (-\n, -\n) grid (\n, \n);
}

\def\rectangulargrid{
\clip[yscale=sqrt(3/4)] (-\n, -\n) rectangle (\n, \n);
\draw [yscale=sqrt(3/4), xslant=0.5] (-2 * \n, -\n) grid (2 * \n, \n);
\draw [yscale=sqrt(3/4), xslant=-0.5] (-2 * \n, -\n) grid (2 * \n, \n);
}

\def\sublattice{
\clip[yscale=sqrt(3/4), xslant=0.5] (0, \n+\nn) -- (\n+\nn, 0) -- (\n+\nn, -\n-\nn) -- (0, -\n-\nn) -- (-\n-\nn, 0) -- (-\n-\nn, \n+\nn) -- cycle;
\foreach \x in {-\n,...,\n}
 \foreach \y in {-\n,...,\n}
 {
  \def\xx{\x * \aa + 0.5 * \x * \bb - 0.5 * \y * \bb + 0.5 * \y * \aa}
  \def\yy{\rr * \x * \bb + \rr * \y * \aa + \rr * \y * \bb}

  \shade[shading=ball, ball color=black] (\xx, \yy) circle (.2);
 }
}

\def\sublatticeinrectangle{
\clip[yscale=sqrt(3/4)] (-\nn, -\nn) rectangle (\nn, \nn);
\foreach \x in {-\n,...,\n}
 \foreach \y in {-\n,...,\n}
 {
  \def\xx{\x * \aa + 0.5 * \x * \bb - 0.5 * \y * \bb + 0.5 * \y * \aa}
  \def\yy{\rr * \x * \bb + \rr * \y * \aa + \rr * \y * \bb}

  \shade[shading=ball, ball color=black] (\xx + \dx + 0.5 * \dy, \yy + \rr * \dy) circle (.3);
 }
}

\def\inclinedgrid{
	\draw [yscale=\bb * sqrt(3/4), xslant=\aa + 0.5 * \bb, xscale=\aa * \aa / \bb + \aa  + \bb , yslant=-\aa / \bb, ultra thick] (-\n, -\n) grid (\n, \n);
}

\def\rectsub-lattice{
\clip[yscale=sqrt(3/4)] (-\nn, -\nn) rectangle (\nn, \nn);
\foreach \x in {-\n,...,\n}
 \foreach \y in {-\n,...,\n}
 {
  \filldraw[yscale=sqrt(3/4), xslant=0.5, \cc] (\ss * \x + \dx, \ss * \y + \dy) rectangle (\ss * \x + \ss - 1 + \dx, \ss * \y + \ss - 1 + \dy);

\definecolor{gray1}{gray}{0.1}
\definecolor{gray2}{gray}{0.2}
\definecolor{gray3}{gray}{0.3}
\definecolor{gray4}{gray}{0.4}
\definecolor{gray5}{gray}{0.5}
\definecolor{gray6}{gray}{0.6}
\definecolor{gray7}{gray}{0.7}
\definecolor{gray8}{gray}{0.8}
\definecolor{gray9}{gray}{0.9}
 }
}

%%%%%%%%%%%%%%%%%%%%%%% El2gFC1i?

\def\FigureA1 %\label{Fig1}
{\begin{figure} \label{Fig1}\centering
(a) \begin{tikzpicture}[scale=0.3]
\clip (2, -0.5) rectangle (17.4, 16.1);
\draw[yscale=sqrt(3/4), xslant=0.5] (-4,-2) grid (29, 12);
\draw[yscale=sqrt(3/4), xslant=-0.5] (-1,-2) grid (29, 12);
% \foreach \pos in {(0, 0), (7, 0), (3.5, 7*\rr),(10.5, 7*\rr)}
% \shade[shading=ball, ball color=black] \pos circle (.4);
\end{tikzpicture}\quad (b)
\begin{tikzpicture}[scale=0.27]
\clip (-0.9, -1.2*\rr) rectangle (17, 13*\rr);
\draw [line width=0.2mm] (0.5,1*\rr) -- (1.5,1*\rr) -- (2,0) -- (1.5,-1*\rr)
-- (0.5,-1*\rr) -- (0,0) -- (0.5,1*\rr);
\draw [line width=0.2mm] (0.5,1*\rr) -- (0,2*\rr) -- (0.5,3*\rr) -- (1.5,3*\rr) 
-- (2,2*\rr) -- (1.5,1*\rr);
\draw [line width=0.2mm] (0.5,3*\rr) -- (0,4*\rr) -- (0.5,5*\rr) -- (1.5,5*\rr) 
-- (2,4*\rr) -- (1.5,3*\rr);
\draw [line width=0.2mm] (0.5,5*\rr) -- (0,6*\rr) -- (0.5,7*\rr) -- (1.5,7*\rr) 
-- (2,6*\rr) -- (1.5,5*\rr);
\draw [line width=0.2mm] (0.5,7*\rr) -- (0,8*\rr) -- (0.5,9*\rr) -- (1.5,9*\rr) 
-- (2,8*\rr) -- (1.5,7*\rr);
\draw [line width=0.2mm] (0.5,9*\rr) -- (0,10*\rr) -- (0.5,11*\rr) -- (1.5,11*\rr) 
-- (2,10*\rr) -- (1.5,9*\rr);
\draw [line width=0.2mm] (0.5,11*\rr) -- (0,12*\rr) -- (0.5,13*\rr) -- (1.5,13*\rr) 
-- (2,12*\rr) -- (1.5,11*\rr);
\draw [line width=0.2mm] (0.5,13*\rr) -- (0,14*\rr) -- (0.5,15*\rr) -- (1.5,15*\rr) 
-- (2,14*\rr) -- (1.5,13*\rr);
\draw [line width=0.2mm] (0.5,15*\rr) -- (0,16*\rr) -- (0.5,17*\rr) -- (1.5,17*\rr) 
-- (2,16*\rr) -- (1.5,15*\rr);
\draw [line width=0.2mm] (0.5,17*\rr) -- (0,18*\rr) -- (0.5,19*\rr) -- (1.5,19*\rr) 
-- (2,18*\rr) -- (1.5,17*\rr);
\draw [line width=0.2mm] (0.5,19*\rr) -- (0,20*\rr) -- (0.5,21*\rr) -- (1.5,21*\rr) 
-- (2,20*\rr) -- (1.5,19*\rr);
%\draw [line width=0.2mm] (0.5,21*\rr) -- (0,22*\rr) -- (0.5,22*\rr) -- (1.5,22*\rr) 
%-- (2,22*\rr) -- (1.5,21*\rr); 1 1010 1212 1313 19191919 1818181818 
% 6666 55555

\draw [line width=0.2mm]  (3,0) -- (3.5,-1*\rr);
\draw [line width=0.2mm] (2,2*\rr) -- (3,2*\rr) -- (3.5,1*\rr) -- (3,0) -- (2,0);
\draw [line width=0.2mm] (2,4*\rr) -- (3,4*\rr) -- (3.5,3*\rr) -- (3,2*\rr);
\draw [line width=0.2mm] (2,6*\rr) -- (3,6*\rr) -- (3.5,5*\rr) -- (3,4*\rr);
\draw [line width=0.2mm] (2,8*\rr) -- (3,8*\rr) -- (3.5,7*\rr) -- (3,6*\rr);
\draw [line width=0.2mm] (2,10*\rr) -- (3,10*\rr) -- (3.5,9*\rr) -- (3,8*\rr);
\draw [line width=0.2mm] (2,12*\rr) -- (3,12*\rr) -- (3.5,11*\rr) -- (3,10*\rr);
\draw [line width=0.2mm] (2,14*\rr) -- (3,14*\rr) -- (3.5,13*\rr) -- (3,12*\rr);
\draw [line width=0.2mm] (2,16*\rr) -- (3,16*\rr) -- (3.5,15*\rr) -- (3,14*\rr);
\draw [line width=0.2mm] (2,18*\rr) -- (3,18*\rr) -- (3.5,17*\rr) -- (3,16*\rr);
\draw [line width=0.2mm] (2,20*\rr) -- (3,20*\rr) -- (3.5,19*\rr) -- (3,18*\rr);
\draw [line width=0.2mm] (3.5,21*\rr) -- (3,20*\rr);

\draw [line width=0.2mm]  (3.5,1*\rr) -- (4.5,1*\rr) -- (5,0) -- (4.5,-1*\rr)
--(3.5,-1*\rr);
\draw [line width=0.2mm]  (3.5,3*\rr) -- (4.5,3*\rr) -- (5,2*\rr) -- (4.5,1*\rr);
\draw [line width=0.2mm]  (3.5,5*\rr) -- (4.5,5*\rr) -- (5,4*\rr) -- (4.5,3*\rr);
\draw [line width=0.2mm]  (3.5,7*\rr) -- (4.5,7*\rr) -- (5,6*\rr) -- (4.5,5*\rr);
\draw [line width=0.2mm]  (3.5,9*\rr) -- (4.5,9*\rr) -- (5,8*\rr) -- (4.5,7*\rr);
\draw [line width=0.2mm]  (3.5,11*\rr) -- (4.5,11*\rr) -- (5,10*\rr) -- (4.5,9*\rr);
\draw [line width=0.2mm]  (3.5,13*\rr) -- (4.5,13*\rr) -- (5,12*\rr) -- (4.5,11*\rr);
\draw [line width=0.2mm]  (3.5,15*\rr) -- (4.5,15*\rr) -- (5,14*\rr) -- (4.5,13*\rr);
\draw [line width=0.2mm]  (3.5,17*\rr) -- (4.5,17*\rr) -- (5,16*\rr) -- (4.5,15*\rr);
\draw [line width=0.2mm]  (3.5,19*\rr) -- (4.5,19*\rr) -- (5,18*\rr) -- (4.5,17*\rr);
\draw [line width=0.2mm]  (3.5,21*\rr) -- (4.5,21*\rr) -- (5,20*\rr) -- (4.5,19*\rr);

\draw [line width=0.2mm]  (5,0) -- (6,0)--(6.5,-1*\rr );
\draw [line width=0.2mm]  (5,2*\rr) -- (6,2*\rr)--(6.5,1*\rr )--(6,0);
\draw [line width=0.2mm]  (5,4*\rr) -- (6,4*\rr)--(6.5,3*\rr )--(6,2*\rr);
\draw [line width=0.2mm]  (5,6*\rr) -- (6,6*\rr)--(6.5,5*\rr )--(6,4*\rr);
\draw [line width=0.2mm]  (5,8*\rr) -- (6,8*\rr)--(6.5,7*\rr )--(6,6*\rr);
\draw [line width=0.2mm]  (5,10*\rr) -- (6,10*\rr)--(6.5,9*\rr )--(6,8*\rr);
\draw [line width=0.2mm]  (5,12*\rr) -- (6,12*\rr)--(6.5,11*\rr )--(6,10*\rr);
\draw [line width=0.2mm]  (5,14*\rr) -- (6,14*\rr)--(6.5,13*\rr )--(6,12*\rr);
\draw [line width=0.2mm]  (5,16*\rr) -- (6,16*\rr)--(6.5,15*\rr )--(6,14*\rr);
\draw [line width=0.2mm]  (5,18*\rr) -- (6,18*\rr)--(6.5,17*\rr )--(6,16*\rr);
\draw [line width=0.2mm]  (5,20*\rr) -- (6,20*\rr)--(6.5,19*\rr )--(6,18*\rr);
\draw [line width=0.2mm]  (6.5,21*\rr )--(6,20*\rr);

\draw [line width=0.2mm]  (6.5,1*\rr) -- (7.5,1*\rr)--(8,0)--(7.5,-1*\rr)
--(6.5,-1*\rr);
\draw [line width=0.2mm]  (6.5,3*\rr) -- (7.5,3*\rr)--(8,2*\rr)--(7.5,1*\rr);
\draw [line width=0.2mm]  (6.5,5*\rr) -- (7.5,5*\rr)--(8,4*\rr)--(7.5,3*\rr);
\draw [line width=0.2mm]  (6.5,7*\rr) -- (7.5,7*\rr)--(8,6*\rr)--(7.5,5*\rr);
\draw [line width=0.2mm]  (6.5,9*\rr) -- (7.5,9*\rr)--(8,8*\rr)--(7.5,7*\rr);
\draw [line width=0.2mm]  (6.5,11*\rr) -- (7.5,11*\rr)--(8,10*\rr)--(7.5,9*\rr);
\draw [line width=0.2mm]  (6.5,13*\rr) -- (7.5,13*\rr)--(8,12*\rr)--(7.5,11*\rr);
\draw [line width=0.2mm]  (6.5,15*\rr) -- (7.5,15*\rr)--(8,14*\rr)--(7.5,13*\rr);
\draw [line width=0.2mm]  (6.5,17*\rr) -- (7.5,17*\rr)--(8,16*\rr)--(7.5,15*\rr);
\draw [line width=0.2mm]  (6.5,19*\rr) -- (7.5,19*\rr)--(8,18*\rr)--(7.5,17*\rr);
\draw [line width=0.2mm]  (6.5,21*\rr) -- (7.5,21*\rr)--(8,20*\rr)--(7.5,19*\rr);

\draw [line width=0.2mm]  (8,0)--(9,0)--(9.5,-1*\rr);
\draw [line width=0.2mm]  (8,2*\rr)--(9,2*\rr)--(9.5,1*\rr)--(9,0);
\draw [line width=0.2mm]  (8,4*\rr)--(9,4*\rr)--(9.5,3*\rr)--(9,2*\rr);
\draw [line width=0.2mm]  (8,6*\rr)--(9,6*\rr)--(9.5,5*\rr)--(9,4*\rr);
\draw [line width=0.2mm]  (8,8*\rr)--(9,8*\rr)--(9.5,7*\rr)--(9,6*\rr);
\draw [line width=0.2mm]  (8,10*\rr)--(9,10*\rr)--(9.5,9*\rr)--(9,8*\rr);
\draw [line width=0.2mm]  (8,12*\rr)--(9,12*\rr)--(9.5,11*\rr)--(9,10*\rr);
\draw [line width=0.2mm]  (8,14*\rr)--(9,14*\rr)--(9.5,13*\rr)--(9,12*\rr);
\draw [line width=0.2mm]  (8,16*\rr)--(9,16*\rr)--(9.5,15*\rr)--(9,14*\rr);
\draw [line width=0.2mm]  (8,18*\rr)--(9,18*\rr)--(9.5,17*\rr)--(9,16*\rr);
\draw [line width=0.2mm]  (8,20*\rr)--(9,20*\rr)--(9.5,19*\rr)--(9,18*\rr);
\draw [line width=0.2mm]  (9.5,21*\rr)--(9,20*\rr);

\draw [line width=0.2mm]  (9.5,1*\rr)--(10.5,1*\rr)--(11,0)--(10.5,-1*\rr)
--(9.5,-1*\rr);
\draw [line width=0.2mm]  (9.5,3*\rr)--(10.5,3*\rr)--(11,2*\rr)--(10.5,1*\rr);
\draw [line width=0.2mm]  (9.5,5*\rr)--(10.5,5*\rr)--(11,4*\rr)--(10.5,3*\rr);
\draw [line width=0.2mm]  (9.5,7*\rr)--(10.5,7*\rr)--(11,6*\rr)--(10.5,5*\rr);
\draw [line width=0.2mm]  (9.5,9*\rr)--(10.5,9*\rr)--(11,8*\rr)--(10.5,7*\rr);
\draw [line width=0.2mm]  (9.5,11*\rr)--(10.5,11*\rr)--(11,10*\rr)--(10.5,9*\rr);
\draw [line width=0.2mm]  (9.5,13*\rr)--(10.5,13*\rr)--(11,12*\rr)--(10.5,11*\rr);
\draw [line width=0.2mm]  (9.5,15*\rr)--(10.5,15*\rr)--(11,14*\rr)--(10.5,13*\rr);
\draw [line width=0.2mm]  (9.5,17*\rr)--(10.5,17*\rr)--(11,16*\rr)--(10.5,15*\rr);
\draw [line width=0.2mm]  (9.5,19*\rr)--(10.5,19*\rr)--(11,18*\rr)--(10.5,17*\rr);
\draw [line width=0.2mm]  (9.5,21*\rr)--(10.5,21*\rr)--(11,20*\rr)--(10.5,19*\rr);

\draw [line width=0.2mm]  (11,0)--(12,0)--(12.5,-1*\rr);
\draw [line width=0.2mm]  (11,2*\rr)--(12,2*\rr)--(12.5,1*\rr)--(12,0);
\draw [line width=0.2mm]  (11,4*\rr)--(12,4*\rr)--(12.5,3*\rr)--(12,2*\rr);
\draw [line width=0.2mm]  (11,6*\rr)--(12,6*\rr)--(12.5,5*\rr)--(12,4*\rr);
\draw [line width=0.2mm]  (11,8*\rr)--(12,8*\rr)--(12.5,7*\rr)--(12,6*\rr);
\draw [line width=0.2mm]  (11,10*\rr)--(12,10*\rr)--(12.5,9*\rr)--(12,8*\rr);
\draw [line width=0.2mm]  (11,12*\rr)--(12,12*\rr)--(12.5,11*\rr)--(12,10*\rr);
\draw [line width=0.2mm]  (11,14*\rr)--(12,14*\rr)--(12.5,13*\rr)--(12,12*\rr);
\draw [line width=0.2mm]  (11,16*\rr)--(12,16*\rr)--(12.5,15*\rr)--(12,14*\rr);
\draw [line width=0.2mm]  (11,18*\rr)--(12,18*\rr)--(12.5,17*\rr)--(12,16*\rr);
\draw [line width=0.2mm]  (11,20*\rr)--(12,20*\rr)--(12.5,19*\rr)--(12,18*\rr);
\draw [line width=0.2mm]  (12,20*\rr)--(12.5,21*\rr);

\draw [line width=0.2mm]  (12.5,1*\rr)--(13.5,1*\rr)--(14,0)--(13.5,-1*\rr)
--(12.5,-1*\rr);
\draw [line width=0.2mm]  (12.5,3*\rr)--(13.5,3*\rr)--(14,2*\rr)--(13.5,1*\rr);
\draw [line width=0.2mm]  (12.5,5*\rr)--(13.5,5*\rr)--(14,4*\rr)--(13.5,3*\rr);
\draw [line width=0.2mm]  (12.5,7*\rr)--(13.5,7*\rr)--(14,6*\rr)--(13.5,5*\rr);
\draw [line width=0.2mm]  (12.5,9*\rr)--(13.5,9*\rr)--(14,8*\rr)--(13.5,7*\rr);
\draw [line width=0.2mm]  (12.5,11*\rr)--(13.5,11*\rr)--(14,10*\rr)--(13.5,9*\rr);
\draw [line width=0.2mm]  (12.5,13*\rr)--(13.5,13*\rr)--(14,12*\rr)--(13.5,11*\rr);
\draw [line width=0.2mm]  (12.5,15*\rr)--(13.5,15*\rr)--(14,14*\rr)--(13.5,13*\rr);
\draw [line width=0.2mm]  (12.5,17*\rr)--(13.5,17*\rr)--(14,16*\rr)--(13.5,15*\rr);
\draw [line width=0.2mm]  (12.5,19*\rr)--(13.5,19*\rr)--(14,18*\rr)--(13.5,17*\rr);
\draw [line width=0.2mm]  (12.5,21*\rr)--(13.5,21*\rr)--(14,20*\rr)--(13.5,19*\rr);

\draw [line width=0.2mm]  (14,0)--(15,0)--(15.5,-1*\rr);
\draw [line width=0.2mm]  (14,2*\rr)--(15,2*\rr)--(15.5,1*\rr)--(15,0);
\draw [line width=0.2mm]  (14,4*\rr)--(15,4*\rr)--(15.5,3*\rr)--(15,2*\rr);
\draw [line width=0.2mm]  (14,6*\rr)--(15,6*\rr)--(15.5,5*\rr)--(15,4*\rr);
\draw [line width=0.2mm]  (14,8*\rr)--(15,8*\rr)--(15.5,7*\rr)--(15,6*\rr);
\draw [line width=0.2mm]  (14,10*\rr)--(15,10*\rr)--(15.5,9*\rr)--(15,8*\rr);
\draw [line width=0.2mm]  (14,12*\rr)--(15,12*\rr)--(15.5,11*\rr)--(15,10*\rr);
\draw [line width=0.2mm]  (14,14*\rr)--(15,14*\rr)--(15.5,13*\rr)--(15,12*\rr);
\draw [line width=0.2mm]  (14,16*\rr)--(15,16*\rr)--(15.5,15*\rr)--(15,14*\rr);
\draw [line width=0.2mm]  (14,18*\rr)--(15,18*\rr)--(15.5,17*\rr)--(15,16*\rr);
\draw [line width=0.2mm]  (14,20*\rr)--(15,20*\rr)--(15.5,19*\rr)--(15,18*\rr);
\draw [line width=0.2mm]  (15,20*\rr)--(15.5,21*\rr);

\draw [line width=0.2mm]  (15.5,1*\rr)--(16.5,1*\rr)--(17,0)--(16.5,-1*\rr)
--(15.5,-1*\rr);
\draw [line width=0.2mm]  (15.5,3*\rr)--(16.5,3*\rr)--(17,2*\rr)--(16.5,1*\rr);
\draw [line width=0.2mm]  (15.5,5*\rr)--(16.5,5*\rr)--(17,4*\rr)--(16.5,3*\rr);
\draw [line width=0.2mm]  (15.5,7*\rr)--(16.5,7*\rr)--(17,6*\rr)--(16.5,5*\rr);
\draw [line width=0.2mm]  (15.5,9*\rr)--(16.5,9*\rr)--(17,8*\rr)--(16.5,7*\rr);
\draw [line width=0.2mm]  (15.5,11*\rr)--(16.5,11*\rr)--(17,10*\rr)--(16.5,9*\rr);
\draw [line width=0.2mm]  (15.5,13*\rr)--(16.5,13*\rr)--(17,12*\rr)--(16.5,11*\rr);
\draw [line width=0.2mm]  (15.5,15*\rr)--(16.5,15*\rr)--(17,14*\rr)--(16.5,13*\rr);
\draw [line width=0.2mm]  (15.5,17*\rr)--(16.5,17*\rr)--(17,16*\rr)--(16.5,15*\rr);
\draw [line width=0.2mm]  (15.5,19*\rr)--(16.5,19*\rr)--(17,18*\rr)--(16.5,17*\rr);
\draw [line width=0.2mm]  (15.5,21*\rr)--(16.5,21*\rr)--(17,20*\rr)--(16.5,19*\rr);

\draw [line width=0.2mm]  (17,0)--(18,0)--(18.5,-1*\rr);
\draw [line width=0.2mm]  (17,2*\rr)--(18,2*\rr)--(18.5,1*\rr)--(18,0);
\draw [line width=0.2mm]  (17,4*\rr)--(18,4*\rr)--(18.5,3*\rr)--(18,2*\rr);
\draw [line width=0.2mm]  (17,6*\rr)--(18,6*\rr)--(18.5,5*\rr)--(18,4*\rr);
\draw [line width=0.2mm]  (17,8*\rr)--(18,8*\rr)--(18.5,7*\rr)--(18,6*\rr);
\draw [line width=0.2mm]  (17,10*\rr)--(18,10*\rr)--(18.5,9*\rr)--(18,8*\rr);
\draw [line width=0.2mm]  (17,12*\rr)--(18,12*\rr)--(18.5,11*\rr)--(18,10*\rr);
\draw [line width=0.2mm]  (17,14*\rr)--(18,14*\rr)--(18.5,13*\rr)--(18,12*\rr);
\draw [line width=0.2mm]  (17,16*\rr)--(18,16*\rr)--(18.5,15*\rr)--(18,14*\rr);
\draw [line width=0.2mm]  (17,18*\rr)--(18,18*\rr)--(18.5,17*\rr)--(18,16*\rr);
\draw [line width=0.2mm]  (17,20*\rr)--(18,20*\rr)--(18.5,19*\rr)--(18,18*\rr);
\draw [line width=0.2mm]  (18,20*\rr)--(18.5,21*\rr);

\draw [line width=0.2mm]  (18.5,1*\rr)--(19.5,1*\rr)--(20,0)--(19.5,-1*\rr)
--(18.5,-1*\rr);
\draw [line width=0.2mm]  (18.5,3*\rr)--(19.5,3*\rr)--(20,2*\rr)--(19.5,1*\rr);
\draw [line width=0.2mm]  (18.5,5*\rr)--(19.5,5*\rr)--(20,4*\rr)--(19.5,3*\rr);
\draw [line width=0.2mm]  (18.5,7*\rr)--(19.5,7*\rr)--(20,6*\rr)--(19.5,5*\rr);
\draw [line width=0.2mm]  (18.5,9*\rr)--(19.5,9*\rr)--(20,8*\rr)--(19.5,7*\rr);
\draw [line width=0.2mm]  (18.5,11*\rr)--(19.5,11*\rr)--(20,10*\rr)--(19.5,9*\rr);
\draw [line width=0.2mm]  (18.5,13*\rr)--(19.5,13*\rr)--(20,12*\rr)--(19.5,11*\rr);
\draw [line width=0.2mm]  (18.5,15*\rr)--(19.5,15*\rr)--(20,14*\rr)--(19.5,13*\rr);
\draw [line width=0.2mm]  (18.5,17*\rr)--(19.5,17*\rr)--(20,16*\rr)--(19.5,15*\rr);
\draw [line width=0.2mm]  (18.5,19*\rr)--(19.5,19*\rr)--(20,18*\rr)--(19.5,17*\rr);
\draw [line width=0.2mm]  (18.5,21*\rr)--(19.5,21*\rr)--(20,20*\rr)--(19.5,19*\rr);

\draw [line width=0.2mm]  (20,0)--(21,0)--(21.5,-1*\rr);
\draw [line width=0.2mm]  (20,2*\rr)--(21,2*\rr)--(21.5,1*\rr)--(21,0);
\draw [line width=0.2mm]  (20,4*\rr)--(21,4*\rr)--(21.5,3*\rr)--(21,2*\rr);
\draw [line width=0.2mm]  (20,6*\rr)--(21,6*\rr)--(21.5,5*\rr)--(21,4*\rr);
\draw [line width=0.2mm]  (20,8*\rr)--(21,8*\rr)--(21.5,7*\rr)--(21,6*\rr);
\draw [line width=0.2mm]  (20,10*\rr)--(21,10*\rr)--(21.5,9*\rr)--(21,8*\rr);
\draw [line width=0.2mm]  (20,12*\rr)--(21,12*\rr)--(21.5,11*\rr)--(21,10*\rr);
\draw [line width=0.2mm]  (20,14*\rr)--(21,14*\rr)--(21.5,13*\rr)--(21,12*\rr);
\draw [line width=0.2mm]  (20,16*\rr)--(21,16*\rr)--(21.5,15*\rr)--(21,14*\rr);
\draw [line width=0.2mm]  (20,18*\rr)--(21,18*\rr)--(21.5,17*\rr)--(21,16*\rr);
\draw [line width=0.2mm]  (20,20*\rr)--(21,20*\rr)--(21.5,19*\rr)--(21,18*\rr);
\draw [line width=0.2mm]  (21,20*\rr)--(21.5,21*\rr);

\draw [line width=0.2mm]  (21.5,1*\rr)--(22.5,1*\rr)--(23,0)--(22.5,-1*\rr)
--(21.5,-1*\rr);
\draw [line width=0.2mm]  (21.5,3*\rr)--(22.5,3*\rr)--(23,2*\rr)--(22.5,1*\rr);
\draw [line width=0.2mm]  (21.5,5*\rr)--(22.5,5*\rr)--(23,4*\rr)--(22.5,3*\rr);
\draw [line width=0.2mm]  (21.5,7*\rr)--(22.5,7*\rr)--(23,6*\rr)--(22.5,5*\rr);
\draw [line width=0.2mm]  (21.5,9*\rr)--(22.5,9*\rr)--(23,8*\rr)--(22.5,7*\rr);
\draw [line width=0.2mm]  (21.5,11*\rr)--(22.5,11*\rr)--(23,10*\rr)--(22.5,9*\rr);
\draw [line width=0.2mm]  (21.5,13*\rr)--(22.5,13*\rr)--(23,12*\rr)--(22.5,11*\rr);
\draw [line width=0.2mm]  (21.5,15*\rr)--(22.5,15*\rr)--(23,14*\rr)--(22.5,13*\rr);
\draw [line width=0.2mm]  (21.5,17*\rr)--(22.5,17*\rr)--(23,16*\rr)--(22.5,15*\rr);
\draw [line width=0.2mm]  (21.5,19*\rr)--(22.5,19*\rr)--(23,18*\rr)--(22.5,17*\rr);
\draw [line width=0.2mm]  (21.5,21*\rr)--(22.5,21*\rr)--(23,20*\rr)--(22.5,19*\rr);

\draw [line width=0.2mm]  (23,0)--(24,0)--(24.5,-1*\rr);
\draw [line width=0.2mm]  (23,2*\rr)--(24,2*\rr)--(24.5,1*\rr)--(24,0);
\draw [line width=0.2mm]  (23,4*\rr)--(24,4*\rr)--(24.5,3*\rr)--(24,2*\rr);
\draw [line width=0.2mm]  (23,6*\rr)--(24,6*\rr)--(24.5,5*\rr)--(24,4*\rr);
\draw [line width=0.2mm]  (23,8*\rr)--(24,8*\rr)--(24.5,7*\rr)--(24,6*\rr);
\draw [line width=0.2mm]  (23,10*\rr)--(24,10*\rr)--(24.5,9*\rr)--(24,8*\rr);
\draw [line width=0.2mm]  (23,12*\rr)--(24,12*\rr)--(24.5,11*\rr)--(24,10*\rr);
\draw [line width=0.2mm]  (23,14*\rr)--(24,14*\rr)--(24.5,13*\rr)--(24,12*\rr);
\draw [line width=0.2mm]  (23,16*\rr)--(24,16*\rr)--(24.5,15*\rr)--(24,14*\rr);
\draw [line width=0.2mm]  (23,18*\rr)--(24,18*\rr)--(24.5,17*\rr)--(24,16*\rr);
\draw [line width=0.2mm]  (23,20*\rr)--(24,20*\rr)--(24.5,19*\rr)--(24,18*\rr);
\draw [line width=0.2mm]  (24,20*\rr)--(24.5,21*\rr); 
\end{tikzpicture} % \vskip .5cm
\caption{Fragments of lattices ${\mathbb A}_2$ and ${\mathbb H}_2$}
\end{figure}} %Fig1
%%%%%%%%%%%%%%%%%%%%%%%%%%%%%%%%%%%%%%%%%%%%%%%%%%

\def\FigureB2 %% \label{Fig2}
{\begin{figure}\label{Fig2} 
\centering
\captionsetup{width=0.8\textwidth} 
(a) \begin{tikzpicture}[scale=0.35]
\clip (2, -0.5) rectangle (17.4, 16.1);
\path [fill=lightgray] (10.5,7*\rr)--(9,8*\rr)--(7.5,7*\rr)
--(7.5,5*\rr)--(9,4*\rr)-- (10.5,5*\rr)--(10.5,7*\rr);
\draw [line width=0.4mm] (10.5,7*\rr)--(9,8*\rr)--(7.5,7*\rr)
--(7.5,5*\rr)--(9,4*\rr)-- (10.5,5*\rr)--(10.5,7*\rr);
\path [fill=lightgray] (12,8*\rr)--(10.5,7*\rr)--(10.5,5*\rr)
--(12,4*\rr)--(13.5,5*\rr)--(13.5,7*\rr)--(12,8*\rr);
\draw [line width=0.4mm] (12,8*\rr)--(10.5,7*\rr)--(10.5,5*\rr)
--(12,4*\rr)--(13.5,5*\rr)--(13.5,7*\rr)--(12,8*\rr);
\path [fill=lightgray] (12,8*\rr)--(13.5,7*\rr)--(15,8*\rr)
--(15,10*\rr)--(13.5,11*\rr)--(12,10*\rr)--(12,8*\rr);
\draw [line width=0.4mm] (12,8*\rr)--(13.5,7*\rr)--(15,8*\rr)
--(15,10*\rr)--(13.5,11*\rr)--(12,10*\rr)--(12,8*\rr);

\draw[yscale=sqrt(3/4), xslant=0.5] (-4,-2) grid (29, 12);
\draw[yscale=sqrt(3/4), xslant=-0.5] (-1,-2) grid (29, 12);

\draw [line width=0.8mm] (4.5, 3* \rr)--(7.5, 3* \rr)--(6,6* \rr)--(4.5, 3* \rr);
\draw [line width=0.8mm] (6,6* \rr)--(9,6* \rr)--(7.5, 3* \rr);
\draw [line width=0.8mm] (9,6* \rr)--(10.5, 3* \rr)--(9, 0)--(7.5, 3* \rr)
--(10.5, 3* \rr);
\foreach \pos in {(3, 0), (6, 0), (9, 0), (12, 0), (15, 0), 
(4.5, 3* \rr), (7.5, 3* \rr), (10.5, 3* \rr), (13.5, 3* \rr), (16.5, 3* \rr), 
(3,6* \rr),(6,6* \rr),(9,6* \rr),(12,6* \rr),(15,6* \rr),
(4.5,9* \rr), (7.5, 9* \rr), (10.5, 9* \rr), (13.5, 9* \rr), (16.5, 9* \rr),
(3,12* \rr),(6,12* \rr),(9,12* \rr),(12,12* \rr),(15,12* \rr)}
\shade[shading=ball, ball color=black] \pos circle (.4);
\end{tikzpicture} \qquad (b) \begin{tikzpicture}[scale=0.35]
\clip (2, -0.5) rectangle (17.4, 16.1);

\path [fill=lightgray] (13,5.33* \rr)--(12.5,7.66* \rr)--(10.5,8.33* \rr)
--(9,6.66* \rr)--(9.5,4.33* \rr)--(11.5,3.66* \rr)--(13,5.33* \rr);
\draw [line width=0.4mm] (13,5.33* \rr)--(12.5,7.66* \rr)--(10.5,8.33* \rr)
--(9,6.66* \rr)--(9.5,4.33* \rr)--(11.5,3.66* \rr)--(13,5.33* \rr);

\path [fill=lightgray] (16.5,6.33* \rr)--(16,8.66* \rr)--(14,9.33* \rr)
--(12.5,7.66* \rr)--(13,5.33* \rr)--(15,4.66* \rr)--(16.5,6.33* \rr);
\draw [line width=0.4mm] (16.5,6.33* \rr)--(16,8.66* \rr)--(14,9.33* \rr)
--(12.5,7.66* \rr)--(13,5.33* \rr)--(15,4.66* \rr)--(16.5,6.33* \rr);

\draw[yscale=sqrt(3/4), xslant=0.5] (-4,-2) grid (29, 12);
\draw[yscale=sqrt(3/4), xslant=-0.5] (-1,-2) grid (29, 12);
\draw [line width=0.8mm] (4,4* \rr)--(7.5,5* \rr)--(6.5,1* \rr)--(4,4* \rr);
\draw [line width=0.8mm] (4,4* \rr)--(5,8* \rr)--(8.5,9* \rr)--(7.5,5* \rr)
--(5,8* \rr);
\draw [line width=0.8mm] (8.5,9* \rr)--(11,6* \rr)--(7.5,5* \rr);
\foreach \pos in {(3, 0), (6.5,1* \rr), (10,2* \rr), (13.5,3* \rr), (17,4* \rr),
(4,4* \rr),(7.5,5* \rr),(11,6* \rr),(14.5,7* \rr),
(5,8* \rr),(8.5,9* \rr),(12,10* \rr),(15.5,11* \rr),
(2.5,11* \rr),(6,12* \rr)}
\shade[shading=ball, ball color=black] \pos circle (.4);
\end{tikzpicture}

\caption{PGSs on ${\mathbb A}_2$ for $D^2=9$ (Class TA1), 
and $D^2=13$ (Class TA2), with V-cells (gray hexagons) and choices of $D$-rhombuses.\\
{\footnotesize{The number of PGSs and EGMs for $D^2=9$ is 9 (frame (a)), and
for $D^2=13$ is 26 (frame (b)). The PGSs for $D^2=9$ are horizontal and for $D^2=13$ are inclined.}}
}
\end{figure}} %{Fig3} %FigureB2

%%%%%%%%%%%%%%%%%%%%%%%%%%%%%%%%%%%%%%%%%%%%%%%%%%

\def\FigureC3 %%\label{Fig3}
{\begin{figure}\label{Fig3} 
\centering
\captionsetup{width=0.8\textwidth} 
(a) \begin{tikzpicture}[scale=0.24]
\clip (2, -0.5) rectangle (24.4, 22.1);
\path [fill=lightgray] (10,4.66* \rr)--(10,9.33* \rr)--(13.5,11.66* \rr)--(17,9.33* \rr)
--(17,4.66* \rr)--(13.5,2.33* \rr)--(10,4.66* \rr);
\draw [line width=0.4mm] (10,4.66* \rr)--(10,9.33* \rr)--(13.5,11.66* \rr)--(17,9.33* \rr)
--(17,4.66* \rr)--(13.5,2.33* \rr)--(10,4.66* \rr);
\path [fill=lightgray] (10,9.33* \rr)--(13.5,11.66* \rr)--(13.5,16.33* \rr)--
(10,18.66* \rr)--(6.5,16.33* \rr)--(6.5,11.66* \rr)--(10,9.33* \rr);
\draw [line width=0.4mm] (10,9.33* \rr)--(13.5,11.66* \rr)--(13.5,16.33* \rr)--
(10,18.66* \rr)--(6.5,16.33* \rr)--(6.5,11.66* \rr)--(10,9.33* \rr);

\draw[yscale=sqrt(3/4), xslant=0.5] (-9,-2) grid (39, 21);
\draw[yscale=sqrt(3/4), xslant=-0.5] (-9,-2) grid (49, 21);

\draw [line width=0.8mm] (10,0)--(6.5, 7* \rr)--(13.5, 7* \rr)
--(10,0)--(17,0)--(13.5, 7* \rr);
\draw [line width=0.8mm] (6.5, 7* \rr)--(10,14* \rr)--(17,14* \rr)
--(13.5, 7* \rr)--(10,14* \rr);
\foreach \pos in {(3,0),(10,0),(17,0),(24,0), 
(6.5, 7* \rr), (13.5, 7* \rr), (20.5, 7* \rr), 
(3,14* \rr),(10,14* \rr),(17,14* \rr),(24,14* \rr),
(6.5,21* \rr), (13.5,21* \rr), (20.5,21* \rr)}
\shade[shading=ball, ball color=black] \pos circle (.4);
\end{tikzpicture} \quad (b) \begin{tikzpicture}[scale=0.24]
\clip (2, -0.5) rectangle (24.4, 22.1);

%\path [fill=lightgray] (8.22,7.33*\rr)--(6.5,11.5*\rr)--(8.75,15.3*\rr)
%--(13,14.7*\rr)--(14.5,10.7*\rr)--(12.2,6.8*\rr);
%\draw [line width=0.4mm] (8.22,7.33*\rr)--(6.5,11.5*\rr)--(8.75,15.3*\rr)
%%--(13,14.7*\rr)--(14.5,10.7*\rr)--(12.2,6.8*\rr)--(8.22,7.33*\rr);14.5

%\draw [line width=0.2mm,color=red] (9.5,3*\rr)--(5,16* \rr);  %(8,7.333*\rr), (6.5,11.666*\rr)
%\draw [line width=0.2mm,color=red] (4,8*\rr)--(11.5,19* \rr); %(6.5,11.666*\rr), (9,15.333*\rr)
%\draw [line width=0.2mm,color=red] (5,16* \rr)--(17,14* \rr); %(9,15.333*\rr), (13,14.666*\rr)
%\draw [line width=0.2mm,color=red] (17,14* \rr)--(9.5,3*\rr); %(12,6.666*\rr),(14.5,10.333*\rr)
%\draw [line width=0.2mm,color=red] (16,6*\rr)--(4,8*\rr);       %(8,7.333*\rr),(12,6.666*\rr)

%\draw [line width=0.2mm,color=red] (11.5,19* \rr)--(23.5,17* \rr); %(15.5,18.333*\rr), (19.5,17.666*\rr)
%\draw [line width=0.2mm,color=red] (23.5,17* \rr)--(16,6*\rr);       %(21,13.333*\rr), (18.5,9.666*\rr)
%\draw [line width=0.2mm,color=red] (16,6*\rr)--(11.5,19* \rr);        %(14.5,10.333*\rr), (13,14.666*\rr)
%\draw [line width=0.2mm,color=red] (10.5,11* \rr)--(22.5,9* \rr);    %(14.5,10.333*\rr), (18.5,9.666*\rr)
%\draw [line width=0.2mm,color=red] (22.5,9* \rr)--(18,22* \rr);       %(21,13.333*\rr), (19.5,17.666*\rr)
%\draw [line width=0.2mm,color=red] (10.5,11* \rr)--(18,22* \rr);     %(13,14.666*\rr), (15.5,18.333*\rr)

\path [fill=lightgray] (14.5,10.333*\rr)--(12,6.666*\rr)--(8,7.333* \rr)--(6.5,11.666* \rr)
-- (9,15.333* \rr)--(13,14.666*\rr);
\draw [line width=0.4mm] (14.5,10.333*\rr)--(12,6.666*\rr)--(8,7.333* \rr)--(6.5,11.666* \rr)
-- (9,15.333* \rr)--(13,14.666*\rr);
\path [fill=lightgray] (15.5,18.333* \rr)--(19.5,17.666* \rr)--(21,13.333* \rr)-- (18.5,9.666* \rr)--
(14.5,10.333* \rr)--(13,14.666* \rr)--(15.5,18.333* \rr);
\draw [line width=0.4mm] (15.5,18.333* \rr)--(19.5,17.666* \rr)--(21,13.333* \rr)-- (18.5,9.666* \rr)--
(14.5,10.333* \rr)--(13,14.666* \rr)--(15.5,18.333* \rr);

\draw[yscale=sqrt(3/4), xslant=0.5] (-9,-2) grid (39, 21);
\draw[yscale=sqrt(3/4), xslant=-0.5] (-9,-2) grid (49, 21);

\draw [line width=0.8mm] (4,8*\rr)--(5,16* \rr)--(11.5,19* \rr)
--(10.5,11* \rr)--(5,16* \rr);
\draw [line width=0.8mm] (9.5,3*\rr)--(4,8*\rr)--(10.5,11* \rr)
--(16,6*\rr)--(9.5,3*\rr)--(10.5,11* \rr);
\foreach \pos in {(3,0),(9.5,3*\rr),(16,6*\rr),(22.5,9* \rr),
(4,8*\rr),(10.5,11* \rr),(17,14* \rr),(23.5,17* \rr),
(5,16* \rr),(11.5,19* \rr)}
\shade[shading=ball, ball color=black] \pos circle (.4);
%\draw (10,-2.5*\rr) node  {\large{\bf (b)}};
\end{tikzpicture}
\caption{PGSs on ${\mathbb A}_2$ for $D^2=49$ (Class  TB),  with V-cells 
(gray hexagons) and choices of $D$-rhombuses.\\
{\footnotesize{There are 49 horizontal PGSs (frame (a)), and 98 inclined 
PGSs (frame (b)). The horizontal PGSs are the only dominant, so there are only 49 EGMs.}}
}
\end{figure}} %Fig3 %FigureC3

%%%%%%%%%%%%%%%%%%%%%%%%%%%%%%%%%%%%%%%%%%%%%%%%%%%%%%%%%%%%%%%%%%%%%%%

\def\FigureD4 %%\label{Fig4}
{\begin{figure}\label{Fig4}
\centering
\captionsetup{width=0.8\textwidth} 
(a) \begin{tikzpicture}[scale=0.25]
\clip (-0.9, -1.2*\rr) rectangle (15.5, 13*\rr);

\path [fill=lightgray] (6,8*\rr)--(10,8*\rr)--(12,4*\rr)--(10,0*\rr)
--(6,0*\rr)--(4,4*\rr)--(6,8*\rr);

\draw [line width=0.2mm] (0.5,1*\rr) -- (1.5,1*\rr) -- (2,0) -- (1.5,-1*\rr)
-- (0.5,-1*\rr) -- (0,0) -- (0.5,1*\rr);
\draw [line width=0.2mm] (0.5,1*\rr) -- (0,2*\rr) -- (0.5,3*\rr) -- (1.5,3*\rr) 
-- (2,2*\rr) -- (1.5,1*\rr);
\draw [line width=0.2mm] (0.5,3*\rr) -- (0,4*\rr) -- (0.5,5*\rr) -- (1.5,5*\rr) 
-- (2,4*\rr) -- (1.5,3*\rr);
\draw [line width=0.2mm] (0.5,5*\rr) -- (0,6*\rr) -- (0.5,7*\rr) -- (1.5,7*\rr) 
-- (2,6*\rr) -- (1.5,5*\rr);
\draw [line width=0.2mm] (0.5,7*\rr) -- (0,8*\rr) -- (0.5,9*\rr) -- (1.5,9*\rr) 
-- (2,8*\rr) -- (1.5,7*\rr);
\draw [line width=0.2mm] (0.5,9*\rr) -- (0,10*\rr) -- (0.5,11*\rr) -- (1.5,11*\rr) 
-- (2,10*\rr) -- (1.5,9*\rr);
\draw [line width=0.2mm] (0.5,11*\rr) -- (0,12*\rr) -- (0.5,13*\rr) -- (1.5,13*\rr) 
-- (2,12*\rr) -- (1.5,11*\rr);
\draw [line width=0.2mm] (0.5,13*\rr) -- (0,14*\rr) -- (0.5,15*\rr) -- (1.5,15*\rr) 
-- (2,14*\rr) -- (1.5,13*\rr);
\draw [line width=0.2mm] (0.5,15*\rr) -- (0,16*\rr) -- (0.5,17*\rr) -- (1.5,17*\rr) 
-- (2,16*\rr) -- (1.5,15*\rr);
\draw [line width=0.2mm] (0.5,17*\rr) -- (0,18*\rr) -- (0.5,19*\rr) -- (1.5,19*\rr) 
-- (2,18*\rr) -- (1.5,17*\rr);
\draw [line width=0.2mm] (0.5,19*\rr) -- (0,20*\rr) -- (0.5,21*\rr) -- (1.5,21*\rr) 
-- (2,20*\rr) -- (1.5,19*\rr);
%\draw [line width=0.2mm] (0.5,21*\rr) -- (0,22*\rr) -- (0.5,22*\rr) -- (1.5,22*\rr) 
%-- (2,22*\rr) -- (1.5,21*\rr);

\draw [line width=0.2mm]  (3,0) -- (3.5,-1*\rr);
\draw [line width=0.2mm] (2,2*\rr) -- (3,2*\rr) -- (3.5,1*\rr) -- (3,0) -- (2,0);
\draw [line width=0.2mm] (2,4*\rr) -- (3,4*\rr) -- (3.5,3*\rr) -- (3,2*\rr);
\draw [line width=0.2mm] (2,6*\rr) -- (3,6*\rr) -- (3.5,5*\rr) -- (3,4*\rr);
\draw [line width=0.2mm] (2,8*\rr) -- (3,8*\rr) -- (3.5,7*\rr) -- (3,6*\rr);
\draw [line width=0.2mm] (2,10*\rr) -- (3,10*\rr) -- (3.5,9*\rr) -- (3,8*\rr);
\draw [line width=0.2mm] (2,12*\rr) -- (3,12*\rr) -- (3.5,11*\rr) -- (3,10*\rr);
\draw [line width=0.2mm] (2,14*\rr) -- (3,14*\rr) -- (3.5,13*\rr) -- (3,12*\rr);
\draw [line width=0.2mm] (2,16*\rr) -- (3,16*\rr) -- (3.5,15*\rr) -- (3,14*\rr);
\draw [line width=0.2mm] (2,18*\rr) -- (3,18*\rr) -- (3.5,17*\rr) -- (3,16*\rr);
\draw [line width=0.2mm] (2,20*\rr) -- (3,20*\rr) -- (3.5,19*\rr) -- (3,18*\rr);
\draw [line width=0.2mm] (3.5,21*\rr) -- (3,20*\rr);

\draw [line width=0.2mm]  (3.5,1*\rr) -- (4.5,1*\rr) -- (5,0) -- (4.5,-1*\rr)
--(3.5,-1*\rr);
\draw [line width=0.2mm]  (3.5,3*\rr) -- (4.5,3*\rr) -- (5,2*\rr) -- (4.5,1*\rr);
\draw [line width=0.2mm]  (3.5,5*\rr) -- (4.5,5*\rr) -- (5,4*\rr) -- (4.5,3*\rr);
\draw [line width=0.2mm]  (3.5,7*\rr) -- (4.5,7*\rr) -- (5,6*\rr) -- (4.5,5*\rr);
\draw [line width=0.2mm]  (3.5,9*\rr) -- (4.5,9*\rr) -- (5,8*\rr) -- (4.5,7*\rr);
\draw [line width=0.2mm]  (3.5,11*\rr) -- (4.5,11*\rr) -- (5,10*\rr) -- (4.5,9*\rr);
\draw [line width=0.2mm]  (3.5,13*\rr) -- (4.5,13*\rr) -- (5,12*\rr) -- (4.5,11*\rr);
\draw [line width=0.2mm]  (3.5,15*\rr) -- (4.5,15*\rr) -- (5,14*\rr) -- (4.5,13*\rr);
\draw [line width=0.2mm]  (3.5,17*\rr) -- (4.5,17*\rr) -- (5,16*\rr) -- (4.5,15*\rr);
\draw [line width=0.2mm]  (3.5,19*\rr) -- (4.5,19*\rr) -- (5,18*\rr) -- (4.5,17*\rr);
\draw [line width=0.2mm]  (3.5,21*\rr) -- (4.5,21*\rr) -- (5,20*\rr) -- (4.5,19*\rr);

\draw [line width=0.2mm]  (5,0) -- (6,0)--(6.5,-1*\rr );
\draw [line width=0.2mm]  (5,2*\rr) -- (6,2*\rr)--(6.5,1*\rr )--(6,0);
\draw [line width=0.2mm]  (5,4*\rr) -- (6,4*\rr)--(6.5,3*\rr )--(6,2*\rr);
\draw [line width=0.2mm]  (5,6*\rr) -- (6,6*\rr)--(6.5,5*\rr )--(6,4*\rr);
\draw [line width=0.2mm]  (5,8*\rr) -- (6,8*\rr)--(6.5,7*\rr )--(6,6*\rr);
\draw [line width=0.2mm]  (5,10*\rr) -- (6,10*\rr)--(6.5,9*\rr )--(6,8*\rr);
\draw [line width=0.2mm]  (5,12*\rr) -- (6,12*\rr)--(6.5,11*\rr )--(6,10*\rr);
\draw [line width=0.2mm]  (5,14*\rr) -- (6,14*\rr)--(6.5,13*\rr )--(6,12*\rr);
\draw [line width=0.2mm]  (5,16*\rr) -- (6,16*\rr)--(6.5,15*\rr )--(6,14*\rr);
\draw [line width=0.2mm]  (5,18*\rr) -- (6,18*\rr)--(6.5,17*\rr )--(6,16*\rr);
\draw [line width=0.2mm]  (5,20*\rr) -- (6,20*\rr)--(6.5,19*\rr )--(6,18*\rr);
\draw [line width=0.2mm]  (6.5,21*\rr )--(6,20*\rr);

\draw [line width=0.2mm]  (6.5,1*\rr) -- (7.5,1*\rr)--(8,0)--(7.5,-1*\rr)
--(6.5,-1*\rr);
\draw [line width=0.2mm]  (6.5,3*\rr) -- (7.5,3*\rr)--(8,2*\rr)--(7.5,1*\rr);
\draw [line width=0.2mm]  (6.5,5*\rr) -- (7.5,5*\rr)--(8,4*\rr)--(7.5,3*\rr);
\draw [line width=0.2mm]  (6.5,7*\rr) -- (7.5,7*\rr)--(8,6*\rr)--(7.5,5*\rr);
\draw [line width=0.2mm]  (6.5,9*\rr) -- (7.5,9*\rr)--(8,8*\rr)--(7.5,7*\rr);
\draw [line width=0.2mm]  (6.5,11*\rr) -- (7.5,11*\rr)--(8,10*\rr)--(7.5,9*\rr);
\draw [line width=0.2mm]  (6.5,13*\rr) -- (7.5,13*\rr)--(8,12*\rr)--(7.5,11*\rr);
\draw [line width=0.2mm]  (6.5,15*\rr) -- (7.5,15*\rr)--(8,14*\rr)--(7.5,13*\rr);
\draw [line width=0.2mm]  (6.5,17*\rr) -- (7.5,17*\rr)--(8,16*\rr)--(7.5,15*\rr);
\draw [line width=0.2mm]  (6.5,19*\rr) -- (7.5,19*\rr)--(8,18*\rr)--(7.5,17*\rr);
\draw [line width=0.2mm]  (6.5,21*\rr) -- (7.5,21*\rr)--(8,20*\rr)--(7.5,19*\rr);

\draw [line width=0.2mm]  (8,0)--(9,0)--(9.5,-1*\rr);
\draw [line width=0.2mm]  (8,2*\rr)--(9,2*\rr)--(9.5,1*\rr)--(9,0);
\draw [line width=0.2mm]  (8,4*\rr)--(9,4*\rr)--(9.5,3*\rr)--(9,2*\rr);
\draw [line width=0.2mm]  (8,6*\rr)--(9,6*\rr)--(9.5,5*\rr)--(9,4*\rr);
\draw [line width=0.2mm]  (8,8*\rr)--(9,8*\rr)--(9.5,7*\rr)--(9,6*\rr);
\draw [line width=0.2mm]  (8,10*\rr)--(9,10*\rr)--(9.5,9*\rr)--(9,8*\rr);
\draw [line width=0.2mm]  (8,12*\rr)--(9,12*\rr)--(9.5,11*\rr)--(9,10*\rr);
\draw [line width=0.2mm]  (8,14*\rr)--(9,14*\rr)--(9.5,13*\rr)--(9,12*\rr);
\draw [line width=0.2mm]  (8,16*\rr)--(9,16*\rr)--(9.5,15*\rr)--(9,14*\rr);
\draw [line width=0.2mm]  (8,18*\rr)--(9,18*\rr)--(9.5,17*\rr)--(9,16*\rr);
\draw [line width=0.2mm]  (8,20*\rr)--(9,20*\rr)--(9.5,19*\rr)--(9,18*\rr);
\draw [line width=0.2mm]  (9.5,21*\rr)--(9,20*\rr);

\draw [line width=0.2mm]  (9.5,1*\rr)--(10.5,1*\rr)--(11,0)--(10.5,-1*\rr)
--(9.5,-1*\rr);
\draw [line width=0.2mm]  (9.5,3*\rr)--(10.5,3*\rr)--(11,2*\rr)--(10.5,1*\rr);
\draw [line width=0.2mm]  (9.5,5*\rr)--(10.5,5*\rr)--(11,4*\rr)--(10.5,3*\rr);
\draw [line width=0.2mm]  (9.5,7*\rr)--(10.5,7*\rr)--(11,6*\rr)--(10.5,5*\rr);
\draw [line width=0.2mm]  (9.5,9*\rr)--(10.5,9*\rr)--(11,8*\rr)--(10.5,7*\rr);
\draw [line width=0.2mm]  (9.5,11*\rr)--(10.5,11*\rr)--(11,10*\rr)--(10.5,9*\rr);
\draw [line width=0.2mm]  (9.5,13*\rr)--(10.5,13*\rr)--(11,12*\rr)--(10.5,11*\rr);
\draw [line width=0.2mm]  (9.5,15*\rr)--(10.5,15*\rr)--(11,14*\rr)--(10.5,13*\rr);
\draw [line width=0.2mm]  (9.5,17*\rr)--(10.5,17*\rr)--(11,16*\rr)--(10.5,15*\rr);
\draw [line width=0.2mm]  (9.5,19*\rr)--(10.5,19*\rr)--(11,18*\rr)--(10.5,17*\rr);
\draw [line width=0.2mm]  (9.5,21*\rr)--(10.5,21*\rr)--(11,20*\rr)--(10.5,19*\rr);

\draw [line width=0.2mm]  (11,0)--(12,0)--(12.5,-1*\rr);
\draw [line width=0.2mm]  (11,2*\rr)--(12,2*\rr)--(12.5,1*\rr)--(12,0);
\draw [line width=0.2mm]  (11,4*\rr)--(12,4*\rr)--(12.5,3*\rr)--(12,2*\rr);
\draw [line width=0.2mm]  (11,6*\rr)--(12,6*\rr)--(12.5,5*\rr)--(12,4*\rr);
\draw [line width=0.2mm]  (11,8*\rr)--(12,8*\rr)--(12.5,7*\rr)--(12,6*\rr);
\draw [line width=0.2mm]  (11,10*\rr)--(12,10*\rr)--(12.5,9*\rr)--(12,8*\rr);
\draw [line width=0.2mm]  (11,12*\rr)--(12,12*\rr)--(12.5,11*\rr)--(12,10*\rr);
\draw [line width=0.2mm]  (11,14*\rr)--(12,14*\rr)--(12.5,13*\rr)--(12,12*\rr);
\draw [line width=0.2mm]  (11,16*\rr)--(12,16*\rr)--(12.5,15*\rr)--(12,14*\rr);
\draw [line width=0.2mm]  (11,18*\rr)--(12,18*\rr)--(12.5,17*\rr)--(12,16*\rr);
\draw [line width=0.2mm]  (11,20*\rr)--(12,20*\rr)--(12.5,19*\rr)--(12,18*\rr);
\draw [line width=0.2mm]  (12,20*\rr)--(12.5,21*\rr);

\draw [line width=0.2mm]  (12.5,1*\rr)--(13.5,1*\rr)--(14,0)--(13.5,-1*\rr)
--(12.5,-1*\rr);
\draw [line width=0.2mm]  (12.5,3*\rr)--(13.5,3*\rr)--(14,2*\rr)--(13.5,1*\rr);
\draw [line width=0.2mm]  (12.5,5*\rr)--(13.5,5*\rr)--(14,4*\rr)--(13.5,3*\rr);
\draw [line width=0.2mm]  (12.5,7*\rr)--(13.5,7*\rr)--(14,6*\rr)--(13.5,5*\rr);
\draw [line width=0.2mm]  (12.5,9*\rr)--(13.5,9*\rr)--(14,8*\rr)--(13.5,7*\rr);
\draw [line width=0.2mm]  (12.5,11*\rr)--(13.5,11*\rr)--(14,10*\rr)--(13.5,9*\rr);
\draw [line width=0.2mm]  (12.5,13*\rr)--(13.5,13*\rr)--(14,12*\rr)--(13.5,11*\rr);
\draw [line width=0.2mm]  (12.5,15*\rr)--(13.5,15*\rr)--(14,14*\rr)--(13.5,13*\rr);
\draw [line width=0.2mm]  (12.5,17*\rr)--(13.5,17*\rr)--(14,16*\rr)--(13.5,15*\rr);
\draw [line width=0.2mm]  (12.5,19*\rr)--(13.5,19*\rr)--(14,18*\rr)--(13.5,17*\rr);
\draw [line width=0.2mm]  (12.5,21*\rr)--(13.5,21*\rr)--(14,20*\rr)--(13.5,19*\rr);

\draw [line width=0.2mm]  (14,0)--(15,0)--(15.5,-1*\rr);
\draw [line width=0.2mm]  (14,2*\rr)--(15,2*\rr)--(15.5,1*\rr)--(15,0);
\draw [line width=0.2mm]  (14,4*\rr)--(15,4*\rr)--(15.5,3*\rr)--(15,2*\rr);
\draw [line width=0.2mm]  (14,6*\rr)--(15,6*\rr)--(15.5,5*\rr)--(15,4*\rr);
\draw [line width=0.2mm]  (14,8*\rr)--(15,8*\rr)--(15.5,7*\rr)--(15,6*\rr);
\draw [line width=0.2mm]  (14,10*\rr)--(15,10*\rr)--(15.5,9*\rr)--(15,8*\rr);
\draw [line width=0.2mm]  (14,12*\rr)--(15,12*\rr)--(15.5,11*\rr)--(15,10*\rr);
\draw [line width=0.2mm]  (14,14*\rr)--(15,14*\rr)--(15.5,13*\rr)--(15,12*\rr);
\draw [line width=0.2mm]  (14,16*\rr)--(15,16*\rr)--(15.5,15*\rr)--(15,14*\rr);
\draw [line width=0.2mm]  (14,18*\rr)--(15,18*\rr)--(15.5,17*\rr)--(15,16*\rr);
\draw [line width=0.2mm]  (14,20*\rr)--(15,20*\rr)--(15.5,19*\rr)--(15,18*\rr);
\draw [line width=0.2mm]  (15,20*\rr)--(15.5,21*\rr);

\draw [line width=0.2mm]  (15.5,1*\rr)--(16.5,1*\rr)--(17,0)--(16.5,-1*\rr)
--(15.5,-1*\rr);
\draw [line width=0.2mm]  (15.5,3*\rr)--(16.5,3*\rr)--(17,2*\rr)--(16.5,1*\rr);
\draw [line width=0.2mm]  (15.5,5*\rr)--(16.5,5*\rr)--(17,4*\rr)--(16.5,3*\rr);
\draw [line width=0.2mm]  (15.5,7*\rr)--(16.5,7*\rr)--(17,6*\rr)--(16.5,5*\rr);
\draw [line width=0.2mm]  (15.5,9*\rr)--(16.5,9*\rr)--(17,8*\rr)--(16.5,7*\rr);
\draw [line width=0.2mm]  (15.5,11*\rr)--(16.5,11*\rr)--(17,10*\rr)--(16.5,9*\rr);
\draw [line width=0.2mm]  (15.5,13*\rr)--(16.5,13*\rr)--(17,12*\rr)--(16.5,11*\rr);
\draw [line width=0.2mm]  (15.5,15*\rr)--(16.5,15*\rr)--(17,14*\rr)--(16.5,13*\rr);
\draw [line width=0.2mm]  (15.5,17*\rr)--(16.5,17*\rr)--(17,16*\rr)--(16.5,15*\rr);
\draw [line width=0.2mm]  (15.5,19*\rr)--(16.5,19*\rr)--(17,18*\rr)--(16.5,17*\rr);
\draw [line width=0.2mm]  (15.5,21*\rr)--(16.5,21*\rr)--(17,20*\rr)--(16.5,19*\rr);

\draw [line width=0.2mm]  (17,0)--(18,0)--(18.5,-1*\rr);
\draw [line width=0.2mm]  (17,2*\rr)--(18,2*\rr)--(18.5,1*\rr)--(18,0);
\draw [line width=0.2mm]  (17,4*\rr)--(18,4*\rr)--(18.5,3*\rr)--(18,2*\rr);
\draw [line width=0.2mm]  (17,6*\rr)--(18,6*\rr)--(18.5,5*\rr)--(18,4*\rr);
\draw [line width=0.2mm]  (17,8*\rr)--(18,8*\rr)--(18.5,7*\rr)--(18,6*\rr);
\draw [line width=0.2mm]  (17,10*\rr)--(18,10*\rr)--(18.5,9*\rr)--(18,8*\rr);
\draw [line width=0.2mm]  (17,12*\rr)--(18,12*\rr)--(18.5,11*\rr)--(18,10*\rr);
\draw [line width=0.2mm]  (17,14*\rr)--(18,14*\rr)--(18.5,13*\rr)--(18,12*\rr);
\draw [line width=0.2mm]  (17,16*\rr)--(18,16*\rr)--(18.5,15*\rr)--(18,14*\rr);
\draw [line width=0.2mm]  (17,18*\rr)--(18,18*\rr)--(18.5,17*\rr)--(18,16*\rr);
\draw [line width=0.2mm]  (17,20*\rr)--(18,20*\rr)--(18.5,19*\rr)--(18,18*\rr);
\draw [line width=0.2mm]  (18,20*\rr)--(18.5,21*\rr);

\draw [line width=0.2mm]  (18.5,1*\rr)--(19.5,1*\rr)--(20,0)--(19.5,-1*\rr)
--(18.5,-1*\rr);
\draw [line width=0.2mm]  (18.5,3*\rr)--(19.5,3*\rr)--(20,2*\rr)--(19.5,1*\rr);
\draw [line width=0.2mm]  (18.5,5*\rr)--(19.5,5*\rr)--(20,4*\rr)--(19.5,3*\rr);
\draw [line width=0.2mm]  (18.5,7*\rr)--(19.5,7*\rr)--(20,6*\rr)--(19.5,5*\rr);
\draw [line width=0.2mm]  (18.5,9*\rr)--(19.5,9*\rr)--(20,8*\rr)--(19.5,7*\rr);
\draw [line width=0.2mm]  (18.5,11*\rr)--(19.5,11*\rr)--(20,10*\rr)--(19.5,9*\rr);
\draw [line width=0.2mm]  (18.5,13*\rr)--(19.5,13*\rr)--(20,12*\rr)--(19.5,11*\rr);
\draw [line width=0.2mm]  (18.5,15*\rr)--(19.5,15*\rr)--(20,14*\rr)--(19.5,13*\rr);
\draw [line width=0.2mm]  (18.5,17*\rr)--(19.5,17*\rr)--(20,16*\rr)--(19.5,15*\rr);
\draw [line width=0.2mm]  (18.5,19*\rr)--(19.5,19*\rr)--(20,18*\rr)--(19.5,17*\rr);
\draw [line width=0.2mm]  (18.5,21*\rr)--(19.5,21*\rr)--(20,20*\rr)--(19.5,19*\rr);

\draw [line width=0.2mm]  (20,0)--(21,0)--(21.5,-1*\rr);
\draw [line width=0.2mm]  (20,2*\rr)--(21,2*\rr)--(21.5,1*\rr)--(21,0);
\draw [line width=0.2mm]  (20,4*\rr)--(21,4*\rr)--(21.5,3*\rr)--(21,2*\rr);
\draw [line width=0.2mm]  (20,6*\rr)--(21,6*\rr)--(21.5,5*\rr)--(21,4*\rr);
\draw [line width=0.2mm]  (20,8*\rr)--(21,8*\rr)--(21.5,7*\rr)--(21,6*\rr);
\draw [line width=0.2mm]  (20,10*\rr)--(21,10*\rr)--(21.5,9*\rr)--(21,8*\rr);
\draw [line width=0.2mm]  (20,12*\rr)--(21,12*\rr)--(21.5,11*\rr)--(21,10*\rr);
\draw [line width=0.2mm]  (20,14*\rr)--(21,14*\rr)--(21.5,13*\rr)--(21,12*\rr);
\draw [line width=0.2mm]  (20,16*\rr)--(21,16*\rr)--(21.5,15*\rr)--(21,14*\rr);
\draw [line width=0.2mm]  (20,18*\rr)--(21,18*\rr)--(21.5,17*\rr)--(21,16*\rr);
\draw [line width=0.2mm]  (20,20*\rr)--(21,20*\rr)--(21.5,19*\rr)--(21,18*\rr);
\draw [line width=0.2mm]  (21,20*\rr)--(21.5,21*\rr);

\draw [line width=0.2mm]  (21.5,1*\rr)--(22.5,1*\rr)--(23,0)--(22.5,-1*\rr)
--(21.5,-1*\rr);
\draw [line width=0.2mm]  (21.5,3*\rr)--(22.5,3*\rr)--(23,2*\rr)--(22.5,1*\rr);
\draw [line width=0.2mm]  (21.5,5*\rr)--(22.5,5*\rr)--(23,4*\rr)--(22.5,3*\rr);
\draw [line width=0.2mm]  (21.5,7*\rr)--(22.5,7*\rr)--(23,6*\rr)--(22.5,5*\rr);
\draw [line width=0.2mm]  (21.5,9*\rr)--(22.5,9*\rr)--(23,8*\rr)--(22.5,7*\rr);
\draw [line width=0.2mm]  (21.5,11*\rr)--(22.5,11*\rr)--(23,10*\rr)--(22.5,9*\rr);
\draw [line width=0.2mm]  (21.5,13*\rr)--(22.5,13*\rr)--(23,12*\rr)--(22.5,11*\rr);
\draw [line width=0.2mm]  (21.5,15*\rr)--(22.5,15*\rr)--(23,14*\rr)--(22.5,13*\rr);
\draw [line width=0.2mm]  (21.5,17*\rr)--(22.5,17*\rr)--(23,16*\rr)--(22.5,15*\rr);
\draw [line width=0.2mm]  (21.5,19*\rr)--(22.5,19*\rr)--(23,18*\rr)--(22.5,17*\rr);
\draw [line width=0.2mm]  (21.5,21*\rr)--(22.5,21*\rr)--(23,20*\rr)--(22.5,19*\rr);

\draw [line width=0.2mm]  (23,0)--(24,0)--(24.5,-1*\rr);
\draw [line width=0.2mm]  (23,2*\rr)--(24,2*\rr)--(24.5,1*\rr)--(24,0);
\draw [line width=0.2mm]  (23,4*\rr)--(24,4*\rr)--(24.5,3*\rr)--(24,2*\rr);
\draw [line width=0.2mm]  (23,6*\rr)--(24,6*\rr)--(24.5,5*\rr)--(24,4*\rr);
\draw [line width=0.2mm]  (23,8*\rr)--(24,8*\rr)--(24.5,7*\rr)--(24,6*\rr);
\draw [line width=0.2mm]  (23,10*\rr)--(24,10*\rr)--(24.5,9*\rr)--(24,8*\rr);
\draw [line width=0.2mm]  (23,12*\rr)--(24,12*\rr)--(24.5,11*\rr)--(24,10*\rr);
\draw [line width=0.2mm]  (23,14*\rr)--(24,14*\rr)--(24.5,13*\rr)--(24,12*\rr);
\draw [line width=0.2mm]  (23,16*\rr)--(24,16*\rr)--(24.5,15*\rr)--(24,14*\rr);
\draw [line width=0.2mm]  (23,18*\rr)--(24,18*\rr)--(24.5,17*\rr)--(24,16*\rr);
\draw [line width=0.2mm]  (23,20*\rr)--(24,20*\rr)--(24.5,19*\rr)--(24,18*\rr);
\draw [line width=0.2mm]  (24,20*\rr)--(24.5,21*\rr);

\draw [line width=0.6mm] (2,0)--(2,8*\rr)--(8,4*\rr)--(2,0);
\draw [line width=0.6mm] (2,8*\rr)--(8,12*\rr)--(8,4*\rr);
\draw [line width=0.6mm] (8,12*\rr)--(14,8*\rr)--(8,4*\rr)--(14,0)--(14,8*\rr);

\foreach \pos in {(2,0), (2,8*\rr),(8,12*\rr),(8,4*\rr),(14,0),(14,8*\rr)}
\shade[shading=ball, ball color=black] \pos circle (.4);

\draw [line width=0.3mm] (6,8*\rr)--(10,8*\rr)--(12,4*\rr)--(10,0*\rr)
--(6,0*\rr)--(4,4*\rr)--(6,8*\rr);
 
\end{tikzpicture}\quad (b) \begin{tikzpicture}[scale=0.25]
\clip (-0.9, -1.2*\rr) rectangle (19, 13*\rr);

\path [fill=lightgray] (13.5,3*\rr)--(9,5)--(9,12*\rr)--(10.5,13*\rr)
--(16.5,13*\rr)--(18,12*\rr)--(18,6*\rr)--(13.5,3*\rr);

\draw [line width=0.2mm] (0.5,1*\rr) -- (1.5,1*\rr) -- (2,0) -- (1.5,-1*\rr)
-- (0.5,-1*\rr) -- (0,0) -- (0.5,1*\rr);
\draw [line width=0.2mm] (0.5,1*\rr) -- (0,2*\rr) -- (0.5,3*\rr) -- (1.5,3*\rr) 
-- (2,2*\rr) -- (1.5,1*\rr);
\draw [line width=0.2mm] (0.5,3*\rr) -- (0,4*\rr) -- (0.5,5*\rr) -- (1.5,5*\rr) 
-- (2,4*\rr) -- (1.5,3*\rr);
\draw [line width=0.2mm] (0.5,5*\rr) -- (0,6*\rr) -- (0.5,7*\rr) -- (1.5,7*\rr) 
-- (2,6*\rr) -- (1.5,5*\rr);
\draw [line width=0.2mm] (0.5,7*\rr) -- (0,8*\rr) -- (0.5,9*\rr) -- (1.5,9*\rr) 
-- (2,8*\rr) -- (1.5,7*\rr);
\draw [line width=0.2mm] (0.5,9*\rr) -- (0,10*\rr) -- (0.5,11*\rr) -- (1.5,11*\rr) 
-- (2,10*\rr) -- (1.5,9*\rr);
\draw [line width=0.2mm] (0.5,11*\rr) -- (0,12*\rr) -- (0.5,13*\rr) -- (1.5,13*\rr) 
-- (2,12*\rr) -- (1.5,11*\rr);
\draw [line width=0.2mm] (0.5,13*\rr) -- (0,14*\rr) -- (0.5,15*\rr) -- (1.5,15*\rr) 
-- (2,14*\rr) -- (1.5,13*\rr);
\draw [line width=0.2mm] (0.5,15*\rr) -- (0,16*\rr) -- (0.5,17*\rr) -- (1.5,17*\rr) 
-- (2,16*\rr) -- (1.5,15*\rr);
\draw [line width=0.2mm] (0.5,17*\rr) -- (0,18*\rr) -- (0.5,19*\rr) -- (1.5,19*\rr) 
-- (2,18*\rr) -- (1.5,17*\rr);
\draw [line width=0.2mm] (0.5,19*\rr) -- (0,20*\rr) -- (0.5,21*\rr) -- (1.5,21*\rr) 
-- (2,20*\rr) -- (1.5,19*\rr);
%\draw [line width=0.2mm] (0.5,21*\rr) -- (0,22*\rr) -- (0.5,22*\rr) -- (1.5,22*\rr) 
%-- (2,22*\rr) -- (1.5,21*\rr);

\draw [line width=0.2mm]  (3,0) -- (3.5,-1*\rr);
\draw [line width=0.2mm] (2,2*\rr) -- (3,2*\rr) -- (3.5,1*\rr) -- (3,0) -- (2,0);
\draw [line width=0.2mm] (2,4*\rr) -- (3,4*\rr) -- (3.5,3*\rr) -- (3,2*\rr);
\draw [line width=0.2mm] (2,6*\rr) -- (3,6*\rr) -- (3.5,5*\rr) -- (3,4*\rr);
\draw [line width=0.2mm] (2,8*\rr) -- (3,8*\rr) -- (3.5,7*\rr) -- (3,6*\rr);
\draw [line width=0.2mm] (2,10*\rr) -- (3,10*\rr) -- (3.5,9*\rr) -- (3,8*\rr);
\draw [line width=0.2mm] (2,12*\rr) -- (3,12*\rr) -- (3.5,11*\rr) -- (3,10*\rr);
\draw [line width=0.2mm] (2,14*\rr) -- (3,14*\rr) -- (3.5,13*\rr) -- (3,12*\rr);
\draw [line width=0.2mm] (2,16*\rr) -- (3,16*\rr) -- (3.5,15*\rr) -- (3,14*\rr);
\draw [line width=0.2mm] (2,18*\rr) -- (3,18*\rr) -- (3.5,17*\rr) -- (3,16*\rr);
\draw [line width=0.2mm] (2,20*\rr) -- (3,20*\rr) -- (3.5,19*\rr) -- (3,18*\rr);
\draw [line width=0.2mm] (3.5,21*\rr) -- (3,20*\rr);

\draw [line width=0.2mm]  (3.5,1*\rr) -- (4.5,1*\rr) -- (5,0) -- (4.5,-1*\rr)
--(3.5,-1*\rr);
\draw [line width=0.2mm]  (3.5,3*\rr) -- (4.5,3*\rr) -- (5,2*\rr) -- (4.5,1*\rr);
\draw [line width=0.2mm]  (3.5,5*\rr) -- (4.5,5*\rr) -- (5,4*\rr) -- (4.5,3*\rr);
\draw [line width=0.2mm]  (3.5,7*\rr) -- (4.5,7*\rr) -- (5,6*\rr) -- (4.5,5*\rr);
\draw [line width=0.2mm]  (3.5,9*\rr) -- (4.5,9*\rr) -- (5,8*\rr) -- (4.5,7*\rr);
\draw [line width=0.2mm]  (3.5,11*\rr) -- (4.5,11*\rr) -- (5,10*\rr) -- (4.5,9*\rr);
\draw [line width=0.2mm]  (3.5,13*\rr) -- (4.5,13*\rr) -- (5,12*\rr) -- (4.5,11*\rr);
\draw [line width=0.2mm]  (3.5,15*\rr) -- (4.5,15*\rr) -- (5,14*\rr) -- (4.5,13*\rr);
\draw [line width=0.2mm]  (3.5,17*\rr) -- (4.5,17*\rr) -- (5,16*\rr) -- (4.5,15*\rr);
\draw [line width=0.2mm]  (3.5,19*\rr) -- (4.5,19*\rr) -- (5,18*\rr) -- (4.5,17*\rr);
\draw [line width=0.2mm]  (3.5,21*\rr) -- (4.5,21*\rr) -- (5,20*\rr) -- (4.5,19*\rr);

\draw [line width=0.2mm]  (5,0) -- (6,0)--(6.5,-1*\rr );
\draw [line width=0.2mm]  (5,2*\rr) -- (6,2*\rr)--(6.5,1*\rr )--(6,0);
\draw [line width=0.2mm]  (5,4*\rr) -- (6,4*\rr)--(6.5,3*\rr )--(6,2*\rr);
\draw [line width=0.2mm]  (5,6*\rr) -- (6,6*\rr)--(6.5,5*\rr )--(6,4*\rr);
\draw [line width=0.2mm]  (5,8*\rr) -- (6,8*\rr)--(6.5,7*\rr )--(6,6*\rr);
\draw [line width=0.2mm]  (5,10*\rr) -- (6,10*\rr)--(6.5,9*\rr )--(6,8*\rr);
\draw [line width=0.2mm]  (5,12*\rr) -- (6,12*\rr)--(6.5,11*\rr )--(6,10*\rr);
\draw [line width=0.2mm]  (5,14*\rr) -- (6,14*\rr)--(6.5,13*\rr )--(6,12*\rr);
\draw [line width=0.2mm]  (5,16*\rr) -- (6,16*\rr)--(6.5,15*\rr )--(6,14*\rr);
\draw [line width=0.2mm]  (5,18*\rr) -- (6,18*\rr)--(6.5,17*\rr )--(6,16*\rr);
\draw [line width=0.2mm]  (5,20*\rr) -- (6,20*\rr)--(6.5,19*\rr )--(6,18*\rr);
\draw [line width=0.2mm]  (6.5,21*\rr )--(6,20*\rr);

\draw [line width=0.2mm]  (6.5,1*\rr) -- (7.5,1*\rr)--(8,0)--(7.5,-1*\rr)
--(6.5,-1*\rr);
\draw [line width=0.2mm]  (6.5,3*\rr) -- (7.5,3*\rr)--(8,2*\rr)--(7.5,1*\rr);
\draw [line width=0.2mm]  (6.5,5*\rr) -- (7.5,5*\rr)--(8,4*\rr)--(7.5,3*\rr);
\draw [line width=0.2mm]  (6.5,7*\rr) -- (7.5,7*\rr)--(8,6*\rr)--(7.5,5*\rr);
\draw [line width=0.2mm]  (6.5,9*\rr) -- (7.5,9*\rr)--(8,8*\rr)--(7.5,7*\rr);
\draw [line width=0.2mm]  (6.5,11*\rr) -- (7.5,11*\rr)--(8,10*\rr)--(7.5,9*\rr);
\draw [line width=0.2mm]  (6.5,13*\rr) -- (7.5,13*\rr)--(8,12*\rr)--(7.5,11*\rr);
\draw [line width=0.2mm]  (6.5,15*\rr) -- (7.5,15*\rr)--(8,14*\rr)--(7.5,13*\rr);
\draw [line width=0.2mm]  (6.5,17*\rr) -- (7.5,17*\rr)--(8,16*\rr)--(7.5,15*\rr);
\draw [line width=0.2mm]  (6.5,19*\rr) -- (7.5,19*\rr)--(8,18*\rr)--(7.5,17*\rr);
\draw [line width=0.2mm]  (6.5,21*\rr) -- (7.5,21*\rr)--(8,20*\rr)--(7.5,19*\rr);

\draw [line width=0.2mm]  (8,0)--(9,0)--(9.5,-1*\rr);
\draw [line width=0.2mm]  (8,2*\rr)--(9,2*\rr)--(9.5,1*\rr)--(9,0);
\draw [line width=0.2mm]  (8,4*\rr)--(9,4*\rr)--(9.5,3*\rr)--(9,2*\rr);
\draw [line width=0.2mm]  (8,6*\rr)--(9,6*\rr)--(9.5,5*\rr)--(9,4*\rr);
\draw [line width=0.2mm]  (8,8*\rr)--(9,8*\rr)--(9.5,7*\rr)--(9,6*\rr);
\draw [line width=0.2mm]  (8,10*\rr)--(9,10*\rr)--(9.5,9*\rr)--(9,8*\rr);
\draw [line width=0.2mm]  (8,12*\rr)--(9,12*\rr)--(9.5,11*\rr)--(9,10*\rr);
\draw [line width=0.2mm]  (8,14*\rr)--(9,14*\rr)--(9.5,13*\rr)--(9,12*\rr);
\draw [line width=0.2mm]  (8,16*\rr)--(9,16*\rr)--(9.5,15*\rr)--(9,14*\rr);
\draw [line width=0.2mm]  (8,18*\rr)--(9,18*\rr)--(9.5,17*\rr)--(9,16*\rr);
\draw [line width=0.2mm]  (8,20*\rr)--(9,20*\rr)--(9.5,19*\rr)--(9,18*\rr);
\draw [line width=0.2mm]  (9.5,21*\rr)--(9,20*\rr);

\draw [line width=0.2mm]  (9.5,1*\rr)--(10.5,1*\rr)--(11,0)--(10.5,-1*\rr)
--(9.5,-1*\rr);
\draw [line width=0.2mm]  (9.5,3*\rr)--(10.5,3*\rr)--(11,2*\rr)--(10.5,1*\rr);
\draw [line width=0.2mm]  (9.5,5*\rr)--(10.5,5*\rr)--(11,4*\rr)--(10.5,3*\rr);
\draw [line width=0.2mm]  (9.5,7*\rr)--(10.5,7*\rr)--(11,6*\rr)--(10.5,5*\rr);
\draw [line width=0.2mm]  (9.5,9*\rr)--(10.5,9*\rr)--(11,8*\rr)--(10.5,7*\rr);
\draw [line width=0.2mm]  (9.5,11*\rr)--(10.5,11*\rr)--(11,10*\rr)--(10.5,9*\rr);
\draw [line width=0.2mm]  (9.5,13*\rr)--(10.5,13*\rr)--(11,12*\rr)--(10.5,11*\rr);
\draw [line width=0.2mm]  (9.5,15*\rr)--(10.5,15*\rr)--(11,14*\rr)--(10.5,13*\rr);
\draw [line width=0.2mm]  (9.5,17*\rr)--(10.5,17*\rr)--(11,16*\rr)--(10.5,15*\rr);
\draw [line width=0.2mm]  (9.5,19*\rr)--(10.5,19*\rr)--(11,18*\rr)--(10.5,17*\rr);
\draw [line width=0.2mm]  (9.5,21*\rr)--(10.5,21*\rr)--(11,20*\rr)--(10.5,19*\rr);

\draw [line width=0.2mm]  (11,0)--(12,0)--(12.5,-1*\rr);
\draw [line width=0.2mm]  (11,2*\rr)--(12,2*\rr)--(12.5,1*\rr)--(12,0);
\draw [line width=0.2mm]  (11,4*\rr)--(12,4*\rr)--(12.5,3*\rr)--(12,2*\rr);
\draw [line width=0.2mm]  (11,6*\rr)--(12,6*\rr)--(12.5,5*\rr)--(12,4*\rr);
\draw [line width=0.2mm]  (11,8*\rr)--(12,8*\rr)--(12.5,7*\rr)--(12,6*\rr);
\draw [line width=0.2mm]  (11,10*\rr)--(12,10*\rr)--(12.5,9*\rr)--(12,8*\rr);
\draw [line width=0.2mm]  (11,12*\rr)--(12,12*\rr)--(12.5,11*\rr)--(12,10*\rr);
\draw [line width=0.2mm]  (11,14*\rr)--(12,14*\rr)--(12.5,13*\rr)--(12,12*\rr);
\draw [line width=0.2mm]  (11,16*\rr)--(12,16*\rr)--(12.5,15*\rr)--(12,14*\rr);
\draw [line width=0.2mm]  (11,18*\rr)--(12,18*\rr)--(12.5,17*\rr)--(12,16*\rr);
\draw [line width=0.2mm]  (11,20*\rr)--(12,20*\rr)--(12.5,19*\rr)--(12,18*\rr);
\draw [line width=0.2mm]  (12,20*\rr)--(12.5,21*\rr);

\draw [line width=0.2mm]  (12.5,1*\rr)--(13.5,1*\rr)--(14,0)--(13.5,-1*\rr)
--(12.5,-1*\rr);
\draw [line width=0.2mm]  (12.5,3*\rr)--(13.5,3*\rr)--(14,2*\rr)--(13.5,1*\rr);
\draw [line width=0.2mm]  (12.5,5*\rr)--(13.5,5*\rr)--(14,4*\rr)--(13.5,3*\rr);
\draw [line width=0.2mm]  (12.5,7*\rr)--(13.5,7*\rr)--(14,6*\rr)--(13.5,5*\rr);
\draw [line width=0.2mm]  (12.5,9*\rr)--(13.5,9*\rr)--(14,8*\rr)--(13.5,7*\rr);
\draw [line width=0.2mm]  (12.5,11*\rr)--(13.5,11*\rr)--(14,10*\rr)--(13.5,9*\rr);
\draw [line width=0.2mm]  (12.5,13*\rr)--(13.5,13*\rr)--(14,12*\rr)--(13.5,11*\rr);
\draw [line width=0.2mm]  (12.5,15*\rr)--(13.5,15*\rr)--(14,14*\rr)--(13.5,13*\rr);
\draw [line width=0.2mm]  (12.5,17*\rr)--(13.5,17*\rr)--(14,16*\rr)--(13.5,15*\rr);
\draw [line width=0.2mm]  (12.5,19*\rr)--(13.5,19*\rr)--(14,18*\rr)--(13.5,17*\rr);
\draw [line width=0.2mm]  (12.5,21*\rr)--(13.5,21*\rr)--(14,20*\rr)--(13.5,19*\rr);

\draw [line width=0.2mm]  (14,0)--(15,0)--(15.5,-1*\rr);
\draw [line width=0.2mm]  (14,2*\rr)--(15,2*\rr)--(15.5,1*\rr)--(15,0);
\draw [line width=0.2mm]  (14,4*\rr)--(15,4*\rr)--(15.5,3*\rr)--(15,2*\rr);
\draw [line width=0.2mm]  (14,6*\rr)--(15,6*\rr)--(15.5,5*\rr)--(15,4*\rr);
\draw [line width=0.2mm]  (14,8*\rr)--(15,8*\rr)--(15.5,7*\rr)--(15,6*\rr);
\draw [line width=0.2mm]  (14,10*\rr)--(15,10*\rr)--(15.5,9*\rr)--(15,8*\rr);
\draw [line width=0.2mm]  (14,12*\rr)--(15,12*\rr)--(15.5,11*\rr)--(15,10*\rr);
\draw [line width=0.2mm]  (14,14*\rr)--(15,14*\rr)--(15.5,13*\rr)--(15,12*\rr);
\draw [line width=0.2mm]  (14,16*\rr)--(15,16*\rr)--(15.5,15*\rr)--(15,14*\rr);
\draw [line width=0.2mm]  (14,18*\rr)--(15,18*\rr)--(15.5,17*\rr)--(15,16*\rr);
\draw [line width=0.2mm]  (14,20*\rr)--(15,20*\rr)--(15.5,19*\rr)--(15,18*\rr);
\draw [line width=0.2mm]  (15,20*\rr)--(15.5,21*\rr);

\draw [line width=0.2mm]  (15.5,1*\rr)--(16.5,1*\rr)--(17,0)--(16.5,-1*\rr)
--(15.5,-1*\rr);
\draw [line width=0.2mm]  (15.5,3*\rr)--(16.5,3*\rr)--(17,2*\rr)--(16.5,1*\rr);
\draw [line width=0.2mm]  (15.5,5*\rr)--(16.5,5*\rr)--(17,4*\rr)--(16.5,3*\rr);
\draw [line width=0.2mm]  (15.5,7*\rr)--(16.5,7*\rr)--(17,6*\rr)--(16.5,5*\rr);
\draw [line width=0.2mm]  (15.5,9*\rr)--(16.5,9*\rr)--(17,8*\rr)--(16.5,7*\rr);
\draw [line width=0.2mm]  (15.5,11*\rr)--(16.5,11*\rr)--(17,10*\rr)--(16.5,9*\rr);
\draw [line width=0.2mm]  (15.5,13*\rr)--(16.5,13*\rr)--(17,12*\rr)--(16.5,11*\rr);
\draw [line width=0.2mm]  (15.5,15*\rr)--(16.5,15*\rr)--(17,14*\rr)--(16.5,13*\rr);
\draw [line width=0.2mm]  (15.5,17*\rr)--(16.5,17*\rr)--(17,16*\rr)--(16.5,15*\rr);
\draw [line width=0.2mm]  (15.5,19*\rr)--(16.5,19*\rr)--(17,18*\rr)--(16.5,17*\rr);
\draw [line width=0.2mm]  (15.5,21*\rr)--(16.5,21*\rr)--(17,20*\rr)--(16.5,19*\rr);

\draw [line width=0.2mm]  (17,0)--(18,0)--(18.5,-1*\rr);
\draw [line width=0.2mm]  (17,2*\rr)--(18,2*\rr)--(18.5,1*\rr)--(18,0);
\draw [line width=0.2mm]  (17,4*\rr)--(18,4*\rr)--(18.5,3*\rr)--(18,2*\rr);
\draw [line width=0.2mm]  (17,6*\rr)--(18,6*\rr)--(18.5,5*\rr)--(18,4*\rr);
\draw [line width=0.2mm]  (17,8*\rr)--(18,8*\rr)--(18.5,7*\rr)--(18,6*\rr);
\draw [line width=0.2mm]  (17,10*\rr)--(18,10*\rr)--(18.5,9*\rr)--(18,8*\rr);
\draw [line width=0.2mm]  (17,12*\rr)--(18,12*\rr)--(18.5,11*\rr)--(18,10*\rr);
\draw [line width=0.2mm]  (17,14*\rr)--(18,14*\rr)--(18.5,13*\rr)--(18,12*\rr);
\draw [line width=0.2mm]  (17,16*\rr)--(18,16*\rr)--(18.5,15*\rr)--(18,14*\rr);
\draw [line width=0.2mm]  (17,18*\rr)--(18,18*\rr)--(18.5,17*\rr)--(18,16*\rr);
\draw [line width=0.2mm]  (17,20*\rr)--(18,20*\rr)--(18.5,19*\rr)--(18,18*\rr);
\draw [line width=0.2mm]  (18,20*\rr)--(18.5,21*\rr);

\draw [line width=0.2mm]  (18.5,1*\rr)--(19.5,1*\rr)--(20,0)--(19.5,-1*\rr)
--(18.5,-1*\rr);
\draw [line width=0.2mm]  (18.5,3*\rr)--(19.5,3*\rr)--(20,2*\rr)--(19.5,1*\rr);
\draw [line width=0.2mm]  (18.5,5*\rr)--(19.5,5*\rr)--(20,4*\rr)--(19.5,3*\rr);
\draw [line width=0.2mm]  (18.5,7*\rr)--(19.5,7*\rr)--(20,6*\rr)--(19.5,5*\rr);
\draw [line width=0.2mm]  (18.5,9*\rr)--(19.5,9*\rr)--(20,8*\rr)--(19.5,7*\rr);
\draw [line width=0.2mm]  (18.5,11*\rr)--(19.5,11*\rr)--(20,10*\rr)--(19.5,9*\rr);
\draw [line width=0.2mm]  (18.5,13*\rr)--(19.5,13*\rr)--(20,12*\rr)--(19.5,11*\rr);
\draw [line width=0.2mm]  (18.5,15*\rr)--(19.5,15*\rr)--(20,14*\rr)--(19.5,13*\rr);
\draw [line width=0.2mm]  (18.5,17*\rr)--(19.5,17*\rr)--(20,16*\rr)--(19.5,15*\rr);
\draw [line width=0.2mm]  (18.5,19*\rr)--(19.5,19*\rr)--(20,18*\rr)--(19.5,17*\rr);
\draw [line width=0.2mm]  (18.5,21*\rr)--(19.5,21*\rr)--(20,20*\rr)--(19.5,19*\rr);

\draw [line width=0.2mm]  (20,0)--(21,0)--(21.5,-1*\rr);
\draw [line width=0.2mm]  (20,2*\rr)--(21,2*\rr)--(21.5,1*\rr)--(21,0);
\draw [line width=0.2mm]  (20,4*\rr)--(21,4*\rr)--(21.5,3*\rr)--(21,2*\rr);
\draw [line width=0.2mm]  (20,6*\rr)--(21,6*\rr)--(21.5,5*\rr)--(21,4*\rr);
\draw [line width=0.2mm]  (20,8*\rr)--(21,8*\rr)--(21.5,7*\rr)--(21,6*\rr);
\draw [line width=0.2mm]  (20,10*\rr)--(21,10*\rr)--(21.5,9*\rr)--(21,8*\rr);
\draw [line width=0.2mm]  (20,12*\rr)--(21,12*\rr)--(21.5,11*\rr)--(21,10*\rr);
\draw [line width=0.2mm]  (20,14*\rr)--(21,14*\rr)--(21.5,13*\rr)--(21,12*\rr);
\draw [line width=0.2mm]  (20,16*\rr)--(21,16*\rr)--(21.5,15*\rr)--(21,14*\rr);
\draw [line width=0.2mm]  (20,18*\rr)--(21,18*\rr)--(21.5,17*\rr)--(21,16*\rr);
\draw [line width=0.2mm]  (20,20*\rr)--(21,20*\rr)--(21.5,19*\rr)--(21,18*\rr);
\draw [line width=0.2mm]  (21,20*\rr)--(21.5,21*\rr);

\draw [line width=0.2mm]  (21.5,1*\rr)--(22.5,1*\rr)--(23,0)--(22.5,-1*\rr)
--(21.5,-1*\rr);
\draw [line width=0.2mm]  (21.5,3*\rr)--(22.5,3*\rr)--(23,2*\rr)--(22.5,1*\rr);
\draw [line width=0.2mm]  (21.5,5*\rr)--(22.5,5*\rr)--(23,4*\rr)--(22.5,3*\rr);
\draw [line width=0.2mm]  (21.5,7*\rr)--(22.5,7*\rr)--(23,6*\rr)--(22.5,5*\rr);
\draw [line width=0.2mm]  (21.5,9*\rr)--(22.5,9*\rr)--(23,8*\rr)--(22.5,7*\rr);
\draw [line width=0.2mm]  (21.5,11*\rr)--(22.5,11*\rr)--(23,10*\rr)--(22.5,9*\rr);
\draw [line width=0.2mm]  (21.5,13*\rr)--(22.5,13*\rr)--(23,12*\rr)--(22.5,11*\rr);
\draw [line width=0.2mm]  (21.5,15*\rr)--(22.5,15*\rr)--(23,14*\rr)--(22.5,13*\rr);
\draw [line width=0.2mm]  (21.5,17*\rr)--(22.5,17*\rr)--(23,16*\rr)--(22.5,15*\rr);
\draw [line width=0.2mm]  (21.5,19*\rr)--(22.5,19*\rr)--(23,18*\rr)--(22.5,17*\rr);
\draw [line width=0.2mm]  (21.5,21*\rr)--(22.5,21*\rr)--(23,20*\rr)--(22.5,19*\rr);

\draw [line width=0.2mm]  (23,0)--(24,0)--(24.5,-1*\rr);
\draw [line width=0.2mm]  (23,2*\rr)--(24,2*\rr)--(24.5,1*\rr)--(24,0);
\draw [line width=0.2mm]  (23,4*\rr)--(24,4*\rr)--(24.5,3*\rr)--(24,2*\rr);
\draw [line width=0.2mm]  (23,6*\rr)--(24,6*\rr)--(24.5,5*\rr)--(24,4*\rr);
\draw [line width=0.2mm]  (23,8*\rr)--(24,8*\rr)--(24.5,7*\rr)--(24,6*\rr);
\draw [line width=0.2mm]  (23,10*\rr)--(24,10*\rr)--(24.5,9*\rr)--(24,8*\rr);
\draw [line width=0.2mm]  (23,12*\rr)--(24,12*\rr)--(24.5,11*\rr)--(24,10*\rr);
\draw [line width=0.2mm]  (23,14*\rr)--(24,14*\rr)--(24.5,13*\rr)--(24,12*\rr);
\draw [line width=0.2mm]  (23,16*\rr)--(24,16*\rr)--(24.5,15*\rr)--(24,14*\rr);
\draw [line width=0.2mm]  (23,18*\rr)--(24,18*\rr)--(24.5,17*\rr)--(24,16*\rr);
\draw [line width=0.2mm]  (23,20*\rr)--(24,20*\rr)--(24.5,19*\rr)--(24,18*\rr);
\draw [line width=0.2mm]  (24,20*\rr)--(24.5,21*\rr);

\draw [line width=0.6mm] (0,0)--(4.5,9*\rr)--(9,0)--(0,0);
\draw [line width=0.6mm] (4.5,9*\rr)--(13.5,9*\rr)--(9,0);
\draw [line width=0.6mm] (13.5,9*\rr)--(18,0)--(9,0);

\foreach \pos in {(0,0), (4.5,9*\rr),(13.5,9*\rr),(9,0),(18,0)}
\shade[shading=ball, ball color=black] \pos circle (.4);

%\draw [line width=0.3mm] (4.5,9*\rr)--(18,0);
%\draw [line width=0.3mm] (0,0)--(13.5,9*\rr);

\draw [line width=0.3mm] (13.5,3*\rr)--(9,5);
\draw [line width=0.3mm] (9,5)--(9,12*\rr);
\draw [line width=0.3mm] (9,12*\rr)--(10.5,13*\rr);
\draw [line width=0.3mm] (16.5,13*\rr)--(18,12*\rr);
\draw [line width=0.3mm] (18,12*\rr)--(18,6*\rr);
\draw [line width=0.3mm] (18,6*\rr)--(13.5,3*\rr);

%\foreach \pos in {(0.5,11*\rr),(3,2*\rr),(10.5,3*\rr),(12,12*\rr),(18.5,5*\rr),(24.5,9*\rr)}
%\shade[shading=ball, ball color=white] \pos circle (.4);
 
\end{tikzpicture}\vskip .5cm

(c) \begin{tikzpicture}[scale=0.25]
\clip (-1, -1.6*\rr) rectangle (20, 13.5*\rr);

%\path [fill=lightgray] (13.5,3*\rr)--(9,5)--(9,12*\rr)--(10.5,13*\rr)
%--(16.5,13*\rr)--(18,12*\rr)--(18,6*\rr)--(13.5,3*\rr);

\path [fill=lightgray] (8.5,5*\rr)--(11,8*\rr)--(14.5,7*\rr)--(15.5,3*\rr)
--(13,0*\rr)--(9.5,1*\rr)-- (8.5,5*\rr);

\draw [line width=0.2mm] (0.5,1*\rr) -- (1.5,1*\rr) -- (2,0) -- (1.5,-1*\rr)
-- (0.5,-1*\rr) -- (0,0) -- (0.5,1*\rr);
\draw [line width=0.2mm] (0.5,1*\rr) -- (0,2*\rr) -- (0.5,3*\rr) -- (1.5,3*\rr) 
-- (2,2*\rr) -- (1.5,1*\rr);
\draw [line width=0.2mm] (0.5,3*\rr) -- (0,4*\rr) -- (0.5,5*\rr) -- (1.5,5*\rr) 
-- (2,4*\rr) -- (1.5,3*\rr);
\draw [line width=0.2mm] (0.5,5*\rr) -- (0,6*\rr) -- (0.5,7*\rr) -- (1.5,7*\rr) 
-- (2,6*\rr) -- (1.5,5*\rr);
\draw [line width=0.2mm] (0.5,7*\rr) -- (0,8*\rr) -- (0.5,9*\rr) -- (1.5,9*\rr) 
-- (2,8*\rr) -- (1.5,7*\rr);
\draw [line width=0.2mm] (0.5,9*\rr) -- (0,10*\rr) -- (0.5,11*\rr) -- (1.5,11*\rr) 
-- (2,10*\rr) -- (1.5,9*\rr);
\draw [line width=0.2mm] (0.5,11*\rr) -- (0,12*\rr) -- (0.5,13*\rr) -- (1.5,13*\rr) 
-- (2,12*\rr) -- (1.5,11*\rr);
\draw [line width=0.2mm] (0.5,13*\rr) -- (0,14*\rr) -- (0.5,15*\rr) -- (1.5,15*\rr) 
-- (2,14*\rr) -- (1.5,13*\rr);
\draw [line width=0.2mm] (0.5,15*\rr) -- (0,16*\rr) -- (0.5,17*\rr) -- (1.5,17*\rr) 
-- (2,16*\rr) -- (1.5,15*\rr);
\draw [line width=0.2mm] (0.5,17*\rr) -- (0,18*\rr) -- (0.5,19*\rr) -- (1.5,19*\rr) 
-- (2,18*\rr) -- (1.5,17*\rr);
\draw [line width=0.2mm] (0.5,19*\rr) -- (0,20*\rr) -- (0.5,21*\rr) -- (1.5,21*\rr) 
-- (2,20*\rr) -- (1.5,19*\rr);
%\draw [line width=0.2mm] (0.5,21*\rr) -- (0,22*\rr) -- (0.5,22*\rr) -- (1.5,22*\rr) 
%-- (2,22*\rr) -- (1.5,21*\rr);

\draw [line width=0.2mm]  (3,0) -- (3.5,-1*\rr);
\draw [line width=0.2mm] (2,2*\rr) -- (3,2*\rr) -- (3.5,1*\rr) -- (3,0) -- (2,0);
\draw [line width=0.2mm] (2,4*\rr) -- (3,4*\rr) -- (3.5,3*\rr) -- (3,2*\rr);
\draw [line width=0.2mm] (2,6*\rr) -- (3,6*\rr) -- (3.5,5*\rr) -- (3,4*\rr);
\draw [line width=0.2mm] (2,8*\rr) -- (3,8*\rr) -- (3.5,7*\rr) -- (3,6*\rr);
\draw [line width=0.2mm] (2,10*\rr) -- (3,10*\rr) -- (3.5,9*\rr) -- (3,8*\rr);
\draw [line width=0.2mm] (2,12*\rr) -- (3,12*\rr) -- (3.5,11*\rr) -- (3,10*\rr);
\draw [line width=0.2mm] (2,14*\rr) -- (3,14*\rr) -- (3.5,13*\rr) -- (3,12*\rr);
\draw [line width=0.2mm] (2,16*\rr) -- (3,16*\rr) -- (3.5,15*\rr) -- (3,14*\rr);
\draw [line width=0.2mm] (2,18*\rr) -- (3,18*\rr) -- (3.5,17*\rr) -- (3,16*\rr);
\draw [line width=0.2mm] (2,20*\rr) -- (3,20*\rr) -- (3.5,19*\rr) -- (3,18*\rr);
\draw [line width=0.2mm] (3.5,21*\rr) -- (3,20*\rr);

\draw [line width=0.2mm]  (3.5,1*\rr) -- (4.5,1*\rr) -- (5,0) -- (4.5,-1*\rr)
--(3.5,-1*\rr);
\draw [line width=0.2mm]  (3.5,3*\rr) -- (4.5,3*\rr) -- (5,2*\rr) -- (4.5,1*\rr);
\draw [line width=0.2mm]  (3.5,5*\rr) -- (4.5,5*\rr) -- (5,4*\rr) -- (4.5,3*\rr);
\draw [line width=0.2mm]  (3.5,7*\rr) -- (4.5,7*\rr) -- (5,6*\rr) -- (4.5,5*\rr);
\draw [line width=0.2mm]  (3.5,9*\rr) -- (4.5,9*\rr) -- (5,8*\rr) -- (4.5,7*\rr);
\draw [line width=0.2mm]  (3.5,11*\rr) -- (4.5,11*\rr) -- (5,10*\rr) -- (4.5,9*\rr);
\draw [line width=0.2mm]  (3.5,13*\rr) -- (4.5,13*\rr) -- (5,12*\rr) -- (4.5,11*\rr);
\draw [line width=0.2mm]  (3.5,15*\rr) -- (4.5,15*\rr) -- (5,14*\rr) -- (4.5,13*\rr);
\draw [line width=0.2mm]  (3.5,17*\rr) -- (4.5,17*\rr) -- (5,16*\rr) -- (4.5,15*\rr);
\draw [line width=0.2mm]  (3.5,19*\rr) -- (4.5,19*\rr) -- (5,18*\rr) -- (4.5,17*\rr);
\draw [line width=0.2mm]  (3.5,21*\rr) -- (4.5,21*\rr) -- (5,20*\rr) -- (4.5,19*\rr);

\draw [line width=0.2mm]  (5,0) -- (6,0)--(6.5,-1*\rr );
\draw [line width=0.2mm]  (5,2*\rr) -- (6,2*\rr)--(6.5,1*\rr )--(6,0);
\draw [line width=0.2mm]  (5,4*\rr) -- (6,4*\rr)--(6.5,3*\rr )--(6,2*\rr);
\draw [line width=0.2mm]  (5,6*\rr) -- (6,6*\rr)--(6.5,5*\rr )--(6,4*\rr);
\draw [line width=0.2mm]  (5,8*\rr) -- (6,8*\rr)--(6.5,7*\rr )--(6,6*\rr);
\draw [line width=0.2mm]  (5,10*\rr) -- (6,10*\rr)--(6.5,9*\rr )--(6,8*\rr);
\draw [line width=0.2mm]  (5,12*\rr) -- (6,12*\rr)--(6.5,11*\rr )--(6,10*\rr);
\draw [line width=0.2mm]  (5,14*\rr) -- (6,14*\rr)--(6.5,13*\rr )--(6,12*\rr);
\draw [line width=0.2mm]  (5,16*\rr) -- (6,16*\rr)--(6.5,15*\rr )--(6,14*\rr);
\draw [line width=0.2mm]  (5,18*\rr) -- (6,18*\rr)--(6.5,17*\rr )--(6,16*\rr);
\draw [line width=0.2mm]  (5,20*\rr) -- (6,20*\rr)--(6.5,19*\rr )--(6,18*\rr);
\draw [line width=0.2mm]  (6.5,21*\rr )--(6,20*\rr);

\draw [line width=0.2mm]  (6.5,1*\rr) -- (7.5,1*\rr)--(8,0)--(7.5,-1*\rr)
--(6.5,-1*\rr);
\draw [line width=0.2mm]  (6.5,3*\rr) -- (7.5,3*\rr)--(8,2*\rr)--(7.5,1*\rr);
\draw [line width=0.2mm]  (6.5,5*\rr) -- (7.5,5*\rr)--(8,4*\rr)--(7.5,3*\rr);
\draw [line width=0.2mm]  (6.5,7*\rr) -- (7.5,7*\rr)--(8,6*\rr)--(7.5,5*\rr);
\draw [line width=0.2mm]  (6.5,9*\rr) -- (7.5,9*\rr)--(8,8*\rr)--(7.5,7*\rr);
\draw [line width=0.2mm]  (6.5,11*\rr) -- (7.5,11*\rr)--(8,10*\rr)--(7.5,9*\rr);
\draw [line width=0.2mm]  (6.5,13*\rr) -- (7.5,13*\rr)--(8,12*\rr)--(7.5,11*\rr);
\draw [line width=0.2mm]  (6.5,15*\rr) -- (7.5,15*\rr)--(8,14*\rr)--(7.5,13*\rr);
\draw [line width=0.2mm]  (6.5,17*\rr) -- (7.5,17*\rr)--(8,16*\rr)--(7.5,15*\rr);
\draw [line width=0.2mm]  (6.5,19*\rr) -- (7.5,19*\rr)--(8,18*\rr)--(7.5,17*\rr);
\draw [line width=0.2mm]  (6.5,21*\rr) -- (7.5,21*\rr)--(8,20*\rr)--(7.5,19*\rr);

\draw [line width=0.2mm]  (8,0)--(9,0)--(9.5,-1*\rr);
\draw [line width=0.2mm]  (8,2*\rr)--(9,2*\rr)--(9.5,1*\rr)--(9,0);
\draw [line width=0.2mm]  (8,4*\rr)--(9,4*\rr)--(9.5,3*\rr)--(9,2*\rr);
\draw [line width=0.2mm]  (8,6*\rr)--(9,6*\rr)--(9.5,5*\rr)--(9,4*\rr);
\draw [line width=0.2mm]  (8,8*\rr)--(9,8*\rr)--(9.5,7*\rr)--(9,6*\rr);
\draw [line width=0.2mm]  (8,10*\rr)--(9,10*\rr)--(9.5,9*\rr)--(9,8*\rr);
\draw [line width=0.2mm]  (8,12*\rr)--(9,12*\rr)--(9.5,11*\rr)--(9,10*\rr);
\draw [line width=0.2mm]  (8,14*\rr)--(9,14*\rr)--(9.5,13*\rr)--(9,12*\rr);
\draw [line width=0.2mm]  (8,16*\rr)--(9,16*\rr)--(9.5,15*\rr)--(9,14*\rr);
\draw [line width=0.2mm]  (8,18*\rr)--(9,18*\rr)--(9.5,17*\rr)--(9,16*\rr);
\draw [line width=0.2mm]  (8,20*\rr)--(9,20*\rr)--(9.5,19*\rr)--(9,18*\rr);
\draw [line width=0.2mm]  (9.5,21*\rr)--(9,20*\rr);

\draw [line width=0.2mm]  (9.5,1*\rr)--(10.5,1*\rr)--(11,0)--(10.5,-1*\rr)
--(9.5,-1*\rr);
\draw [line width=0.2mm]  (9.5,3*\rr)--(10.5,3*\rr)--(11,2*\rr)--(10.5,1*\rr);
\draw [line width=0.2mm]  (9.5,5*\rr)--(10.5,5*\rr)--(11,4*\rr)--(10.5,3*\rr);
\draw [line width=0.2mm]  (9.5,7*\rr)--(10.5,7*\rr)--(11,6*\rr)--(10.5,5*\rr);
\draw [line width=0.2mm]  (9.5,9*\rr)--(10.5,9*\rr)--(11,8*\rr)--(10.5,7*\rr);
\draw [line width=0.2mm]  (9.5,11*\rr)--(10.5,11*\rr)--(11,10*\rr)--(10.5,9*\rr);
\draw [line width=0.2mm]  (9.5,13*\rr)--(10.5,13*\rr)--(11,12*\rr)--(10.5,11*\rr);
\draw [line width=0.2mm]  (9.5,15*\rr)--(10.5,15*\rr)--(11,14*\rr)--(10.5,13*\rr);
\draw [line width=0.2mm]  (9.5,17*\rr)--(10.5,17*\rr)--(11,16*\rr)--(10.5,15*\rr);
\draw [line width=0.2mm]  (9.5,19*\rr)--(10.5,19*\rr)--(11,18*\rr)--(10.5,17*\rr);
\draw [line width=0.2mm]  (9.5,21*\rr)--(10.5,21*\rr)--(11,20*\rr)--(10.5,19*\rr);

\draw [line width=0.2mm]  (11,0)--(12,0)--(12.5,-1*\rr);
\draw [line width=0.2mm]  (11,2*\rr)--(12,2*\rr)--(12.5,1*\rr)--(12,0);
\draw [line width=0.2mm]  (11,4*\rr)--(12,4*\rr)--(12.5,3*\rr)--(12,2*\rr);
\draw [line width=0.2mm]  (11,6*\rr)--(12,6*\rr)--(12.5,5*\rr)--(12,4*\rr);
\draw [line width=0.2mm]  (11,8*\rr)--(12,8*\rr)--(12.5,7*\rr)--(12,6*\rr);
\draw [line width=0.2mm]  (11,10*\rr)--(12,10*\rr)--(12.5,9*\rr)--(12,8*\rr);
\draw [line width=0.2mm]  (11,12*\rr)--(12,12*\rr)--(12.5,11*\rr)--(12,10*\rr);
\draw [line width=0.2mm]  (11,14*\rr)--(12,14*\rr)--(12.5,13*\rr)--(12,12*\rr);
\draw [line width=0.2mm]  (11,16*\rr)--(12,16*\rr)--(12.5,15*\rr)--(12,14*\rr);
\draw [line width=0.2mm]  (11,18*\rr)--(12,18*\rr)--(12.5,17*\rr)--(12,16*\rr);
\draw [line width=0.2mm]  (11,20*\rr)--(12,20*\rr)--(12.5,19*\rr)--(12,18*\rr);
\draw [line width=0.2mm]  (12,20*\rr)--(12.5,21*\rr);

\draw [line width=0.2mm]  (12.5,1*\rr)--(13.5,1*\rr)--(14,0)--(13.5,-1*\rr)
--(12.5,-1*\rr);
\draw [line width=0.2mm]  (12.5,3*\rr)--(13.5,3*\rr)--(14,2*\rr)--(13.5,1*\rr);
\draw [line width=0.2mm]  (12.5,5*\rr)--(13.5,5*\rr)--(14,4*\rr)--(13.5,3*\rr);
\draw [line width=0.2mm]  (12.5,7*\rr)--(13.5,7*\rr)--(14,6*\rr)--(13.5,5*\rr);
\draw [line width=0.2mm]  (12.5,9*\rr)--(13.5,9*\rr)--(14,8*\rr)--(13.5,7*\rr);
\draw [line width=0.2mm]  (12.5,11*\rr)--(13.5,11*\rr)--(14,10*\rr)--(13.5,9*\rr);
\draw [line width=0.2mm]  (12.5,13*\rr)--(13.5,13*\rr)--(14,12*\rr)--(13.5,11*\rr);
\draw [line width=0.2mm]  (12.5,15*\rr)--(13.5,15*\rr)--(14,14*\rr)--(13.5,13*\rr);
\draw [line width=0.2mm]  (12.5,17*\rr)--(13.5,17*\rr)--(14,16*\rr)--(13.5,15*\rr);
\draw [line width=0.2mm]  (12.5,19*\rr)--(13.5,19*\rr)--(14,18*\rr)--(13.5,17*\rr);
\draw [line width=0.2mm]  (12.5,21*\rr)--(13.5,21*\rr)--(14,20*\rr)--(13.5,19*\rr);

\draw [line width=0.2mm]  (14,0)--(15,0)--(15.5,-1*\rr);
\draw [line width=0.2mm]  (14,2*\rr)--(15,2*\rr)--(15.5,1*\rr)--(15,0);
\draw [line width=0.2mm]  (14,4*\rr)--(15,4*\rr)--(15.5,3*\rr)--(15,2*\rr);
\draw [line width=0.2mm]  (14,6*\rr)--(15,6*\rr)--(15.5,5*\rr)--(15,4*\rr);
\draw [line width=0.2mm]  (14,8*\rr)--(15,8*\rr)--(15.5,7*\rr)--(15,6*\rr);
\draw [line width=0.2mm]  (14,10*\rr)--(15,10*\rr)--(15.5,9*\rr)--(15,8*\rr);
\draw [line width=0.2mm]  (14,12*\rr)--(15,12*\rr)--(15.5,11*\rr)--(15,10*\rr);
\draw [line width=0.2mm]  (14,14*\rr)--(15,14*\rr)--(15.5,13*\rr)--(15,12*\rr);
\draw [line width=0.2mm]  (14,16*\rr)--(15,16*\rr)--(15.5,15*\rr)--(15,14*\rr);
\draw [line width=0.2mm]  (14,18*\rr)--(15,18*\rr)--(15.5,17*\rr)--(15,16*\rr);
\draw [line width=0.2mm]  (14,20*\rr)--(15,20*\rr)--(15.5,19*\rr)--(15,18*\rr);
\draw [line width=0.2mm]  (15,20*\rr)--(15.5,21*\rr);

\draw [line width=0.2mm]  (15.5,1*\rr)--(16.5,1*\rr)--(17,0)--(16.5,-1*\rr)
--(15.5,-1*\rr);
\draw [line width=0.2mm]  (15.5,3*\rr)--(16.5,3*\rr)--(17,2*\rr)--(16.5,1*\rr);
\draw [line width=0.2mm]  (15.5,5*\rr)--(16.5,5*\rr)--(17,4*\rr)--(16.5,3*\rr);
\draw [line width=0.2mm]  (15.5,7*\rr)--(16.5,7*\rr)--(17,6*\rr)--(16.5,5*\rr);
\draw [line width=0.2mm]  (15.5,9*\rr)--(16.5,9*\rr)--(17,8*\rr)--(16.5,7*\rr);
\draw [line width=0.2mm]  (15.5,11*\rr)--(16.5,11*\rr)--(17,10*\rr)--(16.5,9*\rr);
\draw [line width=0.2mm]  (15.5,13*\rr)--(16.5,13*\rr)--(17,12*\rr)--(16.5,11*\rr);
\draw [line width=0.2mm]  (15.5,15*\rr)--(16.5,15*\rr)--(17,14*\rr)--(16.5,13*\rr);
\draw [line width=0.2mm]  (15.5,17*\rr)--(16.5,17*\rr)--(17,16*\rr)--(16.5,15*\rr);
\draw [line width=0.2mm]  (15.5,19*\rr)--(16.5,19*\rr)--(17,18*\rr)--(16.5,17*\rr);
\draw [line width=0.2mm]  (15.5,21*\rr)--(16.5,21*\rr)--(17,20*\rr)--(16.5,19*\rr);

\draw [line width=0.2mm]  (17,0)--(18,0)--(18.5,-1*\rr);
\draw [line width=0.2mm]  (17,2*\rr)--(18,2*\rr)--(18.5,1*\rr)--(18,0);
\draw [line width=0.2mm]  (17,4*\rr)--(18,4*\rr)--(18.5,3*\rr)--(18,2*\rr);
\draw [line width=0.2mm]  (17,6*\rr)--(18,6*\rr)--(18.5,5*\rr)--(18,4*\rr);
\draw [line width=0.2mm]  (17,8*\rr)--(18,8*\rr)--(18.5,7*\rr)--(18,6*\rr);
\draw [line width=0.2mm]  (17,10*\rr)--(18,10*\rr)--(18.5,9*\rr)--(18,8*\rr);
\draw [line width=0.2mm]  (17,12*\rr)--(18,12*\rr)--(18.5,11*\rr)--(18,10*\rr);
\draw [line width=0.2mm]  (17,14*\rr)--(18,14*\rr)--(18.5,13*\rr)--(18,12*\rr);
\draw [line width=0.2mm]  (17,16*\rr)--(18,16*\rr)--(18.5,15*\rr)--(18,14*\rr);
\draw [line width=0.2mm]  (17,18*\rr)--(18,18*\rr)--(18.5,17*\rr)--(18,16*\rr);
\draw [line width=0.2mm]  (17,20*\rr)--(18,20*\rr)--(18.5,19*\rr)--(18,18*\rr);
\draw [line width=0.2mm]  (18,20*\rr)--(18.5,21*\rr);

\draw [line width=0.2mm]  (18.5,1*\rr)--(19.5,1*\rr)--(20,0)--(19.5,-1*\rr)
--(18.5,-1*\rr);
\draw [line width=0.2mm]  (18.5,3*\rr)--(19.5,3*\rr)--(20,2*\rr)--(19.5,1*\rr);
\draw [line width=0.2mm]  (18.5,5*\rr)--(19.5,5*\rr)--(20,4*\rr)--(19.5,3*\rr);
\draw [line width=0.2mm]  (18.5,7*\rr)--(19.5,7*\rr)--(20,6*\rr)--(19.5,5*\rr);
\draw [line width=0.2mm]  (18.5,9*\rr)--(19.5,9*\rr)--(20,8*\rr)--(19.5,7*\rr);
\draw [line width=0.2mm]  (18.5,11*\rr)--(19.5,11*\rr)--(20,10*\rr)--(19.5,9*\rr);
\draw [line width=0.2mm]  (18.5,13*\rr)--(19.5,13*\rr)--(20,12*\rr)--(19.5,11*\rr);
\draw [line width=0.2mm]  (18.5,15*\rr)--(19.5,15*\rr)--(20,14*\rr)--(19.5,13*\rr);
\draw [line width=0.2mm]  (18.5,17*\rr)--(19.5,17*\rr)--(20,16*\rr)--(19.5,15*\rr);
\draw [line width=0.2mm]  (18.5,19*\rr)--(19.5,19*\rr)--(20,18*\rr)--(19.5,17*\rr);
\draw [line width=0.2mm]  (18.5,21*\rr)--(19.5,21*\rr)--(20,20*\rr)--(19.5,19*\rr);

\draw [line width=0.2mm]  (20,0)--(21,0)--(21.5,-1*\rr);
\draw [line width=0.2mm]  (20,2*\rr)--(21,2*\rr)--(21.5,1*\rr)--(21,0);
\draw [line width=0.2mm]  (20,4*\rr)--(21,4*\rr)--(21.5,3*\rr)--(21,2*\rr);
\draw [line width=0.2mm]  (20,6*\rr)--(21,6*\rr)--(21.5,5*\rr)--(21,4*\rr);
\draw [line width=0.2mm]  (20,8*\rr)--(21,8*\rr)--(21.5,7*\rr)--(21,6*\rr);
\draw [line width=0.2mm]  (20,10*\rr)--(21,10*\rr)--(21.5,9*\rr)--(21,8*\rr);
\draw [line width=0.2mm]  (20,12*\rr)--(21,12*\rr)--(21.5,11*\rr)--(21,10*\rr);
\draw [line width=0.2mm]  (20,14*\rr)--(21,14*\rr)--(21.5,13*\rr)--(21,12*\rr);
\draw [line width=0.2mm]  (20,16*\rr)--(21,16*\rr)--(21.5,15*\rr)--(21,14*\rr);
\draw [line width=0.2mm]  (20,18*\rr)--(21,18*\rr)--(21.5,17*\rr)--(21,16*\rr);
\draw [line width=0.2mm]  (20,20*\rr)--(21,20*\rr)--(21.5,19*\rr)--(21,18*\rr);
\draw [line width=0.2mm]  (21,20*\rr)--(21.5,21*\rr);

\draw [line width=0.2mm]  (21.5,1*\rr)--(22.5,1*\rr)--(23,0)--(22.5,-1*\rr)
--(21.5,-1*\rr);
\draw [line width=0.2mm]  (21.5,3*\rr)--(22.5,3*\rr)--(23,2*\rr)--(22.5,1*\rr);
\draw [line width=0.2mm]  (21.5,5*\rr)--(22.5,5*\rr)--(23,4*\rr)--(22.5,3*\rr);
\draw [line width=0.2mm]  (21.5,7*\rr)--(22.5,7*\rr)--(23,6*\rr)--(22.5,5*\rr);
\draw [line width=0.2mm]  (21.5,9*\rr)--(22.5,9*\rr)--(23,8*\rr)--(22.5,7*\rr);
\draw [line width=0.2mm]  (21.5,11*\rr)--(22.5,11*\rr)--(23,10*\rr)--(22.5,9*\rr);
\draw [line width=0.2mm]  (21.5,13*\rr)--(22.5,13*\rr)--(23,12*\rr)--(22.5,11*\rr);
\draw [line width=0.2mm]  (21.5,15*\rr)--(22.5,15*\rr)--(23,14*\rr)--(22.5,13*\rr);
\draw [line width=0.2mm]  (21.5,17*\rr)--(22.5,17*\rr)--(23,16*\rr)--(22.5,15*\rr);
\draw [line width=0.2mm]  (21.5,19*\rr)--(22.5,19*\rr)--(23,18*\rr)--(22.5,17*\rr);
\draw [line width=0.2mm]  (21.5,21*\rr)--(22.5,21*\rr)--(23,20*\rr)--(22.5,19*\rr);

\draw [line width=0.2mm]  (23,0)--(24,0)--(24.5,-1*\rr);
\draw [line width=0.2mm]  (23,2*\rr)--(24,2*\rr)--(24.5,1*\rr)--(24,0);
\draw [line width=0.2mm]  (23,4*\rr)--(24,4*\rr)--(24.5,3*\rr)--(24,2*\rr);
\draw [line width=0.2mm]  (23,6*\rr)--(24,6*\rr)--(24.5,5*\rr)--(24,4*\rr);
\draw [line width=0.2mm]  (23,8*\rr)--(24,8*\rr)--(24.5,7*\rr)--(24,6*\rr);
\draw [line width=0.2mm]  (23,10*\rr)--(24,10*\rr)--(24.5,9*\rr)--(24,8*\rr);
\draw [line width=0.2mm]  (23,12*\rr)--(24,12*\rr)--(24.5,11*\rr)--(24,10*\rr);
\draw [line width=0.2mm]  (23,14*\rr)--(24,14*\rr)--(24.5,13*\rr)--(24,12*\rr);
\draw [line width=0.2mm]  (23,16*\rr)--(24,16*\rr)--(24.5,15*\rr)--(24,14*\rr);
\draw [line width=0.2mm]  (23,18*\rr)--(24,18*\rr)--(24.5,17*\rr)--(24,16*\rr);
\draw [line width=0.2mm]  (23,20*\rr)--(24,20*\rr)--(24.5,19*\rr)--(24,18*\rr);
\draw [line width=0.2mm]  (24,20*\rr)--(24.5,21*\rr);

\draw [line width=0.6mm] (6,2*\rr)--(0,0)--(1.5,7*\rr)--(6,2*\rr)--(7.5,9*\rr)
--(1.5,7*\rr);
\draw [line width=0.6mm] (6,2*\rr)--(12,4*\rr)--(7.5,9*\rr)--(13.5,11*\rr)
--(12,4*\rr);
\draw [line width=0.6mm] (13.5,11*\rr)--(18,6*\rr)--(12,4*\rr)
--(16.5,-1*\rr)--(18,6*\rr);
\draw [line width=0.6mm] (13.5,11*\rr)--(19.5,13*\rr)--(18,6*\rr);

\foreach \pos in {(0,0),(1.5,7*\rr),(6,2*\rr),(7.5,9*\rr),(12,4*\rr),(13.5,11*\rr),
(16.5,-1*\rr),(18,6*\rr),(19.5,13*\rr)}
\shade[shading=ball, ball color=black] \pos circle (.4);

%\draw [line width=0.3mm] (7.5,9*\rr)--(18,6*\rr);
%\draw [line width=0.3mm] (13.5,11*\rr)--(16.5,-1*\rr);
%\draw [line width=0.3mm] (13.5,11*\rr)--(6,2*\rr);
%\draw [line width=0.3mm] (7.5,9*\rr)--(10.5,-3*\rr);
%\draw [line width=0.3mm] (16.5,-1*\rr)--(6,2*\rr);
%\draw [line width=0.3mm] (18,6*\rr)--(10.5,-3*\rr);

\draw [line width=0.3mm] (8.5,5*\rr)--(11,8*\rr)--(14.5,7*\rr)--(15.5,3*\rr)
--(13,0*\rr)--(9.5,1*\rr)-- (8.5,5*\rr);

%\foreach \pos in {(0.5,11*\rr),(3,2*\rr),(10.5,3*\rr),(12,12*\rr),(18.5,5*\rr),(24.5,9*\rr)}
%\shade[shading=ball, ball color=white] \pos circle (.4);
 
\end{tikzpicture} 
\caption{$\alpha$-PGSs on ${\mathbb H}_2$,  for $D^2=48,\;81$ (Class HA1, frames (a,b)), 
and  $D^2=39$ (Class HA2, frame (c)), with Voronoi cells (gray hexagons) and choices of 
$D$-rhombuses.\\
{\footnotesize{The number of PGSs and EGMs for $D^2=48,\;81, 39$ equals 32, 
54 and 52, respectively. The PGSs for $D^2=48$ are horizontal, for $D^2=81$ vertical, and
for $D^2=39$ inclined.}}
}
\end{figure}} %FigureD4 Fig4

%%%%%%%%%%%%%%%%%%%%%%%%%%%%%%%%%%%%%%%%%%%%%%%%%%%%%%%%%%%%%%%%%%%%%%%

\def\FigureE5 %%\label{Fig5}
{\begin{figure}\label{Fig5}%\begin{center} 
\centering
\captionsetup{width=0.8\textwidth} 
%\end{center}%\vskip .5cm
\caption{$\alpha$-PGSs on ${\mathbb H}_2$ for $D^2=147$ (Class HB) and $D^2=139$
(Class HC).
{\footnotesize{For $D^2=147$ the vertical PGSs (black rhombuses) are dominant, 
and the number of EGMs equals 98. However, for $D^2=139$ the inclined PGSs (gray
rhombuses) are dominant, and the number of EGMs equals 196.}}
}
\end{figure}}  %FigureE5 Fig5

%%%%%%%%%%%%%%%%%%%%%%%%%%%%%%%%%%%%%%%%%%%%%%%%%%
%%%%%%%%%%%%%%%%%%%%%%%%%%%%%%%%%%%%%%%%%%%%%%%%%%

\def\FigureF6 %%\label{Fig6}
{\begin{figure}[H]\label{Fig6} % \begin{center} 
\centering
\captionsetup{width=0.8\textwidth} 
 %\end{center}%\vskip .5cm

\caption{$\alpha$-PGSs on $\bbH_2$ for a series of values $D^2$ from Class HC, with
${D^*}^2=D^2+2$, and their associated ${D^*}$-triangles.\\
{\footnotesize{(i) For $D^2=19$,  ${D^*}^2=21$ (red). (ii) For $D^2=37$, ${D^*}^2=39$
(turquoise). (iii) For $D^2=61$, ${D^*}^2=63$ (orange). (iv) For $D^2=169$, 
${D^*}^2=171$ (blue). (v) For $D^2=217$, ${D^*}^2=219$ (purple).
The black $\bbH_2$-triangles give the minimal area when the side-lengths 
are at least $D$ and the angles at most $\pi/2$. However, they do not generate PGSs by extension;
they are parts of
$(\beta ,D)$-configurations which are not PGSs; see below. The PGSs are $(\alpha ,D^*)$-configurations
generated from ${D^*}$-triangles of the corresponding color.}}
}
%\centerline{\bf for a straight PGS. The total amount is $25+21=46$}
%\centerline{\bf per a fundamental parallelogram}
\end{figure}} %Fig6 %FigureF6

%%%%%%%%%%%%%%%%%%%%%%%%%%%%%%%%%%%%%%%%%%%%%%%%%%
%%%%%%%%%%%%%%%%%%%%%%%%%%%%%%%%%%%%%%%%%%%%%%%%%%

\def\FigureG7{\begin{figure}\label{Fig7} %\label{D2=13} 
\centering
\captionsetup{width=0.8\textwidth} 
\begin{tikzpicture}[scale=0.28]
\clip (-0.9, -1*\rr) rectangle (32.4, 21*\rr);

%%%\path [fill=lightgray] (4.36,8.3*\rr)--(6.62,9.8*\rr)
%--(7.82,9.57*\rr)--(8.82,6.77*\rr)--(6.7,4.5*\rr)--(5.22,5.02*\rr)--(4.36,8.3*\rr);
%\draw [line width=0.2mm, color=black] (4.36,8.3*\rr)--(6.62,9.8*\rr)
%--(7.82,9.57*\rr)--(8.82,6.77*\rr)--(6.7,4.5*\rr)--(5.22,5.02*\rr)--(4.36,8.3*\rr);

%%%\path [fill=lightgray] (8.82,6.77*\rr)--(10.3,6.25*\rr)
%--(11.25,2.75*\rr)--(8.9,1.27*\rr)--(7.7,1.7*\rr)--(6.7,4.5*\rr)--(8.82,6.75*\rr); 
%\draw [line width=0.2mm, color=black] (8.82,6.77*\rr)--(10.3,6.25*\rr)
%--(11.25,2.75*\rr)--(8.9,1.27*\rr)--(7.7,1.7*\rr)--(6.7,4.5*\rr)--(8.82,6.75*\rr); 

\draw [line width=0.2mm] (0.5,1*\rr) -- (1.5,1*\rr) -- (2,0) -- (1.5,-1*\rr)
-- (0.5,-1*\rr) -- (0,0) -- (0.5,1*\rr);
\draw [line width=0.2mm] (0.5,1*\rr) -- (0,2*\rr) -- (0.5,3*\rr) -- (1.5,3*\rr) 
-- (2,2*\rr) -- (1.5,1*\rr);
\draw [line width=0.2mm] (0.5,3*\rr) -- (0,4*\rr) -- (0.5,5*\rr) -- (1.5,5*\rr) 
-- (2,4*\rr) -- (1.5,3*\rr);
\draw [line width=0.2mm] (0.5,5*\rr) -- (0,6*\rr) -- (0.5,7*\rr) -- (1.5,7*\rr) 
-- (2,6*\rr) -- (1.5,5*\rr);
\draw [line width=0.2mm] (0.5,7*\rr) -- (0,8*\rr) -- (0.5,9*\rr) -- (1.5,9*\rr) 
-- (2,8*\rr) -- (1.5,7*\rr);
\draw [line width=0.2mm] (0.5,9*\rr) -- (0,10*\rr) -- (0.5,11*\rr) -- (1.5,11*\rr) 
-- (2,10*\rr) -- (1.5,9*\rr);
\draw [line width=0.2mm] (0.5,11*\rr) -- (0,12*\rr) -- (0.5,13*\rr) -- (1.5,13*\rr) 
-- (2,12*\rr) -- (1.5,11*\rr);
\draw [line width=0.2mm] (0.5,13*\rr) -- (0,14*\rr) -- (0.5,15*\rr) -- (1.5,15*\rr) 
-- (2,14*\rr) -- (1.5,13*\rr);
\draw [line width=0.2mm] (0.5,15*\rr) -- (0,16*\rr) -- (0.5,17*\rr) -- (1.5,17*\rr) 
-- (2,16*\rr) -- (1.5,15*\rr);
\draw [line width=0.2mm] (0.5,17*\rr) -- (0,18*\rr) -- (0.5,19*\rr) -- (1.5,19*\rr) 
-- (2,18*\rr) -- (1.5,17*\rr);
\draw [line width=0.2mm] (0.5,19*\rr) -- (0,20*\rr) -- (0.5,21*\rr) -- (1.5,21*\rr) 
-- (2,20*\rr) -- (1.5,19*\rr);
%\draw [line width=0.2mm] (0.5,21*\rr) -- (0,22*\rr) -- (0.5,22*\rr) -- (1.5,22*\rr) 
%-- (2,22*\rr) -- (1.5,21*\rr);

\draw [line width=0.2mm]  (3,0) -- (3.5,-1*\rr);
\draw [line width=0.2mm] (2,2*\rr) -- (3,2*\rr) -- (3.5,1*\rr) -- (3,0) -- (2,0);
\draw [line width=0.2mm] (2,4*\rr) -- (3,4*\rr) -- (3.5,3*\rr) -- (3,2*\rr);
\draw [line width=0.2mm] (2,6*\rr) -- (3,6*\rr) -- (3.5,5*\rr) -- (3,4*\rr);
\draw [line width=0.2mm] (2,8*\rr) -- (3,8*\rr) -- (3.5,7*\rr) -- (3,6*\rr);
\draw [line width=0.2mm] (2,10*\rr) -- (3,10*\rr) -- (3.5,9*\rr) -- (3,8*\rr);
\draw [line width=0.2mm] (2,12*\rr) -- (3,12*\rr) -- (3.5,11*\rr) -- (3,10*\rr);
\draw [line width=0.2mm] (2,14*\rr) -- (3,14*\rr) -- (3.5,13*\rr) -- (3,12*\rr);
\draw [line width=0.2mm] (2,16*\rr) -- (3,16*\rr) -- (3.5,15*\rr) -- (3,14*\rr);
\draw [line width=0.2mm] (2,18*\rr) -- (3,18*\rr) -- (3.5,17*\rr) -- (3,16*\rr);
\draw [line width=0.2mm] (2,20*\rr) -- (3,20*\rr) -- (3.5,19*\rr) -- (3,18*\rr);
\draw [line width=0.2mm] (3.5,21*\rr) -- (3,20*\rr);

\draw [line width=0.2mm]  (3.5,1*\rr) -- (4.5,1*\rr) -- (5,0) -- (4.5,-1*\rr)
--(3.5,-1*\rr);
\draw [line width=0.2mm]  (3.5,3*\rr) -- (4.5,3*\rr) -- (5,2*\rr) -- (4.5,1*\rr);
\draw [line width=0.2mm]  (3.5,5*\rr) -- (4.5,5*\rr) -- (5,4*\rr) -- (4.5,3*\rr);
\draw [line width=0.2mm]  (3.5,7*\rr) -- (4.5,7*\rr) -- (5,6*\rr) -- (4.5,5*\rr);
\draw [line width=0.2mm]  (3.5,9*\rr) -- (4.5,9*\rr) -- (5,8*\rr) -- (4.5,7*\rr);
\draw [line width=0.2mm]  (3.5,11*\rr) -- (4.5,11*\rr) -- (5,10*\rr) -- (4.5,9*\rr);
\draw [line width=0.2mm]  (3.5,13*\rr) -- (4.5,13*\rr) -- (5,12*\rr) -- (4.5,11*\rr);
\draw [line width=0.2mm]  (3.5,15*\rr) -- (4.5,15*\rr) -- (5,14*\rr) -- (4.5,13*\rr);
\draw [line width=0.2mm]  (3.5,17*\rr) -- (4.5,17*\rr) -- (5,16*\rr) -- (4.5,15*\rr);
\draw [line width=0.2mm]  (3.5,19*\rr) -- (4.5,19*\rr) -- (5,18*\rr) -- (4.5,17*\rr);
\draw [line width=0.2mm]  (3.5,21*\rr) -- (4.5,21*\rr) -- (5,20*\rr) -- (4.5,19*\rr);

\draw [line width=0.2mm]  (5,0) -- (6,0)--(6.5,-1*\rr );
\draw [line width=0.2mm]  (5,2*\rr) -- (6,2*\rr)--(6.5,1*\rr )--(6,0);
\draw [line width=0.2mm]  (5,4*\rr) -- (6,4*\rr)--(6.5,3*\rr )--(6,2*\rr);
\draw [line width=0.2mm]  (5,6*\rr) -- (6,6*\rr)--(6.5,5*\rr )--(6,4*\rr);
\draw [line width=0.2mm]  (5,8*\rr) -- (6,8*\rr)--(6.5,7*\rr )--(6,6*\rr);
\draw [line width=0.2mm]  (5,10*\rr) -- (6,10*\rr)--(6.5,9*\rr )--(6,8*\rr);
\draw [line width=0.2mm]  (5,12*\rr) -- (6,12*\rr)--(6.5,11*\rr )--(6,10*\rr);
\draw [line width=0.2mm]  (5,14*\rr) -- (6,14*\rr)--(6.5,13*\rr )--(6,12*\rr);
\draw [line width=0.2mm]  (5,16*\rr) -- (6,16*\rr)--(6.5,15*\rr )--(6,14*\rr);
\draw [line width=0.2mm]  (5,18*\rr) -- (6,18*\rr)--(6.5,17*\rr )--(6,16*\rr);
\draw [line width=0.2mm]  (5,20*\rr) -- (6,20*\rr)--(6.5,19*\rr )--(6,18*\rr);
\draw [line width=0.2mm]  (6.5,21*\rr )--(6,20*\rr);

\draw [line width=0.2mm]  (6.5,1*\rr) -- (7.5,1*\rr)--(8,0)--(7.5,-1*\rr)
--(6.5,-1*\rr);
\draw [line width=0.2mm]  (6.5,3*\rr) -- (7.5,3*\rr)--(8,2*\rr)--(7.5,1*\rr);
\draw [line width=0.2mm]  (6.5,5*\rr) -- (7.5,5*\rr)--(8,4*\rr)--(7.5,3*\rr);
\draw [line width=0.2mm]  (6.5,7*\rr) -- (7.5,7*\rr)--(8,6*\rr)--(7.5,5*\rr);
\draw [line width=0.2mm]  (6.5,9*\rr) -- (7.5,9*\rr)--(8,8*\rr)--(7.5,7*\rr);
\draw [line width=0.2mm]  (6.5,11*\rr) -- (7.5,11*\rr)--(8,10*\rr)--(7.5,9*\rr);
\draw [line width=0.2mm]  (6.5,13*\rr) -- (7.5,13*\rr)--(8,12*\rr)--(7.5,11*\rr);
\draw [line width=0.2mm]  (6.5,15*\rr) -- (7.5,15*\rr)--(8,14*\rr)--(7.5,13*\rr);
\draw [line width=0.2mm]  (6.5,17*\rr) -- (7.5,17*\rr)--(8,16*\rr)--(7.5,15*\rr);
\draw [line width=0.2mm]  (6.5,19*\rr) -- (7.5,19*\rr)--(8,18*\rr)--(7.5,17*\rr);
\draw [line width=0.2mm]  (6.5,21*\rr) -- (7.5,21*\rr)--(8,20*\rr)--(7.5,19*\rr);

\draw [line width=0.2mm]  (8,0)--(9,0)--(9.5,-1*\rr);
\draw [line width=0.2mm]  (8,2*\rr)--(9,2*\rr)--(9.5,1*\rr)--(9,0);
\draw [line width=0.2mm]  (8,4*\rr)--(9,4*\rr)--(9.5,3*\rr)--(9,2*\rr);
\draw [line width=0.2mm]  (8,6*\rr)--(9,6*\rr)--(9.5,5*\rr)--(9,4*\rr);
\draw [line width=0.2mm]  (8,8*\rr)--(9,8*\rr)--(9.5,7*\rr)--(9,6*\rr);
\draw [line width=0.2mm]  (8,10*\rr)--(9,10*\rr)--(9.5,9*\rr)--(9,8*\rr);
\draw [line width=0.2mm]  (8,12*\rr)--(9,12*\rr)--(9.5,11*\rr)--(9,10*\rr);
\draw [line width=0.2mm]  (8,14*\rr)--(9,14*\rr)--(9.5,13*\rr)--(9,12*\rr);
\draw [line width=0.2mm]  (8,16*\rr)--(9,16*\rr)--(9.5,15*\rr)--(9,14*\rr);
\draw [line width=0.2mm]  (8,18*\rr)--(9,18*\rr)--(9.5,17*\rr)--(9,16*\rr);
\draw [line width=0.2mm]  (8,20*\rr)--(9,20*\rr)--(9.5,19*\rr)--(9,18*\rr);
\draw [line width=0.2mm]  (9.5,21*\rr)--(9,20*\rr);

\draw [line width=0.2mm]  (9.5,1*\rr)--(10.5,1*\rr)--(11,0)--(10.5,-1*\rr)
--(9.5,-1*\rr);
\draw [line width=0.2mm]  (9.5,3*\rr)--(10.5,3*\rr)--(11,2*\rr)--(10.5,1*\rr);
\draw [line width=0.2mm]  (9.5,5*\rr)--(10.5,5*\rr)--(11,4*\rr)--(10.5,3*\rr);
\draw [line width=0.2mm]  (9.5,7*\rr)--(10.5,7*\rr)--(11,6*\rr)--(10.5,5*\rr);
\draw [line width=0.2mm]  (9.5,9*\rr)--(10.5,9*\rr)--(11,8*\rr)--(10.5,7*\rr);
\draw [line width=0.2mm]  (9.5,11*\rr)--(10.5,11*\rr)--(11,10*\rr)--(10.5,9*\rr);
\draw [line width=0.2mm]  (9.5,13*\rr)--(10.5,13*\rr)--(11,12*\rr)--(10.5,11*\rr);
\draw [line width=0.2mm]  (9.5,15*\rr)--(10.5,15*\rr)--(11,14*\rr)--(10.5,13*\rr);
\draw [line width=0.2mm]  (9.5,17*\rr)--(10.5,17*\rr)--(11,16*\rr)--(10.5,15*\rr);
\draw [line width=0.2mm]  (9.5,19*\rr)--(10.5,19*\rr)--(11,18*\rr)--(10.5,17*\rr);
\draw [line width=0.2mm]  (9.5,21*\rr)--(10.5,21*\rr)--(11,20*\rr)--(10.5,19*\rr);

\draw [line width=0.2mm]  (11,0)--(12,0)--(12.5,-1*\rr);
\draw [line width=0.2mm]  (11,2*\rr)--(12,2*\rr)--(12.5,1*\rr)--(12,0);
\draw [line width=0.2mm]  (11,4*\rr)--(12,4*\rr)--(12.5,3*\rr)--(12,2*\rr);
\draw [line width=0.2mm]  (11,6*\rr)--(12,6*\rr)--(12.5,5*\rr)--(12,4*\rr);
\draw [line width=0.2mm]  (11,8*\rr)--(12,8*\rr)--(12.5,7*\rr)--(12,6*\rr);
\draw [line width=0.2mm]  (11,10*\rr)--(12,10*\rr)--(12.5,9*\rr)--(12,8*\rr);
\draw [line width=0.2mm]  (11,12*\rr)--(12,12*\rr)--(12.5,11*\rr)--(12,10*\rr);
\draw [line width=0.2mm]  (11,14*\rr)--(12,14*\rr)--(12.5,13*\rr)--(12,12*\rr);
\draw [line width=0.2mm]  (11,16*\rr)--(12,16*\rr)--(12.5,15*\rr)--(12,14*\rr);
\draw [line width=0.2mm]  (11,18*\rr)--(12,18*\rr)--(12.5,17*\rr)--(12,16*\rr);
\draw [line width=0.2mm]  (11,20*\rr)--(12,20*\rr)--(12.5,19*\rr)--(12,18*\rr);
\draw [line width=0.2mm]  (12,20*\rr)--(12.5,21*\rr);

\draw [line width=0.2mm]  (12.5,1*\rr)--(13.5,1*\rr)--(14,0)--(13.5,-1*\rr)
--(12.5,-1*\rr);
\draw [line width=0.2mm]  (12.5,3*\rr)--(13.5,3*\rr)--(14,2*\rr)--(13.5,1*\rr);
\draw [line width=0.2mm]  (12.5,5*\rr)--(13.5,5*\rr)--(14,4*\rr)--(13.5,3*\rr);
\draw [line width=0.2mm]  (12.5,7*\rr)--(13.5,7*\rr)--(14,6*\rr)--(13.5,5*\rr);
\draw [line width=0.2mm]  (12.5,9*\rr)--(13.5,9*\rr)--(14,8*\rr)--(13.5,7*\rr);
\draw [line width=0.2mm]  (12.5,11*\rr)--(13.5,11*\rr)--(14,10*\rr)--(13.5,9*\rr);
\draw [line width=0.2mm]  (12.5,13*\rr)--(13.5,13*\rr)--(14,12*\rr)--(13.5,11*\rr);
\draw [line width=0.2mm]  (12.5,15*\rr)--(13.5,15*\rr)--(14,14*\rr)--(13.5,13*\rr);
\draw [line width=0.2mm]  (12.5,17*\rr)--(13.5,17*\rr)--(14,16*\rr)--(13.5,15*\rr);
\draw [line width=0.2mm]  (12.5,19*\rr)--(13.5,19*\rr)--(14,18*\rr)--(13.5,17*\rr);
\draw [line width=0.2mm]  (12.5,21*\rr)--(13.5,21*\rr)--(14,20*\rr)--(13.5,19*\rr);

\draw [line width=0.2mm]  (14,0)--(15,0)--(15.5,-1*\rr);
\draw [line width=0.2mm]  (14,2*\rr)--(15,2*\rr)--(15.5,1*\rr)--(15,0);
\draw [line width=0.2mm]  (14,4*\rr)--(15,4*\rr)--(15.5,3*\rr)--(15,2*\rr);
\draw [line width=0.2mm]  (14,6*\rr)--(15,6*\rr)--(15.5,5*\rr)--(15,4*\rr);
\draw [line width=0.2mm]  (14,8*\rr)--(15,8*\rr)--(15.5,7*\rr)--(15,6*\rr);
\draw [line width=0.2mm]  (14,10*\rr)--(15,10*\rr)--(15.5,9*\rr)--(15,8*\rr);
\draw [line width=0.2mm]  (14,12*\rr)--(15,12*\rr)--(15.5,11*\rr)--(15,10*\rr);
\draw [line width=0.2mm]  (14,14*\rr)--(15,14*\rr)--(15.5,13*\rr)--(15,12*\rr);
\draw [line width=0.2mm]  (14,16*\rr)--(15,16*\rr)--(15.5,15*\rr)--(15,14*\rr);
\draw [line width=0.2mm]  (14,18*\rr)--(15,18*\rr)--(15.5,17*\rr)--(15,16*\rr);
\draw [line width=0.2mm]  (14,20*\rr)--(15,20*\rr)--(15.5,19*\rr)--(15,18*\rr);
\draw [line width=0.2mm]  (15,20*\rr)--(15.5,21*\rr);

\draw [line width=0.2mm]  (15.5,1*\rr)--(16.5,1*\rr)--(17,0)--(16.5,-1*\rr)
--(15.5,-1*\rr);
\draw [line width=0.2mm]  (15.5,3*\rr)--(16.5,3*\rr)--(17,2*\rr)--(16.5,1*\rr);
\draw [line width=0.2mm]  (15.5,5*\rr)--(16.5,5*\rr)--(17,4*\rr)--(16.5,3*\rr);
\draw [line width=0.2mm]  (15.5,7*\rr)--(16.5,7*\rr)--(17,6*\rr)--(16.5,5*\rr);
\draw [line width=0.2mm]  (15.5,9*\rr)--(16.5,9*\rr)--(17,8*\rr)--(16.5,7*\rr);
\draw [line width=0.2mm]  (15.5,11*\rr)--(16.5,11*\rr)--(17,10*\rr)--(16.5,9*\rr);
\draw [line width=0.2mm]  (15.5,13*\rr)--(16.5,13*\rr)--(17,12*\rr)--(16.5,11*\rr);
\draw [line width=0.2mm]  (15.5,15*\rr)--(16.5,15*\rr)--(17,14*\rr)--(16.5,13*\rr);
\draw [line width=0.2mm]  (15.5,17*\rr)--(16.5,17*\rr)--(17,16*\rr)--(16.5,15*\rr);
\draw [line width=0.2mm]  (15.5,19*\rr)--(16.5,19*\rr)--(17,18*\rr)--(16.5,17*\rr);
\draw [line width=0.2mm]  (15.5,21*\rr)--(16.5,21*\rr)--(17,20*\rr)--(16.5,19*\rr);

\draw [line width=0.2mm]  (17,0)--(18,0)--(18.5,-1*\rr);
\draw [line width=0.2mm]  (17,2*\rr)--(18,2*\rr)--(18.5,1*\rr)--(18,0);
\draw [line width=0.2mm]  (17,4*\rr)--(18,4*\rr)--(18.5,3*\rr)--(18,2*\rr);
\draw [line width=0.2mm]  (17,6*\rr)--(18,6*\rr)--(18.5,5*\rr)--(18,4*\rr);
\draw [line width=0.2mm]  (17,8*\rr)--(18,8*\rr)--(18.5,7*\rr)--(18,6*\rr);
\draw [line width=0.2mm]  (17,10*\rr)--(18,10*\rr)--(18.5,9*\rr)--(18,8*\rr);
\draw [line width=0.2mm]  (17,12*\rr)--(18,12*\rr)--(18.5,11*\rr)--(18,10*\rr);
\draw [line width=0.2mm]  (17,14*\rr)--(18,14*\rr)--(18.5,13*\rr)--(18,12*\rr);
\draw [line width=0.2mm]  (17,16*\rr)--(18,16*\rr)--(18.5,15*\rr)--(18,14*\rr);
\draw [line width=0.2mm]  (17,18*\rr)--(18,18*\rr)--(18.5,17*\rr)--(18,16*\rr);
\draw [line width=0.2mm]  (17,20*\rr)--(18,20*\rr)--(18.5,19*\rr)--(18,18*\rr);
\draw [line width=0.2mm]  (18,20*\rr)--(18.5,21*\rr);

\draw [line width=0.2mm]  (18.5,1*\rr)--(19.5,1*\rr)--(20,0)--(19.5,-1*\rr)
--(18.5,-1*\rr);
\draw [line width=0.2mm]  (18.5,3*\rr)--(19.5,3*\rr)--(20,2*\rr)--(19.5,1*\rr);
\draw [line width=0.2mm]  (18.5,5*\rr)--(19.5,5*\rr)--(20,4*\rr)--(19.5,3*\rr);
\draw [line width=0.2mm]  (18.5,7*\rr)--(19.5,7*\rr)--(20,6*\rr)--(19.5,5*\rr);
\draw [line width=0.2mm]  (18.5,9*\rr)--(19.5,9*\rr)--(20,8*\rr)--(19.5,7*\rr);
\draw [line width=0.2mm]  (18.5,11*\rr)--(19.5,11*\rr)--(20,10*\rr)--(19.5,9*\rr);
\draw [line width=0.2mm]  (18.5,13*\rr)--(19.5,13*\rr)--(20,12*\rr)--(19.5,11*\rr);
\draw [line width=0.2mm]  (18.5,15*\rr)--(19.5,15*\rr)--(20,14*\rr)--(19.5,13*\rr);
\draw [line width=0.2mm]  (18.5,17*\rr)--(19.5,17*\rr)--(20,16*\rr)--(19.5,15*\rr);
\draw [line width=0.2mm]  (18.5,19*\rr)--(19.5,19*\rr)--(20,18*\rr)--(19.5,17*\rr);
\draw [line width=0.2mm]  (18.5,21*\rr)--(19.5,21*\rr)--(20,20*\rr)--(19.5,19*\rr);

\draw [line width=0.2mm]  (20,0)--(21,0)--(21.5,-1*\rr);
\draw [line width=0.2mm]  (20,2*\rr)--(21,2*\rr)--(21.5,1*\rr)--(21,0);
\draw [line width=0.2mm]  (20,4*\rr)--(21,4*\rr)--(21.5,3*\rr)--(21,2*\rr);
\draw [line width=0.2mm]  (20,6*\rr)--(21,6*\rr)--(21.5,5*\rr)--(21,4*\rr);
\draw [line width=0.2mm]  (20,8*\rr)--(21,8*\rr)--(21.5,7*\rr)--(21,6*\rr);
\draw [line width=0.2mm]  (20,10*\rr)--(21,10*\rr)--(21.5,9*\rr)--(21,8*\rr);
\draw [line width=0.2mm]  (20,12*\rr)--(21,12*\rr)--(21.5,11*\rr)--(21,10*\rr);
\draw [line width=0.2mm]  (20,14*\rr)--(21,14*\rr)--(21.5,13*\rr)--(21,12*\rr);
\draw [line width=0.2mm]  (20,16*\rr)--(21,16*\rr)--(21.5,15*\rr)--(21,14*\rr);
\draw [line width=0.2mm]  (20,18*\rr)--(21,18*\rr)--(21.5,17*\rr)--(21,16*\rr);
\draw [line width=0.2mm]  (20,20*\rr)--(21,20*\rr)--(21.5,19*\rr)--(21,18*\rr);
\draw [line width=0.2mm]  (21,20*\rr)--(21.5,21*\rr);

\draw [line width=0.2mm]  (21.5,1*\rr)--(22.5,1*\rr)--(23,0)--(22.5,-1*\rr)
--(21.5,-1*\rr);
\draw [line width=0.2mm]  (21.5,3*\rr)--(22.5,3*\rr)--(23,2*\rr)--(22.5,1*\rr);
\draw [line width=0.2mm]  (21.5,5*\rr)--(22.5,5*\rr)--(23,4*\rr)--(22.5,3*\rr);
\draw [line width=0.2mm]  (21.5,7*\rr)--(22.5,7*\rr)--(23,6*\rr)--(22.5,5*\rr);
\draw [line width=0.2mm]  (21.5,9*\rr)--(22.5,9*\rr)--(23,8*\rr)--(22.5,7*\rr);
\draw [line width=0.2mm]  (21.5,11*\rr)--(22.5,11*\rr)--(23,10*\rr)--(22.5,9*\rr);
\draw [line width=0.2mm]  (21.5,13*\rr)--(22.5,13*\rr)--(23,12*\rr)--(22.5,11*\rr);
\draw [line width=0.2mm]  (21.5,15*\rr)--(22.5,15*\rr)--(23,14*\rr)--(22.5,13*\rr);
\draw [line width=0.2mm]  (21.5,17*\rr)--(22.5,17*\rr)--(23,16*\rr)--(22.5,15*\rr);
\draw [line width=0.2mm]  (21.5,19*\rr)--(22.5,19*\rr)--(23,18*\rr)--(22.5,17*\rr);
\draw [line width=0.2mm]  (21.5,21*\rr)--(22.5,21*\rr)--(23,20*\rr)--(22.5,19*\rr);

\draw [line width=0.2mm]  (23,0)--(24,0)--(24.5,-1*\rr);
\draw [line width=0.2mm]  (23,2*\rr)--(24,2*\rr)--(24.5,1*\rr)--(24,0);
\draw [line width=0.2mm]  (23,4*\rr)--(24,4*\rr)--(24.5,3*\rr)--(24,2*\rr);
\draw [line width=0.2mm]  (23,6*\rr)--(24,6*\rr)--(24.5,5*\rr)--(24,4*\rr);
\draw [line width=0.2mm]  (23,8*\rr)--(24,8*\rr)--(24.5,7*\rr)--(24,6*\rr);
\draw [line width=0.2mm]  (23,10*\rr)--(24,10*\rr)--(24.5,9*\rr)--(24,8*\rr);
\draw [line width=0.2mm]  (23,12*\rr)--(24,12*\rr)--(24.5,11*\rr)--(24,10*\rr);
\draw [line width=0.2mm]  (23,14*\rr)--(24,14*\rr)--(24.5,13*\rr)--(24,12*\rr);
\draw [line width=0.2mm]  (23,16*\rr)--(24,16*\rr)--(24.5,15*\rr)--(24,14*\rr);
\draw [line width=0.2mm]  (23,18*\rr)--(24,18*\rr)--(24.5,17*\rr)--(24,16*\rr);
\draw [line width=0.2mm]  (23,20*\rr)--(24,20*\rr)--(24.5,19*\rr)--(24,18*\rr);
\draw [line width=0.2mm]  (24,20*\rr)--(24.5,21*\rr);

%\draw [line width=0.6mm] (0.85,-1*\rr)--(1.5,1*\rr);

\draw [line width=0.6mm] (0.5,9*\rr)--(3,6*\rr);

\draw [line width=0.6mm] (0.5,9*\rr)--(4.5,11*\rr);

\draw [line width=0.6mm] (0,-4*\rr)--(7.5,21*\rr);

\draw [line width=0.6mm] (1.5,1*\rr)--(2.5,-1*\rr);
\draw [line width=0.6mm] (3,6*\rr)--(5,2*\rr);
\draw [line width=0.6mm] (4.5,11*\rr)--(6.5,7*\rr);
\draw [line width=0.6mm] (6,16*\rr)--(8,12*\rr);
\draw [line width=0.6mm] (7.5,21*\rr)--(9.5,17*\rr);

\draw [line width=0.6mm] (1.5,1*\rr)--(5,2*\rr);
\draw [line width=0.6mm] (3,6*\rr)--(6.5,7*\rr);
\draw [line width=0.6mm] (4.5,11*\rr)--(8,12*\rr);
\draw [line width=0.6mm] (6,16*\rr)--(9.5,17*\rr);

\draw [line width=0.6mm] (3.5,-3*\rr)--(11,22*\rr);

\draw [line width=0.6mm] (5,2*\rr)--(7.5,-1*\rr);
\draw [line width=0.6mm] (6.5,7*\rr)--(9,4*\rr);
\draw [line width=0.6mm] (8,12*\rr)--(10.5,9*\rr);
\draw [line width=0.6mm] (9.5,17*\rr)--(12,14*\rr);
\draw [line width=0.6mm] (11,22*\rr)--(13.5,19*\rr);

\draw [line width=0.6mm] (5,2*\rr)--(9,4*\rr);
\draw [line width=0.6mm] (6.5,7*\rr)--(10.5,9*\rr);
\draw [line width=0.6mm] (8,12*\rr)--(12,14*\rr);
\draw [line width=0.6mm] (9.5,17*\rr)--(13.5,19*\rr);

\draw [line width=0.6mm] (7.5,-1*\rr)--(15,24*\rr); 

\draw [line width=0.6mm] (9,4*\rr)--(11,0*\rr);
\draw [line width=0.6mm] (10.5,9*\rr)--(12.5,5*\rr);
\draw [line width=0.6mm] (12,14*\rr)--(14,10*\rr);
\draw [line width=0.6mm] (13.5,19*\rr)--(15.5,15*\rr);
\draw [line width=0.6mm] (16.5,21*\rr)--(17,20*\rr);

\draw [line width=0.6mm] (7.5,-1*\rr)--(11,0*\rr);
\draw [line width=0.6mm] (9,4*\rr)--(12.5,5*\rr);
\draw [line width=0.6mm] (10.5,9*\rr)--(14,10*\rr);
\draw [line width=0.6mm] (12,14*\rr)--(15.5,15*\rr);
\draw [line width=0.6mm] (13.5,19*\rr)--(17,20*\rr);

\draw [line width=0.6mm] (9.5,-5*\rr)--(18.5,25*\rr);

% \draw [line width=0.6mm] (12.5,5*\rr)--(15,2*\rr)--(18.5,3*\rr);
% \draw [line width=0.6mm] (14,10*\rr)--(16.5,7*\rr)--(20,8*\rr);

% \draw [line width=0.6mm] (18.5,3*\rr)--(21,0*\rr)--(24.5,1*\rr)--(22.5,5*\rr)--(18.5,3*\rr);
% \draw [line width=0.6mm] (21,0*\rr)--(22.5,5*\rr);

\draw [line width=0.6mm] (18.5,3*\rr)--(22.5,5*\rr)--(20,8*\rr)--(16.5,7*\rr)-- (18.5,3*\rr);
\draw [line width=0.6mm] (18.5,3*\rr)--(20,8*\rr);
\draw [line width=0.6mm] (16.5,7*\rr)--(14,10*\rr)--(18,12*\rr)--(20,8*\rr);
\draw [line width=0.6mm] (16.5,7*\rr)--(18,12*\rr);
\draw [line width=0.6mm] (12.5,5*\rr)--(16.5,7*\rr)--(15,2*\rr)
--(12.5,5*\rr);
\draw [line width=0.6mm] (15,2*\rr)--(18.5,3*\rr);

\foreach \pos in {(18.5,3*\rr),(22.5,5*\rr),(20,8*\rr),(16.5,7*\rr),(14,10*\rr),
(12.5,5*\rr),(15,2*\rr),(18,12*\rr)}
\shade[shading=ball, ball color=black] \pos circle (.4);

\draw (14.7,7.2*\rr) node   {\Large {${\mbox{\boldmath$O$}}$}};
\draw (12.2,11.2*\rr) node   {\Large {${\mbox{\boldmath$A$}}$}};
\draw (18,14.2*\rr) node   {\Large {${\mbox{\boldmath$B$}}$}};
\draw (21,10.5*\rr) node   {\Large {${\mbox{\boldmath$C$}}$}};
\draw (23,3*\rr) node   {\Large {${\mbox{\boldmath$D$}}$}};
\draw (18.5,1*\rr) node   {\Large {${\mbox{\boldmath$E$}}$}};
\draw (15,0.5*\rr) node   {\Large {${\mbox{\boldmath$F$}}$}}; 
\draw (13,2.3*\rr) node   {\Large {${\mbox{\boldmath$H$}}$}}; 
\end{tikzpicture}
\caption{A  $\beta$-PGS on $\bbH_2$ for $D^2=13$, with $a=1$, $b=3$ 
(Class HD1).\\
{\footnotesize{Here $|AO|^2=|OC|^2=|CD|^2=|HF|^2=|FE|^2=13$, $|BC|^2=|OE|^2=16$, 
$|AB|^2=|ED|^2=|OH|^2=19$, $|BO|^2=|OF|^2=|CE|^2=21$. Triangles $OBC$, $FOE$ and 
$OCE$ have area $4{\sqrt 3}$. Triangles $OAB$, $FHO$ and $CED$ have area $17{\sqrt 3}/4$. 
The number of PGSs equals 66; they are obtained from each other by $\bbH_2$-shifts and rotations by 
$\pm 2\pi /3$. Every PGS generates an EGM, 
and every EGM is generated by a PGS. Thus,  there are 66 EGMs which inherit the symmetries 
between their generating PGSs.}}}
\end{figure}} %Fig7 %FigureG7

%%%%%%%%%%%%%%%%%%%%%%%%%%%%%%%%%%%%%%%%%%%%%%%%%%
%%%%%%%%%%%%%%%%%%%%%%%%%%%%%%%%%%%%%%%%%%%%%%%%%%

\def\FigureH8  %%\label{Fig8}
{\begin{figure} \label{Fig8} %subClass \tE3
%\begin{center} 
\centering
\captionsetup{width=0.8\textwidth} 

\caption{PGSs ($\gamma$-configurations) on $\bbH_2$  for $D^2=16$, (Class HD2).\\
{\footnotesize{Inner equilateral $(D^2+D+1)$-triangles are colored gray.}}
}
%\vskip .2cm

\end{figure}} %FigureH8     %Fig8
%%%%%%%%%%%%%%%%%%%%%%%%%%%%%%%%%%%%%%%%%%%%%%%%%%%%%%%%%%%%%%%%%%%%%%% \captionsetup{width=0.8\textwidth}

\def\FigureI9  %%\label{Fig9} %\label{D2=67Triangle75}
{\begin{figure}\label{Fig9}\centering
\captionsetup{width=0.8\textwidth} 
(a)%
 
\caption{PGSs on $\bbH_2$  for $D^2=67$ (Class HE): (a)  type $\alpha$
 (vertical), (b) type $\beta$ (inclined).\\
{\footnotesize{The $\alpha$-PGSs are conjectured to be dominant.}}
}
%\vskip .2cm

\end{figure}} %\FigureI9      Fig9  

%%%%%%%%%%%%%%%%%%%%%%%%%%%%%%%%%%%%%%%%%%%%%%
%%%%%%%%%%%%%%%%%%%%%%%%%

\def\FigureJ10 %\figurea %%Fig10
{\begin{figure}[H]\label{Fig10}
\captionsetup{width=0.8\textwidth}
\centering
\begin{tikzpicture}[scale=0.45]
\filldraw [lightgray, yscale=sqrt(3/4), xslant=0.5] (0, 0) rectangle (7, 7);

\foreach \x in {0,...,6}
 \foreach \y in {0,...,6}
 {
  \filldraw[black] (\x + 0.5 * \y, +{sqrt(3/4)} * \y) circle (.15);
 }

\clip[yscale=sqrt(3/4), xslant=0.5] (-1, -1) rectangle (8, 8);
\draw[yscale=sqrt(3/4), xslant=0.5] (-1,-1) grid (8, 8);
\draw[yscale=sqrt(3/4), xslant=-0.5] (-1,-1) grid (15, 8);

\filldraw[black] (0,0) circle (.15);

\begin{scope}
\def\n{5}
\def\aa{2}
\def\bb{1}
\inclinedgrid
\end{scope}
\end{tikzpicture} \begin{tikzpicture}[scale=0.45]
\filldraw [lightgray, yscale=sqrt(3/4), xslant=0.5] (0, 0) rectangle (7, 7);

\foreach \x in {0,...,6}
 \foreach \y in {0,...,6}
 {
  \filldraw[black] (\x + 0.5 * \y, +{sqrt(3/4)} * \y) circle (.15);
 }

\clip[yscale=sqrt(3/4), xslant=0.5] (-1, -1) rectangle (8, 8);
\draw[yscale=sqrt(3/4), xslant=0.5] (-1,-1) grid (8, 8);
\draw[yscale=sqrt(3/4), xslant=-0.5] (-1,-1) grid (15, 8);

\filldraw[black] (0,0) circle (.15);

\begin{scope}
\def\n{5}
\def\aa{1}
\def\bb{2}
\inclinedgrid
\end{scope}
\end{tikzpicture}
\caption{The template (black circles) and $D$-rhombuses (thick lines) on $\bbA_2$, for $D^2=7$.\\
{\footnotesize{The gray rhombus represents a fundamental parallelogram 
for sub-lattice $\bbA_2(D^2)$.}}
}
%Templates and fundamental  parallelograms (thick lines) on ${\mathbb A}_2$ for $D^2=49$. 
%Marked sites (black circles) belong to the template. The larger mark identifies the origin.}
\end{figure}} %FigureJ10%\figurea %%Fig10

%%%%%%%%%%%%%%%%%%%%%%%%%%%%%%%%%%%%%%%%%%%%%%%%%%%%%%%%%%%%%%%%%%%%%%%

\def\FigureK11 %\figureb %%\label{Fig11}
{\begin{figure}[H]\label{Fig11}
\captionsetup{width=0.8\textwidth}
\centering

\begin{tikzpicture}[scale=0.28]

\begin{scope}
\def\n{4} \def\nn{18} \def\ss{7} \def\dx{0} \def\dy{0} \def\cc{lightgray}{0.2}
\rectsub-lattice
\end{scope} % the background color

\begin{scope}
\def\n{1} \def\nn{18} \def\ss{7} \def\dx{-7} \def\dy{0} \def\cc{gray}{0.5}
\rectsub-lattice
\end{scope}

%\begin{scope}
%\def\n{1} \def\nn{18} \def\ss{7} \def\dx{0} \def\dy{-7} \def\cc{gray}{0.5}
%\rectsub-lattice
%\end{scope}

\begin{scope}
\definecolor{gray6}{gray}{0.6}
\definecolor{gray5}{gray}{0.5}
\filldraw[yscale=sqrt(3/4), xslant=0.5, gray6] (-14,7) rectangle (-8,13);
\filldraw[yscale=sqrt(3/4), xslant=0.5, gray6] (-7,7) rectangle (-1,13);
\filldraw[yscale=sqrt(3/4), xslant=0.5, gray6] (0,7) rectangle (6,13);
\filldraw[yscale=sqrt(3/4), xslant=0.5, gray6] (-14,0) rectangle (-8,6);
\filldraw[yscale=sqrt(3/4), xslant=0.5, gray6] (0,0) rectangle (6,6);
\filldraw[yscale=sqrt(3/4), xslant=0.5, gray6] (7,0) rectangle (13,6);
\filldraw[yscale=sqrt(3/4), xslant=0.5, gray6] (-14,-7) rectangle (-8,-1);
\filldraw[yscale=sqrt(3/4), xslant=0.5, gray6] (-7,-7) rectangle (-1,-1);
\filldraw[yscale=sqrt(3/4), xslant=0.5, gray6] (7,-7) rectangle (13,-1);
\filldraw[yscale=sqrt(3/4), xslant=0.5, gray6] (-7,-14) rectangle (-1,-8);
\filldraw[yscale=sqrt(3/4), xslant=0.5, gray6] (0,-14) rectangle (6,-8);
\filldraw[yscale=sqrt(3/4), xslant=0.5, gray6] (7,-14) rectangle (13,-8);
\end{scope}

\begin{scope}
\definecolor{gray4}{gray}{0.4}
\filldraw[yscale=sqrt(3/4), xslant=0.5, gray4] (-7,0) rectangle (-1,6);
\filldraw[yscale=sqrt(3/4), xslant=0.5, gray4] (0,-7) rectangle (6,-1);
\end{scope}

\begin{scope}
\def\n{18}
\rectangulargrid(6, 8 *\rr)
\end{scope} % the co-ordinate mesh

\begin{scope}
\def\nn{18} \def\n{9} \def\dx{1} \def\dy{1} \def\aa{2} \def\bb{1}
\sublatticeinrectangle
\end{scope}   %the medium black balls

\begin{scope}
\def\n{4} \def\aa{0} \def\bb{7} \def\nn{18}
\clip[yscale=sqrt(3/4)] (-\nn, -\nn) rectangle (\nn, \nn);
\inclinedgrid
\end{scope}   %border lines of the templates

\begin{scope}
\def\xx{-4} \def\yy{2}
\shade[shading=ball, ball color=white] (\xx + 0.5 * \yy, \rr * \yy) circle (.5);
\def\xx{-6} \def\yy{1}
\shade[shading=ball, ball color=white] (\xx + 0.5 * \yy, \rr * \yy) circle (.5);
\def\xx{-7} \def\yy{4}
\shade[shading=ball, ball color=white] (\xx + 0.5 * \yy, \rr * \yy) circle (.5);
\def\xx{2} \def\yy{-2}
\shade[shading=ball, ball color=white] (\xx + 0.5 * \yy, \rr * \yy) circle (.5);
\def\xx{3} \def\yy{-5}
\shade[shading=ball, ball color=white] (\xx + 0.5 * \yy, \rr * \yy) circle (.5);
\def\xx{-5} \def\yy{2}
\shade[shading=ball, ball color=black] (\xx + 0.5 * \yy, \rr * \yy) circle (.5);
\end{scope}

\end{tikzpicture}

\caption{Templates: $\vphi$-correct (light-gray) and non-$\vphi$-correct (medium- and dark-gray), 
on  $\bbA_2$, in a $D$-AC for $D^2=7$. \\
{\footnotesize{Figure 11 demonstrates various types of templates. The small black balls indicate 
occupied sites in a $D$-PGS $\vphi$. The large black ball together with the small black balls 
indicate occupied sites $D$-AC $\phi$. The white balls mark vacant sites in $\phi$ which would be 
occupied in $\vphi$. The light-gray color indicates $\vphi$-correct templates. The dark-gray color 
indicates non-$\vphi$-correct templates containing some defects (where $\phi (\bx)\neq \vphi (\bx)$ 
for some site $\bx$). Over the medium-gray templates, configurations $\phi$ and $\vphi$ coincide.
However, these templates are still not $\vphi$-correct as they have neighboring templates with defects.}}
}
\end{figure}} %FigureK11 %\figureb %%Fig11

%%%%%%%%%%%%%%%%%%%%%%%%%%%%%%%%%%%%%%%%%%%%%%%%%%%%%%%%%%%%%%%%%%%%%%%

\def\FigureL12 %\figurec %%\label{Fig12}
{\begin{figure}[H] \label{Fig12}
\captionsetup{width=0.8\textwidth}
\centering
\begin{tikzpicture}[scale=0.25]
\clip (4, 0.5) rectangle (30, 24);

\begin{scope}[yscale=sqrt(3/4), xslant=1/2]
\definecolor{gray3}{gray}{0.3} %\definecolor{gray3}{lightgray} %
\definecolor{gray5}{gray}{0.7}
\path [fill=gray5, draw=black, very thick] (2, 6) -- (5, 6) -- (5, 9) -- (4, 9) -- (4, 11) -- (1, 11) -- (1, 8) -- (2, 8) -- cycle;
\path [fill=gray3, draw=black, thick] (3,7) -- (4, 7) -- (4, 8) -- (3, 8) -- cycle;
\path [fill=gray3, draw=black, thick] (2,9) -- (3, 9) -- (3, 10) -- (2, 10) -- cycle;
\path [fill=gray5, draw=black, very thick] (5,3) -- (8,3) -- (8,4) -- (11,4) -- (11,6) -- (14,6) -- (14,9) -- (15,9) --
(15,12) -- (16,12) -- (16,13) -- (18,13) -- (18,12) -- (19,12) -- (19,6) -- (14,6) -- (14,3) -- (15,3) -- (15,2) -- (20,2) --
(20,3) -- (23,3) -- (23,6) -- (22, 6) -- (22,14) -- (21,14) -- (21,15) -- (20,15) -- (20,16) -- (16,16) -- (16,21) -- (15,21) --
(15,22) -- (14,22) -- (14,23) -- (13,23) -- (13,25) -- (12,25) -- (12,26) -- (2,26) -- (2, 16) -- (6,16) -- (6,14) -- (4,14) --
(4,11) -- (7,11) -- (7,10) -- (8,10) -- (8,9) -- (9,9) -- (9,8) -- (10,8) -- (10,7) -- (8,7) -- (8,6) -- (5,6) -- cycle;
\path [fill=gray3, draw=black, thick] (6,4) -- (7,4) -- (7,5) -- (6,5) -- cycle;
\path [fill=gray3, draw=black, thick] (9,5) -- (10,5) -- (10,6) -- (9,6) -- cycle;
\path [fill=gray3, draw=black, thick] (5,12) -- (6,12) -- (6,13) -- (5,13) -- cycle;
\path [fill=gray3, draw=black, thick] (15,4) -- (16,4) -- (16,3) -- (19,3) -- (19,4) -- (22,4) -- (22,5) -- (21,5) -- (21,13)--
(20,13) -- (20,14) -- (19,14) -- (19,15) -- (15,15) -- (15,20) -- (14,20) -- (14,21) -- (13,21) -- (13,22) -- (12,22) -- (12,24)--
(11,24) -- (11,25) -- (3,25) -- (3,17) -- (7,17) -- (7,12) -- (8,12) -- (8,11) -- (9,11) -- (9,10) -- (10,10) -- (10,9) --
(11,9) -- (11,8) -- (12,8) --(12,7) -- (13,7) -- (13,10) -- (14,10) -- (14,11) -- (12,11) -- (12,10) -- (10,10) -- (10,11) --
(9,11) -- (9,12) -- (8,12) -- (8,13) -- (8,18) -- (4,18) -- (4,24) -- (9,24) -- (9,21) -- (13,21) -- (13,20) -- (14,20) --
(14,16) -- (13,16) -- (13,13) -- (15,13) -- (15,14) -- (19,14) -- (19,13) -- (20,13) -- (20,5) -- (17,5) -- (17,4) -- (16,4)--
(16,5) -- (15,5) -- cycle;
\path [fill=white, draw=black, very thick] (10,12) -- (12,12) -- (12,17) -- (13,17) -- (13,19) -- (12,19) -- (12,20) -- (8,20) --
(8,23) -- (5,23) -- (5,19) -- (9,19) -- (9,13) -- (10,13) -- cycle;

\draw (-10,0) grid (34,28);
\end{scope}

\path [draw=black, line width=0.5mm] (9,17) -- (19,14);
\path [draw=black, line width=0.5mm] (8,17) -- (22,8);
\path [draw=black, line width=0.5mm] (7,17) -- (11,7);

\path [fill=white, draw=white] (6,21) -- (6.5,21.5) -- (9.5,21.5) -- (10,21) -- (9.5,20.5) -- (6.5,20.5) -- (6,21) -- cycle;
\draw (8,21) node   {\Large {${\bf{Ext}}({\mbox{\boldmath$\Gamma$}})$}};
\path [fill=white, draw=white] (6,18) -- (6.5,18.5) -- (9.5,18.5) -- (10,18) -- (9.5,17.5) -- (6.5,17.5) -- (6,18) -- cycle;
\draw (8,18) node   {\Large {${\bf{Int}}({\mbox{\boldmath$\Gamma$}})$}};
\end{tikzpicture}

\caption{A contour support (the union of gray and dark-gray templates). \\
{\footnotesize{Here the internal area ${\rm{Int}}(\Gam )$ includes 
three components ${\rm{Int}}_{\vphi_i}(\Gam )$. The boundary layers are shown 
as the union of gray templates.}}
}
\end{figure}} %FigureL12 %\figurec %%\label{Fig12}

%%%%%%%%%%%%%%%%%%%%%%%%%%%%%%%%%%%%%%%%%%%%%%%%%%%%%%%%%%%%%%%%%%%%%%%
\def\FigureM13  %\figuree %%\label{Fig13}
{\begin{figure}\label{Fig13}%\begin{center}
\centering
\captionsetup{width=0.8\textwidth}

(a)\begin{tikzpicture}[scale=0.12]\label{Comp1}
\clip (-7.6, -3.4) rectangle (30.1, 26.6);

\filldraw [lightgray] (16,16)--(16,6)--(7,6)--(5,10)--(5,14)
--(9,18)--(14,18)--(16,16);

\path [draw=black, line width=0.5mm] (12,10)--(1,21);
\path [draw=black, line width=0.5mm] (12,10)--(12,26);
\path [draw=black, line width=0.5mm] (12,10)--(22,20);
\path [draw=black, line width=0.5mm] (12,10)--(20, 10);
\path [draw=black, line width=0.5mm] (12,10)--(12, 2);
\path [draw=black, line width=0.5mm] (12,10)--(-2,10);
\path [draw=black, line width=0.5mm] (12,10)--(0.666,4.333);
\path [draw=black, line width=0.5mm] (12,2)--(20,10);

\path [draw=black, line width=0.5mm] (1,21)--(12,26)--(22,20)--(20,10)--(21.5,4.5)
--(12,2)--(0.666,4.333)--(-2,10)--(1,21);

\filldraw[black] (12,26) circle (.5); \filldraw[black] (1, 21) circle (.5);
\filldraw[black] (22, 20) circle (.5); \filldraw[black] (-2,10) circle (.5);
\filldraw[black] (12,10) circle (.5); \filldraw[black] (20, 10) circle (.5);
\filldraw[black] (12,2) circle (.5); \filldraw[black] (0.666,4.333) circle (.5);

\filldraw[black] (28,23) circle (.5); \filldraw[black] (27,14) circle (.5);
\filldraw[black] (21.5,4.5) circle (.5); \filldraw[black] (-3,0) circle (.5);

\path [draw=black, line width=0.3mm] (16,16) circle (7.22);
\path [draw=black, line width=0.3mm] (16,6) circle (5.75);
\path [draw=black, line width=0.3mm] (7,6) circle (6.43);
\end{tikzpicture} \quad (b)\;\begin{tikzpicture}[scale=0.18]\label{Comp2}
\clip (7.9, 5.9) rectangle (26.1, 24);

\filldraw[lightgray] (12, 10) circle (4); \filldraw[lightgray] (22,20) circle (4);
\filldraw[lightgray] (20,10) circle (4); 

\filldraw[black] (12,10) circle (.5); \filldraw[black] (22, 20) circle (.5); 
\filldraw[black] (20, 10) circle (.5);

%\filldraw arc
\path [fill=gray] (12,10)--(16,10) arc (0:47:4)--(12,10);  
\path [fill=gray] (20,10)--(20.7845,13.9223) arc (80:180:4)--(20,10);
\path [fill=gray] (22,20)--(19.1716,17.1716) arc (225:261:4)--(22,20);

\filldraw[red] (16,16) circle (.5);

\path [draw=black, line width=0.5mm] (12,10)--(22,20)--(20,10)--(12,10);
\path [draw=black, line width=0.3mm] (16,16) circle (7.22);
\end{tikzpicture} %\end{center}
\caption{(a) A V-cell (a gray polygon), V-circles and \rC-triangles. (b) Set $\bbS (\triangle )$ 
for a \rC-triangle (the union of three dark-gray sectors; see Lemma 4.1).\\
{\footnotesize{Frame (a): vertices of a V-cell are centers of V-circles; all points of  a 
$D$-AC $\phi$ lie on V-circles, not inside. Frame (b): The radius of a V-circle around 
a \rC-triangle is $\leq D$. Consequently, one can't add a particle at the center (a red spot). In 
a saturated $D$-AC $\phi$, the radii of V-circles are $\leq D+1$: cf. Lemma 4.2.}} 
}
\end{figure}} %FigureM13 %Fig13

%%%%%%%%%%%%%%%%%%%%%%%%%%%%%%%%%%%%%%%%%%%%%%%%%%%%%%%%%%%%%%%%%%%%%%%
\def\FigureN14 %\figuree %%\label{Fig14}
{\begin{figure} \label{Fig14}\centering
\captionsetup{width=0.8\textwidth}
(a)\begin{tikzpicture}[scale=0.3]\label{Comp5}
\clip (-2.1,-4.5) rectangle (8.1, 10.1);

\filldraw [lightgray] (0,0)--(0,6)--(5,8)--(0,0);

%\draw[yscale=1, xslant=0] (-2,-3) grid (11, 11);

\path [draw=black, line width=0.5mm] (0,0)--(0,6)--(5,8)--(7,-1)--(0,0)
--(5,8);

\filldraw[black] (0,0) circle (.25);
\filldraw[black] (0,6) circle (.25);
\filldraw[black] (5,8) circle (.25);
\filldraw[black] (7,-1) circle (.25);

\draw (-1,-0.5) node   {\Large {${\mbox{\boldmath$A$}}$}};
\draw (-1.1,6) node   {\Large {${\mbox{\boldmath$B$}}$}};
\draw (5,9.3) node   {\Large {${\mbox{\boldmath$C$}}$}};
\draw (7,-2.2) node   {\Large {${\mbox{\boldmath$E$}}$}};
\end{tikzpicture}\quad (b)\begin{tikzpicture}[scale=0.3]\label{Comp6}
\clip (-2.1,-4.5) rectangle (11.1, 10.1);

\filldraw [lightgray] (0,0)--(0,6)--(5,8)--(0,0);
\filldraw [lightgray] (5,8)--(9,4)--(7,-1)--(5,8);

%\draw[yscale=1, xslant=0] (-2,-3) grid (11, 11);

\path [draw=black, line width=0.5mm] (0,0)--(0,6)--(5,8)--(9,4)--(7,-1)
--(5,8)--(0,0)--(7,-1);

\filldraw[black] (0,0) circle (.25);
\filldraw[black] (0,6) circle (.25);
\filldraw[black] (5,8) circle (.25);
\filldraw[black] (9,4) circle (.25);
\filldraw[black] (7,-1) circle (.25);

\draw (-1,-0.5) node   {\Large {${\mbox{\boldmath$A$}}$}};
\draw (-1.1,6) node   {\Large {${\mbox{\boldmath$B$}}$}};
\draw (5,9.3) node   {\Large {${\mbox{\boldmath$C$}}$}};
\draw (10.1,4) node   {\Large {${\mbox{\boldmath$D$}}$}};
\draw (7,-2.2) node   {\Large {${\mbox{\boldmath$E$}}$}};
\end{tikzpicture}
\quad (c)\begin{tikzpicture}[scale=0.3]\label{Comp4}
\clip (-2.1,-4.5) rectangle (11.1, 10.1);

\filldraw [lightgray] (0,0)--(0,6)--(5,8)--(0,0);
\filldraw [lightgray] (5,8)--(9,4)--(7,-1)--(5,8);
\filldraw [lightgray] (0,0)--(7,-1)--(3,-2.5)--(0,0);
%\draw[yscale=1, xslant=0] (-2,-3) grid (11, 11);

\path [draw=black, line width=0.5mm] (0,0)--(0,6)--(5,8)--(9,4)--(7,-1)--(5,8)--(0,0)--(7,-1)
--(3,-2.5)--(0,0);

\filldraw[black] (0,0) circle (.25);
\filldraw[black] (0,6) circle (.25);
\filldraw[black] (5,8) circle (.25);
\filldraw[black] (9,4) circle (.25);
\filldraw[black] (7,-1) circle (.25);
\filldraw[black] (3,-2.5) circle (.25);

\draw (-1,-0.5) node   {\Large {${\mbox{\boldmath$A$}}$}};
\draw (-1.1,6) node   {\Large {${\mbox{\boldmath$B$}}$}};
\draw (5,9.3) node   {\Large {${\mbox{\boldmath$C$}}$}};
\draw (10.1,4) node   {\Large {${\mbox{\boldmath$D$}}$}};
\draw (7,-2.2) node   {\Large {${\mbox{\boldmath$E$}}$}};
\draw (3,-3.7) node   {\Large {${\mbox{\boldmath$F$}}$}};
\end{tikzpicture}
\caption{Obtuse C-triangles (gray) with an adjacent acute tringle of a large area.\\
{\footnotesize{Frame (a): a 2-triangle group, frame (b): a 3-triangle group, frame (c): a 4-triangle
group. A quantitative estimation of the re-distributed area of the involved triangles is provided in Lemmas 4.5.1 - 4.5.3.}}
}
\end{figure}} %FigureN14 %\figuree %%Fig14

%%%%%%%%%%%%%%%%%%%%%%%%%%%%%%%%%%%%%%%%%%%%%%%%%%%%%%%%%%%
%%%%%%%%%%%%

\def\FigureO15{\begin{figure}[H]\label{Fig15} %\label{Fig15}  
\centering
\begin{tikzpicture}[scale=0.2]
%\clip (-0.9, -1*\rr) 
\clip (0, -1.5*\rr) rectangle (25, 21*\rr);

%\path [fill=lightgray] (1,-1*\rr)-- (4,-1*\rr)--(10.7,21*\rr)--(7.5,21*\rr);
%\path [fill=cyan] (4.2,-1*\rr)--(7.4,-1*\rr)--(14,21*\rr)--(10.8,21*\rr);
%\path [fill=lightgray] (7.5,-1*\rr)--(10.6,-1*\rr)--(17.3,21*\rr)--(14.2,21*\rr);

\draw [line width=0.2mm] (0.5,1*\rr) -- (1.5,1*\rr) -- (2,0) -- (1.5,-1*\rr)
-- (0.5,-1*\rr) -- (0,0) -- (0.5,1*\rr);
\draw [line width=0.2mm] (0.5,1*\rr) -- (0,2*\rr) -- (0.5,3*\rr) -- (1.5,3*\rr) 
-- (2,2*\rr) -- (1.5,1*\rr);
\draw [line width=0.2mm] (0.5,3*\rr) -- (0,4*\rr) -- (0.5,5*\rr) -- (1.5,5*\rr) 
-- (2,4*\rr) -- (1.5,3*\rr);
\draw [line width=0.2mm] (0.5,5*\rr) -- (0,6*\rr) -- (0.5,7*\rr) -- (1.5,7*\rr) 
-- (2,6*\rr) -- (1.5,5*\rr);
\draw [line width=0.2mm] (0.5,7*\rr) -- (0,8*\rr) -- (0.5,9*\rr) -- (1.5,9*\rr) 
-- (2,8*\rr) -- (1.5,7*\rr);
\draw [line width=0.2mm] (0.5,9*\rr) -- (0,10*\rr) -- (0.5,11*\rr) -- (1.5,11*\rr) 
-- (2,10*\rr) -- (1.5,9*\rr);
\draw [line width=0.2mm] (0.5,11*\rr) -- (0,12*\rr) -- (0.5,13*\rr) -- (1.5,13*\rr) 
-- (2,12*\rr) -- (1.5,11*\rr);
\draw [line width=0.2mm] (0.5,13*\rr) -- (0,14*\rr) -- (0.5,15*\rr) -- (1.5,15*\rr) 
-- (2,14*\rr) -- (1.5,13*\rr);
\draw [line width=0.2mm] (0.5,15*\rr) -- (0,16*\rr) -- (0.5,17*\rr) -- (1.5,17*\rr) 
-- (2,16*\rr) -- (1.5,15*\rr);
\draw [line width=0.2mm] (0.5,17*\rr) -- (0,18*\rr) -- (0.5,19*\rr) -- (1.5,19*\rr) 
-- (2,18*\rr) -- (1.5,17*\rr);
\draw [line width=0.2mm] (0.5,19*\rr) -- (0,20*\rr) -- (0.5,21*\rr) -- (1.5,21*\rr) 
-- (2,20*\rr) -- (1.5,19*\rr);
%\draw [line width=0.2mm] (0.5,21*\rr) -- (0,22*\rr) -- (0.5,22*\rr) -- (1.5,22*\rr) 
%-- (2,22*\rr) -- (1.5,21*\rr);

\draw [line width=0.2mm]  (3,0) -- (3.5,-1*\rr);
\draw [line width=0.2mm] (2,2*\rr) -- (3,2*\rr) -- (3.5,1*\rr) -- (3,0) -- (2,0);
\draw [line width=0.2mm] (2,4*\rr) -- (3,4*\rr) -- (3.5,3*\rr) -- (3,2*\rr);
\draw [line width=0.2mm] (2,6*\rr) -- (3,6*\rr) -- (3.5,5*\rr) -- (3,4*\rr);
\draw [line width=0.2mm] (2,8*\rr) -- (3,8*\rr) -- (3.5,7*\rr) -- (3,6*\rr);
\draw [line width=0.2mm] (2,10*\rr) -- (3,10*\rr) -- (3.5,9*\rr) -- (3,8*\rr);
\draw [line width=0.2mm] (2,12*\rr) -- (3,12*\rr) -- (3.5,11*\rr) -- (3,10*\rr);
\draw [line width=0.2mm] (2,14*\rr) -- (3,14*\rr) -- (3.5,13*\rr) -- (3,12*\rr);
\draw [line width=0.2mm] (2,16*\rr) -- (3,16*\rr) -- (3.5,15*\rr) -- (3,14*\rr);
\draw [line width=0.2mm] (2,18*\rr) -- (3,18*\rr) -- (3.5,17*\rr) -- (3,16*\rr);
\draw [line width=0.2mm] (2,20*\rr) -- (3,20*\rr) -- (3.5,19*\rr) -- (3,18*\rr);
\draw [line width=0.2mm] (3.5,21*\rr) -- (3,20*\rr);

\draw [line width=0.2mm]  (3.5,1*\rr) -- (4.5,1*\rr) -- (5,0) -- (4.5,-1*\rr)
--(3.5,-1*\rr);
\draw [line width=0.2mm]  (3.5,3*\rr) -- (4.5,3*\rr) -- (5,2*\rr) -- (4.5,1*\rr);
\draw [line width=0.2mm]  (3.5,5*\rr) -- (4.5,5*\rr) -- (5,4*\rr) -- (4.5,3*\rr);
\draw [line width=0.2mm]  (3.5,7*\rr) -- (4.5,7*\rr) -- (5,6*\rr) -- (4.5,5*\rr);
\draw [line width=0.2mm]  (3.5,9*\rr) -- (4.5,9*\rr) -- (5,8*\rr) -- (4.5,7*\rr);
\draw [line width=0.2mm]  (3.5,11*\rr) -- (4.5,11*\rr) -- (5,10*\rr) -- (4.5,9*\rr);
\draw [line width=0.2mm]  (3.5,13*\rr) -- (4.5,13*\rr) -- (5,12*\rr) -- (4.5,11*\rr);
\draw [line width=0.2mm]  (3.5,15*\rr) -- (4.5,15*\rr) -- (5,14*\rr) -- (4.5,13*\rr);
\draw [line width=0.2mm]  (3.5,17*\rr) -- (4.5,17*\rr) -- (5,16*\rr) -- (4.5,15*\rr);
\draw [line width=0.2mm]  (3.5,19*\rr) -- (4.5,19*\rr) -- (5,18*\rr) -- (4.5,17*\rr);
\draw [line width=0.2mm]  (3.5,21*\rr) -- (4.5,21*\rr) -- (5,20*\rr) -- (4.5,19*\rr);

\draw [line width=0.2mm]  (5,0) -- (6,0)--(6.5,-1*\rr );
\draw [line width=0.2mm]  (5,2*\rr) -- (6,2*\rr)--(6.5,1*\rr )--(6,0);
\draw [line width=0.2mm]  (5,4*\rr) -- (6,4*\rr)--(6.5,3*\rr )--(6,2*\rr);
\draw [line width=0.2mm]  (5,6*\rr) -- (6,6*\rr)--(6.5,5*\rr )--(6,4*\rr);
\draw [line width=0.2mm]  (5,8*\rr) -- (6,8*\rr)--(6.5,7*\rr )--(6,6*\rr);
\draw [line width=0.2mm]  (5,10*\rr) -- (6,10*\rr)--(6.5,9*\rr )--(6,8*\rr);
\draw [line width=0.2mm]  (5,12*\rr) -- (6,12*\rr)--(6.5,11*\rr )--(6,10*\rr);
\draw [line width=0.2mm]  (5,14*\rr) -- (6,14*\rr)--(6.5,13*\rr )--(6,12*\rr);
\draw [line width=0.2mm]  (5,16*\rr) -- (6,16*\rr)--(6.5,15*\rr )--(6,14*\rr);
\draw [line width=0.2mm]  (5,18*\rr) -- (6,18*\rr)--(6.5,17*\rr )--(6,16*\rr);
\draw [line width=0.2mm]  (5,20*\rr) -- (6,20*\rr)--(6.5,19*\rr )--(6,18*\rr);
\draw [line width=0.2mm]  (6.5,21*\rr )--(6,20*\rr);

\draw [line width=0.2mm]  (6.5,1*\rr) -- (7.5,1*\rr)--(8,0)--(7.5,-1*\rr)
--(6.5,-1*\rr);
\draw [line width=0.2mm]  (6.5,3*\rr) -- (7.5,3*\rr)--(8,2*\rr)--(7.5,1*\rr);
\draw [line width=0.2mm]  (6.5,5*\rr) -- (7.5,5*\rr)--(8,4*\rr)--(7.5,3*\rr);
\draw [line width=0.2mm]  (6.5,7*\rr) -- (7.5,7*\rr)--(8,6*\rr)--(7.5,5*\rr);
\draw [line width=0.2mm]  (6.5,9*\rr) -- (7.5,9*\rr)--(8,8*\rr)--(7.5,7*\rr);
\draw [line width=0.2mm]  (6.5,11*\rr) -- (7.5,11*\rr)--(8,10*\rr)--(7.5,9*\rr);
\draw [line width=0.2mm]  (6.5,13*\rr) -- (7.5,13*\rr)--(8,12*\rr)--(7.5,11*\rr);
\draw [line width=0.2mm]  (6.5,15*\rr) -- (7.5,15*\rr)--(8,14*\rr)--(7.5,13*\rr);
\draw [line width=0.2mm]  (6.5,17*\rr) -- (7.5,17*\rr)--(8,16*\rr)--(7.5,15*\rr);
\draw [line width=0.2mm]  (6.5,19*\rr) -- (7.5,19*\rr)--(8,18*\rr)--(7.5,17*\rr);
\draw [line width=0.2mm]  (6.5,21*\rr) -- (7.5,21*\rr)--(8,20*\rr)--(7.5,19*\rr);

\draw [line width=0.2mm]  (8,0)--(9,0)--(9.5,-1*\rr);
\draw [line width=0.2mm]  (8,2*\rr)--(9,2*\rr)--(9.5,1*\rr)--(9,0);
\draw [line width=0.2mm]  (8,4*\rr)--(9,4*\rr)--(9.5,3*\rr)--(9,2*\rr);
\draw [line width=0.2mm]  (8,6*\rr)--(9,6*\rr)--(9.5,5*\rr)--(9,4*\rr);
\draw [line width=0.2mm]  (8,8*\rr)--(9,8*\rr)--(9.5,7*\rr)--(9,6*\rr);
\draw [line width=0.2mm]  (8,10*\rr)--(9,10*\rr)--(9.5,9*\rr)--(9,8*\rr);
\draw [line width=0.2mm]  (8,12*\rr)--(9,12*\rr)--(9.5,11*\rr)--(9,10*\rr);
\draw [line width=0.2mm]  (8,14*\rr)--(9,14*\rr)--(9.5,13*\rr)--(9,12*\rr);
\draw [line width=0.2mm]  (8,16*\rr)--(9,16*\rr)--(9.5,15*\rr)--(9,14*\rr);
\draw [line width=0.2mm]  (8,18*\rr)--(9,18*\rr)--(9.5,17*\rr)--(9,16*\rr);
\draw [line width=0.2mm]  (8,20*\rr)--(9,20*\rr)--(9.5,19*\rr)--(9,18*\rr);
\draw [line width=0.2mm]  (9.5,21*\rr)--(9,20*\rr);

\draw [line width=0.2mm]  (9.5,1*\rr)--(10.5,1*\rr)--(11,0)--(10.5,-1*\rr)
--(9.5,-1*\rr);
\draw [line width=0.2mm]  (9.5,3*\rr)--(10.5,3*\rr)--(11,2*\rr)--(10.5,1*\rr);
\draw [line width=0.2mm]  (9.5,5*\rr)--(10.5,5*\rr)--(11,4*\rr)--(10.5,3*\rr);
\draw [line width=0.2mm]  (9.5,7*\rr)--(10.5,7*\rr)--(11,6*\rr)--(10.5,5*\rr);
\draw [line width=0.2mm]  (9.5,9*\rr)--(10.5,9*\rr)--(11,8*\rr)--(10.5,7*\rr);
\draw [line width=0.2mm]  (9.5,11*\rr)--(10.5,11*\rr)--(11,10*\rr)--(10.5,9*\rr);
\draw [line width=0.2mm]  (9.5,13*\rr)--(10.5,13*\rr)--(11,12*\rr)--(10.5,11*\rr);
\draw [line width=0.2mm]  (9.5,15*\rr)--(10.5,15*\rr)--(11,14*\rr)--(10.5,13*\rr);
\draw [line width=0.2mm]  (9.5,17*\rr)--(10.5,17*\rr)--(11,16*\rr)--(10.5,15*\rr);
\draw [line width=0.2mm]  (9.5,19*\rr)--(10.5,19*\rr)--(11,18*\rr)--(10.5,17*\rr);
\draw [line width=0.2mm]  (9.5,21*\rr)--(10.5,21*\rr)--(11,20*\rr)--(10.5,19*\rr);

\draw [line width=0.2mm]  (11,0)--(12,0)--(12.5,-1*\rr);
\draw [line width=0.2mm]  (11,2*\rr)--(12,2*\rr)--(12.5,1*\rr)--(12,0);
\draw [line width=0.2mm]  (11,4*\rr)--(12,4*\rr)--(12.5,3*\rr)--(12,2*\rr);
\draw [line width=0.2mm]  (11,6*\rr)--(12,6*\rr)--(12.5,5*\rr)--(12,4*\rr);
\draw [line width=0.2mm]  (11,8*\rr)--(12,8*\rr)--(12.5,7*\rr)--(12,6*\rr);
\draw [line width=0.2mm]  (11,10*\rr)--(12,10*\rr)--(12.5,9*\rr)--(12,8*\rr);
\draw [line width=0.2mm]  (11,12*\rr)--(12,12*\rr)--(12.5,11*\rr)--(12,10*\rr);
\draw [line width=0.2mm]  (11,14*\rr)--(12,14*\rr)--(12.5,13*\rr)--(12,12*\rr);
\draw [line width=0.2mm]  (11,16*\rr)--(12,16*\rr)--(12.5,15*\rr)--(12,14*\rr);
\draw [line width=0.2mm]  (11,18*\rr)--(12,18*\rr)--(12.5,17*\rr)--(12,16*\rr);
\draw [line width=0.2mm]  (11,20*\rr)--(12,20*\rr)--(12.5,19*\rr)--(12,18*\rr);
\draw [line width=0.2mm]  (12,20*\rr)--(12.5,21*\rr);

\draw [line width=0.2mm]  (12.5,1*\rr)--(13.5,1*\rr)--(14,0)--(13.5,-1*\rr)
--(12.5,-1*\rr);
\draw [line width=0.2mm]  (12.5,3*\rr)--(13.5,3*\rr)--(14,2*\rr)--(13.5,1*\rr);
\draw [line width=0.2mm]  (12.5,5*\rr)--(13.5,5*\rr)--(14,4*\rr)--(13.5,3*\rr);
\draw [line width=0.2mm]  (12.5,7*\rr)--(13.5,7*\rr)--(14,6*\rr)--(13.5,5*\rr);
\draw [line width=0.2mm]  (12.5,9*\rr)--(13.5,9*\rr)--(14,8*\rr)--(13.5,7*\rr);
\draw [line width=0.2mm]  (12.5,11*\rr)--(13.5,11*\rr)--(14,10*\rr)--(13.5,9*\rr);
\draw [line width=0.2mm]  (12.5,13*\rr)--(13.5,13*\rr)--(14,12*\rr)--(13.5,11*\rr);
\draw [line width=0.2mm]  (12.5,15*\rr)--(13.5,15*\rr)--(14,14*\rr)--(13.5,13*\rr);
\draw [line width=0.2mm]  (12.5,17*\rr)--(13.5,17*\rr)--(14,16*\rr)--(13.5,15*\rr);
\draw [line width=0.2mm]  (12.5,19*\rr)--(13.5,19*\rr)--(14,18*\rr)--(13.5,17*\rr);
\draw [line width=0.2mm]  (12.5,21*\rr)--(13.5,21*\rr)--(14,20*\rr)--(13.5,19*\rr);

\draw [line width=0.2mm]  (14,0)--(15,0)--(15.5,-1*\rr);
\draw [line width=0.2mm]  (14,2*\rr)--(15,2*\rr)--(15.5,1*\rr)--(15,0);
\draw [line width=0.2mm]  (14,4*\rr)--(15,4*\rr)--(15.5,3*\rr)--(15,2*\rr);
\draw [line width=0.2mm]  (14,6*\rr)--(15,6*\rr)--(15.5,5*\rr)--(15,4*\rr);
\draw [line width=0.2mm]  (14,8*\rr)--(15,8*\rr)--(15.5,7*\rr)--(15,6*\rr);
\draw [line width=0.2mm]  (14,10*\rr)--(15,10*\rr)--(15.5,9*\rr)--(15,8*\rr);
\draw [line width=0.2mm]  (14,12*\rr)--(15,12*\rr)--(15.5,11*\rr)--(15,10*\rr);
\draw [line width=0.2mm]  (14,14*\rr)--(15,14*\rr)--(15.5,13*\rr)--(15,12*\rr);
\draw [line width=0.2mm]  (14,16*\rr)--(15,16*\rr)--(15.5,15*\rr)--(15,14*\rr);
\draw [line width=0.2mm]  (14,18*\rr)--(15,18*\rr)--(15.5,17*\rr)--(15,16*\rr);
\draw [line width=0.2mm]  (14,20*\rr)--(15,20*\rr)--(15.5,19*\rr)--(15,18*\rr);
\draw [line width=0.2mm]  (15,20*\rr)--(15.5,21*\rr);

\draw [line width=0.2mm]  (15.5,1*\rr)--(16.5,1*\rr)--(17,0)--(16.5,-1*\rr)
--(15.5,-1*\rr);
\draw [line width=0.2mm]  (15.5,3*\rr)--(16.5,3*\rr)--(17,2*\rr)--(16.5,1*\rr);
\draw [line width=0.2mm]  (15.5,5*\rr)--(16.5,5*\rr)--(17,4*\rr)--(16.5,3*\rr);
\draw [line width=0.2mm]  (15.5,7*\rr)--(16.5,7*\rr)--(17,6*\rr)--(16.5,5*\rr);
\draw [line width=0.2mm]  (15.5,9*\rr)--(16.5,9*\rr)--(17,8*\rr)--(16.5,7*\rr);
\draw [line width=0.2mm]  (15.5,11*\rr)--(16.5,11*\rr)--(17,10*\rr)--(16.5,9*\rr);
\draw [line width=0.2mm]  (15.5,13*\rr)--(16.5,13*\rr)--(17,12*\rr)--(16.5,11*\rr);
\draw [line width=0.2mm]  (15.5,15*\rr)--(16.5,15*\rr)--(17,14*\rr)--(16.5,13*\rr);
\draw [line width=0.2mm]  (15.5,17*\rr)--(16.5,17*\rr)--(17,16*\rr)--(16.5,15*\rr);
\draw [line width=0.2mm]  (15.5,19*\rr)--(16.5,19*\rr)--(17,18*\rr)--(16.5,17*\rr);
\draw [line width=0.2mm]  (15.5,21*\rr)--(16.5,21*\rr)--(17,20*\rr)--(16.5,19*\rr);

\draw [line width=0.2mm]  (17,0)--(18,0)--(18.5,-1*\rr);
\draw [line width=0.2mm]  (17,2*\rr)--(18,2*\rr)--(18.5,1*\rr)--(18,0);
\draw [line width=0.2mm]  (17,4*\rr)--(18,4*\rr)--(18.5,3*\rr)--(18,2*\rr);
\draw [line width=0.2mm]  (17,6*\rr)--(18,6*\rr)--(18.5,5*\rr)--(18,4*\rr);
\draw [line width=0.2mm]  (17,8*\rr)--(18,8*\rr)--(18.5,7*\rr)--(18,6*\rr);
\draw [line width=0.2mm]  (17,10*\rr)--(18,10*\rr)--(18.5,9*\rr)--(18,8*\rr);
\draw [line width=0.2mm]  (17,12*\rr)--(18,12*\rr)--(18.5,11*\rr)--(18,10*\rr);
\draw [line width=0.2mm]  (17,14*\rr)--(18,14*\rr)--(18.5,13*\rr)--(18,12*\rr);
\draw [line width=0.2mm]  (17,16*\rr)--(18,16*\rr)--(18.5,15*\rr)--(18,14*\rr);
\draw [line width=0.2mm]  (17,18*\rr)--(18,18*\rr)--(18.5,17*\rr)--(18,16*\rr);
\draw [line width=0.2mm]  (17,20*\rr)--(18,20*\rr)--(18.5,19*\rr)--(18,18*\rr);
\draw [line width=0.2mm]  (18,20*\rr)--(18.5,21*\rr);

\draw [line width=0.2mm]  (18.5,1*\rr)--(19.5,1*\rr)--(20,0)--(19.5,-1*\rr)
--(18.5,-1*\rr);
\draw [line width=0.2mm]  (18.5,3*\rr)--(19.5,3*\rr)--(20,2*\rr)--(19.5,1*\rr);
\draw [line width=0.2mm]  (18.5,5*\rr)--(19.5,5*\rr)--(20,4*\rr)--(19.5,3*\rr);
\draw [line width=0.2mm]  (18.5,7*\rr)--(19.5,7*\rr)--(20,6*\rr)--(19.5,5*\rr);
\draw [line width=0.2mm]  (18.5,9*\rr)--(19.5,9*\rr)--(20,8*\rr)--(19.5,7*\rr);
\draw [line width=0.2mm]  (18.5,11*\rr)--(19.5,11*\rr)--(20,10*\rr)--(19.5,9*\rr);
\draw [line width=0.2mm]  (18.5,13*\rr)--(19.5,13*\rr)--(20,12*\rr)--(19.5,11*\rr);
\draw [line width=0.2mm]  (18.5,15*\rr)--(19.5,15*\rr)--(20,14*\rr)--(19.5,13*\rr);
\draw [line width=0.2mm]  (18.5,17*\rr)--(19.5,17*\rr)--(20,16*\rr)--(19.5,15*\rr);
\draw [line width=0.2mm]  (18.5,19*\rr)--(19.5,19*\rr)--(20,18*\rr)--(19.5,17*\rr);
\draw [line width=0.2mm]  (18.5,21*\rr)--(19.5,21*\rr)--(20,20*\rr)--(19.5,19*\rr);

\draw [line width=0.2mm]  (20,0)--(21,0)--(21.5,-1*\rr);
\draw [line width=0.2mm]  (20,2*\rr)--(21,2*\rr)--(21.5,1*\rr)--(21,0);
\draw [line width=0.2mm]  (20,4*\rr)--(21,4*\rr)--(21.5,3*\rr)--(21,2*\rr);
\draw [line width=0.2mm]  (20,6*\rr)--(21,6*\rr)--(21.5,5*\rr)--(21,4*\rr);
\draw [line width=0.2mm]  (20,8*\rr)--(21,8*\rr)--(21.5,7*\rr)--(21,6*\rr);
\draw [line width=0.2mm]  (20,10*\rr)--(21,10*\rr)--(21.5,9*\rr)--(21,8*\rr);
\draw [line width=0.2mm]  (20,12*\rr)--(21,12*\rr)--(21.5,11*\rr)--(21,10*\rr);
\draw [line width=0.2mm]  (20,14*\rr)--(21,14*\rr)--(21.5,13*\rr)--(21,12*\rr);
\draw [line width=0.2mm]  (20,16*\rr)--(21,16*\rr)--(21.5,15*\rr)--(21,14*\rr);
\draw [line width=0.2mm]  (20,18*\rr)--(21,18*\rr)--(21.5,17*\rr)--(21,16*\rr);
\draw [line width=0.2mm]  (20,20*\rr)--(21,20*\rr)--(21.5,19*\rr)--(21,18*\rr);
\draw [line width=0.2mm]  (21,20*\rr)--(21.5,21*\rr);

\draw [line width=0.2mm]  (21.5,1*\rr)--(22.5,1*\rr)--(23,0)--(22.5,-1*\rr)
--(21.5,-1*\rr);
\draw [line width=0.2mm]  (21.5,3*\rr)--(22.5,3*\rr)--(23,2*\rr)--(22.5,1*\rr);
\draw [line width=0.2mm]  (21.5,5*\rr)--(22.5,5*\rr)--(23,4*\rr)--(22.5,3*\rr);
\draw [line width=0.2mm]  (21.5,7*\rr)--(22.5,7*\rr)--(23,6*\rr)--(22.5,5*\rr);
\draw [line width=0.2mm]  (21.5,9*\rr)--(22.5,9*\rr)--(23,8*\rr)--(22.5,7*\rr);
\draw [line width=0.2mm]  (21.5,11*\rr)--(22.5,11*\rr)--(23,10*\rr)--(22.5,9*\rr);
\draw [line width=0.2mm]  (21.5,13*\rr)--(22.5,13*\rr)--(23,12*\rr)--(22.5,11*\rr);
\draw [line width=0.2mm]  (21.5,15*\rr)--(22.5,15*\rr)--(23,14*\rr)--(22.5,13*\rr);
\draw [line width=0.2mm]  (21.5,17*\rr)--(22.5,17*\rr)--(23,16*\rr)--(22.5,15*\rr);
\draw [line width=0.2mm]  (21.5,19*\rr)--(22.5,19*\rr)--(23,18*\rr)--(22.5,17*\rr);
\draw [line width=0.2mm]  (21.5,21*\rr)--(22.5,21*\rr)--(23,20*\rr)--(22.5,19*\rr);

\draw [line width=0.2mm]  (23,0)--(24,0)--(24.5,-1*\rr);
\draw [line width=0.2mm]  (23,2*\rr)--(24,2*\rr)--(24.5,1*\rr)--(24,0);
\draw [line width=0.2mm]  (23,4*\rr)--(24,4*\rr)--(24.5,3*\rr)--(24,2*\rr);
\draw [line width=0.2mm]  (23,6*\rr)--(24,6*\rr)--(24.5,5*\rr)--(24,4*\rr);
\draw [line width=0.2mm]  (23,8*\rr)--(24,8*\rr)--(24.5,7*\rr)--(24,6*\rr);
\draw [line width=0.2mm]  (23,10*\rr)--(24,10*\rr)--(24.5,9*\rr)--(24,8*\rr);
\draw [line width=0.2mm]  (23,12*\rr)--(24,12*\rr)--(24.5,11*\rr)--(24,10*\rr);
\draw [line width=0.2mm]  (23,14*\rr)--(24,14*\rr)--(24.5,13*\rr)--(24,12*\rr);
\draw [line width=0.2mm]  (23,16*\rr)--(24,16*\rr)--(24.5,15*\rr)--(24,14*\rr);
\draw [line width=0.2mm]  (23,18*\rr)--(24,18*\rr)--(24.5,17*\rr)--(24,16*\rr);
\draw [line width=0.2mm]  (23,20*\rr)--(24,20*\rr)--(24.5,19*\rr)--(24,18*\rr);
\draw [line width=0.2mm]  (24,20*\rr)--(24.5,21*\rr);

%\foreach \pos in {(13,18*\rr)}
%\shade[shading=ball, ball color=red] \pos circle (.4);

\draw [line width=0.5mm] (15,6*\rr)--(20,14*\rr);
\draw [line width=0.5mm] (6.5,7*\rr)--(16.5,23*\rr);
\draw [line width=0.5mm] (15,6*\rr)--(9,18*\rr)--(20,14*\rr);
\draw [line width=0.5mm] (15,6*\rr)--(11.5,15*\rr)--(20,14*\rr);
\draw [line width=0.5mm] (11.5,15*\rr)--(9,18*\rr);

\foreach \pos in {(11.5,15*\rr)}
\shade[shading=ball, ball color=red] \pos circle (.4);
\foreach \pos in {(6.5,7*\rr),(9,18*\rr),(15,6*\rr),(20,14*\rr)}
% (2,8*\rr),(0.5,17*\rr)}
\shade[shading=ball, ball color=black] \pos circle (.4);

% \foreach \pos in {(9,4*\rr),(15.5,1*\rr),(21,6*\rr),(3.5,-1*\rr)}
% \shade[shading=ball, ball color=black] \pos circle (.6);

%\foreach \pos in {(10.5, 7*\rr),(0, 14*\rr),(12.5, 17*\rr)}
 %\shade[shading=ball, ball color=white] \pos circle (.4);

 \draw (14.5,4*\rr) node   {\Large {${\mbox{\boldmath$O$}}$}};
 %\draw (5.7,13.2*\rr) node   {\Large {${\mbox{\boldmath$A$}}$}};
\draw (7,19*\rr) node   {\Large {${\mbox{\boldmath$B''$}}$}};
\draw (14,13*\rr) node   {\Large {${\mbox{\boldmath$B'$}}$}};
\draw (22,15.5*\rr) node   {\Large {${\mbox{\boldmath$C$}}$}};
%\draw (18.5,1*\rr) node   {\Large {${\mbox{\boldmath$E$}}$}};
 
\end{tikzpicture}
\caption{}
\end{figure}}  %Fig15 %\FigureO15

\def\FigureP16  %\figuree %%\label{Fig16}
{\begin{figure} \label{Fig16}
\captionsetup{width=0.8\textwidth} %, labelfont=bf}
\centering
(a)\quad \begin{tikzpicture}[scale=0.37]

\path [fill=lightgray, draw=black] (6.5, 3*\rr)  arc (21.79:81.79:7) arc (201.79:261.79:7);

\clip (-0.9, -1) rectangle (8.4, 10.6);

%\path [draw=black, thick] (8, 2*\rr) arc (0:98:7);
%\path [draw=black, thick] (-0.5, 9*\rr) arc (180:278:7);

\draw[yscale=sqrt(3/4), xslant=0.5] (-7,-2) grid (12, 13);
\draw[yscale=sqrt(3/4), xslant=-0.5] (-4,-2) grid (17, 13);

\foreach \pos in {(0, 0), (6.5, 3*\rr),(1,8* \rr), (7.5, 11*\rr)}
 \shade[shading=ball, ball color=black] \pos circle (.3);

\draw [line width=0.6mm] (0, 0)--(6.5,3* \rr)--(7.5, 11*\rr)--(1, 8*\rr)--(0, 0);
\draw [line width=0.6mm] (6.5,3* \rr)--(1, 8*\rr);

%\path [draw=orange, very thick] (2, 2*\rr) -- (5.5, 9*\rr) -- (1, 2*\rr) -- (6.5, 9*\rr) -- (2, 2*\rr)
%(1.5, 3*\rr) -- (6.5, 9*\rr)   (1, 2*\rr) -- (6, 8*\rr);

\foreach \pos in {(1,2*\rr),(2,2*\rr),(1.5,3*\rr),(2.5,3*\rr),(3.5,3*\rr),(4.5,3*\rr),
(2,4*\rr),(3,4*\rr),(4,4*\rr),(1.5,5*\rr),(2.5,5*\rr),(2,6*\rr)}
 \shade[shading=ball, ball color=orange] \pos circle (.25);
\foreach \pos in {(5.5,9*\rr),(6.5,9*\rr),(3,8*\rr),(4,8*\rr),(5,8*\rr),(6,8*\rr),
(3.5,7*\rr),(4.5,7*\rr),(5.5,7*\rr),(5,6*\rr),(6,6*\rr),(5.5,5*\rr)}
 \shade[shading=ball, ball color=orange] \pos circle (.25);
  
\draw [line width=0.4mm] (0, 0) arc (321.79:381.79:7) arc (201.79:261.79:7) arc (81.79:141.79:7);
\draw [line width=0.4mm] (6.5, 3*\rr) arc (21.79:81.79:7) arc (261.79:321.79:7) arc (141.79:201.79:7);

\end{tikzpicture}\qquad (b)\quad \begin{tikzpicture}[scale=0.38]

\clip (-0.9, -1) rectangle (8.4, 10.6);

\draw[yscale=sqrt(3/4), xslant=0.5] (-7,-2) grid (12, 13);
\draw[yscale=sqrt(3/4), xslant=-0.5] (-4,-2) grid (17, 13);

\foreach \pos in {(0, 0), (6.5, 3*\rr),(1, 8*\rr),(7.5, 11*\rr)}
 \shade[shading=ball, ball color=black] \pos circle (.3);

\draw [line width=0.6mm] (0, 0)--(6.5,3* \rr)--(7.5, 11*\rr)--(1, 8*\rr)--(0, 0);
\draw [line width=0.6mm] (6.5,3* \rr)--(1, 8*\rr);

\draw [line width=0.6mm, color=red] (2, 2*\rr)--(5.5,9*\rr);
\draw [line width=0.6mm, color=red] (1.5, 3*\rr)--(6.5,9*\rr);
%\path [draw=orange, very thick] (2, 2*\rr) -- (5.5, 9*\rr);
 %-- (1, 2*\rr) -- (6.5, 9*\rr) -- (2, 2*\rr)(1.5, 3*\rr) -- (6.5, 9*\rr)   (1, 2*\rr) -- (6, 8*\rr);

\foreach \pos in {(2,2*\rr),(1.5, 3*\rr),(5.5,9*\rr),(6.5,9*\rr)}
 \shade[shading=ball, ball color=red] \pos circle (.25);

%\foreach \pos in {(1, 2*\rr),(6.5, 9*\rr) }
 %\shade[shading=ball, ball color=lightgray] \pos circle (.3);
\end{tikzpicture}

(c) \begin{tikzpicture}[scale=0.37]

\path [fill=lightgray, draw=black] (7, 0)  arc (0:60:7) arc (180:240:7);

%\path [draw=black, thick] (8.5, \rr) arc (0:60:7);
%\path [draw=black, thick] (2, 6*\rr) arc (180:240:7);

\clip (-0.9, -1) rectangle (11.4, 7.1);
\draw[yscale=sqrt(3/4), xslant=0.5] (-4,-2) grid (12, 9);
\draw[yscale=sqrt(3/4), xslant=-0.5] (-1,-2) grid (15, 9);

\draw [line width=0.6mm] (0,0)--(7,0)--(10.5,7*\rr)--(3.5,7*\rr)--(0,0);
\draw [line width=0.6mm] (7,0)--(3.5,7*\rr);

%\path [draw=orange, very thick] (1.5, \rr) -- (8, 6*\rr) -- (2.5, \rr) -- (9, 6*\rr) -- (2, 2*\rr)
%-- (8.5, 5*\rr) -- (1.5, \rr) -- (9, 6*\rr);

\foreach \pos in {(5,6*\rr),(6,6*\rr),(8,6*\rr),(9,6*\rr),(5.5,5*\rr),(6.5,5*\rr),
(7.5,5*\rr),(8.5,5*\rr),(7,4*\rr),(6.5,3*\rr),(7.5,3*\rr),(7,2*\rr)}
 \shade[shading=ball, ball color=orange] \pos circle (.25);

\foreach \pos in {(1.5,1*\rr),(2.5,1*\rr),(4.5,1*\rr),(5.5,1*\rr),(2,2*\rr),(3,2*\rr),
(4,2*\rr),(5,2*\rr),(3.5,3*\rr),(3,4*\rr),(4,4*\rr),(3.5,5*\rr)}
\shade[shading=ball, ball color=orange] \pos circle (.25);
 %\filldraw[white] \pos circle (.25);
 
\foreach \pos in {(0, 0), (7, 0), (3.5,7*\rr),(10.5,7*\rr)}
\shade[shading=ball, ball color=black] \pos circle (.3);

%\draw (8,-0.3*\rr) node {\large {${\mbox{\boldmath$O$}}$}};
%\draw (-0.3,1.2*\rr) node {\large {${\mbox{\boldmath$A$}}$}};
%\draw (2.3,7*\rr) node {\large {${\mbox{\boldmath$B$}}$}};
%\draw (10.7,5.6*\rr) node {\large {${\mbox{\boldmath$C$}}$}};
\draw [line width=0.4mm] (0, 0) arc (300:360:7) arc (180:240:7) arc (60:120:7);
\draw [line width=0.4mm] (7, 0) arc (0:60:7) arc (240:300:7) arc (120:180:7);

\end{tikzpicture}\;(d) \begin{tikzpicture}[scale=0.3]
\clip (-0.9, -1) rectangle (11.4, 7.1);
\draw[yscale=sqrt(3/4), xslant=0.5] (-4,-2) grid (12, 9);
\draw[yscale=sqrt(3/4), xslant=-0.5] (-1,-2) grid (15, 9);

\draw [line width=0.6mm] (0,0)--(7,0)--(10.5,7*\rr)--(3.5,7*\rr)--(0,0);
\draw [line width=0.6mm] (7,0)--(3.5,7*\rr);

\draw [line width=0.6mm, color=red] (2.5,1*\rr)--(9,6*\rr);
\draw [line width=0.6mm, color=red] (2,2*\rr)--(9,6*\rr);
% \draw [line width=0.6mm, color=orange] (2,2*\rr)--(8.5,5* \rr);
%\path [draw=orange, very thick] (2.5,1*\rr)--(9, 6*\rr);

% (1.5,1*\rr)--(8,6*\rr)--(2.5,1*\rr)--(9,6*\rr)--(2,2*\rr)
% -- (8.5,5* \rr)-- (1.5,1*\rr)--(9, 6*\rr);

\foreach \pos in {(9,6*\rr)}
\shade[shading=ball, ball color=red] \pos circle (.25);

\foreach \pos in {(2.5,1*\rr),(2,2*\rr)} 
\shade[shading=ball, ball color=red] \pos circle (.25);
 
\foreach \pos in {(0, 0), (7, 0), (3.5,7*\rr),(10.5,7*\rr)}
\shade[shading=ball, ball color=black] \pos circle (.3);
\end{tikzpicture}
\caption{Single and double $u^{-2}$-insertions for  $D^2=49$ on $\bbA_2$.\\
{\footnotesize{For each PGS there are 12 sites inside a $D$-triangle (orange balls) where
a single insertion repels 3 particles at  the vertices of the triangle (black balls).
Examples of double insertions repelling 4 particles at the vertices of 
$D$-rhombus are shown in red. There are 6 double insertions in a
$(5, 3)$-rhombus and 7 in a $(7, 0)$-one. As there is no other $u^{-2}$-insertions,  the 
$(7, 0)$-class is dominant.}}
}
\end{figure}} %FigureP16  %\figuree %%\label{Fig16}

%%%%%%%%%%%%%%%%%%%%%%%%%%%%%%%%%%%%%%%%%%%%%%%%%%%%%%%%%%%%%%%%%%%%%%%

\def\FigureQ17  %%\figuref%%\label{Fig17}
{\begin{figure}[H] \label{Fig17}\centering
\captionsetup{width=0.8\textwidth}  %, labelfont=bf}
\centering
(a) \begin{tikzpicture}[scale=0.25]
\clip (-0.9, -1) rectangle (12.9, 20.1);

\path [fill=lightgray, draw=black] (11.5, 7*\rr) arc (27.8:87.8:13) arc (207.8:267.8:13);

%\draw [line width=0.9mm,color=orange] (10,18*\rr)--(2,6*\rr);
%\draw [line width=0.9mm,color=orange] (10,14*\rr)--(1,2*\rr);

\draw [line width=0.6mm] (0,0) -- (0.5, 15*\rr) -- (12, 22*\rr) -- (11.5, 7*\rr) -- (0,0);
\draw [line width=0.6mm] (0.5, 15*\rr) -- (11.5, 7*\rr);

\draw[yscale=sqrt(3/4), xslant=0.5] (-12,-2) grid (21, 24);
\draw[yscale=sqrt(3/4), xslant=-0.5] (-1,-2) grid (31, 24);

\foreach \pos in {(0, 0), (11.5, 7*\rr),(0.5, 15*\rr),(12, 22*\rr) }
 \shade[shading=ball, ball color=black] \pos circle (.4);

\foreach \pos in {(0.5,1*\rr),(1,2*\rr),(1.5,3*\rr), % 1+1+1
(2,4*\rr),(3,4*\rr),  % 2
(2.5,5*\rr),(3.5,5*\rr),(4.5,5*\rr),  % 3
(2,6*\rr),(3,6*\rr),(4,6*\rr),(5,6*\rr),(6,6*\rr), % 5
(2.5,7*\rr),(3.5,7*\rr),(4.5,7*\rr),(5.5,7*\rr),(6.5,7*\rr),
(7.5,7*\rr), (8.5,7*\rr),(9.5,7*\rr),(10.5,7*\rr), % 9
(2,8*\rr),(3,8*\rr),(4,8*\rr),(5,8*\rr),(6,8*\rr),(7,8*\rr), % 6
(2.5,9*\rr),(3.5,9*\rr),(4.5,9*\rr),(5.5,9*\rr), % 4
(2,10*\rr),(3,10*\rr),(4,10*\rr), % 3
(2.5,11*\rr),(2,12*\rr), (1.5,13*\rr),(1,14*\rr),  % 1+1+1+1
(11.5,21*\rr),(11,20*\rr),(10.5,19*\rr), % 1+1+1
(10,18*\rr),(9,18*\rr),  % 2
(9.5,17*\rr),(8.5,17*\rr),(7.5,17*\rr), % 3
(10,16*\rr),(9,16*\rr),(8,16*\rr),(7,16*\rr),(6,16*\rr), % 5
(9.5,15*\rr),(8.5,15*\rr),(7.5,15*\rr),(6.5,15*\rr),(5.5,15*\rr),
(4.5,15*\rr),(3.5,15*\rr),(2.5,15*\rr),(1.5,15*\rr), % 9
(10,14*\rr),(9,14*\rr),(8,14*\rr),(7,14*\rr),(6,14*\rr),(5,14*\rr), % 6
(9.5,13*\rr),(8.5,13*\rr),(7.5,13*\rr),(6.5,13*\rr), % 4
(10,12*\rr),(9,12*\rr),(8,12*\rr), % 3
(9.5,11*\rr),(10,10*\rr),(10.5,9*\rr),(11,8*\rr)} % 1+1+1+1
\shade[shading=ball, ball color=orange] \pos circle (.25);
% Total:  2(1+1+1+2+3+5+9+6+4+3+1+1+1+1)=2x39=78

%\foreach \pos in {(10,18*\rr),(2,6*\rr)}
%\shade[shading=ball, ball color=orange] \pos circle (.3);

%\foreach \pos in {(10,14*\rr),(1,2*\rr)}
%\shade[shading=ball, ball color=orange] \pos circle (.3);

%\foreach \pos in {(0.5,\rr),(11.5, 21*\rr)}
%\shade[shading=ball, ball color=white] \pos circle (.4);
\path [draw=black,line width=0.4mm] (0, 0) arc (327.8:387.8:13) arc (207.8:267.8:13) arc (87.8:147.8:13);
\path [draw=black,line width=0.4mm] (11.5, 7*\rr) arc (27.8:87.8:13) arc (267.8:327.8:13) arc (147.8:207.8:13);

\end{tikzpicture}\quad (b)\begin{tikzpicture}[scale=0.25]
\clip (-0.9, -1) rectangle (20.4, 12.5);

\path [fill=lightgray, draw=black] (13,0) arc (0:60:13) arc (180:240:13);

%\draw [line width=0.9mm, color=red] (4,2*\rr)--(15.5,9*\rr);
%\draw [line width=0.9mm, color=red] (3.5,3*\rr)--(16.5,11*\rr);

% \draw [line width=0.9mm,color=orange] (10,18*\rr)--(2,6*\rr);
% \draw [line width=0.9mm,color=orange] (10,14*\rr)--(1,2*\rr);

\draw [line width=0.6mm]  (0,0) -- (6.5, 13*\rr) -- (19.5, 13*\rr) -- (13,0) -- (0,0);
\draw [line width=0.6mm]  (6.5, 13*\rr) -- (13,0);

\draw[yscale=sqrt(3/4), xslant=0.5] (-8,-2) grid (21, 15);
\draw[yscale=sqrt(3/4), xslant=-0.5] (-1,-2) grid (27, 15);

\foreach \pos in {(0, 0), (13, 0), (6.5, 13*\rr),(19.5, 13*\rr) }
 \shade[shading=ball, ball color=black] \pos circle (.4);

% \foreach \pos in {(1.5, \rr),(18, 12*\rr) }
 %  \shade[shading=ball, ball color=red] \pos circle (.3);

\foreach \xx in {1,2,3,4,5,6}
 \foreach \yy in {1,...,\xx}
   \shade[shading=ball, ball color=orange] (5 - 0.5 * \xx + \yy, 8 * \rr - \xx * \rr ) circle (.25);

\foreach \xx in {1,2,3,4,5,6}
 \foreach \yy in {1,...,\xx}
   \shade[shading=ball, ball color=orange] (13.5 - 0.5 * \xx + \yy, 5 * \rr + \xx * \rr ) circle (.25);

\foreach \pos in {(1.5,1*\rr),(11.5,1*\rr),(9,2*\rr),(10,2*\rr),(8.5,3*\rr),
(9.5,3*\rr),(8,4*\rr),(9,4*\rr),(7.5,5*\rr),(8.5,5*\rr),(7,6*\rr),(8,6*\rr),(6.5,7*\rr),(7.5,7*\rr),
(6,8*\rr),(7,8*\rr),(6.5,9*\rr),(6.5,11*\rr)}
\shade[shading=ball, ball color=orange] \pos circle (.25);

\foreach \pos in {(8,12*\rr),(18,12*\rr),(9.5,11*\rr),(10.5,11*\rr),
(10,10*\rr),(11,10*\rr),(10.5,9*\rr),(11.5,9*\rr),(11,8*\rr),(12,8*\rr),
(11.5,7*\rr),(12.5,7*\rr),(12,6*\rr),(13,6*\rr),(12.5,5*\rr),(13.5,5*\rr),
(13,4*\rr),(13,2*\rr)}
\shade[shading=ball, ball color=orange] \pos circle (.25);
 %  Total: 2(2+8+7+6+5+4+3+2+1+1)=2x39=78

%\foreach \pos in {(4,2*\rr),(15.5,9*\rr)}
%\shade[shading=ball, ball color=red] \pos circle (.3);

%\foreach \pos in {(3.5,3*\rr),(16.5,11*\rr)}
%\shade[shading=ball, ball color=red] \pos circle (.3);

\path [draw=black,line width=0.4mm] (0, 0) arc (300:360:13) arc (180:240:13) arc (60:120:13);
\path [draw=black,line width=0.4mm] (13, 0) arc (0:60:13) arc (240:300:13) arc (120:180:13);

\end{tikzpicture}
\caption{Single $u^{-2}$-insertions for  $D^2=169$  on $\bbA_2$.\\
{\footnotesize{We again use orange balls and circular triangles for marking the positions 
where an inserted particle repels 3 vertices of a $D$-triangle in a PGS. In both 
PGS types the number of single $u^{-2}$-insertions equals 39 per triangle or 78 per a
$D$-rhombus.}}
}

% \caption{The structure of pair defects for a distinguishing contour of weight $u^{-2}$ used
% in the proof of Theorem~5,
% with $D^2=169$. The contour exhibits 4 removed and 2 added
% particles, affecting the vertices of a $D$-rhombus. Both inclined (left) and horizontal (right)
%ground states are shown: white circles indicate the removed ground-state particles. The
% origin ${\mathbf 0}=(0,0)$ is at the position of the bottom left white ball. The
% light-gray areas (closed circular triangles) cover the lattice sites which repel 3 out of 4
% white balls but do not repel any other ground-state particles. The gray area
% (an open bi-convex lense) includes the sites that repel all 4 white balls but no other ground-state
% particles. Each black circle from one triangular stack, together
% with at least one black circle from the opposite stack, can be selected to form
% a pair defect. There are 113 pairs that can be formed for the inclined $D$-sub-lattice
% and 78 pairs for the horizontal one (out of 625 and 484 nominal pairs, respectively). Again,
% a larger black ball can be paired with any ball from the opposite stack. }
\end{figure}} %FigureQ17  %%Fig17

%%%%%%%%%%%%%%%%%%%%%%%%%%%%%%%%%%%%%%%%%%%%%%%%%%%%%%%%%%%%%%%%%%%%%%%

\def\FigureR18 %\fgri %\label{Fig18}
{\begin{figure}[H] \label{Fig18}
\centering
\captionsetup{width=0.8\textwidth} %, labelfont=bf}
\def\rr{0.86602540378443864676372317075294} %\sqrt{3/4}

(a) \begin{tikzpicture}[scale=0.2]
\clip (-1.2, -1.2*\rr) rectangle (24.2, 32.2*\rr);

\draw [line width=0.6mm, color=red] (2.5, 15*\rr)--(16, 16*\rr);
\draw [line width=0.6mm, color=red] (10.5,19*\rr)--(4,6*\rr);
\draw [line width=0.6mm, color=blue] (1,2*\rr)--(9.5,17*\rr)--(1.5, 29*\rr);
\draw [line width=0.6mm, color=green] (1,2*\rr)--(8,16*\rr)--(1.5, 29*\rr);
\draw [line width=0.6mm, color=green] (8,16*\rr)--(21, 14*\rr);

\draw [line width=0.6mm] (0,0) -- (1, 30*\rr) -- (23,14*\rr) -- cycle
(0.5, 15*\rr) -- (11.5, 7*\rr) -- (12, 22*\rr) -- cycle;

\draw[yscale=sqrt(3/4), xslant=0.5] (-19,-2) grid (90, 60);
\draw[yscale=sqrt(3/4), xslant=-0.5] (-1,-2) grid (90, 60);

\foreach \pos in {(0, 0), (11.5, 7*\rr),(23, 14*\rr),(0.5, 15*\rr),(1, 30*\rr),(12, 22*\rr)}
 \shade[shading=ball, ball color=black] \pos circle (.4);

\foreach \pos in {(10.5,19*\rr),(4,6*\rr),(2.5, 15*\rr),(16, 16*\rr)}
\shade[shading=ball, ball color=red] \pos circle (.5);
 
 \foreach \pos in {(9.5,17*\rr),(1,2*\rr),(1.5, 29*\rr)}
\shade[shading=ball, ball color=blue] \pos circle (.5);

\foreach \pos in {(8,16*\rr),(21, 14*\rr)}
 \shade[shading=ball, ball color=green] \pos circle (.5);
\foreach \pos in {(1,2*\rr),(1.5, 29*\rr)}
 \shade[shading=ball, ball color=green] \pos circle (.25);

\end{tikzpicture}\quad
(b) \begin{tikzpicture}[scale=0.2]
\clip (-1.2, -1.2*\rr) rectangle (27.2, 27.5*\rr);

% \path [draw=black, fill=lightgray]
% (0,0) arc (120:60:13) arc (240:180:13) arc (360:300:13)
% (6.5,13*\rr) arc (120:60:13) arc (240:180:13) arc (360:300:13)
% (13,0) arc (120:60:13) arc (240:180:13) arc (360:300:13)
% (13,0) arc (180:120:13) arc (300:240:13) arc (60:0:13);

\path [draw=red, line width=0.6mm] (5,4*\rr)--(18, 12*\rr);
\path [draw=red, line width=0.6mm] (8,12*\rr)--(21, 5*\rr);
\path [draw=blue, line width=0.6mm] (1.5, \rr)--(15,8*\rr)--(13, 24*\rr);
\path [draw=green, line width=0.6mm] (1.5, \rr)--(13,8*\rr)--(13, 24*\rr);
\path [draw=green, line width=0.6mm] (13,8*\rr)--(24.5, \rr);

\path [draw=black, line width=0.6mm] (0,0) -- (13, 26*\rr) -- (26, 0) -- cycle
(13,0) -- (6.5,13*\rr) -- (19.5, 13*\rr) -- cycle;

\draw[yscale=sqrt(3/4), xslant=0.5] (-90,-60) grid (90, 60);
\draw[yscale=sqrt(3/4), xslant=-0.5] (-90,-60) grid (90, 60);

\foreach \pos in {(0,0),(6.5,13*\rr),(13, 26*\rr),(13,0),(26, 0), (19.5, 13*\rr)}
 \shade[shading=ball, ball color=black] \pos circle (.4);

\foreach \pos in {(8,12*\rr),(21, 5*\rr),(18, 12*\rr),(5,4*\rr)}
\shade[shading=ball, ball color=red] \pos circle (.5);

\foreach \pos in {(1.5, \rr), (13, 24*\rr),(24.5, \rr), (13,8*\rr)}
 \shade[shading=ball, ball color=green] \pos circle (.5);

\foreach \pos in {(1.5, \rr), (13, 24*\rr),(15,8*\rr)}
 \shade[shading=ball, ball color=blue] \pos circle (.3);
\foreach \pos in {(15,8*\rr)}
 \shade[shading=ball, ball color=blue] \pos circle (.5);

\end{tikzpicture}
\caption{Double, triple and quadruple $u^{-2}$-insertions for $D^2=169$  on $\bbA_2$.\\
{\footnotesize{Examples of double $u^{-2}$-insertions repelling 4 vertices of a
$D$-rhombus are shown in red. There are 113 double insertions in an (8, 7)-rhombus 
and 86 in a (13, 0)-rhombus.
Examples of triple $u^{-2}$-insertions are shown in blue: they 
repel 5 sites on the boundary of a trapeze. Quadruple $u^{-2}$-insertions (green) 
repel 6 sites on the boundary of a $2D$-triangle.
The number of triple insertions is 61 in an (8, 7)- and 20 in a (13, 0)-triangle,
while the number of quadruple insertions is 39 in an (8, 7)-PGS and 3 in an (13 ,0)-PGS.}}
}
%\caption{ Examples of double, triple and quadruple defects of statistical weight $u^{-2}$
% for $D^2=169$.
% Pairs of inserted sites (black balls) remove four repelled sites (the vertices of a
% $D$-rhombus). Triples of inserted sites remove five repelled sites (white balls
% forming a trapezoid). Quadruples of inserted sites remove six repelled sites (white balls
% forming a $2D$-triangle).}
\end{figure}}  %FigureR18 %\fgri %\label{Fig18}%% 

%%%%%%%%%%%%%%%%%%%%%%%%%%%%%%%%%%%%%%%%%%%
%%%%%%%%%%%%%%%%%%%%%%%%%%%%

\def\FigureS19  %%\label{Fig19} %%\fgrj
{\begin{figure}[H]\label{Fig19}
\centering
\captionsetup{width=0.8\textwidth} %, labelfont=bf}

(a) \begin{tikzpicture}[scale=0.27]
\clip (-1.2, -1.2*\rr) rectangle (12, 22.5*\rr);

\path [fill=lightgray, draw=black] (0, 14*\rr) arc (210:270:14*\rr) arc (30:90:14*\rr);
 
\draw [line width=0.7mm] (0, 0) -- (0, 14*\rr) -- (10.5, 21*\rr) -- (10.5, 7*\rr) -- cycle
(0, 14*\rr) -- (10.5, 7*\rr) -- cycle;

\draw[yscale=sqrt(3/4), xslant=0.5] (-19,-2) grid (90, 60);
\draw[yscale=sqrt(3/4), xslant=-0.5] (-1,-2) grid (90, 60);

\foreach \pos in {(0, 0), (10.5, 7*\rr),(0, 14*\rr),(10.5, 21*\rr)}
 \shade[shading=ball, ball color=black] \pos circle (.4);

\foreach \pos in {(0.5,1*\rr),(1,2*\rr),(1.5,3*\rr), % 1+1+1
(2,4*\rr),(3,4*\rr),  % 2 
(1.5,5*\rr),(2.5,5*\rr),(3.5,5*\rr),  % 3
(2,6*\rr),(3,6*\rr),(4,6*\rr),(5,6*\rr),(6,6*\rr),  % 5
(2.5,7*\rr),(3.5,7*\rr),(4.5,7*\rr),(5.5,7*\rr),(6.5,7*\rr), 
(7.5,7*\rr),(8.5,7*\rr),(9.5,7*\rr),  %  8
(2,8*\rr),(3,8*\rr),(4,8*\rr),(5,8*\rr),(6,8*\rr), %  5
(1.5,9*\rr),(2.5,9*\rr),(3.5,9*\rr),  %  3 
(2,10*\rr),(3,10*\rr),  %  2 
(1.5,11*\rr),(1,12*\rr),(0.5,13*\rr),  %  1+1+1
(10,8*\rr),(9.5,9*\rr),(9,10*\rr),  %  1+1+1
(7.5,11*\rr),(8.5,11*\rr),  %  2 
(7,12*\rr),(8,12*\rr),(9,12*\rr),  %  3
(4.5,13*\rr),(5.5,13*\rr),(6.5,13*\rr),(7.5,13*\rr),(8.5,13*\rr), %  5
(1,14*\rr),(2,14*\rr),(3,14*\rr),(4,14*\rr),(5,14*\rr),(6,14*\rr),(7,14*\rr),(8,14*\rr), %  8
(4.5,15*\rr),(5.5,15*\rr),(6.5,15*\rr),(7.5,15*\rr),(8.5,15*\rr),  %  5
(7,16*\rr),(8,16*\rr),(9,16*\rr),  %  3
(7.5,17*\rr),(8.5,17*\rr), %  2 
(9,18*\rr),(9.5,19*\rr),(10,20*\rr)} % 1+1+1
\shade[shading=ball, ball color=orange] \pos circle (.25);
%  Total: 2(1+1+1+2+3+5+8+5+3+2+1+1+1)=2x34=68
\path [draw=black,line width=0.4mm]
(0, 0) arc (-30:30:14*\rr) arc (210:270:14*\rr) arc (90:150:14*\rr)
(10.5, 7*\rr) arc (30:90:14*\rr) arc (270:330:14*\rr) arc (150:210:14*\rr);

\end{tikzpicture}\quad (b) \begin{tikzpicture}[scale=0.27]
\clip (-1.2, -1.2*\rr) rectangle (18, 17*\rr);

\path [fill=lightgray, draw=black] (12, 2*\rr) arc (248.21:188.21:14*\rr) arc (68.21:8.21:14*\rr);

\path [draw=black, line width=0.6mm]
(0,0) -- (4.5,13*\rr) -- (16.5, 15*\rr) -- (12, 2*\rr) -- cycle
(12, 2*\rr) -- (4.5,13*\rr) -- cycle;

\draw[yscale=sqrt(3/4), xslant=0.5] (-90,-60) grid (90, 60);
\draw[yscale=sqrt(3/4), xslant=-0.5] (-90,-60) grid (90, 60);

\foreach \pos in {(0,0), (4.5,13*\rr),(12, 2*\rr),(16.5, 15*\rr)}
\shade[shading=ball, ball color=black] \pos circle (.4);

\foreach \pos in {(2, 2*\rr),(3, 2*\rr),   %  2
(2.5, 3*\rr),(3.5, 3*\rr),(4.5, 3*\rr),(5.5, 3*\rr),(6.5, 3*\rr),(7.5, 3*\rr),(8.5, 3*\rr),(9.5, 3*\rr), % 8
(3,4*\rr),(4,4*\rr),(5,4*\rr),(6,4*\rr),(7,4*\rr),(8,4*\rr),(9,4*\rr), %  7
(3.5,5*\rr),(4.5,5*\rr),(5.5,5*\rr),(6.5,5*\rr),(7.5,5*\rr), %  5
(4,6*\rr),(5,6*\rr),(6,6*\rr),(7,6*\rr), %  4
(4.5,7*\rr),(5.5,7*\rr),(6.5,7*\rr), %  3
(5,8*\rr),(6,8*\rr),  %  2 
(4.5,9*\rr),(5.5,9*\rr),(5,10*\rr),  %  2+1
(11.5,5*\rr),(11,6*\rr),  (12,6*\rr),  %  2+1
(10.5,7*\rr),(11.5,7*\rr), %  2
(10,8*\rr),(11,8*\rr),(12,8*\rr),  %  3
(9.5,9*\rr),(10.5,9*\rr),(11.5,9*\rr),(12.5,9*\rr),  %  4
(9,10*\rr),(10,10*\rr),(11,10*\rr),(12,10*\rr),(13,10*\rr),  %  5
(7.5,11*\rr),(8.5,11*\rr),(9.5,11*\rr),(10.5,11*\rr),(11.5,11*\rr),(12.5,11*\rr),(13.5,11*\rr), %  7
(7,12*\rr),(8,12*\rr),(9,12*\rr),(10,12*\rr),(11,12*\rr),(12,12*\rr),(13,12*\rr),(14,12*\rr),  %  8
(13.5, 13*\rr),(14.5, 13*\rr)}  %  2
 \shade[shading=ball, ball color=orange] \pos circle (.25);
%  Total: 2(2+8+7+5+4+3+2+2+1) =  2x34=68

\path [draw=black,line width=0.4mm]
(0,0) arc (128.21:68.21:14*\rr) arc (248.21:188.21:14*\rr) arc (368.21:308.21:14*\rr)
(12, 2*\rr) arc (188.21:128.21:14*\rr) arc (308.21:248.21:14*\rr) arc (68.21:8.21:14*\rr);

\end{tikzpicture} 

\caption{Single $u^{-2}$-insertions for $D^2=147$ on $\bbA_2$: (a) in a (7, 7)-PGS and 
(b) in an (11 ,2)-PGS.\\
{\footnotesize{As before, circular triangles  and orange balls mark the positions 
where an inserted particle repels 3 vertices  of a $D$-triangle in a PGS. In both PGS types 
the number of single $u^{-2}$-insertions 
equals 34 per a $D$-triangle or 68 per a $D$-rhombus.}}
}
%The structure of pair defects for a distinguishing contour of weight $u^{-2}$ used
%in the proof of Theorem~6,
%with $D^2=147$. The contour exhibits 4 removed and 2 added
%particles, affecting the vertices of inclined $D$-rhombus. Both $(11,2)$ (left) and $(7,7)$ (right)
%ground states are shown: white circles indicate the removed ground-state particles. The
%origin ${\mathbf 0}=(0,0)$ is at the position of the bottom left white ball. The
%light-gray areas (closed circular triangles) cover the lattice sites which repel 3 out of 4
%white balls but do not repel any other ground-state particles. The gray area
%(an open bi-convex lense) includes the sites that repel all 4 white balls but no other ground-state
%particles. Each black circle from one triangular stack, together
%with at least one black circle from the opposite stack, can be selected to form
%a pair defect. There are $51$ pairs that can be formed for the $(11,2)$ $D$-sub-lattice
%and $86$ pairs for the vertical $(7,7)$ one (out of .$144$ and $441$ nominal pairs, respectively).
%As before, a larger black ball can be paired with any ball from the opposite stack.}

\end{figure}} %Fig19 $FigureS19

%%%%%%%%%%%%%%%%%%%%%%%%%%%%%%%%%%%%%%%%%%%
%%%%%%%%%%%%%%%%%%%%%%%%%%%%

\def\FigureT20 %%\label{Fig20} %%\fgrk
{\begin{figure}[H]\label{Fig20}
\centering
\captionsetup{width=0.8\textwidth}  %, labelfont=bf}

(a)\begin{tikzpicture}[scale=0.22]
\clip (-1.2, -1.2*\rr) rectangle (23.2, 30.2*\rr);

\draw[yscale=sqrt(3/4), xslant=0.5] (-19,-2) grid (90, 60);
\draw[yscale=sqrt(3/4), xslant=-0.5] (-1,-2) grid (90, 60);

\draw [line width=0.6mm, color=red] (3,4*\rr)--(7.5,17*\rr);
\draw [line width=0.6mm, color=red] (2,14*\rr)--(14.5,13*\rr);
\foreach \pos in {(3,4*\rr),(7.5,17*\rr)}
\shade[shading=ball, ball color=red] \pos circle (.45);
\foreach \pos in {(2,14*\rr),(14.5,13*\rr)}
\shade[shading=ball, ball color=red] \pos circle (.45);

\draw [line width=0.6mm, color=blue] (0.5, 1*\rr)--(7.5,15*\rr)--(0.5, 27*\rr);

\draw [line width=0.6mm, color=green] (0.5,1*\rr)--(6,14*\rr)--(0.5,27*\rr);
\draw [line width=0.8mm, color=green] (6,14*\rr)--(19,14*\rr);

\draw [line width=0.6mm] (0,0) -- (0, 28*\rr) -- (21, 14*\rr) -- cycle
(0, 14*\rr) -- (10.5, 7*\rr) -- (10.5, 21*\rr) -- cycle;

\foreach \pos in {(0, 0), (10.5, 7*\rr),(21, 14*\rr),(0, 14*\rr),(0, 28*\rr),(10.5, 21*\rr)}
 \shade[shading=ball, ball color=black] \pos circle (.5);

 \foreach \pos in {(0.5, 1*\rr),(6, 14*\rr),(0.5, 27*\rr),(19, 14*\rr)}
 \shade[shading=ball, ball color=green] \pos circle (.45);

\foreach \pos in {(0.5, 1*\rr),(0.5, 27*\rr)}
 \shade[shading=ball, ball color=blue] \pos circle (.3);
\foreach \pos in {(7.5,15*\rr)}
 \shade[shading=ball, ball color=blue] \pos circle (.45); 
\end{tikzpicture}\quad (b) \begin{tikzpicture}[scale=0.22]
\clip (-1.2, -1.2*\rr) rectangle (26.2, 27.2*\rr);

\draw[yscale=sqrt(3/4), xslant=0.5] (-90,-60) grid (90, 60);
\draw[yscale=sqrt(3/4), xslant=-0.5] (-90,-60) grid (90, 60);

\path [draw=black, line width=0.6mm, color=red] (9.5,21*\rr)--(12,6*\rr);
\path [draw=black, line width=0.6mm, color=red] (3,4*\rr)--(14.5,13*\rr);
\path [draw=black, line width=0.6mm, color=blue] (2,2*\rr)--(11,12*\rr)--(21.5,5*\rr);
%\path [draw=black, line width=0.6mm, color=blue] (9.5,23*\rr)--(12.5,9*\rr)--(2,2*\rr);

\path [draw=black, line width=0.6mm] (0,0) -- (9, 26*\rr) -- (24, 4*\rr) -- cycle
(12, 2*\rr) -- (4.5,13*\rr) -- (16.5, 15*\rr) -- cycle;

\foreach \pos in {(0,0),(4.5,13*\rr),(9, 26*\rr),(12, 2*\rr),(24, 4*\rr),(16.5, 15*\rr)}
 \shade[shading=ball, ball color=black] \pos circle (.5);

\foreach \pos in {(3,4*\rr),(14.5,13*\rr)}
 \shade[shading=ball, ball color=red] \pos circle (.45);
\foreach \pos in {(9.5,21*\rr),(12,6*\rr)}
\shade[shading=ball, ball color=red] \pos circle (.45);

\foreach \pos in {(2,2*\rr),(21.5,5*\rr),(11,12*\rr)}
\shade[shading=ball, ball color=blue] \pos circle (.5);

%\foreach \pos in {(9.5,23*\rr),(12.5,9*\rr)}
%\shade[shading=ball, ball color=blue] \pos circle (.45);
 
\end{tikzpicture}

\caption{Double, triple and quadruple $u^{-2}$-insertions for $D^2=147$ on $\bbA_2$, in a
(7, 7)-PGS (a) and an (11, 2)-PGS (b).\\
{\footnotesize{Double $u^{-2}$-insertions (red) again repel 4 vertices
of a rhombus; triple insertions (blue) repel 5 sites on the boundary of a
trapeze; quadruple insertions (green) repel 6 sites on the boundary of a $2D$-triangle. 
There are 86  double insertions in (7, 7)- and 51 in an (11, 2)-rhombus. Next, we have 
39 triple insertions in a (7, 7)- and 1 in an (11, 2)-triangle. Finally, the number 
of quadruple insertions is 1 in a (7, 7)- PGS and 0 in an (11, 2)-PGS.}}
}
%Examples of double, triple and quadruple defects of statistical weight $u^{-2}$
%for $D^2=147$. As before,
%pairs of inserted sites (black balls) remove four repelled sites (the vertices of a
%$D$-rhombus). Triples of inserted sites remove five repelled sites (white balls
%forming a trapezoid). A quadruple of inserted sites removes six repelled sites from the
%vertical $(7,7)$ ground state (white balls forming a $2D$-triangle); for the $(11,2)$ such
%a quadruple is impossible }

\end{figure}} %{Fig20} %FigureT20

%%%%%%%%%%%%%%%%%%%%%%%%%%%%%%%%%%%%%%%%%%%%
%%%%%%%%%%%%%%%%%%%%%%%%%%%

\def\FigureU21 %\label{Fig21} % {D2=147Stfraight}
{\begin{figure}\label{Fig21}
\centering
\captionsetup{width=0.8\textwidth}  %, labelfont=bf}
(a)
\vskip .2cm
%\end{center}%\vskip .5cm
\caption{Admissible $u^{-2}$-insertions in a vertical PGS for $D^2=147$ on $\bbH_2$.\\
{\footnotesize{(a) Single insertions are again marked by orange balls: they repel 3 
vertices of the covering $D$-triangle. 
(b)  A double insertion is marked by a red bar; here we remove 4 vertices of the
covering $D$-rhombus. Triple insertions are represented by a pair of blue bars and a triple of
green bars. Such insertions remove 5 vertices in a PGS,  on the boundary of the covering
trapeze. Quadruple insertions repel 6 vertices of a  $2D$-triangle (triangles $AFD$ and $GBE$), 
with a single insertion in each of $4$ involved $D$-triangles. Such an insertion is represented 
by a triple of green bars. As in the previous examples, there is no other admissible $u^{-2}$-insertion.}}
}
\end{figure}} %%{Fig21} %FigureU21 {D2=147Stfraight}

%%%%%%%%%%%%%%%%%%%%%%%%%%%%%%%%%%%%%%%%%%%%
%%%%%%%%%%%%%%%%%%%%%%%%%%%

\def\FigureV22 %%\label{Fig22} {D2=147Inclined}
{\begin{figure}\label{Fig22} 
\centering 
\captionsetup{width=0.8\textwidth}  %, labelfont=bf}
(a)
%\end{center}
\caption{Admissible $u^{-2}$-insertions in an inclined PGS for  $D^2=147$ on $\bbH_2$.\\
{\footnotesize{The marking here is as in the previous figure. (a) As in Figure 21 (a), 
single insertions with a triple repulsion are marked by orange balls: they repel 3 
vertices of the covering $D$-triangle. (b) Again, double insertions 
repel 4 vertices only when the latter ones are the vertices of the covering $D$-rhombus. Such 
insertions are marked by red bars. Triple insertions are marked by pairs of blue bars joining 
triples of occupied sites.}}
}
\end{figure}} %Fig22 FigureV22 {D2=147Inclined}

%%%%%%%%%%%%%%%%%%%%%%%%%%%%%%%%%%%%%%%%%
%%%%%%%%%%%%%%%%%%%%%%%%%%%%%%

%%%%%%%%%%%%%%%%%%%%%%%%%%%%%%%%%%%%%%%%%%%%
%%%%%%%%%%%%%%%%%%%%%%%%%%%

\def\FigureW23 %\label{Fig23} % {D2=139Stfraight}
{\begin{figure}\label{Fig23}
\centering 
\captionsetup{width=0.8\textwidth}  %, labelfont=bf}
(a)

%\vskip .2cm
%\end{center}%\vskip .5cm
\caption{Admissible $u^{-2}$-insertions in a vertical PGS for $D^2=139$ on $\bbH_2$, with 
$(D^*)^2=147$.\\
{\footnotesize{Frame (a): dotted orange balls indicate the additional single 
insertions allowed for $D^2=139$ compared to $D^2=147$ in Figure 18. Frame 
(b) shows examples of double (red), triple (blue) and quadruple (green) admissible 
$u^{-2}$-insertions.}}
}
%\centerline{\bf for a straight PGS. The total amount is $25+21=46$}
%\centerline{\bf per a fundamental parallelogram}
%\vskip .2cm

\end{figure}} %FigureW23 %Fig23 % {D2=139Stfraight}

%%%%%%%%%%%%%%%%%%%%%%%%%%%%%%%%%%%%%%%%%%%%
%%%%%%%%%%%%%%%%%%%%%%%%%%%

\def\FigureX24 %%\label{Fig24} {D2=139Inclined}
{\begin{figure}\label{Fig24} 
\centering 
\captionsetup{width=0.8\textwidth}  %, labelfont=bf}

(a)

%\end{center}
\caption{Admissible $u^{-2}$-insertions in an inclined PGS for $D^2=139$ on $\bbH_2$, with 
$(D^*)^2=147$.\\
{\footnotesize{As before, dotted orange balls in frame (a) indicate the additional 
single insertions compared to $D^2=147$ in Figure 19. Frame (b) again shows 
examples of double, triple and quadruple admissible insertions (red, blue, and green).}}
}
\end{figure}} %Fig24 FigureX24  {D2=139Inclined}

%%%%%%%%%%%%%%%%%%%%%%%%%%%%%%%%%%%%%%%%%%

\def\Fig15 %\label{D2=67Triangle75}
{\begin{figure}\begin{center} 
%
\caption{A deletable vertex of type $\pi/3$ (a large black  ball), in a horizontal PGS (a) 
and in an inclined PGS (b), for $D^2=49$. The white balls show sites of insertion.
The same meaning is assigned in Figure 25.\\
%\indent\hskip .6cm{\footnotesize{Thick lines in Figures 26--27 indicate the boundary $\partial\Delta$ 
%of polygon $\Delta$, and white}}\\ 
%\indent\hskip .6cm{\footnotesize{balls mark inserted sites.}}\\
{\footnotesize{Medium-size black balls in Figures 23, 24 mark the positions
of removed sites of type $\pi$ in $\partial\Delta$.}}
}
%A deletable vertex of the $\pi/3$ type (a large white ball), for horizontal
% (left) and inclined (right) sub-lattices. The thick lines in Figures 7--9 mark the
% boundary $\partial\Delta$ of polygon $\Delta$.}
\end{figure}} %%Fig{25} FigureY26

%%%%%%%%%%%%%%%%%%%%%%%%%%%%%%%%%%%%%%%%%%%%
%%%%%%%%%%%%%%%%%%%%%%%%%%%

\def\FigureZ27 %%\label{Fig26} %\fgrab
{\begin{figure}[H]\label{Fig26}
\centering
\captionsetup{width=0.8\textwidth} %, labelfont=bf}

(a)\;\begin{tikzpicture}[scale=0.35]
\clip (-0.9, -1.3) rectangle (11.4, 8.1);

% \path [fill=lightgray, draw=black] (0, 0) arc (180:60:7)  arc (360:240:7);
% \path [fill=lightgray, draw=black] (11.5, 5*\rr) arc (141.79:21.79:7)  arc (321.79:201.79:7);

\path [fill=lightgray, draw=black]
(0, 0) -- (3.5, 7*\rr) -- (10.5, 7*\rr) arc (360:240:7);
\path [fill=gray, draw=black]
(0, 0) -- (3.5, 7*\rr) -- (10.5, 7*\rr) -- (7,0) -- (0,0);

\draw [line width=1mm]
(-0.5, -\rr) -- (3.5, 7*\rr) -- (11.5, 7*\rr);

\draw[yscale=sqrt(3/4), xslant=0.5] (-5,-2) grid (25, 10);
\draw[yscale=sqrt(3/4), xslant=-0.5] (-1,-2) grid (29, 10);

\foreach \pos in {(3.5, 7*\rr)}
\shade[shading=ball, ball color=black] \pos circle (.4);

\foreach \pos in {(0, 0), (10.5, 7*\rr)}
\shade[shading=ball, ball color=black] \pos circle (.3);

\foreach \pos in {(7, 0)}
 \shade[shading=ball, ball color=black] \pos circle (.2);

\foreach \pos in {(9.5, 3*\rr)}
 \shade[shading=ball, ball color=white] \pos circle (.4);

\end{tikzpicture}\quad (b)\;\begin{tikzpicture}[scale=0.35]
\clip (-0.9, -1.3) rectangle (13.4, 8.1);

% \path [fill=lightgray, draw=black] (0, 0) arc (180:60:7)  arc (360:240:7);
% \path [fill=lightgray, draw=black] (11.5, 5*\rr) arc (141.79:21.79:7)  arc (321.79:201.79:7);

\path [fill=lightgray, draw=black]
(0.5, 5*\rr) -- (7, 8*\rr) -- (12.5, 3*\rr) arc (321.79:201.79:7);
\path [fill=gray, draw=black]
(0.5, 5*\rr) -- (7, 8*\rr) -- (12.5, 3*\rr) -- (6,0) -- (0.5, 5*\rr);

\draw [line width=1mm] (-6, 2*\rr) -- (7, 8*\rr) -- (18, -2*\rr);

\draw[yscale=sqrt(3/4), xslant=0.5] (-5,-2) grid (25, 10);
\draw[yscale=sqrt(3/4), xslant=-0.5] (-1,-2) grid (29, 10);

\foreach \pos in {(7, 8*\rr)}
 \shade[shading=ball, ball color=black] \pos circle (.4);

\foreach \pos in {(0.5, 5*\rr),(12.5, 3*\rr)}
 \shade[shading=ball, ball color=black] \pos circle (.3);

\foreach \pos in {(6, 0)}
 \shade[shading=ball, ball color=black] \pos circle (.2);

\foreach \pos in {(9.5, 3*\rr)}
 \shade[shading=ball, ball color=white] \pos circle (.4);
\end{tikzpicture}

\caption{A deletable vertex of type $2\pi/3$ (a large black ball),
in a horizontal PGS (a) and in an inclined PGS (b), for $D^2=49$.} 

%A deletable vertex of the $2\pi/3$ type (a large white ball) for horizontal
%(left) and inclined (right) sub-lattices. The medium-size white balls mark the positions
%of type $\pi$ repelled sites in $\partial\Delta$.}
\end{figure}}  %% Fig26 %FigureZ27

%%%%%%%%%%%%%%%%%%%%%%%%%%%%%%%%%%%%%%%%%%%%
%%%%%%%%%%%%%%%%%%%%%%%%%%%
\def\Figurea1{\begin{figure} %sliding/cost-free demarcation: D^2=4,7
\centering
\captionsetup{width=0.8\textwidth} %, labelfont=bf}

(a) 
%\end{center}%\vskip .5cm
\caption{PGSs for $D^2=4$ (frame (a)) and $D^2=7$ (frame (b)).\\
{\footnotesize{Dotted gray lines mark the `cost-free' boundaries between PGSs.}}\\
}
%\vskip .2cm

\end{figure}} %\Figurea1
%%%%%%%%%%%%%%%%%%%%%%%%%%%%%%%%%%%%%%%%%
%%%%%%%%%%%%%%%%%%%%%%%%%%%%%%%%%%%%%%%%%

\def\Figureb2{\begin{figure} %sliding/cost-free demarcation: D^2=133
%\begin{center} 
\centering
\captionsetup{width=0.8\textwidth} %, labelfont=bf}
(a)%
 

%\end{center}%\vskip .5cm
\caption{PGSs for $D^2=31$ (frame (a)) and $D^2=133$ (frame (b)).\\
{\footnotesize{As before, dotted gray lines mark the `cost-free' boundaries
between PGSs.}}
}
%\vskip .2cm
\end{figure}} %\Figureb2 

%%%%%%%%%%%%%%%%%%%%%%%%%%%%%%%%%%%%%%%%%
%%%%%%%%%%%%%%%%%%%%%%%%%%%%%%

\def\WFigure25   %%\label{Fig24A} \fgraa
{\begin{figure}\label{Fig24A}
\centering 
\captionsetup{width=0.8\textwidth} 
%\captionsetup{width=0.8\textwidth, labelfont=bf}

\begin{tikzpicture}[scale=0.28]

\clip (-0.9, -8) rectangle (39.4, 7.1);

\draw[gray, line width=2.5] (121/14, 43/7*\rr)--(3.5, 7*\rr);
\draw[gray, line width=2.5] (121/14, 43/7*\rr)--(10.5, 7*\rr);
\draw[gray, line width=2.5] (121/14, 43/7*\rr)--(7, 0);

\draw[gray, line width=2.5] (2, 2*\rr)--(3.5, 7*\rr);
\draw[gray, line width=2.5] (2, 2*\rr)--(0, 0);
\draw[gray, line width=2.5] (2, 2*\rr)--(7, 0);

\draw[gray, line width=2.5] (51/14, -41/7*\rr)--(3.5, -7*\rr);
\draw[gray, line width=2.5] (51/14, -41/7*\rr)--(10.5, -7*\rr);
\draw[gray, line width=2.5] (51/14, -41/7*\rr)--(7, 0);
\draw[gray, line width=2.5] (51/14, -41/7*\rr)--(0, 0);

\draw[gray, line width=2.5] (219/14, 43/7*\rr)--(10.5, 7*\rr);
\draw[gray, line width=2.5] (219/14, 43/7*\rr)--(17.5, 7*\rr);
\draw[gray, line width=2.5] (219/14, 43/7*\rr)--(14, 0);

\path [draw=black, thick] (7, 0) -- (10.5, 7*\rr) -- (17.5, 7*\rr) -- (14, 0) -- (7, 0);
\path [draw=black, thick] (28, 0) -- (35, 0) -- (38.5, -7*\rr) -- (31.5, -7*\rr) -- (28, 0);

\path [draw=black, thick] (0,0) -- (7, 0) -- (10.5, -7*\rr) -- (3.5, -7*\rr) -- (0, 0)
(17.5, 7*\rr) -- (24.5, 7*\rr) -- (28, 0) -- (21, 0) -- (17.5, 7*\rr);

\path [draw=black, thick] (0, 0) -- (3.5, 7*\rr) -- (10.5, 7*\rr) -- (7, 0) -- (0, 0)
(24.5, 7*\rr) -- (31.5, 7*\rr) -- (35, 0) -- (28, 0) -- (24.5, 7*\rr);

\draw [line width=1mm] (12.5, -7*\rr) -- (3.5, -7*\rr) -- (0,0) -- (3.5, 7*\rr) -- (31.5, 7*\rr) -- (39, -8*\rr) ;

\path [draw=black, thick] (17.5, 7*\rr) -- (17.5, 7*\rr) -- (14, 0)  -- (0, 0);
\path [draw=black, thick] (7,0) -- (10.5, 7*\rr) -- (17.5, 7*\rr) -- (14, 0);

\begin{scope}[rotate=21.79]
\draw[yscale=sqrt(3/4), xslant=0.5] (-8,-26) grid (46, 9);
\draw[yscale=sqrt(3/4), xslant=-0.5] (-8,-26) grid (46, 9);
\end{scope}

\foreach \pos in {(0, 0), (7, 0),(14, 0),(21,0),(28, 0), (35,0),
(3.5, 7*\rr),(10.5, 7*\rr),(17.5, 7*\rr), (24.5, 7*\rr),
(3.5, -7*\rr), (10.5, -7*\rr),(31.5, 7*\rr), (31.5, -7*\rr),(38.5, -7*\rr)}
 \shade[shading=ball, ball color=black] \pos circle (.35);

\foreach \pos in {(51/14, -41/7*\rr),(2, 2*\rr),(121/14, 43/7*\rr),(219/14, 43/7*\rr),(296/14, 8/7*\rr),(378/14, 42/7*\rr),
(233/7, 2/7*\rr),(515/14, -47/7*\rr)}
 \shade[shading=ball, ball color=lightgray] \pos circle (.45);

\end{tikzpicture}
\caption{A fragment of set $\Delta$, for an inclined PGS $\vphi$, for $D^2=169$.
\footnotesize{Black balls in Figures 23--25 mark removed sites $\bx\in\vphi$ while 
gray balls mark inserted sites $\by\in\phi$. Gray lines indicate which sites 
are removed by a given inserted site. Thick black lines indicate the boundary $\partial\Delta$.}
}
\end{figure}} %WFigure25   %%\label{Fig24A} \fgraa

\def\XFigure24A %\label{Fig24A} 
{\begin{figure}\label{Fig24A} 
%\begin{center}
(a)\begin{tikzpicture}[scale=0.18]
\clip (0, -1.5*\rr) rectangle (25, 21*\rr);

%\path [fill=lightgray] (1,-1*\rr)-- (4,-1*\rr)--(10.7,21*\rr)--(7.5,21*\rr);
%\path [fill=cyan] (4.2,-1*\rr)--(7.4,-1*\rr)--(14,21*\rr)--(10.8,21*\rr);
%\path [fill=lightgray] (7.5,-1*\rr)--(10.6,-1*\rr)--(17.3,21*\rr)--(14.2,21*\rr);

\draw [line width=0.2mm] (0.5,1*\rr) -- (1.5,1*\rr) -- (2,0) -- (1.5,-1*\rr)
-- (0.5,-1*\rr) -- (0,0) -- (0.5,1*\rr);
\draw [line width=0.2mm] (0.5,1*\rr) -- (0,2*\rr) -- (0.5,3*\rr) -- (1.5,3*\rr) 
-- (2,2*\rr) -- (1.5,1*\rr);
\draw [line width=0.2mm] (0.5,3*\rr) -- (0,4*\rr) -- (0.5,5*\rr) -- (1.5,5*\rr) 
-- (2,4*\rr) -- (1.5,3*\rr);
\draw [line width=0.2mm] (0.5,5*\rr) -- (0,6*\rr) -- (0.5,7*\rr) -- (1.5,7*\rr) 
-- (2,6*\rr) -- (1.5,5*\rr);
\draw [line width=0.2mm] (0.5,7*\rr) -- (0,8*\rr) -- (0.5,9*\rr) -- (1.5,9*\rr) 
-- (2,8*\rr) -- (1.5,7*\rr);
\draw [line width=0.2mm] (0.5,9*\rr) -- (0,10*\rr) -- (0.5,11*\rr) -- (1.5,11*\rr) 
-- (2,10*\rr) -- (1.5,9*\rr);
\draw [line width=0.2mm] (0.5,11*\rr) -- (0,12*\rr) -- (0.5,13*\rr) -- (1.5,13*\rr) 
-- (2,12*\rr) -- (1.5,11*\rr);
\draw [line width=0.2mm] (0.5,13*\rr) -- (0,14*\rr) -- (0.5,15*\rr) -- (1.5,15*\rr) 
-- (2,14*\rr) -- (1.5,13*\rr);
\draw [line width=0.2mm] (0.5,15*\rr) -- (0,16*\rr) -- (0.5,17*\rr) -- (1.5,17*\rr) 
-- (2,16*\rr) -- (1.5,15*\rr);
\draw [line width=0.2mm] (0.5,17*\rr) -- (0,18*\rr) -- (0.5,19*\rr) -- (1.5,19*\rr) 
-- (2,18*\rr) -- (1.5,17*\rr);
\draw [line width=0.2mm] (0.5,19*\rr) -- (0,20*\rr) -- (0.5,21*\rr) -- (1.5,21*\rr) 
-- (2,20*\rr) -- (1.5,19*\rr);
%\draw [line width=0.2mm] (0.5,21*\rr) -- (0,22*\rr) -- (0.5,22*\rr) -- (1.5,22*\rr) 
%-- (2,22*\rr) -- (1.5,21*\rr);

\draw [line width=0.2mm]  (3,0) -- (3.5,-1*\rr);
\draw [line width=0.2mm] (2,2*\rr) -- (3,2*\rr) -- (3.5,1*\rr) -- (3,0) -- (2,0);
\draw [line width=0.2mm] (2,4*\rr) -- (3,4*\rr) -- (3.5,3*\rr) -- (3,2*\rr);
\draw [line width=0.2mm] (2,6*\rr) -- (3,6*\rr) -- (3.5,5*\rr) -- (3,4*\rr);
\draw [line width=0.2mm] (2,8*\rr) -- (3,8*\rr) -- (3.5,7*\rr) -- (3,6*\rr);
\draw [line width=0.2mm] (2,10*\rr) -- (3,10*\rr) -- (3.5,9*\rr) -- (3,8*\rr);
\draw [line width=0.2mm] (2,12*\rr) -- (3,12*\rr) -- (3.5,11*\rr) -- (3,10*\rr);
\draw [line width=0.2mm] (2,14*\rr) -- (3,14*\rr) -- (3.5,13*\rr) -- (3,12*\rr);
\draw [line width=0.2mm] (2,16*\rr) -- (3,16*\rr) -- (3.5,15*\rr) -- (3,14*\rr);
\draw [line width=0.2mm] (2,18*\rr) -- (3,18*\rr) -- (3.5,17*\rr) -- (3,16*\rr);
\draw [line width=0.2mm] (2,20*\rr) -- (3,20*\rr) -- (3.5,19*\rr) -- (3,18*\rr);
\draw [line width=0.2mm] (3.5,21*\rr) -- (3,20*\rr);

\draw [line width=0.2mm]  (3.5,1*\rr) -- (4.5,1*\rr) -- (5,0) -- (4.5,-1*\rr)
--(3.5,-1*\rr);
\draw [line width=0.2mm]  (3.5,3*\rr) -- (4.5,3*\rr) -- (5,2*\rr) -- (4.5,1*\rr);
\draw [line width=0.2mm]  (3.5,5*\rr) -- (4.5,5*\rr) -- (5,4*\rr) -- (4.5,3*\rr);
\draw [line width=0.2mm]  (3.5,7*\rr) -- (4.5,7*\rr) -- (5,6*\rr) -- (4.5,5*\rr);
\draw [line width=0.2mm]  (3.5,9*\rr) -- (4.5,9*\rr) -- (5,8*\rr) -- (4.5,7*\rr);
\draw [line width=0.2mm]  (3.5,11*\rr) -- (4.5,11*\rr) -- (5,10*\rr) -- (4.5,9*\rr);
\draw [line width=0.2mm]  (3.5,13*\rr) -- (4.5,13*\rr) -- (5,12*\rr) -- (4.5,11*\rr);
\draw [line width=0.2mm]  (3.5,15*\rr) -- (4.5,15*\rr) -- (5,14*\rr) -- (4.5,13*\rr);
\draw [line width=0.2mm]  (3.5,17*\rr) -- (4.5,17*\rr) -- (5,16*\rr) -- (4.5,15*\rr);
\draw [line width=0.2mm]  (3.5,19*\rr) -- (4.5,19*\rr) -- (5,18*\rr) -- (4.5,17*\rr);
\draw [line width=0.2mm]  (3.5,21*\rr) -- (4.5,21*\rr) -- (5,20*\rr) -- (4.5,19*\rr);

\draw [line width=0.2mm]  (5,0) -- (6,0)--(6.5,-1*\rr );
\draw [line width=0.2mm]  (5,2*\rr) -- (6,2*\rr)--(6.5,1*\rr )--(6,0);
\draw [line width=0.2mm]  (5,4*\rr) -- (6,4*\rr)--(6.5,3*\rr )--(6,2*\rr);
\draw [line width=0.2mm]  (5,6*\rr) -- (6,6*\rr)--(6.5,5*\rr )--(6,4*\rr);
\draw [line width=0.2mm]  (5,8*\rr) -- (6,8*\rr)--(6.5,7*\rr )--(6,6*\rr);
\draw [line width=0.2mm]  (5,10*\rr) -- (6,10*\rr)--(6.5,9*\rr )--(6,8*\rr);
\draw [line width=0.2mm]  (5,12*\rr) -- (6,12*\rr)--(6.5,11*\rr )--(6,10*\rr);
\draw [line width=0.2mm]  (5,14*\rr) -- (6,14*\rr)--(6.5,13*\rr )--(6,12*\rr);
\draw [line width=0.2mm]  (5,16*\rr) -- (6,16*\rr)--(6.5,15*\rr )--(6,14*\rr);
\draw [line width=0.2mm]  (5,18*\rr) -- (6,18*\rr)--(6.5,17*\rr )--(6,16*\rr);
\draw [line width=0.2mm]  (5,20*\rr) -- (6,20*\rr)--(6.5,19*\rr )--(6,18*\rr);
\draw [line width=0.2mm]  (6.5,21*\rr )--(6,20*\rr);

\draw [line width=0.2mm]  (6.5,1*\rr) -- (7.5,1*\rr)--(8,0)--(7.5,-1*\rr)
--(6.5,-1*\rr);
\draw [line width=0.2mm]  (6.5,3*\rr) -- (7.5,3*\rr)--(8,2*\rr)--(7.5,1*\rr);
\draw [line width=0.2mm]  (6.5,5*\rr) -- (7.5,5*\rr)--(8,4*\rr)--(7.5,3*\rr);
\draw [line width=0.2mm]  (6.5,7*\rr) -- (7.5,7*\rr)--(8,6*\rr)--(7.5,5*\rr);
\draw [line width=0.2mm]  (6.5,9*\rr) -- (7.5,9*\rr)--(8,8*\rr)--(7.5,7*\rr);
\draw [line width=0.2mm]  (6.5,11*\rr) -- (7.5,11*\rr)--(8,10*\rr)--(7.5,9*\rr);
\draw [line width=0.2mm]  (6.5,13*\rr) -- (7.5,13*\rr)--(8,12*\rr)--(7.5,11*\rr);
\draw [line width=0.2mm]  (6.5,15*\rr) -- (7.5,15*\rr)--(8,14*\rr)--(7.5,13*\rr);
\draw [line width=0.2mm]  (6.5,17*\rr) -- (7.5,17*\rr)--(8,16*\rr)--(7.5,15*\rr);
\draw [line width=0.2mm]  (6.5,19*\rr) -- (7.5,19*\rr)--(8,18*\rr)--(7.5,17*\rr);
\draw [line width=0.2mm]  (6.5,21*\rr) -- (7.5,21*\rr)--(8,20*\rr)--(7.5,19*\rr);

\draw [line width=0.2mm]  (8,0)--(9,0)--(9.5,-1*\rr);
\draw [line width=0.2mm]  (8,2*\rr)--(9,2*\rr)--(9.5,1*\rr)--(9,0);
\draw [line width=0.2mm]  (8,4*\rr)--(9,4*\rr)--(9.5,3*\rr)--(9,2*\rr);
\draw [line width=0.2mm]  (8,6*\rr)--(9,6*\rr)--(9.5,5*\rr)--(9,4*\rr);
\draw [line width=0.2mm]  (8,8*\rr)--(9,8*\rr)--(9.5,7*\rr)--(9,6*\rr);
\draw [line width=0.2mm]  (8,10*\rr)--(9,10*\rr)--(9.5,9*\rr)--(9,8*\rr);
\draw [line width=0.2mm]  (8,12*\rr)--(9,12*\rr)--(9.5,11*\rr)--(9,10*\rr);
\draw [line width=0.2mm]  (8,14*\rr)--(9,14*\rr)--(9.5,13*\rr)--(9,12*\rr);
\draw [line width=0.2mm]  (8,16*\rr)--(9,16*\rr)--(9.5,15*\rr)--(9,14*\rr);
\draw [line width=0.2mm]  (8,18*\rr)--(9,18*\rr)--(9.5,17*\rr)--(9,16*\rr);
\draw [line width=0.2mm]  (8,20*\rr)--(9,20*\rr)--(9.5,19*\rr)--(9,18*\rr);
\draw [line width=0.2mm]  (9.5,21*\rr)--(9,20*\rr);

\draw [line width=0.2mm]  (9.5,1*\rr)--(10.5,1*\rr)--(11,0)--(10.5,-1*\rr)
--(9.5,-1*\rr);
\draw [line width=0.2mm]  (9.5,3*\rr)--(10.5,3*\rr)--(11,2*\rr)--(10.5,1*\rr);
\draw [line width=0.2mm]  (9.5,5*\rr)--(10.5,5*\rr)--(11,4*\rr)--(10.5,3*\rr);
\draw [line width=0.2mm]  (9.5,7*\rr)--(10.5,7*\rr)--(11,6*\rr)--(10.5,5*\rr);
\draw [line width=0.2mm]  (9.5,9*\rr)--(10.5,9*\rr)--(11,8*\rr)--(10.5,7*\rr);
\draw [line width=0.2mm]  (9.5,11*\rr)--(10.5,11*\rr)--(11,10*\rr)--(10.5,9*\rr);
\draw [line width=0.2mm]  (9.5,13*\rr)--(10.5,13*\rr)--(11,12*\rr)--(10.5,11*\rr);
\draw [line width=0.2mm]  (9.5,15*\rr)--(10.5,15*\rr)--(11,14*\rr)--(10.5,13*\rr);
\draw [line width=0.2mm]  (9.5,17*\rr)--(10.5,17*\rr)--(11,16*\rr)--(10.5,15*\rr);
\draw [line width=0.2mm]  (9.5,19*\rr)--(10.5,19*\rr)--(11,18*\rr)--(10.5,17*\rr);
\draw [line width=0.2mm]  (9.5,21*\rr)--(10.5,21*\rr)--(11,20*\rr)--(10.5,19*\rr);

\draw [line width=0.2mm]  (11,0)--(12,0)--(12.5,-1*\rr);
\draw [line width=0.2mm]  (11,2*\rr)--(12,2*\rr)--(12.5,1*\rr)--(12,0);
\draw [line width=0.2mm]  (11,4*\rr)--(12,4*\rr)--(12.5,3*\rr)--(12,2*\rr);
\draw [line width=0.2mm]  (11,6*\rr)--(12,6*\rr)--(12.5,5*\rr)--(12,4*\rr);
\draw [line width=0.2mm]  (11,8*\rr)--(12,8*\rr)--(12.5,7*\rr)--(12,6*\rr);
\draw [line width=0.2mm]  (11,10*\rr)--(12,10*\rr)--(12.5,9*\rr)--(12,8*\rr);
\draw [line width=0.2mm]  (11,12*\rr)--(12,12*\rr)--(12.5,11*\rr)--(12,10*\rr);
\draw [line width=0.2mm]  (11,14*\rr)--(12,14*\rr)--(12.5,13*\rr)--(12,12*\rr);
\draw [line width=0.2mm]  (11,16*\rr)--(12,16*\rr)--(12.5,15*\rr)--(12,14*\rr);
\draw [line width=0.2mm]  (11,18*\rr)--(12,18*\rr)--(12.5,17*\rr)--(12,16*\rr);
\draw [line width=0.2mm]  (11,20*\rr)--(12,20*\rr)--(12.5,19*\rr)--(12,18*\rr);
\draw [line width=0.2mm]  (12,20*\rr)--(12.5,21*\rr);

\draw [line width=0.2mm]  (12.5,1*\rr)--(13.5,1*\rr)--(14,0)--(13.5,-1*\rr)
--(12.5,-1*\rr);
\draw [line width=0.2mm]  (12.5,3*\rr)--(13.5,3*\rr)--(14,2*\rr)--(13.5,1*\rr);
\draw [line width=0.2mm]  (12.5,5*\rr)--(13.5,5*\rr)--(14,4*\rr)--(13.5,3*\rr);
\draw [line width=0.2mm]  (12.5,7*\rr)--(13.5,7*\rr)--(14,6*\rr)--(13.5,5*\rr);
\draw [line width=0.2mm]  (12.5,9*\rr)--(13.5,9*\rr)--(14,8*\rr)--(13.5,7*\rr);
\draw [line width=0.2mm]  (12.5,11*\rr)--(13.5,11*\rr)--(14,10*\rr)--(13.5,9*\rr);
\draw [line width=0.2mm]  (12.5,13*\rr)--(13.5,13*\rr)--(14,12*\rr)--(13.5,11*\rr);
\draw [line width=0.2mm]  (12.5,15*\rr)--(13.5,15*\rr)--(14,14*\rr)--(13.5,13*\rr);
\draw [line width=0.2mm]  (12.5,17*\rr)--(13.5,17*\rr)--(14,16*\rr)--(13.5,15*\rr);
\draw [line width=0.2mm]  (12.5,19*\rr)--(13.5,19*\rr)--(14,18*\rr)--(13.5,17*\rr);
\draw [line width=0.2mm]  (12.5,21*\rr)--(13.5,21*\rr)--(14,20*\rr)--(13.5,19*\rr);

\draw [line width=0.2mm]  (14,0)--(15,0)--(15.5,-1*\rr);
\draw [line width=0.2mm]  (14,2*\rr)--(15,2*\rr)--(15.5,1*\rr)--(15,0);
\draw [line width=0.2mm]  (14,4*\rr)--(15,4*\rr)--(15.5,3*\rr)--(15,2*\rr);
\draw [line width=0.2mm]  (14,6*\rr)--(15,6*\rr)--(15.5,5*\rr)--(15,4*\rr);
\draw [line width=0.2mm]  (14,8*\rr)--(15,8*\rr)--(15.5,7*\rr)--(15,6*\rr);
\draw [line width=0.2mm]  (14,10*\rr)--(15,10*\rr)--(15.5,9*\rr)--(15,8*\rr);
\draw [line width=0.2mm]  (14,12*\rr)--(15,12*\rr)--(15.5,11*\rr)--(15,10*\rr);
\draw [line width=0.2mm]  (14,14*\rr)--(15,14*\rr)--(15.5,13*\rr)--(15,12*\rr);
\draw [line width=0.2mm]  (14,16*\rr)--(15,16*\rr)--(15.5,15*\rr)--(15,14*\rr);
\draw [line width=0.2mm]  (14,18*\rr)--(15,18*\rr)--(15.5,17*\rr)--(15,16*\rr);
\draw [line width=0.2mm]  (14,20*\rr)--(15,20*\rr)--(15.5,19*\rr)--(15,18*\rr);
\draw [line width=0.2mm]  (15,20*\rr)--(15.5,21*\rr);

\draw [line width=0.2mm]  (15.5,1*\rr)--(16.5,1*\rr)--(17,0)--(16.5,-1*\rr)
--(15.5,-1*\rr);
\draw [line width=0.2mm]  (15.5,3*\rr)--(16.5,3*\rr)--(17,2*\rr)--(16.5,1*\rr);
\draw [line width=0.2mm]  (15.5,5*\rr)--(16.5,5*\rr)--(17,4*\rr)--(16.5,3*\rr);
\draw [line width=0.2mm]  (15.5,7*\rr)--(16.5,7*\rr)--(17,6*\rr)--(16.5,5*\rr);
\draw [line width=0.2mm]  (15.5,9*\rr)--(16.5,9*\rr)--(17,8*\rr)--(16.5,7*\rr);
\draw [line width=0.2mm]  (15.5,11*\rr)--(16.5,11*\rr)--(17,10*\rr)--(16.5,9*\rr);
\draw [line width=0.2mm]  (15.5,13*\rr)--(16.5,13*\rr)--(17,12*\rr)--(16.5,11*\rr);
\draw [line width=0.2mm]  (15.5,15*\rr)--(16.5,15*\rr)--(17,14*\rr)--(16.5,13*\rr);
\draw [line width=0.2mm]  (15.5,17*\rr)--(16.5,17*\rr)--(17,16*\rr)--(16.5,15*\rr);
\draw [line width=0.2mm]  (15.5,19*\rr)--(16.5,19*\rr)--(17,18*\rr)--(16.5,17*\rr);
\draw [line width=0.2mm]  (15.5,21*\rr)--(16.5,21*\rr)--(17,20*\rr)--(16.5,19*\rr);

\draw [line width=0.2mm]  (17,0)--(18,0)--(18.5,-1*\rr);
\draw [line width=0.2mm]  (17,2*\rr)--(18,2*\rr)--(18.5,1*\rr)--(18,0);
\draw [line width=0.2mm]  (17,4*\rr)--(18,4*\rr)--(18.5,3*\rr)--(18,2*\rr);
\draw [line width=0.2mm]  (17,6*\rr)--(18,6*\rr)--(18.5,5*\rr)--(18,4*\rr);
\draw [line width=0.2mm]  (17,8*\rr)--(18,8*\rr)--(18.5,7*\rr)--(18,6*\rr);
\draw [line width=0.2mm]  (17,10*\rr)--(18,10*\rr)--(18.5,9*\rr)--(18,8*\rr);
\draw [line width=0.2mm]  (17,12*\rr)--(18,12*\rr)--(18.5,11*\rr)--(18,10*\rr);
\draw [line width=0.2mm]  (17,14*\rr)--(18,14*\rr)--(18.5,13*\rr)--(18,12*\rr);
\draw [line width=0.2mm]  (17,16*\rr)--(18,16*\rr)--(18.5,15*\rr)--(18,14*\rr);
\draw [line width=0.2mm]  (17,18*\rr)--(18,18*\rr)--(18.5,17*\rr)--(18,16*\rr);
\draw [line width=0.2mm]  (17,20*\rr)--(18,20*\rr)--(18.5,19*\rr)--(18,18*\rr);
\draw [line width=0.2mm]  (18,20*\rr)--(18.5,21*\rr);

\draw [line width=0.2mm]  (18.5,1*\rr)--(19.5,1*\rr)--(20,0)--(19.5,-1*\rr)
--(18.5,-1*\rr);
\draw [line width=0.2mm]  (18.5,3*\rr)--(19.5,3*\rr)--(20,2*\rr)--(19.5,1*\rr);
\draw [line width=0.2mm]  (18.5,5*\rr)--(19.5,5*\rr)--(20,4*\rr)--(19.5,3*\rr);
\draw [line width=0.2mm]  (18.5,7*\rr)--(19.5,7*\rr)--(20,6*\rr)--(19.5,5*\rr);
\draw [line width=0.2mm]  (18.5,9*\rr)--(19.5,9*\rr)--(20,8*\rr)--(19.5,7*\rr);
\draw [line width=0.2mm]  (18.5,11*\rr)--(19.5,11*\rr)--(20,10*\rr)--(19.5,9*\rr);
\draw [line width=0.2mm]  (18.5,13*\rr)--(19.5,13*\rr)--(20,12*\rr)--(19.5,11*\rr);
\draw [line width=0.2mm]  (18.5,15*\rr)--(19.5,15*\rr)--(20,14*\rr)--(19.5,13*\rr);
\draw [line width=0.2mm]  (18.5,17*\rr)--(19.5,17*\rr)--(20,16*\rr)--(19.5,15*\rr);
\draw [line width=0.2mm]  (18.5,19*\rr)--(19.5,19*\rr)--(20,18*\rr)--(19.5,17*\rr);
\draw [line width=0.2mm]  (18.5,21*\rr)--(19.5,21*\rr)--(20,20*\rr)--(19.5,19*\rr);

\draw [line width=0.2mm]  (20,0)--(21,0)--(21.5,-1*\rr);
\draw [line width=0.2mm]  (20,2*\rr)--(21,2*\rr)--(21.5,1*\rr)--(21,0);
\draw [line width=0.2mm]  (20,4*\rr)--(21,4*\rr)--(21.5,3*\rr)--(21,2*\rr);
\draw [line width=0.2mm]  (20,6*\rr)--(21,6*\rr)--(21.5,5*\rr)--(21,4*\rr);
\draw [line width=0.2mm]  (20,8*\rr)--(21,8*\rr)--(21.5,7*\rr)--(21,6*\rr);
\draw [line width=0.2mm]  (20,10*\rr)--(21,10*\rr)--(21.5,9*\rr)--(21,8*\rr);
\draw [line width=0.2mm]  (20,12*\rr)--(21,12*\rr)--(21.5,11*\rr)--(21,10*\rr);
\draw [line width=0.2mm]  (20,14*\rr)--(21,14*\rr)--(21.5,13*\rr)--(21,12*\rr);
\draw [line width=0.2mm]  (20,16*\rr)--(21,16*\rr)--(21.5,15*\rr)--(21,14*\rr);
\draw [line width=0.2mm]  (20,18*\rr)--(21,18*\rr)--(21.5,17*\rr)--(21,16*\rr);
\draw [line width=0.2mm]  (20,20*\rr)--(21,20*\rr)--(21.5,19*\rr)--(21,18*\rr);
\draw [line width=0.2mm]  (21,20*\rr)--(21.5,21*\rr);

\draw [line width=0.2mm]  (21.5,1*\rr)--(22.5,1*\rr)--(23,0)--(22.5,-1*\rr)
--(21.5,-1*\rr);
\draw [line width=0.2mm]  (21.5,3*\rr)--(22.5,3*\rr)--(23,2*\rr)--(22.5,1*\rr);
\draw [line width=0.2mm]  (21.5,5*\rr)--(22.5,5*\rr)--(23,4*\rr)--(22.5,3*\rr);
\draw [line width=0.2mm]  (21.5,7*\rr)--(22.5,7*\rr)--(23,6*\rr)--(22.5,5*\rr);
\draw [line width=0.2mm]  (21.5,9*\rr)--(22.5,9*\rr)--(23,8*\rr)--(22.5,7*\rr);
\draw [line width=0.2mm]  (21.5,11*\rr)--(22.5,11*\rr)--(23,10*\rr)--(22.5,9*\rr);
\draw [line width=0.2mm]  (21.5,13*\rr)--(22.5,13*\rr)--(23,12*\rr)--(22.5,11*\rr);
\draw [line width=0.2mm]  (21.5,15*\rr)--(22.5,15*\rr)--(23,14*\rr)--(22.5,13*\rr);
\draw [line width=0.2mm]  (21.5,17*\rr)--(22.5,17*\rr)--(23,16*\rr)--(22.5,15*\rr);
\draw [line width=0.2mm]  (21.5,19*\rr)--(22.5,19*\rr)--(23,18*\rr)--(22.5,17*\rr);
\draw [line width=0.2mm]  (21.5,21*\rr)--(22.5,21*\rr)--(23,20*\rr)--(22.5,19*\rr);

\draw [line width=0.2mm]  (23,0)--(24,0)--(24.5,-1*\rr);
\draw [line width=0.2mm]  (23,2*\rr)--(24,2*\rr)--(24.5,1*\rr)--(24,0);
\draw [line width=0.2mm]  (23,4*\rr)--(24,4*\rr)--(24.5,3*\rr)--(24,2*\rr);
\draw [line width=0.2mm]  (23,6*\rr)--(24,6*\rr)--(24.5,5*\rr)--(24,4*\rr);
\draw [line width=0.2mm]  (23,8*\rr)--(24,8*\rr)--(24.5,7*\rr)--(24,6*\rr);
\draw [line width=0.2mm]  (23,10*\rr)--(24,10*\rr)--(24.5,9*\rr)--(24,8*\rr);
\draw [line width=0.2mm]  (23,12*\rr)--(24,12*\rr)--(24.5,11*\rr)--(24,10*\rr);
\draw [line width=0.2mm]  (23,14*\rr)--(24,14*\rr)--(24.5,13*\rr)--(24,12*\rr);
\draw [line width=0.2mm]  (23,16*\rr)--(24,16*\rr)--(24.5,15*\rr)--(24,14*\rr);
\draw [line width=0.2mm]  (23,18*\rr)--(24,18*\rr)--(24.5,17*\rr)--(24,16*\rr);
\draw [line width=0.2mm]  (23,20*\rr)--(24,20*\rr)--(24.5,19*\rr)--(24,18*\rr);
\draw [line width=0.2mm]  (24,20*\rr)--(24.5,21*\rr);

\foreach \pos in {(13,18*\rr)}
\shade[shading=ball, ball color=red] \pos circle (.4);

\draw [line width=0.6mm] (12.5,19*\rr)--(7.5,13*\rr)--(14,10*\rr)--(12.5,19*\rr)
--(19.5,15*\rr)--(14,10*\rr);

\foreach \pos in {(7.5,13*\rr),(14,10*\rr),(19.5,15*\rr),(12.5,19*\rr)}
% (2,8*\rr),(0.5,17*\rr)}
\shade[shading=ball, ball color=black] \pos circle (.6);

% \foreach \pos in {(9,4*\rr),(15.5,1*\rr),(21,6*\rr),(3.5,-1*\rr)}
% \shade[shading=ball, ball color=black] \pos circle (.6);

%\foreach \pos in {(10.5, 7*\rr),(0, 14*\rr),(12.5, 17*\rr)}
 %\shade[shading=ball, ball color=white] \pos circle (.4);

\draw (14.5,8*\rr) node   {\Large {${\mbox{\boldmath$O$}}$}};
\draw (5.7,13.2*\rr) node   {\Large {${\mbox{\boldmath$A$}}$}};
\draw (10,19*\rr) node   {\Large {${\mbox{\boldmath$B$}}$}};
\draw (22,15.5*\rr) node   {\Large {${\mbox{\boldmath$C$}}$}};
%\draw (23,3*\rr) node   {\Large {${\mbox{\boldmath$D$}}$}};
%\draw (18.5,1*\rr) node   {\Large {${\mbox{\boldmath$E$}}$}};
 
\end{tikzpicture}B\begin{tikzpicture}[scale=0.18]
\clip (0, -1.5*\rr) rectangle (24.5, 21*\rr);

%\path [fill=lightgray] (1,-1*\rr)-- (4,-1*\rr)--(10.7,21*\rr)--(7.5,21*\rr);
%\path [fill=cyan] (4.2,-1*\rr)--(7.4,-1*\rr)--(14,21*\rr)--(10.8,21*\rr);
%\path [fill=lightgray] (7.5,-1*\rr)--(10.6,-1*\rr)--(17.3,21*\rr)--(14.2,21*\rr);

\draw [line width=0.2mm] (0.5,1*\rr) -- (1.5,1*\rr) -- (2,0) -- (1.5,-1*\rr)
-- (0.5,-1*\rr) -- (0,0) -- (0.5,1*\rr);
\draw [line width=0.2mm] (0.5,1*\rr) -- (0,2*\rr) -- (0.5,3*\rr) -- (1.5,3*\rr) 
-- (2,2*\rr) -- (1.5,1*\rr);
\draw [line width=0.2mm] (0.5,3*\rr) -- (0,4*\rr) -- (0.5,5*\rr) -- (1.5,5*\rr) 
-- (2,4*\rr) -- (1.5,3*\rr);
\draw [line width=0.2mm] (0.5,5*\rr) -- (0,6*\rr) -- (0.5,7*\rr) -- (1.5,7*\rr) 
-- (2,6*\rr) -- (1.5,5*\rr);
\draw [line width=0.2mm] (0.5,7*\rr) -- (0,8*\rr) -- (0.5,9*\rr) -- (1.5,9*\rr) 
-- (2,8*\rr) -- (1.5,7*\rr);
\draw [line width=0.2mm] (0.5,9*\rr) -- (0,10*\rr) -- (0.5,11*\rr) -- (1.5,11*\rr) 
-- (2,10*\rr) -- (1.5,9*\rr);
\draw [line width=0.2mm] (0.5,11*\rr) -- (0,12*\rr) -- (0.5,13*\rr) -- (1.5,13*\rr) 
-- (2,12*\rr) -- (1.5,11*\rr);
\draw [line width=0.2mm] (0.5,13*\rr) -- (0,14*\rr) -- (0.5,15*\rr) -- (1.5,15*\rr) 
-- (2,14*\rr) -- (1.5,13*\rr);
\draw [line width=0.2mm] (0.5,15*\rr) -- (0,16*\rr) -- (0.5,17*\rr) -- (1.5,17*\rr) 
-- (2,16*\rr) -- (1.5,15*\rr);
\draw [line width=0.2mm] (0.5,17*\rr) -- (0,18*\rr) -- (0.5,19*\rr) -- (1.5,19*\rr) 
-- (2,18*\rr) -- (1.5,17*\rr);
\draw [line width=0.2mm] (0.5,19*\rr) -- (0,20*\rr) -- (0.5,21*\rr) -- (1.5,21*\rr) 
-- (2,20*\rr) -- (1.5,19*\rr);
%\draw [line width=0.2mm] (0.5,21*\rr) -- (0,22*\rr) -- (0.5,22*\rr) -- (1.5,22*\rr) 
%-- (2,22*\rr) -- (1.5,21*\rr);

\draw [line width=0.2mm]  (3,0) -- (3.5,-1*\rr);
\draw [line width=0.2mm] (2,2*\rr) -- (3,2*\rr) -- (3.5,1*\rr) -- (3,0) -- (2,0);
\draw [line width=0.2mm] (2,4*\rr) -- (3,4*\rr) -- (3.5,3*\rr) -- (3,2*\rr);
\draw [line width=0.2mm] (2,6*\rr) -- (3,6*\rr) -- (3.5,5*\rr) -- (3,4*\rr);
\draw [line width=0.2mm] (2,8*\rr) -- (3,8*\rr) -- (3.5,7*\rr) -- (3,6*\rr);
\draw [line width=0.2mm] (2,10*\rr) -- (3,10*\rr) -- (3.5,9*\rr) -- (3,8*\rr);
\draw [line width=0.2mm] (2,12*\rr) -- (3,12*\rr) -- (3.5,11*\rr) -- (3,10*\rr);
\draw [line width=0.2mm] (2,14*\rr) -- (3,14*\rr) -- (3.5,13*\rr) -- (3,12*\rr);
\draw [line width=0.2mm] (2,16*\rr) -- (3,16*\rr) -- (3.5,15*\rr) -- (3,14*\rr);
\draw [line width=0.2mm] (2,18*\rr) -- (3,18*\rr) -- (3.5,17*\rr) -- (3,16*\rr);
\draw [line width=0.2mm] (2,20*\rr) -- (3,20*\rr) -- (3.5,19*\rr) -- (3,18*\rr);
\draw [line width=0.2mm] (3.5,21*\rr) -- (3,20*\rr);

\draw [line width=0.2mm]  (3.5,1*\rr) -- (4.5,1*\rr) -- (5,0) -- (4.5,-1*\rr)
--(3.5,-1*\rr);
\draw [line width=0.2mm]  (3.5,3*\rr) -- (4.5,3*\rr) -- (5,2*\rr) -- (4.5,1*\rr);
\draw [line width=0.2mm]  (3.5,5*\rr) -- (4.5,5*\rr) -- (5,4*\rr) -- (4.5,3*\rr);
\draw [line width=0.2mm]  (3.5,7*\rr) -- (4.5,7*\rr) -- (5,6*\rr) -- (4.5,5*\rr);
\draw [line width=0.2mm]  (3.5,9*\rr) -- (4.5,9*\rr) -- (5,8*\rr) -- (4.5,7*\rr);
\draw [line width=0.2mm]  (3.5,11*\rr) -- (4.5,11*\rr) -- (5,10*\rr) -- (4.5,9*\rr);
\draw [line width=0.2mm]  (3.5,13*\rr) -- (4.5,13*\rr) -- (5,12*\rr) -- (4.5,11*\rr);
\draw [line width=0.2mm]  (3.5,15*\rr) -- (4.5,15*\rr) -- (5,14*\rr) -- (4.5,13*\rr);
\draw [line width=0.2mm]  (3.5,17*\rr) -- (4.5,17*\rr) -- (5,16*\rr) -- (4.5,15*\rr);
\draw [line width=0.2mm]  (3.5,19*\rr) -- (4.5,19*\rr) -- (5,18*\rr) -- (4.5,17*\rr);
\draw [line width=0.2mm]  (3.5,21*\rr) -- (4.5,21*\rr) -- (5,20*\rr) -- (4.5,19*\rr);

\draw [line width=0.2mm]  (5,0) -- (6,0)--(6.5,-1*\rr );
\draw [line width=0.2mm]  (5,2*\rr) -- (6,2*\rr)--(6.5,1*\rr )--(6,0);
\draw [line width=0.2mm]  (5,4*\rr) -- (6,4*\rr)--(6.5,3*\rr )--(6,2*\rr);
\draw [line width=0.2mm]  (5,6*\rr) -- (6,6*\rr)--(6.5,5*\rr )--(6,4*\rr);
\draw [line width=0.2mm]  (5,8*\rr) -- (6,8*\rr)--(6.5,7*\rr )--(6,6*\rr);
\draw [line width=0.2mm]  (5,10*\rr) -- (6,10*\rr)--(6.5,9*\rr )--(6,8*\rr);
\draw [line width=0.2mm]  (5,12*\rr) -- (6,12*\rr)--(6.5,11*\rr )--(6,10*\rr);
\draw [line width=0.2mm]  (5,14*\rr) -- (6,14*\rr)--(6.5,13*\rr )--(6,12*\rr);
\draw [line width=0.2mm]  (5,16*\rr) -- (6,16*\rr)--(6.5,15*\rr )--(6,14*\rr);
\draw [line width=0.2mm]  (5,18*\rr) -- (6,18*\rr)--(6.5,17*\rr )--(6,16*\rr);
\draw [line width=0.2mm]  (5,20*\rr) -- (6,20*\rr)--(6.5,19*\rr )--(6,18*\rr);
\draw [line width=0.2mm]  (6.5,21*\rr )--(6,20*\rr);

\draw [line width=0.2mm]  (6.5,1*\rr) -- (7.5,1*\rr)--(8,0)--(7.5,-1*\rr)
--(6.5,-1*\rr);
\draw [line width=0.2mm]  (6.5,3*\rr) -- (7.5,3*\rr)--(8,2*\rr)--(7.5,1*\rr);
\draw [line width=0.2mm]  (6.5,5*\rr) -- (7.5,5*\rr)--(8,4*\rr)--(7.5,3*\rr);
\draw [line width=0.2mm]  (6.5,7*\rr) -- (7.5,7*\rr)--(8,6*\rr)--(7.5,5*\rr);
\draw [line width=0.2mm]  (6.5,9*\rr) -- (7.5,9*\rr)--(8,8*\rr)--(7.5,7*\rr);
\draw [line width=0.2mm]  (6.5,11*\rr) -- (7.5,11*\rr)--(8,10*\rr)--(7.5,9*\rr);
\draw [line width=0.2mm]  (6.5,13*\rr) -- (7.5,13*\rr)--(8,12*\rr)--(7.5,11*\rr);
\draw [line width=0.2mm]  (6.5,15*\rr) -- (7.5,15*\rr)--(8,14*\rr)--(7.5,13*\rr);
\draw [line width=0.2mm]  (6.5,17*\rr) -- (7.5,17*\rr)--(8,16*\rr)--(7.5,15*\rr);
\draw [line width=0.2mm]  (6.5,19*\rr) -- (7.5,19*\rr)--(8,18*\rr)--(7.5,17*\rr);
\draw [line width=0.2mm]  (6.5,21*\rr) -- (7.5,21*\rr)--(8,20*\rr)--(7.5,19*\rr);

\draw [line width=0.2mm]  (8,0)--(9,0)--(9.5,-1*\rr);
\draw [line width=0.2mm]  (8,2*\rr)--(9,2*\rr)--(9.5,1*\rr)--(9,0);
\draw [line width=0.2mm]  (8,4*\rr)--(9,4*\rr)--(9.5,3*\rr)--(9,2*\rr);
\draw [line width=0.2mm]  (8,6*\rr)--(9,6*\rr)--(9.5,5*\rr)--(9,4*\rr);
\draw [line width=0.2mm]  (8,8*\rr)--(9,8*\rr)--(9.5,7*\rr)--(9,6*\rr);
\draw [line width=0.2mm]  (8,10*\rr)--(9,10*\rr)--(9.5,9*\rr)--(9,8*\rr);
\draw [line width=0.2mm]  (8,12*\rr)--(9,12*\rr)--(9.5,11*\rr)--(9,10*\rr);
\draw [line width=0.2mm]  (8,14*\rr)--(9,14*\rr)--(9.5,13*\rr)--(9,12*\rr);
\draw [line width=0.2mm]  (8,16*\rr)--(9,16*\rr)--(9.5,15*\rr)--(9,14*\rr);
\draw [line width=0.2mm]  (8,18*\rr)--(9,18*\rr)--(9.5,17*\rr)--(9,16*\rr);
\draw [line width=0.2mm]  (8,20*\rr)--(9,20*\rr)--(9.5,19*\rr)--(9,18*\rr);
\draw [line width=0.2mm]  (9.5,21*\rr)--(9,20*\rr);

\draw [line width=0.2mm]  (9.5,1*\rr)--(10.5,1*\rr)--(11,0)--(10.5,-1*\rr)
--(9.5,-1*\rr);
\draw [line width=0.2mm]  (9.5,3*\rr)--(10.5,3*\rr)--(11,2*\rr)--(10.5,1*\rr);
\draw [line width=0.2mm]  (9.5,5*\rr)--(10.5,5*\rr)--(11,4*\rr)--(10.5,3*\rr);
\draw [line width=0.2mm]  (9.5,7*\rr)--(10.5,7*\rr)--(11,6*\rr)--(10.5,5*\rr);
\draw [line width=0.2mm]  (9.5,9*\rr)--(10.5,9*\rr)--(11,8*\rr)--(10.5,7*\rr);
\draw [line width=0.2mm]  (9.5,11*\rr)--(10.5,11*\rr)--(11,10*\rr)--(10.5,9*\rr);
\draw [line width=0.2mm]  (9.5,13*\rr)--(10.5,13*\rr)--(11,12*\rr)--(10.5,11*\rr);
\draw [line width=0.2mm]  (9.5,15*\rr)--(10.5,15*\rr)--(11,14*\rr)--(10.5,13*\rr);
\draw [line width=0.2mm]  (9.5,17*\rr)--(10.5,17*\rr)--(11,16*\rr)--(10.5,15*\rr);
\draw [line width=0.2mm]  (9.5,19*\rr)--(10.5,19*\rr)--(11,18*\rr)--(10.5,17*\rr);
\draw [line width=0.2mm]  (9.5,21*\rr)--(10.5,21*\rr)--(11,20*\rr)--(10.5,19*\rr);

\draw [line width=0.2mm]  (11,0)--(12,0)--(12.5,-1*\rr);
\draw [line width=0.2mm]  (11,2*\rr)--(12,2*\rr)--(12.5,1*\rr)--(12,0);
\draw [line width=0.2mm]  (11,4*\rr)--(12,4*\rr)--(12.5,3*\rr)--(12,2*\rr);
\draw [line width=0.2mm]  (11,6*\rr)--(12,6*\rr)--(12.5,5*\rr)--(12,4*\rr);
\draw [line width=0.2mm]  (11,8*\rr)--(12,8*\rr)--(12.5,7*\rr)--(12,6*\rr);
\draw [line width=0.2mm]  (11,10*\rr)--(12,10*\rr)--(12.5,9*\rr)--(12,8*\rr);
\draw [line width=0.2mm]  (11,12*\rr)--(12,12*\rr)--(12.5,11*\rr)--(12,10*\rr);
\draw [line width=0.2mm]  (11,14*\rr)--(12,14*\rr)--(12.5,13*\rr)--(12,12*\rr);
\draw [line width=0.2mm]  (11,16*\rr)--(12,16*\rr)--(12.5,15*\rr)--(12,14*\rr);
\draw [line width=0.2mm]  (11,18*\rr)--(12,18*\rr)--(12.5,17*\rr)--(12,16*\rr);
\draw [line width=0.2mm]  (11,20*\rr)--(12,20*\rr)--(12.5,19*\rr)--(12,18*\rr);
\draw [line width=0.2mm]  (12,20*\rr)--(12.5,21*\rr);

\draw [line width=0.2mm]  (12.5,1*\rr)--(13.5,1*\rr)--(14,0)--(13.5,-1*\rr)
--(12.5,-1*\rr);
\draw [line width=0.2mm]  (12.5,3*\rr)--(13.5,3*\rr)--(14,2*\rr)--(13.5,1*\rr);
\draw [line width=0.2mm]  (12.5,5*\rr)--(13.5,5*\rr)--(14,4*\rr)--(13.5,3*\rr);
\draw [line width=0.2mm]  (12.5,7*\rr)--(13.5,7*\rr)--(14,6*\rr)--(13.5,5*\rr);
\draw [line width=0.2mm]  (12.5,9*\rr)--(13.5,9*\rr)--(14,8*\rr)--(13.5,7*\rr);
\draw [line width=0.2mm]  (12.5,11*\rr)--(13.5,11*\rr)--(14,10*\rr)--(13.5,9*\rr);
\draw [line width=0.2mm]  (12.5,13*\rr)--(13.5,13*\rr)--(14,12*\rr)--(13.5,11*\rr);
\draw [line width=0.2mm]  (12.5,15*\rr)--(13.5,15*\rr)--(14,14*\rr)--(13.5,13*\rr);
\draw [line width=0.2mm]  (12.5,17*\rr)--(13.5,17*\rr)--(14,16*\rr)--(13.5,15*\rr);
\draw [line width=0.2mm]  (12.5,19*\rr)--(13.5,19*\rr)--(14,18*\rr)--(13.5,17*\rr);
\draw [line width=0.2mm]  (12.5,21*\rr)--(13.5,21*\rr)--(14,20*\rr)--(13.5,19*\rr);

\draw [line width=0.2mm]  (14,0)--(15,0)--(15.5,-1*\rr);
\draw [line width=0.2mm]  (14,2*\rr)--(15,2*\rr)--(15.5,1*\rr)--(15,0);
\draw [line width=0.2mm]  (14,4*\rr)--(15,4*\rr)--(15.5,3*\rr)--(15,2*\rr);
\draw [line width=0.2mm]  (14,6*\rr)--(15,6*\rr)--(15.5,5*\rr)--(15,4*\rr);
\draw [line width=0.2mm]  (14,8*\rr)--(15,8*\rr)--(15.5,7*\rr)--(15,6*\rr);
\draw [line width=0.2mm]  (14,10*\rr)--(15,10*\rr)--(15.5,9*\rr)--(15,8*\rr);
\draw [line width=0.2mm]  (14,12*\rr)--(15,12*\rr)--(15.5,11*\rr)--(15,10*\rr);
\draw [line width=0.2mm]  (14,14*\rr)--(15,14*\rr)--(15.5,13*\rr)--(15,12*\rr);
\draw [line width=0.2mm]  (14,16*\rr)--(15,16*\rr)--(15.5,15*\rr)--(15,14*\rr);
\draw [line width=0.2mm]  (14,18*\rr)--(15,18*\rr)--(15.5,17*\rr)--(15,16*\rr);
\draw [line width=0.2mm]  (14,20*\rr)--(15,20*\rr)--(15.5,19*\rr)--(15,18*\rr);
\draw [line width=0.2mm]  (15,20*\rr)--(15.5,21*\rr);

\draw [line width=0.2mm]  (15.5,1*\rr)--(16.5,1*\rr)--(17,0)--(16.5,-1*\rr)
--(15.5,-1*\rr);
\draw [line width=0.2mm]  (15.5,3*\rr)--(16.5,3*\rr)--(17,2*\rr)--(16.5,1*\rr);
\draw [line width=0.2mm]  (15.5,5*\rr)--(16.5,5*\rr)--(17,4*\rr)--(16.5,3*\rr);
\draw [line width=0.2mm]  (15.5,7*\rr)--(16.5,7*\rr)--(17,6*\rr)--(16.5,5*\rr);
\draw [line width=0.2mm]  (15.5,9*\rr)--(16.5,9*\rr)--(17,8*\rr)--(16.5,7*\rr);
\draw [line width=0.2mm]  (15.5,11*\rr)--(16.5,11*\rr)--(17,10*\rr)--(16.5,9*\rr);
\draw [line width=0.2mm]  (15.5,13*\rr)--(16.5,13*\rr)--(17,12*\rr)--(16.5,11*\rr);
\draw [line width=0.2mm]  (15.5,15*\rr)--(16.5,15*\rr)--(17,14*\rr)--(16.5,13*\rr);
\draw [line width=0.2mm]  (15.5,17*\rr)--(16.5,17*\rr)--(17,16*\rr)--(16.5,15*\rr);
\draw [line width=0.2mm]  (15.5,19*\rr)--(16.5,19*\rr)--(17,18*\rr)--(16.5,17*\rr);
\draw [line width=0.2mm]  (15.5,21*\rr)--(16.5,21*\rr)--(17,20*\rr)--(16.5,19*\rr);

\draw [line width=0.2mm]  (17,0)--(18,0)--(18.5,-1*\rr);
\draw [line width=0.2mm]  (17,2*\rr)--(18,2*\rr)--(18.5,1*\rr)--(18,0);
\draw [line width=0.2mm]  (17,4*\rr)--(18,4*\rr)--(18.5,3*\rr)--(18,2*\rr);
\draw [line width=0.2mm]  (17,6*\rr)--(18,6*\rr)--(18.5,5*\rr)--(18,4*\rr);
\draw [line width=0.2mm]  (17,8*\rr)--(18,8*\rr)--(18.5,7*\rr)--(18,6*\rr);
\draw [line width=0.2mm]  (17,10*\rr)--(18,10*\rr)--(18.5,9*\rr)--(18,8*\rr);
\draw [line width=0.2mm]  (17,12*\rr)--(18,12*\rr)--(18.5,11*\rr)--(18,10*\rr);
\draw [line width=0.2mm]  (17,14*\rr)--(18,14*\rr)--(18.5,13*\rr)--(18,12*\rr);
\draw [line width=0.2mm]  (17,16*\rr)--(18,16*\rr)--(18.5,15*\rr)--(18,14*\rr);
\draw [line width=0.2mm]  (17,18*\rr)--(18,18*\rr)--(18.5,17*\rr)--(18,16*\rr);
\draw [line width=0.2mm]  (17,20*\rr)--(18,20*\rr)--(18.5,19*\rr)--(18,18*\rr);
\draw [line width=0.2mm]  (18,20*\rr)--(18.5,21*\rr);

\draw [line width=0.2mm]  (18.5,1*\rr)--(19.5,1*\rr)--(20,0)--(19.5,-1*\rr)
--(18.5,-1*\rr);
\draw [line width=0.2mm]  (18.5,3*\rr)--(19.5,3*\rr)--(20,2*\rr)--(19.5,1*\rr);
\draw [line width=0.2mm]  (18.5,5*\rr)--(19.5,5*\rr)--(20,4*\rr)--(19.5,3*\rr);
\draw [line width=0.2mm]  (18.5,7*\rr)--(19.5,7*\rr)--(20,6*\rr)--(19.5,5*\rr);
\draw [line width=0.2mm]  (18.5,9*\rr)--(19.5,9*\rr)--(20,8*\rr)--(19.5,7*\rr);
\draw [line width=0.2mm]  (18.5,11*\rr)--(19.5,11*\rr)--(20,10*\rr)--(19.5,9*\rr);
\draw [line width=0.2mm]  (18.5,13*\rr)--(19.5,13*\rr)--(20,12*\rr)--(19.5,11*\rr);
\draw [line width=0.2mm]  (18.5,15*\rr)--(19.5,15*\rr)--(20,14*\rr)--(19.5,13*\rr);
\draw [line width=0.2mm]  (18.5,17*\rr)--(19.5,17*\rr)--(20,16*\rr)--(19.5,15*\rr);
\draw [line width=0.2mm]  (18.5,19*\rr)--(19.5,19*\rr)--(20,18*\rr)--(19.5,17*\rr);
\draw [line width=0.2mm]  (18.5,21*\rr)--(19.5,21*\rr)--(20,20*\rr)--(19.5,19*\rr);

\draw [line width=0.2mm]  (20,0)--(21,0)--(21.5,-1*\rr);
\draw [line width=0.2mm]  (20,2*\rr)--(21,2*\rr)--(21.5,1*\rr)--(21,0);
\draw [line width=0.2mm]  (20,4*\rr)--(21,4*\rr)--(21.5,3*\rr)--(21,2*\rr);
\draw [line width=0.2mm]  (20,6*\rr)--(21,6*\rr)--(21.5,5*\rr)--(21,4*\rr);
\draw [line width=0.2mm]  (20,8*\rr)--(21,8*\rr)--(21.5,7*\rr)--(21,6*\rr);
\draw [line width=0.2mm]  (20,10*\rr)--(21,10*\rr)--(21.5,9*\rr)--(21,8*\rr);
\draw [line width=0.2mm]  (20,12*\rr)--(21,12*\rr)--(21.5,11*\rr)--(21,10*\rr);
\draw [line width=0.2mm]  (20,14*\rr)--(21,14*\rr)--(21.5,13*\rr)--(21,12*\rr);
\draw [line width=0.2mm]  (20,16*\rr)--(21,16*\rr)--(21.5,15*\rr)--(21,14*\rr);
\draw [line width=0.2mm]  (20,18*\rr)--(21,18*\rr)--(21.5,17*\rr)--(21,16*\rr);
\draw [line width=0.2mm]  (20,20*\rr)--(21,20*\rr)--(21.5,19*\rr)--(21,18*\rr);
\draw [line width=0.2mm]  (21,20*\rr)--(21.5,21*\rr);

\draw [line width=0.2mm]  (21.5,1*\rr)--(22.5,1*\rr)--(23,0)--(22.5,-1*\rr)
--(21.5,-1*\rr);
\draw [line width=0.2mm]  (21.5,3*\rr)--(22.5,3*\rr)--(23,2*\rr)--(22.5,1*\rr);
\draw [line width=0.2mm]  (21.5,5*\rr)--(22.5,5*\rr)--(23,4*\rr)--(22.5,3*\rr);
\draw [line width=0.2mm]  (21.5,7*\rr)--(22.5,7*\rr)--(23,6*\rr)--(22.5,5*\rr);
\draw [line width=0.2mm]  (21.5,9*\rr)--(22.5,9*\rr)--(23,8*\rr)--(22.5,7*\rr);
\draw [line width=0.2mm]  (21.5,11*\rr)--(22.5,11*\rr)--(23,10*\rr)--(22.5,9*\rr);
\draw [line width=0.2mm]  (21.5,13*\rr)--(22.5,13*\rr)--(23,12*\rr)--(22.5,11*\rr);
\draw [line width=0.2mm]  (21.5,15*\rr)--(22.5,15*\rr)--(23,14*\rr)--(22.5,13*\rr);
\draw [line width=0.2mm]  (21.5,17*\rr)--(22.5,17*\rr)--(23,16*\rr)--(22.5,15*\rr);
\draw [line width=0.2mm]  (21.5,19*\rr)--(22.5,19*\rr)--(23,18*\rr)--(22.5,17*\rr);
\draw [line width=0.2mm]  (21.5,21*\rr)--(22.5,21*\rr)--(23,20*\rr)--(22.5,19*\rr);

\draw [line width=0.2mm]  (23,0)--(24,0)--(24.5,-1*\rr);
\draw [line width=0.2mm]  (23,2*\rr)--(24,2*\rr)--(24.5,1*\rr)--(24,0);
\draw [line width=0.2mm]  (23,4*\rr)--(24,4*\rr)--(24.5,3*\rr)--(24,2*\rr);
\draw [line width=0.2mm]  (23,6*\rr)--(24,6*\rr)--(24.5,5*\rr)--(24,4*\rr);
\draw [line width=0.2mm]  (23,8*\rr)--(24,8*\rr)--(24.5,7*\rr)--(24,6*\rr);
\draw [line width=0.2mm]  (23,10*\rr)--(24,10*\rr)--(24.5,9*\rr)--(24,8*\rr);
\draw [line width=0.2mm]  (23,12*\rr)--(24,12*\rr)--(24.5,11*\rr)--(24,10*\rr);
\draw [line width=0.2mm]  (23,14*\rr)--(24,14*\rr)--(24.5,13*\rr)--(24,12*\rr);
\draw [line width=0.2mm]  (23,16*\rr)--(24,16*\rr)--(24.5,15*\rr)--(24,14*\rr);
\draw [line width=0.2mm]  (23,18*\rr)--(24,18*\rr)--(24.5,17*\rr)--(24,16*\rr);
\draw [line width=0.2mm]  (23,20*\rr)--(24,20*\rr)--(24.5,19*\rr)--(24,18*\rr);
\draw [line width=0.2mm]  (24,20*\rr)--(24.5,21*\rr);

\draw [line width=0.6mm] (12.5,19*\rr)--(7.5,13*\rr)--(14,10*\rr)--(12.5,19*\rr)
--(19.5,15*\rr)--(14,10*\rr)--(15.5,1*\rr);

\draw [line width=0.6mm] (7.5,13*\rr)--(9,4*\rr)--(14,10*\rr);
\draw [line width=0.6mm] (9,4*\rr)--(15.5,1*\rr);
\draw [line width=0.6mm] (19.5,15*\rr)--(21,6*\rr)--(14,10*\rr);
\draw [line width=0.6mm] (21,6*\rr)--(15.5,1*\rr);

\draw [line width=0.6mm] (21,6*\rr)--(22.5,-3*\rr);
\draw [line width=0.6mm] (15.5,1*\rr)--(17,-8*\rr);
\draw [line width=0.6mm] (9,4*\rr)--(10.5,-5*\rr);
\draw [line width=0.6mm] (0.5,17*\rr)--(2,8*\rr)--(3.5,-1*\rr);

\draw [line width=0.6mm]  (18,24*\rr)--(19.5,15*\rr);
\draw [line width=0.6mm]  (6,22*\rr)--(7.5,13*\rr);
\draw [line width=0.6mm] (12.5,19*\rr)--(11,28*\rr);
\draw [line width=0.6mm] (3.5,-1*\rr)--(9,4*\rr)--(2,8*\rr)--(7.5,13*\rr)
--(0.5,17*\rr)--(6,22*\rr)--(12.5,19*\rr)--(18,24*\rr);
\draw [line width=0.6mm] (22.5,-3*\rr)--(15.5,1*\rr)--(10.5,-5*\rr);
\draw [line width=0.6mm] (27.5,3*\rr)--(21,6*\rr)--(26,12*\rr)--(19.5,15*\rr)--(24.5,21*\rr);
\draw [line width=0.6mm] (3.5,-1*\rr)--(-3,2*\rr)--(2,8*\rr)--(-4.5,11*\rr);

% \draw [line width=0.6mm] (14,10*\rr)--(16.5,7*\rr)--(20,8*\rr);

\foreach \pos in {(7.5,13*\rr),(14,10*\rr),(19.5,15*\rr),(12.5,19*\rr),
(2,8*\rr),(0.5,17*\rr)}
\shade[shading=ball, ball color=black] \pos circle (.6);
% \foreach \pos in {(13,18*\rr)}
% \shade[shading=ball, ball color=red] \pos circle (.3);

\foreach \pos in {(9,4*\rr),(15.5,1*\rr),(21,6*\rr),(3.5,-1*\rr)}
\shade[shading=ball, ball color=black] \pos circle (.6);

%\foreach \pos in {(10.5, 7*\rr),(0, 14*\rr),(12.5, 17*\rr)}
 %\shade[shading=ball, ball color=white] \pos circle (.4);

\draw (15,13*\rr) node   {\Large {${\mbox{\boldmath$O$}}$}};
%\draw (5.7,13.2*\rr) node   {\Large {${\mbox{\boldmath$A$}}$}};
%\draw (10,19*\rr) node   {\Large {${\mbox{\boldmath$B$}}$}};
%\draw (22,15.5*\rr) node   {\Large {${\mbox{\boldmath$C$}}$}};
%\draw (23,3*\rr) node   {\Large {${\mbox{\boldmath$D$}}$}};
%\draw (18.5,1*\rr) node   {\Large {${\mbox{\boldmath$E$}}$}};
 
\end{tikzpicture} C\begin{tikzpicture}[scale=0.18]
\clip (0, -1.5*\rr) rectangle (24.5, 21*\rr);

%\path [fill=lightgray] (1,-1*\rr)-- (4,-1*\rr)--(10.7,21*\rr)--(7.5,21*\rr);
%\path [fill=cyan] (4.2,-1*\rr)--(7.4,-1*\rr)--(14,21*\rr)--(10.8,21*\rr);
%\path [fill=lightgray] (7.5,-1*\rr)--(10.6,-1*\rr)--(17.3,21*\rr)--(14.2,21*\rr);

\draw [line width=0.2mm] (0.5,1*\rr) -- (1.5,1*\rr) -- (2,0) -- (1.5,-1*\rr)
-- (0.5,-1*\rr) -- (0,0) -- (0.5,1*\rr);
\draw [line width=0.2mm] (0.5,1*\rr) -- (0,2*\rr) -- (0.5,3*\rr) -- (1.5,3*\rr) 
-- (2,2*\rr) -- (1.5,1*\rr);
\draw [line width=0.2mm] (0.5,3*\rr) -- (0,4*\rr) -- (0.5,5*\rr) -- (1.5,5*\rr) 
-- (2,4*\rr) -- (1.5,3*\rr);
\draw [line width=0.2mm] (0.5,5*\rr) -- (0,6*\rr) -- (0.5,7*\rr) -- (1.5,7*\rr) 
-- (2,6*\rr) -- (1.5,5*\rr);
\draw [line width=0.2mm] (0.5,7*\rr) -- (0,8*\rr) -- (0.5,9*\rr) -- (1.5,9*\rr) 
-- (2,8*\rr) -- (1.5,7*\rr);
\draw [line width=0.2mm] (0.5,9*\rr) -- (0,10*\rr) -- (0.5,11*\rr) -- (1.5,11*\rr) 
-- (2,10*\rr) -- (1.5,9*\rr);
\draw [line width=0.2mm] (0.5,11*\rr) -- (0,12*\rr) -- (0.5,13*\rr) -- (1.5,13*\rr) 
-- (2,12*\rr) -- (1.5,11*\rr);
\draw [line width=0.2mm] (0.5,13*\rr) -- (0,14*\rr) -- (0.5,15*\rr) -- (1.5,15*\rr) 
-- (2,14*\rr) -- (1.5,13*\rr);
\draw [line width=0.2mm] (0.5,15*\rr) -- (0,16*\rr) -- (0.5,17*\rr) -- (1.5,17*\rr) 
-- (2,16*\rr) -- (1.5,15*\rr);
\draw [line width=0.2mm] (0.5,17*\rr) -- (0,18*\rr) -- (0.5,19*\rr) -- (1.5,19*\rr) 
-- (2,18*\rr) -- (1.5,17*\rr);
\draw [line width=0.2mm] (0.5,19*\rr) -- (0,20*\rr) -- (0.5,21*\rr) -- (1.5,21*\rr) 
-- (2,20*\rr) -- (1.5,19*\rr);
%\draw [line width=0.2mm] (0.5,21*\rr) -- (0,22*\rr) -- (0.5,22*\rr) -- (1.5,22*\rr) 
%-- (2,22*\rr) -- (1.5,21*\rr);

\draw [line width=0.2mm]  (3,0) -- (3.5,-1*\rr);
\draw [line width=0.2mm] (2,2*\rr) -- (3,2*\rr) -- (3.5,1*\rr) -- (3,0) -- (2,0);
\draw [line width=0.2mm] (2,4*\rr) -- (3,4*\rr) -- (3.5,3*\rr) -- (3,2*\rr);
\draw [line width=0.2mm] (2,6*\rr) -- (3,6*\rr) -- (3.5,5*\rr) -- (3,4*\rr);
\draw [line width=0.2mm] (2,8*\rr) -- (3,8*\rr) -- (3.5,7*\rr) -- (3,6*\rr);
\draw [line width=0.2mm] (2,10*\rr) -- (3,10*\rr) -- (3.5,9*\rr) -- (3,8*\rr);
\draw [line width=0.2mm] (2,12*\rr) -- (3,12*\rr) -- (3.5,11*\rr) -- (3,10*\rr);
\draw [line width=0.2mm] (2,14*\rr) -- (3,14*\rr) -- (3.5,13*\rr) -- (3,12*\rr);
\draw [line width=0.2mm] (2,16*\rr) -- (3,16*\rr) -- (3.5,15*\rr) -- (3,14*\rr);
\draw [line width=0.2mm] (2,18*\rr) -- (3,18*\rr) -- (3.5,17*\rr) -- (3,16*\rr);
\draw [line width=0.2mm] (2,20*\rr) -- (3,20*\rr) -- (3.5,19*\rr) -- (3,18*\rr);
\draw [line width=0.2mm] (3.5,21*\rr) -- (3,20*\rr);

\draw [line width=0.2mm]  (3.5,1*\rr) -- (4.5,1*\rr) -- (5,0) -- (4.5,-1*\rr)
--(3.5,-1*\rr);
\draw [line width=0.2mm]  (3.5,3*\rr) -- (4.5,3*\rr) -- (5,2*\rr) -- (4.5,1*\rr);
\draw [line width=0.2mm]  (3.5,5*\rr) -- (4.5,5*\rr) -- (5,4*\rr) -- (4.5,3*\rr);
\draw [line width=0.2mm]  (3.5,7*\rr) -- (4.5,7*\rr) -- (5,6*\rr) -- (4.5,5*\rr);
\draw [line width=0.2mm]  (3.5,9*\rr) -- (4.5,9*\rr) -- (5,8*\rr) -- (4.5,7*\rr);
\draw [line width=0.2mm]  (3.5,11*\rr) -- (4.5,11*\rr) -- (5,10*\rr) -- (4.5,9*\rr);
\draw [line width=0.2mm]  (3.5,13*\rr) -- (4.5,13*\rr) -- (5,12*\rr) -- (4.5,11*\rr);
\draw [line width=0.2mm]  (3.5,15*\rr) -- (4.5,15*\rr) -- (5,14*\rr) -- (4.5,13*\rr);
\draw [line width=0.2mm]  (3.5,17*\rr) -- (4.5,17*\rr) -- (5,16*\rr) -- (4.5,15*\rr);
\draw [line width=0.2mm]  (3.5,19*\rr) -- (4.5,19*\rr) -- (5,18*\rr) -- (4.5,17*\rr);
\draw [line width=0.2mm]  (3.5,21*\rr) -- (4.5,21*\rr) -- (5,20*\rr) -- (4.5,19*\rr);

\draw [line width=0.2mm]  (5,0) -- (6,0)--(6.5,-1*\rr );
\draw [line width=0.2mm]  (5,2*\rr) -- (6,2*\rr)--(6.5,1*\rr )--(6,0);
\draw [line width=0.2mm]  (5,4*\rr) -- (6,4*\rr)--(6.5,3*\rr )--(6,2*\rr);
\draw [line width=0.2mm]  (5,6*\rr) -- (6,6*\rr)--(6.5,5*\rr )--(6,4*\rr);
\draw [line width=0.2mm]  (5,8*\rr) -- (6,8*\rr)--(6.5,7*\rr )--(6,6*\rr);
\draw [line width=0.2mm]  (5,10*\rr) -- (6,10*\rr)--(6.5,9*\rr )--(6,8*\rr);
\draw [line width=0.2mm]  (5,12*\rr) -- (6,12*\rr)--(6.5,11*\rr )--(6,10*\rr);
\draw [line width=0.2mm]  (5,14*\rr) -- (6,14*\rr)--(6.5,13*\rr )--(6,12*\rr);
\draw [line width=0.2mm]  (5,16*\rr) -- (6,16*\rr)--(6.5,15*\rr )--(6,14*\rr);
\draw [line width=0.2mm]  (5,18*\rr) -- (6,18*\rr)--(6.5,17*\rr )--(6,16*\rr);
\draw [line width=0.2mm]  (5,20*\rr) -- (6,20*\rr)--(6.5,19*\rr )--(6,18*\rr);
\draw [line width=0.2mm]  (6.5,21*\rr )--(6,20*\rr);

\draw [line width=0.2mm]  (6.5,1*\rr) -- (7.5,1*\rr)--(8,0)--(7.5,-1*\rr)
--(6.5,-1*\rr);
\draw [line width=0.2mm]  (6.5,3*\rr) -- (7.5,3*\rr)--(8,2*\rr)--(7.5,1*\rr);
\draw [line width=0.2mm]  (6.5,5*\rr) -- (7.5,5*\rr)--(8,4*\rr)--(7.5,3*\rr);
\draw [line width=0.2mm]  (6.5,7*\rr) -- (7.5,7*\rr)--(8,6*\rr)--(7.5,5*\rr);
\draw [line width=0.2mm]  (6.5,9*\rr) -- (7.5,9*\rr)--(8,8*\rr)--(7.5,7*\rr);
\draw [line width=0.2mm]  (6.5,11*\rr) -- (7.5,11*\rr)--(8,10*\rr)--(7.5,9*\rr);
\draw [line width=0.2mm]  (6.5,13*\rr) -- (7.5,13*\rr)--(8,12*\rr)--(7.5,11*\rr);
\draw [line width=0.2mm]  (6.5,15*\rr) -- (7.5,15*\rr)--(8,14*\rr)--(7.5,13*\rr);
\draw [line width=0.2mm]  (6.5,17*\rr) -- (7.5,17*\rr)--(8,16*\rr)--(7.5,15*\rr);
\draw [line width=0.2mm]  (6.5,19*\rr) -- (7.5,19*\rr)--(8,18*\rr)--(7.5,17*\rr);
\draw [line width=0.2mm]  (6.5,21*\rr) -- (7.5,21*\rr)--(8,20*\rr)--(7.5,19*\rr);

\draw [line width=0.2mm]  (8,0)--(9,0)--(9.5,-1*\rr);
\draw [line width=0.2mm]  (8,2*\rr)--(9,2*\rr)--(9.5,1*\rr)--(9,0);
\draw [line width=0.2mm]  (8,4*\rr)--(9,4*\rr)--(9.5,3*\rr)--(9,2*\rr);
\draw [line width=0.2mm]  (8,6*\rr)--(9,6*\rr)--(9.5,5*\rr)--(9,4*\rr);
\draw [line width=0.2mm]  (8,8*\rr)--(9,8*\rr)--(9.5,7*\rr)--(9,6*\rr);
\draw [line width=0.2mm]  (8,10*\rr)--(9,10*\rr)--(9.5,9*\rr)--(9,8*\rr);
\draw [line width=0.2mm]  (8,12*\rr)--(9,12*\rr)--(9.5,11*\rr)--(9,10*\rr);
\draw [line width=0.2mm]  (8,14*\rr)--(9,14*\rr)--(9.5,13*\rr)--(9,12*\rr);
\draw [line width=0.2mm]  (8,16*\rr)--(9,16*\rr)--(9.5,15*\rr)--(9,14*\rr);
\draw [line width=0.2mm]  (8,18*\rr)--(9,18*\rr)--(9.5,17*\rr)--(9,16*\rr);
\draw [line width=0.2mm]  (8,20*\rr)--(9,20*\rr)--(9.5,19*\rr)--(9,18*\rr);
\draw [line width=0.2mm]  (9.5,21*\rr)--(9,20*\rr);

\draw [line width=0.2mm]  (9.5,1*\rr)--(10.5,1*\rr)--(11,0)--(10.5,-1*\rr)
--(9.5,-1*\rr);
\draw [line width=0.2mm]  (9.5,3*\rr)--(10.5,3*\rr)--(11,2*\rr)--(10.5,1*\rr);
\draw [line width=0.2mm]  (9.5,5*\rr)--(10.5,5*\rr)--(11,4*\rr)--(10.5,3*\rr);
\draw [line width=0.2mm]  (9.5,7*\rr)--(10.5,7*\rr)--(11,6*\rr)--(10.5,5*\rr);
\draw [line width=0.2mm]  (9.5,9*\rr)--(10.5,9*\rr)--(11,8*\rr)--(10.5,7*\rr);
\draw [line width=0.2mm]  (9.5,11*\rr)--(10.5,11*\rr)--(11,10*\rr)--(10.5,9*\rr);
\draw [line width=0.2mm]  (9.5,13*\rr)--(10.5,13*\rr)--(11,12*\rr)--(10.5,11*\rr);
\draw [line width=0.2mm]  (9.5,15*\rr)--(10.5,15*\rr)--(11,14*\rr)--(10.5,13*\rr);
\draw [line width=0.2mm]  (9.5,17*\rr)--(10.5,17*\rr)--(11,16*\rr)--(10.5,15*\rr);
\draw [line width=0.2mm]  (9.5,19*\rr)--(10.5,19*\rr)--(11,18*\rr)--(10.5,17*\rr);
\draw [line width=0.2mm]  (9.5,21*\rr)--(10.5,21*\rr)--(11,20*\rr)--(10.5,19*\rr);

\draw [line width=0.2mm]  (11,0)--(12,0)--(12.5,-1*\rr);
\draw [line width=0.2mm]  (11,2*\rr)--(12,2*\rr)--(12.5,1*\rr)--(12,0);
\draw [line width=0.2mm]  (11,4*\rr)--(12,4*\rr)--(12.5,3*\rr)--(12,2*\rr);
\draw [line width=0.2mm]  (11,6*\rr)--(12,6*\rr)--(12.5,5*\rr)--(12,4*\rr);
\draw [line width=0.2mm]  (11,8*\rr)--(12,8*\rr)--(12.5,7*\rr)--(12,6*\rr);
\draw [line width=0.2mm]  (11,10*\rr)--(12,10*\rr)--(12.5,9*\rr)--(12,8*\rr);
\draw [line width=0.2mm]  (11,12*\rr)--(12,12*\rr)--(12.5,11*\rr)--(12,10*\rr);
\draw [line width=0.2mm]  (11,14*\rr)--(12,14*\rr)--(12.5,13*\rr)--(12,12*\rr);
\draw [line width=0.2mm]  (11,16*\rr)--(12,16*\rr)--(12.5,15*\rr)--(12,14*\rr);
\draw [line width=0.2mm]  (11,18*\rr)--(12,18*\rr)--(12.5,17*\rr)--(12,16*\rr);
\draw [line width=0.2mm]  (11,20*\rr)--(12,20*\rr)--(12.5,19*\rr)--(12,18*\rr);
\draw [line width=0.2mm]  (12,20*\rr)--(12.5,21*\rr);

\draw [line width=0.2mm]  (12.5,1*\rr)--(13.5,1*\rr)--(14,0)--(13.5,-1*\rr)
--(12.5,-1*\rr);
\draw [line width=0.2mm]  (12.5,3*\rr)--(13.5,3*\rr)--(14,2*\rr)--(13.5,1*\rr);
\draw [line width=0.2mm]  (12.5,5*\rr)--(13.5,5*\rr)--(14,4*\rr)--(13.5,3*\rr);
\draw [line width=0.2mm]  (12.5,7*\rr)--(13.5,7*\rr)--(14,6*\rr)--(13.5,5*\rr);
\draw [line width=0.2mm]  (12.5,9*\rr)--(13.5,9*\rr)--(14,8*\rr)--(13.5,7*\rr);
\draw [line width=0.2mm]  (12.5,11*\rr)--(13.5,11*\rr)--(14,10*\rr)--(13.5,9*\rr);
\draw [line width=0.2mm]  (12.5,13*\rr)--(13.5,13*\rr)--(14,12*\rr)--(13.5,11*\rr);
\draw [line width=0.2mm]  (12.5,15*\rr)--(13.5,15*\rr)--(14,14*\rr)--(13.5,13*\rr);
\draw [line width=0.2mm]  (12.5,17*\rr)--(13.5,17*\rr)--(14,16*\rr)--(13.5,15*\rr);
\draw [line width=0.2mm]  (12.5,19*\rr)--(13.5,19*\rr)--(14,18*\rr)--(13.5,17*\rr);
\draw [line width=0.2mm]  (12.5,21*\rr)--(13.5,21*\rr)--(14,20*\rr)--(13.5,19*\rr);

\draw [line width=0.2mm]  (14,0)--(15,0)--(15.5,-1*\rr);
\draw [line width=0.2mm]  (14,2*\rr)--(15,2*\rr)--(15.5,1*\rr)--(15,0);
\draw [line width=0.2mm]  (14,4*\rr)--(15,4*\rr)--(15.5,3*\rr)--(15,2*\rr);
\draw [line width=0.2mm]  (14,6*\rr)--(15,6*\rr)--(15.5,5*\rr)--(15,4*\rr);
\draw [line width=0.2mm]  (14,8*\rr)--(15,8*\rr)--(15.5,7*\rr)--(15,6*\rr);
\draw [line width=0.2mm]  (14,10*\rr)--(15,10*\rr)--(15.5,9*\rr)--(15,8*\rr);
\draw [line width=0.2mm]  (14,12*\rr)--(15,12*\rr)--(15.5,11*\rr)--(15,10*\rr);
\draw [line width=0.2mm]  (14,14*\rr)--(15,14*\rr)--(15.5,13*\rr)--(15,12*\rr);
\draw [line width=0.2mm]  (14,16*\rr)--(15,16*\rr)--(15.5,15*\rr)--(15,14*\rr);
\draw [line width=0.2mm]  (14,18*\rr)--(15,18*\rr)--(15.5,17*\rr)--(15,16*\rr);
\draw [line width=0.2mm]  (14,20*\rr)--(15,20*\rr)--(15.5,19*\rr)--(15,18*\rr);
\draw [line width=0.2mm]  (15,20*\rr)--(15.5,21*\rr);

\draw [line width=0.2mm]  (15.5,1*\rr)--(16.5,1*\rr)--(17,0)--(16.5,-1*\rr)
--(15.5,-1*\rr);
\draw [line width=0.2mm]  (15.5,3*\rr)--(16.5,3*\rr)--(17,2*\rr)--(16.5,1*\rr);
\draw [line width=0.2mm]  (15.5,5*\rr)--(16.5,5*\rr)--(17,4*\rr)--(16.5,3*\rr);
\draw [line width=0.2mm]  (15.5,7*\rr)--(16.5,7*\rr)--(17,6*\rr)--(16.5,5*\rr);
\draw [line width=0.2mm]  (15.5,9*\rr)--(16.5,9*\rr)--(17,8*\rr)--(16.5,7*\rr);
\draw [line width=0.2mm]  (15.5,11*\rr)--(16.5,11*\rr)--(17,10*\rr)--(16.5,9*\rr);
\draw [line width=0.2mm]  (15.5,13*\rr)--(16.5,13*\rr)--(17,12*\rr)--(16.5,11*\rr);
\draw [line width=0.2mm]  (15.5,15*\rr)--(16.5,15*\rr)--(17,14*\rr)--(16.5,13*\rr);
\draw [line width=0.2mm]  (15.5,17*\rr)--(16.5,17*\rr)--(17,16*\rr)--(16.5,15*\rr);
\draw [line width=0.2mm]  (15.5,19*\rr)--(16.5,19*\rr)--(17,18*\rr)--(16.5,17*\rr);
\draw [line width=0.2mm]  (15.5,21*\rr)--(16.5,21*\rr)--(17,20*\rr)--(16.5,19*\rr);

\draw [line width=0.2mm]  (17,0)--(18,0)--(18.5,-1*\rr);
\draw [line width=0.2mm]  (17,2*\rr)--(18,2*\rr)--(18.5,1*\rr)--(18,0);
\draw [line width=0.2mm]  (17,4*\rr)--(18,4*\rr)--(18.5,3*\rr)--(18,2*\rr);
\draw [line width=0.2mm]  (17,6*\rr)--(18,6*\rr)--(18.5,5*\rr)--(18,4*\rr);
\draw [line width=0.2mm]  (17,8*\rr)--(18,8*\rr)--(18.5,7*\rr)--(18,6*\rr);
\draw [line width=0.2mm]  (17,10*\rr)--(18,10*\rr)--(18.5,9*\rr)--(18,8*\rr);
\draw [line width=0.2mm]  (17,12*\rr)--(18,12*\rr)--(18.5,11*\rr)--(18,10*\rr);
\draw [line width=0.2mm]  (17,14*\rr)--(18,14*\rr)--(18.5,13*\rr)--(18,12*\rr);
\draw [line width=0.2mm]  (17,16*\rr)--(18,16*\rr)--(18.5,15*\rr)--(18,14*\rr);
\draw [line width=0.2mm]  (17,18*\rr)--(18,18*\rr)--(18.5,17*\rr)--(18,16*\rr);
\draw [line width=0.2mm]  (17,20*\rr)--(18,20*\rr)--(18.5,19*\rr)--(18,18*\rr);
\draw [line width=0.2mm]  (18,20*\rr)--(18.5,21*\rr);

\draw [line width=0.2mm]  (18.5,1*\rr)--(19.5,1*\rr)--(20,0)--(19.5,-1*\rr)
--(18.5,-1*\rr);
\draw [line width=0.2mm]  (18.5,3*\rr)--(19.5,3*\rr)--(20,2*\rr)--(19.5,1*\rr);
\draw [line width=0.2mm]  (18.5,5*\rr)--(19.5,5*\rr)--(20,4*\rr)--(19.5,3*\rr);
\draw [line width=0.2mm]  (18.5,7*\rr)--(19.5,7*\rr)--(20,6*\rr)--(19.5,5*\rr);
\draw [line width=0.2mm]  (18.5,9*\rr)--(19.5,9*\rr)--(20,8*\rr)--(19.5,7*\rr);
\draw [line width=0.2mm]  (18.5,11*\rr)--(19.5,11*\rr)--(20,10*\rr)--(19.5,9*\rr);
\draw [line width=0.2mm]  (18.5,13*\rr)--(19.5,13*\rr)--(20,12*\rr)--(19.5,11*\rr);
\draw [line width=0.2mm]  (18.5,15*\rr)--(19.5,15*\rr)--(20,14*\rr)--(19.5,13*\rr);
\draw [line width=0.2mm]  (18.5,17*\rr)--(19.5,17*\rr)--(20,16*\rr)--(19.5,15*\rr);
\draw [line width=0.2mm]  (18.5,19*\rr)--(19.5,19*\rr)--(20,18*\rr)--(19.5,17*\rr);
\draw [line width=0.2mm]  (18.5,21*\rr)--(19.5,21*\rr)--(20,20*\rr)--(19.5,19*\rr);

\draw [line width=0.2mm]  (20,0)--(21,0)--(21.5,-1*\rr);
\draw [line width=0.2mm]  (20,2*\rr)--(21,2*\rr)--(21.5,1*\rr)--(21,0);
\draw [line width=0.2mm]  (20,4*\rr)--(21,4*\rr)--(21.5,3*\rr)--(21,2*\rr);
\draw [line width=0.2mm]  (20,6*\rr)--(21,6*\rr)--(21.5,5*\rr)--(21,4*\rr);
\draw [line width=0.2mm]  (20,8*\rr)--(21,8*\rr)--(21.5,7*\rr)--(21,6*\rr);
\draw [line width=0.2mm]  (20,10*\rr)--(21,10*\rr)--(21.5,9*\rr)--(21,8*\rr);
\draw [line width=0.2mm]  (20,12*\rr)--(21,12*\rr)--(21.5,11*\rr)--(21,10*\rr);
\draw [line width=0.2mm]  (20,14*\rr)--(21,14*\rr)--(21.5,13*\rr)--(21,12*\rr);
\draw [line width=0.2mm]  (20,16*\rr)--(21,16*\rr)--(21.5,15*\rr)--(21,14*\rr);
\draw [line width=0.2mm]  (20,18*\rr)--(21,18*\rr)--(21.5,17*\rr)--(21,16*\rr);
\draw [line width=0.2mm]  (20,20*\rr)--(21,20*\rr)--(21.5,19*\rr)--(21,18*\rr);
\draw [line width=0.2mm]  (21,20*\rr)--(21.5,21*\rr);

\draw [line width=0.2mm]  (21.5,1*\rr)--(22.5,1*\rr)--(23,0)--(22.5,-1*\rr)
--(21.5,-1*\rr);
\draw [line width=0.2mm]  (21.5,3*\rr)--(22.5,3*\rr)--(23,2*\rr)--(22.5,1*\rr);
\draw [line width=0.2mm]  (21.5,5*\rr)--(22.5,5*\rr)--(23,4*\rr)--(22.5,3*\rr);
\draw [line width=0.2mm]  (21.5,7*\rr)--(22.5,7*\rr)--(23,6*\rr)--(22.5,5*\rr);
\draw [line width=0.2mm]  (21.5,9*\rr)--(22.5,9*\rr)--(23,8*\rr)--(22.5,7*\rr);
\draw [line width=0.2mm]  (21.5,11*\rr)--(22.5,11*\rr)--(23,10*\rr)--(22.5,9*\rr);
\draw [line width=0.2mm]  (21.5,13*\rr)--(22.5,13*\rr)--(23,12*\rr)--(22.5,11*\rr);
\draw [line width=0.2mm]  (21.5,15*\rr)--(22.5,15*\rr)--(23,14*\rr)--(22.5,13*\rr);
\draw [line width=0.2mm]  (21.5,17*\rr)--(22.5,17*\rr)--(23,16*\rr)--(22.5,15*\rr);
\draw [line width=0.2mm]  (21.5,19*\rr)--(22.5,19*\rr)--(23,18*\rr)--(22.5,17*\rr);
\draw [line width=0.2mm]  (21.5,21*\rr)--(22.5,21*\rr)--(23,20*\rr)--(22.5,19*\rr);

\draw [line width=0.2mm]  (23,0)--(24,0)--(24.5,-1*\rr);
\draw [line width=0.2mm]  (23,2*\rr)--(24,2*\rr)--(24.5,1*\rr)--(24,0);
\draw [line width=0.2mm]  (23,4*\rr)--(24,4*\rr)--(24.5,3*\rr)--(24,2*\rr);
\draw [line width=0.2mm]  (23,6*\rr)--(24,6*\rr)--(24.5,5*\rr)--(24,4*\rr);
\draw [line width=0.2mm]  (23,8*\rr)--(24,8*\rr)--(24.5,7*\rr)--(24,6*\rr);
\draw [line width=0.2mm]  (23,10*\rr)--(24,10*\rr)--(24.5,9*\rr)--(24,8*\rr);
\draw [line width=0.2mm]  (23,12*\rr)--(24,12*\rr)--(24.5,11*\rr)--(24,10*\rr);
\draw [line width=0.2mm]  (23,14*\rr)--(24,14*\rr)--(24.5,13*\rr)--(24,12*\rr);
\draw [line width=0.2mm]  (23,16*\rr)--(24,16*\rr)--(24.5,15*\rr)--(24,14*\rr);
\draw [line width=0.2mm]  (23,18*\rr)--(24,18*\rr)--(24.5,17*\rr)--(24,16*\rr);
\draw [line width=0.2mm]  (23,20*\rr)--(24,20*\rr)--(24.5,19*\rr)--(24,18*\rr);
\draw [line width=0.2mm]  (24,20*\rr)--(24.5,21*\rr);

%\draw [line width=0.6mm] (12.5,19*\rr)--(7.5,13*\rr)--(14,10*\rr)--(12.5,19*\rr)
%--(19.5,15*\rr)--(14,10*\rr);

%\draw [line width=0.6mm] (9,4*\rr)--(14,10*\rr);
%\draw [line width=0.6mm] (9,4*\rr)--(15.5,1*\rr);
%\draw [line width=0.6mm] (21,6*\rr)--(14,10*\rr);
%\draw [line width=0.6mm] (21,6*\rr)--(15.5,1*\rr);

%\draw [line width=0.6mm] (21,6*\rr)--(22.5,-3*\rr);
%\draw [line width=0.6mm] (15.5,1*\rr)--(17,-8*\rr);
%\draw [line width=0.6mm] (9,4*\rr)--(10.5,-5*\rr);
%\draw [line width=0.6mm] (0.5,17*\rr)--(2,8*\rr)--(3.5,-1*\rr);

%\draw [line width=0.6mm]  (18,24*\rr)--(19.5,15*\rr);
%\draw [line width=0.6mm]  (6,22*\rr)--(7.5,13*\rr);
%\draw [line width=0.6mm] (12.5,19*\rr)--(11,28*\rr);
%\draw [line width=0.6mm] (3.5,-1*\rr)--(9,4*\rr)--(2,8*\rr)--(7.5,13*\rr)
%--(0.5,17*\rr)--(6,22*\rr)--(12.5,19*\rr)--(18,24*\rr);
%\draw [line width=0.6mm] (22.5,-3*\rr)--(15.5,1*\rr)--(10.5,-5*\rr);
%\draw [line width=0.6mm] (27.5,3*\rr)--(21,6*\rr)--(26,12*\rr)--(19.5,15*\rr)--(24.5,21*\rr);
%\draw [line width=0.6mm] (3.5,-1*\rr)--(-3,2*\rr)--(2,8*\rr)--(-4.5,11*\rr);

\draw [line width=0.6mm, color=green] (7.5,13*\rr)--(19.5,15*\rr)--(21,6*\rr);
\draw [line width=0.6mm, color=green] (7.5,13*\rr)--(9,4*\rr)--(21,6*\rr);
\foreach \pos in {(7.5,13*\rr),(19.5,15*\rr),(21,6*\rr),(9,4*\rr)}
\shade[shading=ball, ball color=green] \pos circle (.6);

\draw [line width=0.6mm, color=blue] (3.5,-1*\rr)--(15.5,1*\rr)--(14,10*\rr);
\draw [line width=0.6mm, color=blue] (3.5,-1*\rr)--(2,8*\rr)--(14,10*\rr);

\foreach \pos in {(14,10*\rr),(12.5,19*\rr),
(2,8*\rr),(0.5,17*\rr)}
\shade[shading=ball, ball color=blue] \pos circle (.6);
% \foreach \pos in {(13,18*\rr)}
% \shade[shading=ball, ball color=red] \pos circle (.3);

\foreach \pos in {(15.5,1*\rr),(3.5,-1*\rr)}
\shade[shading=ball, ball color=black] \pos circle (.6);

%\foreach \pos in {(10.5, 7*\rr),(0, 14*\rr),(12.5, 17*\rr)}
 %\shade[shading=ball, ball color=white] \pos circle (.4);

\draw (15,13*\rr) node   {\Large {${\mbox{\boldmath$O$}}$}};
% \draw (5.7,13.2*\rr) node   {\Large {${\mbox{\boldmath$A$}}$}};
% \draw (10,19*\rr) node   {\Large {${\mbox{\boldmath$B$}}$}};
%\draw (22,15.5*\rr) node   {\Large {${\mbox{\boldmath$C$}}$}};
%\draw (23,3*\rr) node   {\Large {${\mbox{\boldmath$D$}}$}};
%\draw (18.5,1*\rr) node   {\Large {${\mbox{\boldmath$E$}}$}};
\end{tikzpicture} 
%\end{center}
%\caption{}
\end{figure}} % Fig24A %XFigure24A

%%%%%%%%%%%%%%%%%%%%%%%%%%%%%%%%%%%%%%%%%%%%
%%%%%%%%%%%%%%%%%%%%%%%%%%%

\def\YFigure25{\begin{figure}\label{Fig25A} %\label{D2=25} 
\begin{center}
\begin{tikzpicture}[scale=0.2]
%\clip (-0.9, -1*\rr) 
\clip (0, -1.5*\rr) rectangle (27, 21*\rr);

%\path [fill=lightgray] (1,-1*\rr)-- (4,-1*\rr)--(10.7,21*\rr)--(7.5,21*\rr);
%\path [fill=cyan] (4.2,-1*\rr)--(7.4,-1*\rr)--(14,21*\rr)--(10.8,21*\rr);
%\path [fill=lightgray] (7.5,-1*\rr)--(10.6,-1*\rr)--(17.3,21*\rr)--(14.2,21*\rr);

\path [fill=lightgray]  (5,14*\rr)--(10.5,13*\rr)--(13.5,7*\rr)--(11,2*\rr)--(6,2*\rr)
--(3,8*\rr) ;
\draw [line width=0.6mm, color=gray] (5,14*\rr)--(3,8*\rr)--(8,8*\rr)
--(5,14*\rr)--(10.5,13*\rr)--(8,8*\rr)--(11,2*\rr)--(13.5,7*\rr)--(8,8*\rr)
--(6,2*\rr)--(11,2*\rr);

\draw [line width=0.6mm, color=gray] (6,2*\rr)--(3,8*\rr);

% \draw [line width=0.6mm, color=gray] (13.5,13*\rr)--(16.5,7*\rr);
\draw [line width=0.6mm, color=gray] (10.5,13*\rr)--(13.5,7*\rr);

\draw [line width=0.2mm] (0.5,1*\rr) -- (1.5,1*\rr) -- (2,0) -- (1.5,-1*\rr)
-- (0.5,-1*\rr) -- (0,0) -- (0.5,1*\rr);
\draw [line width=0.2mm] (0.5,1*\rr) -- (0,2*\rr) -- (0.5,3*\rr) -- (1.5,3*\rr) 
-- (2,2*\rr) -- (1.5,1*\rr);
\draw [line width=0.2mm] (0.5,3*\rr) -- (0,4*\rr) -- (0.5,5*\rr) -- (1.5,5*\rr) 
-- (2,4*\rr) -- (1.5,3*\rr);
\draw [line width=0.2mm] (0.5,5*\rr) -- (0,6*\rr) -- (0.5,7*\rr) -- (1.5,7*\rr) 
-- (2,6*\rr) -- (1.5,5*\rr);
\draw [line width=0.2mm] (0.5,7*\rr) -- (0,8*\rr) -- (0.5,9*\rr) -- (1.5,9*\rr) 
-- (2,8*\rr) -- (1.5,7*\rr);
\draw [line width=0.2mm] (0.5,9*\rr) -- (0,10*\rr) -- (0.5,11*\rr) -- (1.5,11*\rr) 
-- (2,10*\rr) -- (1.5,9*\rr);
\draw [line width=0.2mm] (0.5,11*\rr) -- (0,12*\rr) -- (0.5,13*\rr) -- (1.5,13*\rr) 
-- (2,12*\rr) -- (1.5,11*\rr);
\draw [line width=0.2mm] (0.5,13*\rr) -- (0,14*\rr) -- (0.5,15*\rr) -- (1.5,15*\rr) 
-- (2,14*\rr) -- (1.5,13*\rr);
\draw [line width=0.2mm] (0.5,15*\rr) -- (0,16*\rr) -- (0.5,17*\rr) -- (1.5,17*\rr) 
-- (2,16*\rr) -- (1.5,15*\rr);
\draw [line width=0.2mm] (0.5,17*\rr) -- (0,18*\rr) -- (0.5,19*\rr) -- (1.5,19*\rr) 
-- (2,18*\rr) -- (1.5,17*\rr);
\draw [line width=0.2mm] (0.5,19*\rr) -- (0,20*\rr) -- (0.5,21*\rr) -- (1.5,21*\rr) 
-- (2,20*\rr) -- (1.5,19*\rr);
%\draw [line width=0.2mm] (0.5,21*\rr) -- (0,22*\rr) -- (0.5,22*\rr) -- (1.5,22*\rr) 
%-- (2,22*\rr) -- (1.5,21*\rr);

\draw [line width=0.2mm]  (3,0) -- (3.5,-1*\rr);
\draw [line width=0.2mm] (2,2*\rr) -- (3,2*\rr) -- (3.5,1*\rr) -- (3,0) -- (2,0);
\draw [line width=0.2mm] (2,4*\rr) -- (3,4*\rr) -- (3.5,3*\rr) -- (3,2*\rr);
\draw [line width=0.2mm] (2,6*\rr) -- (3,6*\rr) -- (3.5,5*\rr) -- (3,4*\rr);
\draw [line width=0.2mm] (2,8*\rr) -- (3,8*\rr) -- (3.5,7*\rr) -- (3,6*\rr);
\draw [line width=0.2mm] (2,10*\rr) -- (3,10*\rr) -- (3.5,9*\rr) -- (3,8*\rr);
\draw [line width=0.2mm] (2,12*\rr) -- (3,12*\rr) -- (3.5,11*\rr) -- (3,10*\rr);
\draw [line width=0.2mm] (2,14*\rr) -- (3,14*\rr) -- (3.5,13*\rr) -- (3,12*\rr);
\draw [line width=0.2mm] (2,16*\rr) -- (3,16*\rr) -- (3.5,15*\rr) -- (3,14*\rr);
\draw [line width=0.2mm] (2,18*\rr) -- (3,18*\rr) -- (3.5,17*\rr) -- (3,16*\rr);
\draw [line width=0.2mm] (2,20*\rr) -- (3,20*\rr) -- (3.5,19*\rr) -- (3,18*\rr);
\draw [line width=0.2mm] (3.5,21*\rr) -- (3,20*\rr);

\draw [line width=0.2mm]  (3.5,1*\rr) -- (4.5,1*\rr) -- (5,0) -- (4.5,-1*\rr)
--(3.5,-1*\rr);
\draw [line width=0.2mm]  (3.5,3*\rr) -- (4.5,3*\rr) -- (5,2*\rr) -- (4.5,1*\rr);
\draw [line width=0.2mm]  (3.5,5*\rr) -- (4.5,5*\rr) -- (5,4*\rr) -- (4.5,3*\rr);
\draw [line width=0.2mm]  (3.5,7*\rr) -- (4.5,7*\rr) -- (5,6*\rr) -- (4.5,5*\rr);
\draw [line width=0.2mm]  (3.5,9*\rr) -- (4.5,9*\rr) -- (5,8*\rr) -- (4.5,7*\rr);
\draw [line width=0.2mm]  (3.5,11*\rr) -- (4.5,11*\rr) -- (5,10*\rr) -- (4.5,9*\rr);
\draw [line width=0.2mm]  (3.5,13*\rr) -- (4.5,13*\rr) -- (5,12*\rr) -- (4.5,11*\rr);
\draw [line width=0.2mm]  (3.5,15*\rr) -- (4.5,15*\rr) -- (5,14*\rr) -- (4.5,13*\rr);
\draw [line width=0.2mm]  (3.5,17*\rr) -- (4.5,17*\rr) -- (5,16*\rr) -- (4.5,15*\rr);
\draw [line width=0.2mm]  (3.5,19*\rr) -- (4.5,19*\rr) -- (5,18*\rr) -- (4.5,17*\rr);
\draw [line width=0.2mm]  (3.5,21*\rr) -- (4.5,21*\rr) -- (5,20*\rr) -- (4.5,19*\rr);

\draw [line width=0.2mm]  (5,0) -- (6,0)--(6.5,-1*\rr );
\draw [line width=0.2mm]  (5,2*\rr) -- (6,2*\rr)--(6.5,1*\rr )--(6,0);
\draw [line width=0.2mm]  (5,4*\rr) -- (6,4*\rr)--(6.5,3*\rr )--(6,2*\rr);
\draw [line width=0.2mm]  (5,6*\rr) -- (6,6*\rr)--(6.5,5*\rr )--(6,4*\rr);
\draw [line width=0.2mm]  (5,8*\rr) -- (6,8*\rr)--(6.5,7*\rr )--(6,6*\rr);
\draw [line width=0.2mm]  (5,10*\rr) -- (6,10*\rr)--(6.5,9*\rr )--(6,8*\rr);
\draw [line width=0.2mm]  (5,12*\rr) -- (6,12*\rr)--(6.5,11*\rr )--(6,10*\rr);
\draw [line width=0.2mm]  (5,14*\rr) -- (6,14*\rr)--(6.5,13*\rr )--(6,12*\rr);
\draw [line width=0.2mm]  (5,16*\rr) -- (6,16*\rr)--(6.5,15*\rr )--(6,14*\rr);
\draw [line width=0.2mm]  (5,18*\rr) -- (6,18*\rr)--(6.5,17*\rr )--(6,16*\rr);
\draw [line width=0.2mm]  (5,20*\rr) -- (6,20*\rr)--(6.5,19*\rr )--(6,18*\rr);
\draw [line width=0.2mm]  (6.5,21*\rr )--(6,20*\rr);

\draw [line width=0.2mm]  (6.5,1*\rr) -- (7.5,1*\rr)--(8,0)--(7.5,-1*\rr)
--(6.5,-1*\rr);
\draw [line width=0.2mm]  (6.5,3*\rr) -- (7.5,3*\rr)--(8,2*\rr)--(7.5,1*\rr);
\draw [line width=0.2mm]  (6.5,5*\rr) -- (7.5,5*\rr)--(8,4*\rr)--(7.5,3*\rr);
\draw [line width=0.2mm]  (6.5,7*\rr) -- (7.5,7*\rr)--(8,6*\rr)--(7.5,5*\rr);
\draw [line width=0.2mm]  (6.5,9*\rr) -- (7.5,9*\rr)--(8,8*\rr)--(7.5,7*\rr);
\draw [line width=0.2mm]  (6.5,11*\rr) -- (7.5,11*\rr)--(8,10*\rr)--(7.5,9*\rr);
\draw [line width=0.2mm]  (6.5,13*\rr) -- (7.5,13*\rr)--(8,12*\rr)--(7.5,11*\rr);
\draw [line width=0.2mm]  (6.5,15*\rr) -- (7.5,15*\rr)--(8,14*\rr)--(7.5,13*\rr);
\draw [line width=0.2mm]  (6.5,17*\rr) -- (7.5,17*\rr)--(8,16*\rr)--(7.5,15*\rr);
\draw [line width=0.2mm]  (6.5,19*\rr) -- (7.5,19*\rr)--(8,18*\rr)--(7.5,17*\rr);
\draw [line width=0.2mm]  (6.5,21*\rr) -- (7.5,21*\rr)--(8,20*\rr)--(7.5,19*\rr);

\draw [line width=0.2mm]  (8,0)--(9,0)--(9.5,-1*\rr);
\draw [line width=0.2mm]  (8,2*\rr)--(9,2*\rr)--(9.5,1*\rr)--(9,0);
\draw [line width=0.2mm]  (8,4*\rr)--(9,4*\rr)--(9.5,3*\rr)--(9,2*\rr);
\draw [line width=0.2mm]  (8,6*\rr)--(9,6*\rr)--(9.5,5*\rr)--(9,4*\rr);
\draw [line width=0.2mm]  (8,8*\rr)--(9,8*\rr)--(9.5,7*\rr)--(9,6*\rr);
\draw [line width=0.2mm]  (8,10*\rr)--(9,10*\rr)--(9.5,9*\rr)--(9,8*\rr);
\draw [line width=0.2mm]  (8,12*\rr)--(9,12*\rr)--(9.5,11*\rr)--(9,10*\rr);
\draw [line width=0.2mm]  (8,14*\rr)--(9,14*\rr)--(9.5,13*\rr)--(9,12*\rr);
\draw [line width=0.2mm]  (8,16*\rr)--(9,16*\rr)--(9.5,15*\rr)--(9,14*\rr);
\draw [line width=0.2mm]  (8,18*\rr)--(9,18*\rr)--(9.5,17*\rr)--(9,16*\rr);
\draw [line width=0.2mm]  (8,20*\rr)--(9,20*\rr)--(9.5,19*\rr)--(9,18*\rr);
\draw [line width=0.2mm]  (9.5,21*\rr)--(9,20*\rr);

\draw [line width=0.2mm]  (9.5,1*\rr)--(10.5,1*\rr)--(11,0)--(10.5,-1*\rr)
--(9.5,-1*\rr);
\draw [line width=0.2mm]  (9.5,3*\rr)--(10.5,3*\rr)--(11,2*\rr)--(10.5,1*\rr);
\draw [line width=0.2mm]  (9.5,5*\rr)--(10.5,5*\rr)--(11,4*\rr)--(10.5,3*\rr);
\draw [line width=0.2mm]  (9.5,7*\rr)--(10.5,7*\rr)--(11,6*\rr)--(10.5,5*\rr);
\draw [line width=0.2mm]  (9.5,9*\rr)--(10.5,9*\rr)--(11,8*\rr)--(10.5,7*\rr);
\draw [line width=0.2mm]  (9.5,11*\rr)--(10.5,11*\rr)--(11,10*\rr)--(10.5,9*\rr);
\draw [line width=0.2mm]  (9.5,13*\rr)--(10.5,13*\rr)--(11,12*\rr)--(10.5,11*\rr);
\draw [line width=0.2mm]  (9.5,15*\rr)--(10.5,15*\rr)--(11,14*\rr)--(10.5,13*\rr);
\draw [line width=0.2mm]  (9.5,17*\rr)--(10.5,17*\rr)--(11,16*\rr)--(10.5,15*\rr);
\draw [line width=0.2mm]  (9.5,19*\rr)--(10.5,19*\rr)--(11,18*\rr)--(10.5,17*\rr);
\draw [line width=0.2mm]  (9.5,21*\rr)--(10.5,21*\rr)--(11,20*\rr)--(10.5,19*\rr);

\draw [line width=0.2mm]  (11,0)--(12,0)--(12.5,-1*\rr);
\draw [line width=0.2mm]  (11,2*\rr)--(12,2*\rr)--(12.5,1*\rr)--(12,0);
\draw [line width=0.2mm]  (11,4*\rr)--(12,4*\rr)--(12.5,3*\rr)--(12,2*\rr);
\draw [line width=0.2mm]  (11,6*\rr)--(12,6*\rr)--(12.5,5*\rr)--(12,4*\rr);
\draw [line width=0.2mm]  (11,8*\rr)--(12,8*\rr)--(12.5,7*\rr)--(12,6*\rr);
\draw [line width=0.2mm]  (11,10*\rr)--(12,10*\rr)--(12.5,9*\rr)--(12,8*\rr);
\draw [line width=0.2mm]  (11,12*\rr)--(12,12*\rr)--(12.5,11*\rr)--(12,10*\rr);
\draw [line width=0.2mm]  (11,14*\rr)--(12,14*\rr)--(12.5,13*\rr)--(12,12*\rr);
\draw [line width=0.2mm]  (11,16*\rr)--(12,16*\rr)--(12.5,15*\rr)--(12,14*\rr);
\draw [line width=0.2mm]  (11,18*\rr)--(12,18*\rr)--(12.5,17*\rr)--(12,16*\rr);
\draw [line width=0.2mm]  (11,20*\rr)--(12,20*\rr)--(12.5,19*\rr)--(12,18*\rr);
\draw [line width=0.2mm]  (12,20*\rr)--(12.5,21*\rr);

\draw [line width=0.2mm]  (12.5,1*\rr)--(13.5,1*\rr)--(14,0)--(13.5,-1*\rr)
--(12.5,-1*\rr);
\draw [line width=0.2mm]  (12.5,3*\rr)--(13.5,3*\rr)--(14,2*\rr)--(13.5,1*\rr);
\draw [line width=0.2mm]  (12.5,5*\rr)--(13.5,5*\rr)--(14,4*\rr)--(13.5,3*\rr);
\draw [line width=0.2mm]  (12.5,7*\rr)--(13.5,7*\rr)--(14,6*\rr)--(13.5,5*\rr);
\draw [line width=0.2mm]  (12.5,9*\rr)--(13.5,9*\rr)--(14,8*\rr)--(13.5,7*\rr);
\draw [line width=0.2mm]  (12.5,11*\rr)--(13.5,11*\rr)--(14,10*\rr)--(13.5,9*\rr);
\draw [line width=0.2mm]  (12.5,13*\rr)--(13.5,13*\rr)--(14,12*\rr)--(13.5,11*\rr);
\draw [line width=0.2mm]  (12.5,15*\rr)--(13.5,15*\rr)--(14,14*\rr)--(13.5,13*\rr);
\draw [line width=0.2mm]  (12.5,17*\rr)--(13.5,17*\rr)--(14,16*\rr)--(13.5,15*\rr);
\draw [line width=0.2mm]  (12.5,19*\rr)--(13.5,19*\rr)--(14,18*\rr)--(13.5,17*\rr);
\draw [line width=0.2mm]  (12.5,21*\rr)--(13.5,21*\rr)--(14,20*\rr)--(13.5,19*\rr);

\draw [line width=0.2mm]  (14,0)--(15,0)--(15.5,-1*\rr);
\draw [line width=0.2mm]  (14,2*\rr)--(15,2*\rr)--(15.5,1*\rr)--(15,0);
\draw [line width=0.2mm]  (14,4*\rr)--(15,4*\rr)--(15.5,3*\rr)--(15,2*\rr);
\draw [line width=0.2mm]  (14,6*\rr)--(15,6*\rr)--(15.5,5*\rr)--(15,4*\rr);
\draw [line width=0.2mm]  (14,8*\rr)--(15,8*\rr)--(15.5,7*\rr)--(15,6*\rr);
\draw [line width=0.2mm]  (14,10*\rr)--(15,10*\rr)--(15.5,9*\rr)--(15,8*\rr);
\draw [line width=0.2mm]  (14,12*\rr)--(15,12*\rr)--(15.5,11*\rr)--(15,10*\rr);
\draw [line width=0.2mm]  (14,14*\rr)--(15,14*\rr)--(15.5,13*\rr)--(15,12*\rr);
\draw [line width=0.2mm]  (14,16*\rr)--(15,16*\rr)--(15.5,15*\rr)--(15,14*\rr);
\draw [line width=0.2mm]  (14,18*\rr)--(15,18*\rr)--(15.5,17*\rr)--(15,16*\rr);
\draw [line width=0.2mm]  (14,20*\rr)--(15,20*\rr)--(15.5,19*\rr)--(15,18*\rr);
\draw [line width=0.2mm]  (15,20*\rr)--(15.5,21*\rr);

\draw [line width=0.2mm]  (15.5,1*\rr)--(16.5,1*\rr)--(17,0)--(16.5,-1*\rr)
--(15.5,-1*\rr);
\draw [line width=0.2mm]  (15.5,3*\rr)--(16.5,3*\rr)--(17,2*\rr)--(16.5,1*\rr);
\draw [line width=0.2mm]  (15.5,5*\rr)--(16.5,5*\rr)--(17,4*\rr)--(16.5,3*\rr);
\draw [line width=0.2mm]  (15.5,7*\rr)--(16.5,7*\rr)--(17,6*\rr)--(16.5,5*\rr);
\draw [line width=0.2mm]  (15.5,9*\rr)--(16.5,9*\rr)--(17,8*\rr)--(16.5,7*\rr);
\draw [line width=0.2mm]  (15.5,11*\rr)--(16.5,11*\rr)--(17,10*\rr)--(16.5,9*\rr);
\draw [line width=0.2mm]  (15.5,13*\rr)--(16.5,13*\rr)--(17,12*\rr)--(16.5,11*\rr);
\draw [line width=0.2mm]  (15.5,15*\rr)--(16.5,15*\rr)--(17,14*\rr)--(16.5,13*\rr);
\draw [line width=0.2mm]  (15.5,17*\rr)--(16.5,17*\rr)--(17,16*\rr)--(16.5,15*\rr);
\draw [line width=0.2mm]  (15.5,19*\rr)--(16.5,19*\rr)--(17,18*\rr)--(16.5,17*\rr);
\draw [line width=0.2mm]  (15.5,21*\rr)--(16.5,21*\rr)--(17,20*\rr)--(16.5,19*\rr);

\draw [line width=0.2mm]  (17,0)--(18,0)--(18.5,-1*\rr);
\draw [line width=0.2mm]  (17,2*\rr)--(18,2*\rr)--(18.5,1*\rr)--(18,0);
\draw [line width=0.2mm]  (17,4*\rr)--(18,4*\rr)--(18.5,3*\rr)--(18,2*\rr);
\draw [line width=0.2mm]  (17,6*\rr)--(18,6*\rr)--(18.5,5*\rr)--(18,4*\rr);
\draw [line width=0.2mm]  (17,8*\rr)--(18,8*\rr)--(18.5,7*\rr)--(18,6*\rr);
\draw [line width=0.2mm]  (17,10*\rr)--(18,10*\rr)--(18.5,9*\rr)--(18,8*\rr);
\draw [line width=0.2mm]  (17,12*\rr)--(18,12*\rr)--(18.5,11*\rr)--(18,10*\rr);
\draw [line width=0.2mm]  (17,14*\rr)--(18,14*\rr)--(18.5,13*\rr)--(18,12*\rr);
\draw [line width=0.2mm]  (17,16*\rr)--(18,16*\rr)--(18.5,15*\rr)--(18,14*\rr);
\draw [line width=0.2mm]  (17,18*\rr)--(18,18*\rr)--(18.5,17*\rr)--(18,16*\rr);
\draw [line width=0.2mm]  (17,20*\rr)--(18,20*\rr)--(18.5,19*\rr)--(18,18*\rr);
\draw [line width=0.2mm]  (18,20*\rr)--(18.5,21*\rr);

\draw [line width=0.2mm]  (18.5,1*\rr)--(19.5,1*\rr)--(20,0)--(19.5,-1*\rr)
--(18.5,-1*\rr);
\draw [line width=0.2mm]  (18.5,3*\rr)--(19.5,3*\rr)--(20,2*\rr)--(19.5,1*\rr);
\draw [line width=0.2mm]  (18.5,5*\rr)--(19.5,5*\rr)--(20,4*\rr)--(19.5,3*\rr);
\draw [line width=0.2mm]  (18.5,7*\rr)--(19.5,7*\rr)--(20,6*\rr)--(19.5,5*\rr);
\draw [line width=0.2mm]  (18.5,9*\rr)--(19.5,9*\rr)--(20,8*\rr)--(19.5,7*\rr);
\draw [line width=0.2mm]  (18.5,11*\rr)--(19.5,11*\rr)--(20,10*\rr)--(19.5,9*\rr);
\draw [line width=0.2mm]  (18.5,13*\rr)--(19.5,13*\rr)--(20,12*\rr)--(19.5,11*\rr);
\draw [line width=0.2mm]  (18.5,15*\rr)--(19.5,15*\rr)--(20,14*\rr)--(19.5,13*\rr);
\draw [line width=0.2mm]  (18.5,17*\rr)--(19.5,17*\rr)--(20,16*\rr)--(19.5,15*\rr);
\draw [line width=0.2mm]  (18.5,19*\rr)--(19.5,19*\rr)--(20,18*\rr)--(19.5,17*\rr);
\draw [line width=0.2mm]  (18.5,21*\rr)--(19.5,21*\rr)--(20,20*\rr)--(19.5,19*\rr);

\draw [line width=0.2mm]  (20,0)--(21,0)--(21.5,-1*\rr);
\draw [line width=0.2mm]  (20,2*\rr)--(21,2*\rr)--(21.5,1*\rr)--(21,0);
\draw [line width=0.2mm]  (20,4*\rr)--(21,4*\rr)--(21.5,3*\rr)--(21,2*\rr);
\draw [line width=0.2mm]  (20,6*\rr)--(21,6*\rr)--(21.5,5*\rr)--(21,4*\rr);
\draw [line width=0.2mm]  (20,8*\rr)--(21,8*\rr)--(21.5,7*\rr)--(21,6*\rr);
\draw [line width=0.2mm]  (20,10*\rr)--(21,10*\rr)--(21.5,9*\rr)--(21,8*\rr);
\draw [line width=0.2mm]  (20,12*\rr)--(21,12*\rr)--(21.5,11*\rr)--(21,10*\rr);
\draw [line width=0.2mm]  (20,14*\rr)--(21,14*\rr)--(21.5,13*\rr)--(21,12*\rr);
\draw [line width=0.2mm]  (20,16*\rr)--(21,16*\rr)--(21.5,15*\rr)--(21,14*\rr);
\draw [line width=0.2mm]  (20,18*\rr)--(21,18*\rr)--(21.5,17*\rr)--(21,16*\rr);
\draw [line width=0.2mm]  (20,20*\rr)--(21,20*\rr)--(21.5,19*\rr)--(21,18*\rr);
\draw [line width=0.2mm]  (21,20*\rr)--(21.5,21*\rr);

\draw [line width=0.2mm]  (21.5,1*\rr)--(22.5,1*\rr)--(23,0)--(22.5,-1*\rr)
--(21.5,-1*\rr);
\draw [line width=0.2mm]  (21.5,3*\rr)--(22.5,3*\rr)--(23,2*\rr)--(22.5,1*\rr);
\draw [line width=0.2mm]  (21.5,5*\rr)--(22.5,5*\rr)--(23,4*\rr)--(22.5,3*\rr);
\draw [line width=0.2mm]  (21.5,7*\rr)--(22.5,7*\rr)--(23,6*\rr)--(22.5,5*\rr);
\draw [line width=0.2mm]  (21.5,9*\rr)--(22.5,9*\rr)--(23,8*\rr)--(22.5,7*\rr);
\draw [line width=0.2mm]  (21.5,11*\rr)--(22.5,11*\rr)--(23,10*\rr)--(22.5,9*\rr);
\draw [line width=0.2mm]  (21.5,13*\rr)--(22.5,13*\rr)--(23,12*\rr)--(22.5,11*\rr);
\draw [line width=0.2mm]  (21.5,15*\rr)--(22.5,15*\rr)--(23,14*\rr)--(22.5,13*\rr);
\draw [line width=0.2mm]  (21.5,17*\rr)--(22.5,17*\rr)--(23,16*\rr)--(22.5,15*\rr);
\draw [line width=0.2mm]  (21.5,19*\rr)--(22.5,19*\rr)--(23,18*\rr)--(22.5,17*\rr);
\draw [line width=0.2mm]  (21.5,21*\rr)--(22.5,21*\rr)--(23,20*\rr)--(22.5,19*\rr);

\draw [line width=0.2mm]  (23,0)--(24,0)--(24.5,-1*\rr);
\draw [line width=0.2mm]  (23,2*\rr)--(24,2*\rr)--(24.5,1*\rr)--(24,0);
\draw [line width=0.2mm]  (23,4*\rr)--(24,4*\rr)--(24.5,3*\rr)--(24,2*\rr);
\draw [line width=0.2mm]  (23,6*\rr)--(24,6*\rr)--(24.5,5*\rr)--(24,4*\rr);
\draw [line width=0.2mm]  (23,8*\rr)--(24,8*\rr)--(24.5,7*\rr)--(24,6*\rr);
\draw [line width=0.2mm]  (23,10*\rr)--(24,10*\rr)--(24.5,9*\rr)--(24,8*\rr);
\draw [line width=0.2mm]  (23,12*\rr)--(24,12*\rr)--(24.5,11*\rr)--(24,10*\rr);
\draw [line width=0.2mm]  (23,14*\rr)--(24,14*\rr)--(24.5,13*\rr)--(24,12*\rr);
\draw [line width=0.2mm]  (23,16*\rr)--(24,16*\rr)--(24.5,15*\rr)--(24,14*\rr);
\draw [line width=0.2mm]  (23,18*\rr)--(24,18*\rr)--(24.5,17*\rr)--(24,16*\rr);
\draw [line width=0.2mm]  (23,20*\rr)--(24,20*\rr)--(24.5,19*\rr)--(24,18*\rr);
\draw [line width=0.2mm]  (24,20*\rr)--(24.5,21*\rr);

\draw [line width=0.5mm] (17,18*\rr)--(21.5,15*\rr)--(21.5,9*\rr)
--(17,6*\rr)--(12.5,9*\rr)--(12.5,15*\rr)--(17,18*\rr);
\draw [line width=0.5mm] (17,18*\rr)--(17,6*\rr);
\draw [line width=0.5mm] (21.5,15*\rr)--(12.5,9*\rr);
\draw [line width=0.5mm] (21.5,9*\rr)--(12.5,15*\rr);

\foreach \pos in {(5.5,13*\rr)}
\shade[shading=ball, ball color=red] \pos circle (.4);

\foreach \pos in {(3,8*\rr),(8,8*\rr),(5,14*\rr),(10.5,13*\rr),(13.5,7*\rr),
(11,2*\rr),(6,2*\rr)}
\shade[shading=ball, ball color=gray] \pos circle (.6);

\foreach \pos in {(3.5,3*\rr),(3.5,9*\rr),(3.5,15*\rr),(3.5,21*\rr),
(8,0*\rr),(8,6*\rr),(8,12*\rr),(8,18*\rr),
(12.5,3*\rr),(12.5,9*\rr),(12.5,15*\rr),(12.5,21*\rr),
(17,0*\rr),(17,6*\rr),(17,12*\rr),(17,18*\rr),(21.5,3*\rr),(21.5,9*\rr),
(21.5,15*\rr),(21.5,21*\rr)}
\shade[shading=ball, ball color=black] \pos circle (.6);

% \foreach \pos in {(9,4*\rr),(15.5,1*\rr),(21,6*\rr),(3.5,-1*\rr)}
% \shade[shading=ball, ball color=black] \pos circle (.6);

%\foreach \pos in {(10.5, 7*\rr),(0, 14*\rr),(12.5, 17*\rr)}
 %\shade[shading=ball, ball color=white] \pos circle (.4);

% \draw (12,6*\rr) node   {\Large {${\mbox{\boldmath$O$}}$}};
% \draw (5.7,6*\rr) node   {\Large {${\mbox{\boldmath$A$}}$}};
% \draw (7,16*\rr) node   {\Large {${\mbox{\boldmath$B$}}$}};
% \draw (16,15*\rr) node   {\Large {${\mbox{\boldmath$C$}}$}};
 
\end{tikzpicture} 
\end{center}
%\caption{}
\end{figure}} % Fig25 YFigure25

%%%%%%%%%%%%%%%%%%%%%%%%%%%%%%%%%%%%%%%%%%%%
%%%%%%%%%%%%%%%%%%%%%%%%%%%

\def\ZFigure26  %%\label{Fig26A}
{\begin{figure} \label{Fig26A} %subClass \tE3
\begin{center} 



\end{center}%\vskip .5cm
%\caption{}
\vskip .2cm

\end{figure}} %% {Fig26A} ZFigure26
%%%%%%%%%%%%%%%%%%%%%%%%%%%%%%%%%%%%%%%%%%%%%%%%%%%%%%%%%%%%%%%%%%%%%%%

%%%%%%%%%%%%%%%%%%%%%%%%%%%%%%%%%%%%%%%%%%%%
%%%%%%%%%%%%%%%%%%%%%%%%%%%

\title{\bf High-density hard-core model\\
on triangular and hexagonal lattices}

\author{\bf A. Mazel$^1$, I. Stuhl$^2$, Y. Suhov$^{2-4}$}

\date{}
\footnotetext{2010 {\em Mathematics Subject Classification:\; primary 60G60, 82B20, 82B26}}
\footnotetext{{\em Key words and phrases:} triangular lattice, hexagonal lattice, hard-core 
configuration, disk-packing, extreme Gibbs measure, high-density/large fugacity, periodic ground state, 
Delaunay triangulation, minimal re-distributed area of a triangle, maximally-dense sub-lattice, 
maximally-dense non-sub-lattice configuration, contour representation of the partition function, 
Peierls bound, Pirogov-Sinai theory, dominance, local repelling forces, computer-assisted 
enumeration, sliding

\noindent
$^1$ AMC Health, New York, NY, USA;\;\;
$^2$ Math Dept, Penn State University, PA, USA;,\;\;
$^3$ DPMMS, University of Cambridge and St John's College, Cambridge, UK,\;\;
$^4$ IITP RAS, Moscow, RF}

\maketitle

\begin{abstract}
We perform a rigorous study of the Gibbs statistics of high-density hard-core 
random configurations on a unit triangular lattice $\bbA_2$ and a unit
honeycomb graph $\bbH_2$, for any value of the (Euclidean) repulsion 
diameter $D>0$. Only attainable values of $D$ are relevant,
for which $D^2=a^2+b^2+ab$, $a, b \in\mathbb{Z}$ (L\"oschian numbers).
Depending on arithmetic properties of $D^2$, we identify, 
for large fugacities, the pure phases (extreme Gibbs measures) and
specify their symmetries. The answers depend on the way(s) an
equilateral triangle of side-length $D$ can be inscribed in $\bbA_2$
or $\bbH_2$. On $\bbA_2$, our approach works for all
attainable $D^2$; 
on $\bbH_2$ we have to exclude $D^2 = 4, 7, 31, 133$, where a sliding 
phenomenon occurs, similar to that on a unit square lattice $\bbZ^2$.  
For all values $D^2$ apart from the excluded ones we prove the existence
of a first-order phase transition where the number of co-existing pure phases
grows at least as $O(D^2)$.

The proof is based on the Pirogov--Sinai theory which requires non-trivial verifications of 
key assumptions: finiteness of the set of periodic ground states and the Peierls bound. 
To establish the Peierls bound, we develop a general method based on the concept of a
re-distributed area for Delaunay triangles. Some of the presented proofs are 
computer-assisted. 

As a by-product of the ground state identification, we solve the disk-packing problem on  
$\bbA_2$ and $\bbH_2$ for any value of the disk diameter $D$.
\end{abstract}

\section{A summary of results}\label{Sec1}

\subsection{Introduction}\label{SubSec1.1} We analyze properties of random configurations 
of hard disks of a 
given diameter $D$, with centers in a unit triangular lattice $\bbA_2$ and 
a unit honeycomb lattice$^{1)}$\footnote{$^{1)}$Strictly speaking,  $\bbH_2$ is 
not a lattice in the 
algebraic sense. However, we follow a physical tradition and refer to $\bbH_2$
as a lattice.} $\bbH_2$, both embedded in $\bbR^2$. It is also convenient to consider
$\bbH_2$ as a subset in $\bbA_2$. Cf. Figure 1. Together
with a unit square lattice $\bbZ^2$, these are popular examples of `regular' planar graphs
for a number of probabilistic models, including percolation and phase transitions. A separate 
place belongs to a model of hard disks in $\bbR^2$. Historically,
the hard-core model emerged about a 150 years ago in an attempt to describe
a system of atoms, molecules or granules, as represented by rigid spheres of a given
diameter;
a famous example of its application was the Boltzmann equation. Since then, the model
proliferated in a number of pure and applied mathematical disciplines and generated
a substantial literature. A comprehensive discussion of various aspects of the hard-core
model and its applications (including elements of criticism) can be found, e.g.,
in \cite{CKPUZ}, \cite{PaZ}, \cite{PeS}, \cite{BKZZ}, \cite{KMRTSZ}.

The study of lattice hard-core (H-C) models started with the result by Dobrushin \cite{Dob}
about non-uniqueness of pure phases on $\bbZ^d$, $d\geq 2$, with a nearest-neighbor exclusion in a high-density/large fugacity regime.
The paper \cite{HeP} established non-uniqueness of a pure phase for a particular
sequence of exclusion distances
on $\bbA_2$, without specifying the pure phases. We also note the (rather remarkable) result of paper
\cite{Ba} where a critical value of fugacity  $\diy u_{\rm{cr}}=\frac{1}{2}(11+5{\sqrt 5})$ has been calculated,
for the H-C diameter $D={\sqrt 3}$ on $\bbA_2$. This, apparently, indicates an upper limiting value for 
fugacity $u$ for the low-density regime where a pure phase is unique and given via a polymer expansion
around an empty configuration. The paper \cite{JaL} establishes existence of order-disorder phase transitions for a class of `non-sliding'
H-C lattice particle systems on a lattice in two or more dimensions. 

\vskip-35pt
\FigureA1

%\indent\hskip .6cm{\bf Figure 1. Fragments of lattices ${\mathbb A}_2$ and ${\mathbb H}_2$}

The present paper continues and extends the works \cite{Dob} and \cite{HeP} in a general setting.
A detailed study of the H-C model on $\bbZ^2$ has been performed in \cite{MSS1}.
We analyze the {\it ground states} and {\it Gibbs} or {\it DLR} measures for the (H-C) model
on $\bbA_2/\bbH_2$, in a regime of high-density/large-fugacity. The assumption that
the fugacity $u$ is large is adopted throughout the paper without stressing it every time again. The analysis of
Gibbs measures is reduced to {\it extreme Gibbs measure} (EGM, or $D$-EGM 
when dependence on $D$ is emphasized). An EGM is interpreted as a {\it pure phase} in the 
phase diagram of the model. Formal definitions of the notions used in Introduction
are provided in Sections \ref{Sec2}.  The H-C exclusion is imposed in the Euclidean 
$\bbR^2$-metric and is defined by the H-C {\it exclusion diameter} $D$: the shortest 
allowed distance between two occupied sites. Without loss of generality we assume 
throughout the paper that the value 
$D$ (or $D^2$) is {\it attainable}, i.e., there are pairs of sites in $\bbA_2$ and $\bbH_2$ with the 
distance exactly $D$ between them. The set of attainable values is the same for
$\bbA_2$ and $\bbH_2$ and is characterized through a {\it L\"oschian decomposition} of 
the number $D^2$. Referring to $D^2$ rather than to $D$ is 
more convenient since $D^2$ is a positive integer. 

The problem of identification of the EGM structure is reduced -- via the {\it Pirogov--Sinai} 
(PS) {\it theory} \cite{PiS}, \cite{Za} -- to an analysis of {\it periodic ground states} 
(PGSs or $D$-PGSs), including a verification of the {\it Peierls bound}. 
Informally, speaking, the outcome of the PS theory is that every EGM is generated
by a PGS. The inverse is not always true: there may be PGSs that do not generate 
EGMs. The PGSs that generate EGMs are referred to as {\it dominant} (stable 
in the terminology of \cite{Za}). 

\medskip
Our results can be briefly summarized as follows. 
\begin{description}
\item[{\rm{(i)}}] On both $\bbA_2$ and $\bbH_2$ we describe the grounds states 
(periodic and non-periodic) for all values of $D$ and fugacity $u>1$.  See Remark 4.1. 
This solves the disk-packing problem on 
$\bbA_2$ and $\bbH_2$ for any disk diameter $D$;  cf.\cite{CD}. The PGSs are 
naturally partitioned
into equivalence classes defined by lattice symmetries (shifts and reflections).
Apart from 13 values of $D$ on $\bbH_2$, all $D$-PGSs are constructed from
sub-lattices.

\item[{\rm{(ii)}}] For the PGSs we establish a Peierls bound where the Peierls constant 
grows with $u$ and decreases with $D$. See Lemmas \ref{Lem5.1} and
\ref{Lem5.4}.

\item[{\rm{(iii)}}] The structure of $D$-EGMs (the phase diagram) inherits that of $D$-PGSs.
First, for at least one PGS-equivalence class, each PGS from the class generates 
a distinct EGM. That is, we have a first-order phase transition.
This fact is proven for every $D$, except for 4 values of $D$ on $\bbH_2$. See Theorem III in Section \ref{SubSec3.2}.

\item[{\rm{(iv)}}] In the case where the PGS-equivalence class is unique, we obtain a 
complete phase diagram. The sets of values $D$ with a unique  
PGS-equivalence class on $\bbA_2$ and $\bbH_2$ are infinite and explicitly
described. The number of $D$-EGMs in this case grows as $O(D^2)$ and
is further specified. See Theorems 1, 2, 7, 8, 11, 12 in Section \ref{Sec3}.

\item[{\rm{(v)}}] The sets of values $D$  for which the PGS-equivalence class is non-unique
is also infinite and explicitly described, on both $\bbA_2$ and $\bbH_2$. 
The number of $D$-EGMs in this case grows at least as $O(D^2)$. 
The question which classes are dominant, i.e., generate EGMs requires an additional 
analysis. We conduct such analysis on a number of values $D$ on $\bbA_2$ and $\bbH_2$
exploring various emerging possibilities. See Theorems  4, 5, 6, 10 in Section \ref{Sec3}.

\item[{\rm{(vi)}}] We establish that a phenomenon of {\it sliding} occurs only on $\bbH_2$,
for $D^2= 4, 7, 31, 133$. See Lemmas 4.7 - 4.10. For these values, 
the structure of the phase diagram remains
open. Cf. Section \ref{Sec8}. The phenomenon of sliding was first discovered by Dobrushin 
(1968) on $\bbZ^2$. Cf. \cite{MSS1}.
\end{description}

\subsection{The PGSs and EGMs}\label{SubSec1.2} 

The structure of PGSs on both $\bbA_2$ and $\bbH_2$ (and also on $\bbZ^2$: cf. \cite{MSS1}) 
depends on arithmetic properties of the number $D^2$. Moreover, the image of a $D$-PGS 
under a {\it symmetry} (of $\bbA_2$ or $\bbH_2$), 
i.e., a lattice shift or reflection, is also a $D$-PGS, and one can speak about the corresponding 
{\it equivalence classes} of PGSs with respect to $\bbA_2/\bbH_2$-symmetries. This is 
important since dominance is a class property: if an equivalence class contains a dominant
PGS then all PGSs from the class are dominant.  

Referring to the formation of PGS-equivalence classes, the entire set of attainable values of
$D$ (or $D^2$) is divided into disjoint subsets. For lattice $\bbA_2$ we consider
two subsets of values $D^2$ (both infinite), called Classes TA and TB. On $\bbH_2$ we deal 
with six subsets of values $D^2$ called Classes HA, HB, HC (infinite) and 
HD, HE, HS (finite). Here T stands for triangular and H for honeycomb. These subsets 
are further divided, regarding specific aspects of the structure of PGSs and EGMs. 
Cf. Sections \ref{SubSec1.3}, \ref{SubSec1.4}.

Physically speaking, the above subsets are characterized by a possibility (or possibilities or a lack 
of them) to inscribe an equilateral triangle of side-length $D$ in $\bbA_2$ or $\bbH_2$. In the case 
of $\bbA_2$ this is always possible, but the inscription may be non-unique. For $\bbH_2$ it is not 
always possible, which leads to a more complicated partition of the values $D^2$.
 
On  $\bbA_2$, the $D$-PGSs are constructed from $D$-{\it sub-lattices}, i.e., sub-lattices
for which a fundamental parallelogram is a $D$-rhombus formed by $2$ equilateral triangles 
of side-length $D$ ($D$-{\it triangles}, for short). If $D^2$ is from Class TA, the
$D$-sub-lattice is unique, and the PGSs  
form a single equivalence class. Consequently, for Class TA we establish a complete phase 
diagram in the large-fugacity regime. In this regime each PGS applied as a boundary condition 
generates a distinct $D$-EGM, and all $D$-EGMs are obtained in this way; cf. Theorems 1, 2 in
Section \ref{SubSec3.3}. In other words, for values $D$ from Class TA all PGSs are dominant. Also, the 
EGMs inherit  symmetries between their generating PGSs. Similar properties hold true for 
Class HA on $\bbH_2$; cf. Theorems 7, 8 in Section \ref{SubSec3.5}. 
Classes TA and HA yield the simplest cases of the PGS/EGM analysis on $\bbA_2/\bbH_2$.

For $D^2$ from Classes TB or HB, the $D$-PGSs are still constructed from $D$-sub-lattices, 
but there are multiple PGS-equivalence classes. There is always a dominant 
equivalence class (we conjecture that it is unique), but the problem of identifying which 
classes are dominant is more involved. Here we solve it for some specific values of $D$, 
indicating various emerging possibilities for a plausible general answer. 

Arithmetically, Classes HA and HB are formed by the values $D^2$ from Classes TA and TB 
divisible by 3.

Next, Classes HC, HD, HE, and HS on $\bbH_2$ consist of values $D^2$ non-divisible by $3$.
These classes stem from the above-mentioned features of 
$\bbH_2$, about a luck of possibility to inscribe a $D$-triangle for some values of $D$. 
Class HC (which covers a bulk of values of $D$ on $\bbH_2$) is determined by the condition that
$D^2$ is not divisible by 3 and is non-exceptional in a sense made precise below. For values $D$ from 
Class HC we look for
the attainable $D^*>D$ which (i) has $(D^*)^2$ divisible by 3 (i.e., falls in Class HA or HB) and (ii) 
$D^*$ is nearest to $D$
with this property. Then the PGSs and EGMs for $D^2$ are the same as for $(D^*)^2$, i.e., 
are constructed from $D^*$-sub-lattices. 

Classes HD, HE and HS are deemed exceptional and are dealt with on a case-by-case basis
(with the help of a computer). Cf. Section \ref{SubSec1.4}. For these classes not all PGSs are  
constructed from $D$- or $D^*$-sub-lattices. (For 
$D^2$ from Class HD none of the PGSs is a sub-lattice.) In particular, Class HS consists
of values $D^2= 4, 7, 31, 133$ which exhibit a phenomenon of  sliding on $\bbH_2$. A similar
phenomenon occurs on lattice $\bbZ^2$ as well; cf. \cite{MSS1}. For the values of $D^2$ with 
sliding, the PS theory is not applicable since the number of PGSs is
infinite and -- more importantly -- the Peierls bound does not hold.
Our conjecture is that the EGM for these 
values of $D$ is unique when fugacity $u$ is large enough (and, indeed, for all values $u>0$). 
We briefly comment on sliding on $\bbH_2$ in Section \ref{Sec8}.
Cf. Section 2.2 in \cite{MSS1} where the similar problem is treated on lattice $\bbZ^2$.
 
%As a by-product, the above results on PGSs solve the problem of close-packing of disks 
%of diameter $D$ on both $\bbA_2$ (cf. [Con]) and $\bbH_2$, for any given $D>1$.
%$\bbA_2$ these are always triangular $D$-sub-lattices and their shifts; ]. On $\bbH_2$
%these configurations are triangular sub-lattices when $D^2$ is from Classes HA, HB, HC.
%For other classes they are described for each value $D^2$ individually.} 

As was mentioned before, an attainable value $D^2$ admits a L\"oschian decomposition. 
It means that $D^2$ is a positive integer of the form $a^2 + b^2 + ab$ where $a$ and $b$ are integers. 
L\"oschian numbers arise naturally in this context as they are the norms of Eisenstein integers that form $\mathbb{A}_2$. 
The L\"{o}schian numbers are the sequence A003136 in OEIS, the 
on-line Encyclopedia of integer sequences and their initial list is:
$$\beac 1,3,4,7,9,12,13,16,19,21,25,27,28,31,36,37,39,43,48,49,52,57,61,63,64,67,73,\\
75,76,79,81,84,91,93,97,100,103,108,109,111,112,117,121,124,127,129,133,139,\\
144,147,148,151,156,157,163,169,171,172,175,181,183,189,192,193,196,199,201,\\
208,211,217,219,223,225,228,229,237,241,243,244,247,252,256,259,268,271,273, \\
277,279,283,289,291,292,300,301,304,307,309,313,316,324,325,327,331,333,336,\\ 
337,343,351,361,363,364,367,372,379,387,388,397,399,400,403,409,412,417,421,\\ 
427,432,433,436,439,441,444,448,453,457,463,468,469,471,475,481,484,487,489. \ena$$

%,496,499,507,508,511,513,516,523,525,529,532,541,
%543,547,549,553,556,559,\\ %567,571,576,577,579,588,589,%592,597,601,
%603,604,607,613,619,624,625,628,631,\\ %633,637,643,651,652,657,661,
%669,673,675,676,679,684,687,688,691,700,703,709,711,721,\\ %723,724,727,
%729,732,733,739,%741,751,756,757,763,768,769,772,775,777,784,787,793,796,804,811,\\
%813,817,819,823,829,831,832,837,841,844,847,849,853,859,867,868,871,873,876,877,\\
%883,889,892,900. \ena$$

An equivalent characterization of a L\"oschian number is that its rational prime 
factorization must contain primes of the form $3v+2$, $v\in\bbZ_+$ in even 
powers (there is no restriction on factor 3 or primes of the form $3v+1$).
We use some classical results regarding these integers which are presented in a 
convenient form in \cite{M, N}. A compendium of the related theory is given in the 
monograph \cite{CS}. 

%https://www.calculatorsoup.com/calculators/math/prime-factors.php prime decompositions

%The structure of $D$-EGMs on both $\bbA_2$ and $\bbH_2$ is
%connected to arithmetic properties of $D^2$. 
%The H-C model on $\bbZ^2$ also 
%exhibits deep ties with arithmetic properties of $D^2$; see \cite{MSS1}.
%For $\bbA_2$ the connection is rather straightforward whereas for $\bbH_2$ it 
%involves additional technicalities. As was mentioned, by virtue of the Pirogov--Sinai (PS) theory,
%the problem of analysis of the EGM structure is reduced to an identification of dominant
%PGSs and verification of a Peierls bound between them. For the H-C model, the PGSs are 
%represented by periodic admissible configurations of maximal density. 

The PGS identification and a specification of the Peierls-bound is done with the help of {\it Voronoi 
cells} (V-cells) or through the construction of Delaunay $\bbA_2/\bbH_2$-triangles which {\it minimize the 
re-distributed area} (MRA-triangles). As was said, the Peierls bound is given in 
Theorem II from Section \ref{SubSec3.2}. In this paper we employ the approach based on MRA-triangles (also 
used in \cite{MSS1}), but for completeness provide a brief account of the V-cell method as well. Cf. 
Sections \ref{Sec4}, \ref{Sec5}. 

%A concise statement about PGSs on $\bbA_2/\bbH_2$ is presented in Theorem I from 
%Section 3.2. It summarizes the material from Sections 1.2 and 1.3.

In Sections \ref{SubSec1.3} and \ref{SubSec1.4} we give a description of our results on PGSs for the H-C model 
on $\bbA_2$ and $\bbH_2$.

\subsection{PGSs and EGMs on $\bbA_2$}\label{SubSec1.3}

On $\bbA_2$ the situation is made easier by the above-mentioned fact that the PGSs are constructed 
from $D$-sub-lattices. That is, a PGS-equivalence class is determined either by 
a $D$-sub-lattice -- if it is reflection-invariant -- or by a pair of $D$-sub-lattices taken to 
each other by a reflection. Consequently, a PGS-class contains $D^2$ or $2D^2$ PGSs obtained
from each other by lattice shifts and reflections. The value $D^2$ represents the number of 
sites in a $D$-rhombus which gives the number of different lattice shifts for PGSs.

The number of PGS-equivalence classes is related to the number and structure of non-negative solutions to the equation
$D^2=a^2+b^2+ab$. Accordingly, it is natural to extract the following classes of values of $D^2$. 

$\bullet$ Class TA1: $D^2$ is an integer whose prime decomposition contains (i) a factor
 $3$ in any power, (ii)  primes of the form $3v+2$, in even powers, possibly zero, 
 and (iii) no prime of the form $3v+1$. This happens iff $D^2=a^2$ or $3a^2$ where $a\in\bbN$
 has only primes $3v+2$ in its prime decomposition. 
%Either $D$ or $D/\sqrt{3}$ is an integer whose prime decomposition 
%does not contain factors of the form $3v+1$. 
%The list of primes of the form $3v+1$ up to 449: 7,13,19,31,37,43,61,67,73,79,91,97,103,109,127,
%139,151,157,%163,181,193,199,211,217,223,229,241,247,271,277,283,307,313,319,325,331,337,
%343,349,355,361,%367,373,379,391,397,409,421,433,439,451,457
%In other words, each prime factor of $D^2$ is either $3$, entering in any power, or of 
%the form $3v+2$, entering in an even power.
The first 40 values of $D^2$ falling in this category are $D^2 = $ 
1, 3, 4, 9, 12, 16, 25, 27, 36, 48, 
64, 75, 81, 100, 108, 121, 144, 192, 225, 243, 
256, 289, 300, 324, 363, 400, 432, 484, 529, 576, 
625, 675, 729, 768, 841, 867, 900, 972, 1024, 1089. 
%1156, 1200, 1296, 1452, 1587, 1600, 1681, 1728, 1875, 1936.

$\bullet$ Class TA2: $D^2$ is an integer whose prime decomposition contains (i) a factor 
3 in any power, (ii) primes of the form $3v+2$, in even powers, possibly zero, and (iii)
a single prime of the form $3v+1$ (entering in power 1). This happens iff $D^2$ admits a unique
decomposition as $a^2+b^2+ab$ (modulo the permutation of $a$ and $b$) and we have $a,b\in\bbN$,
$a\neq b$. 
%(The factor $3$ can enter with any power, including zero.) 
The first 40 values of
$D^2$ from Class TA2 are 7, 13, 19, 21, 28, 31, 37, 39, 43, 52, 57, 61, 63, 67, 73, 76, 
79, 84, 93, 97, 103, 109, 111, 112, 117, 124, 127, 129, 139, 148, 151, 156, 157, 163, 
171, 172, 175, 181, 183, 189.  %193, 199, 201, 208, 211, 219, 223, 228, 229, 237.

$\bullet$ Class TA: the union of Classes TA1 and TA2. 

$\bullet$ Class TB: All remaining attainable values of $D$ (or $D^2$). Class TB consists of 
positive integers $D^2$ that contain (i) a factor 3 in any power, (ii)  primes of the form $3v+2$, 
in even powers, possibly zero, and (iii) at least two primes of the form $3v+1$
(possibly, identical).  It occurs iff $D^2$ admits a non-unique L\"oschian decomposition, i.e., 
there are more than 1 solutions $(a,b)$ to the Diophantine equation $D^2=a^2+b^2+ab$, with
$a, b$ non-negative integers, again, modulo the permutation of $a$, $b$. The first 40 values of $D^2$ from
Class TB are
49, 91, 133, 147, 169, 196, 217, 247, 259, 273, 
301, 343, 361, 364, 399, 403, 427, 441, 469, 481, 
507, 511, 532, 553, 559, 588, 589, 637, 651, 676, 
679, 703, 721, 741, 763, 777, 784, 793, 817, 819.
%868, 871, 889, 903, 931. % 903, 931, 949, 961, 973, 988. 

%Correspondingly, we refer to the value of $D$ of type TA1, TA2, A and TB.
%\vskip .3cm

For $D^2$ from TA there is a single PGS-equivalence class, while for $D^2$
from TB the number of  PGS-equivalence classes is greater than one. 

\vskip-40pt
\FigureB2

%\indent\hskip .6cm{\bf Figure 2. PGSs on ${\mathbb A}_2$ for $D^2=9$ (Class TA1), 
%and $D^2=13$ (Class TA2),}
%\indent\hskip .6cm{\bf  with V-cells (gray hexagons) and choices of $D$-rhombuses.}\\
%\indent\hskip .6cm{\footnotesize{The number of PGSs and EGMs for $D^2=9$ is 9 (frame (a)), and
%for $D^2=13$ is 26 (frame (b))}}
%\indent\hskip .6cm{\footnotesize{The PGSs for $D^2=9$ are horizontal and for $D^2=13$ inclined.}}

It is convenient to refer to triangles with vertices in $\bbA_2$ as $\bbA_2$-triangles. 
As we saw earlier, an important role is played by $D$-triangles. We will 
distinguish between 3 types of $D$-triangles: {\it horizontal}, with sides fitting $\bbA_2$,
{\it vertical}, with sides perpendicular to constituent lines of $\bbA_2$, and {\it inclined},
covering the remaining cases. We will also use $D$-triangles in the plane $\bbR^2$, for 
general diameters $D>0$. The above terminology is extended to the
$D$-sub-lattices and $D$-PGSs: we speak of horizontal PGSs, vertical PGSs and 
inclined PGSs, respectively, on both $\bbA_2/\bbH_2$.

The PGSs for Class TA1 are all horizontal when $D^2=a^2$ and all vertical when 
$D^2=3a^2$; for Class TA2 they are all inclined. 

As was said before, in Theorems 1, 2 from Section \ref{SubSec3.3} we prove that in the large-fugacity
regime on $\bbA_2$, there are exactly $D^2$ 
EGMs if $D$ is from Class TA1 and exactly $2D^2$ EGMs if $D$ is from Class TA2. 
In both cases, there is a single PGS-equivalence class which is dominant.

\FigureC3

%\indent\hskip .6cm{\bf Figure 3. PGSs on ${\mathbb A}_2$ for $D^2=49$ (Class  TB),  with V-cells }\\
%\indent\hskip .6cm{\bf (gray hexagons) and choices of $D$-rhombuses.}\\ 
%\indent\hskip .6cm{\footnotesize{There are 49 horizontal PGSs (frame (a)), and 98 inclined 
%PGSs (frame (b)).} }\\
%\indent\hskip .6cm{\footnotesize{The horizontal PGSs are the only dominant, so there are only 49 EGMs.}}

For every $D^2$ from Class TB we prove that at least one PGS-equivalence 
class generates $D$-EGMs. 
%that a $D$-EGM is non-unique, and 
%provide a list of possibilities for the set of $D$-EGMs, including its structure. 
See Theorem 3 in Section \ref{SubSec3.4}. As we mentioned before, the structure of EGMs is 
defined by the property of dominance of PGSs. For the specific values 
$D^2 = 49, 147, 169$ we present a new technique that allows us to determine 
which PGS-class is dominant, via a specific count of density of {\it local excitations}. 
In the terminology from \cite{Za}, it is a specific analysis of {\it small contours}. See Theorems 
4, 5, and 6 in Section \ref{SubSec3.4}.

\subsection{PGSs and EGMs on $\bbH_2$}\label{SubSec1.4} On lattice $\bbH_2$ we identify the 
following pair-wise disjoint sets of values of $D$:  HA1, HA2, HA (the union of HA1 and HA2), HB, HC (all infinite), HD, HE, HS (all  finite).

$\bullet$ Class HS: 4 values  with {\it sliding} $D^2= 4, 7, 31, 133$; see Section \ref{Sec8}.

$\bullet$ Class HA1: the values $D$ from the above Class TA1 such that $D^2$ is divisible 
by $3$. That is, $D^2=9b^2$ or $D^2=3b^2$ where $b\in\bbN$. The initial list of 30 such values has 
$D^2 =$ 3, 9, 12, 27, 36, 48, 75, 81, 108, 144, 
192, 225, 243, 300, 324, 363, 432, 576, 675, 729, 
768, 867, 900, 972, 1089, 1137, 1200, 1296, 1389, 1452.

$\bullet$ Class HA2: the values $D$ from the above Class TA2 such that $D^2$ is divisible 
by $3$. The initial list of 30 such values has 
$D^2 =$ 21, 39, 57, 63, 84, 93, 111, 117, 129, 156, 
171, 183, 189, 201, 219, 237, 252, 279, 291, 309, 
327, 333, 336, 351, 372, 381, 387, 417, 444, 453.

$\bullet$ Class HB: the values $D$ from the above Class TB such that $D^2$ is divisible 
by $3$. The initial list of 30 such values has $D^2 =$ 
147, 273, 399, 441, 507, 588, 651, 741, 777, 819, 
903, 1029, 1083, 1092, 1197, 1209, 1281, 1323, 1407, 1443,
1521, 1533, 1596, 1659, 1677, 1764, 1767, 1911, 1953, 2028.

$\bullet$ Class HC: the remaining values of $D$, except for the values from Classes HD
and HE below. Here the initial list of 30 values has $D^2=$ 
19, 25, 37, 38, 43, 52, 61, 73, 76, 79, 
84, 91, 100, 103, 109, 121, 124, 127, 139, 148, 
151, 157, 163, 169, 172, 175, 181, 193, 196, 199. 
% 208, 211, 217, 223, 229, 241, 243.
%289,  400.

$\bullet$ Class HD: 9 values where $D^2=$ 1, 13, 16, 28, 49, 64,  97,
157, 256.

$\bullet$ Class HE: 1 value $D^2=$ 67.

In the analysis of the EGMs, the values $D$ from Class HS are disregarded. As was said 
before, the PS theory does not apply for such $D$. 
\FigureD4

%\indent\hskip .6cm{\bf Figure 4. $\alpha$-PGSs on ${\mathbb H}_2$,  
%for $D^2=48,\;81$ (Class HA1, frames (a,b)),}\\ 
%\indent\hskip .6cm{\bf and  $D^2=39$ (Class HA2, frame (c)), with Voronoi cells (gray }\\
%\indent\hskip .6cm{\bf hexagons) and choices of $D$-rhombuses.}\\
%\indent\hskip .6cm{\footnotesize{The number of PGSs and EGMs for $D^2=48,\;81, 39$ equals 32, 
%54 and 52, respectively.}}\\
%\indent\hskip .6cm{\footnotesize{The PGSs for $D^2=48$ are horizontal, for $D^2=81$ vertical, and
%for $D^2=39$ inclined.}}

%and to those with vertices in $\bbH_2$ as $\bbH_2$-triangles.
%Also, on $\bbH_2$ we will need an additional classification of admissible configurations
%(types $\alpha$, $\beta$ and $\gam$). 

Now, suppose $D$ is from Class HA. Then the model on $\bbH_2$ with a large fugacity 
has $2D^2/3$ EGMs if $D$ falls in Class TA1 and $4D^2/3$ if $D$ falls in Class TA2. See 
Theorems 7, 8 in Section \ref{SubSec3.5}. In Class HA, the PGSs stem from $D$-sub-lattices in $\bbA_2$. 
For HA1, the PGS are all horizontal if $D^2=9b^2$ or all vertical, if $D^2=3b^2$; for HA2
the PGSs are all inclined.

We refer to $D$-admissible configurations on $\bbH_2$ constructed from a $D$-sub-lattice 
as $\alpha$-configurations or 
configurations of type $\alpha$, or -- when the value $D$ should be highlighted -- as 
 $(D,\alpha )$-configurations (in short: $\alpha$-ACs or 
 $(D,\alpha)$-ACs). We will also use the terms an $\alpha$-PGS and a $(D,\alpha )$-PGS.
  
Summarizing, for Class HA we obtain a situation similar to Class TA on $\bbA_2$. Cf. Figure 4.

The picture for Class HB is analogous to that for Class TB. That is, only the dominant 
$\alpha$-PGSs give rise to EGMs, and the issue of dominance is resolved by counting 
local excitations. Cf. Theorem 9 in Section \ref{SubSec3.5}. As an example, we analyze the case 
$D^2=147$ and find out that on $\bbH_2$ there are 98 dominant vertical 
PGSs and 196 non-dominant inclined PGSs. Cf. Figure 5 and Theorem 10.
Consequently, the number of $D$-EGMs for a large fugacity $u$ also equals $98$, and these 
EGMs inherit symmetries between the vertical PGSs. %As above, the decisive part
%is a count of local $u^{-2}$-excitations. 
\FigureE5

A new situation arises for $D$ from Class HC. Here the PGSs stem from $D^*$-sub-lattices
where $D^*>D$ is the nearest L\"oschian number divisible by $3$ (i.e., from Classes HA or HB).
The minimal value for the difference $(D^*)^2-D^2$  equals 
2 and is achieved when $D^2=3b^2+3b+1$ and $(D^*)^2=3b^2+3b+3$, for integer $n\geq 2$. 
(Here, for $n=1$ we obtain $D^2=3b^2+3b+1=7$ which yields a value with sliding.)
If $D^*$ has type HA1, the number of $D$-EGMs in $\bbH_2$ equals 
$2(D^*)^2/3$ while if $D^*$ has type HA2, the number of $D$-EGMs equals $4(D^*)^2/3$.  
Moreover, the PGSs are $\alpha$-configurations and are obtained from each other by $\bbH_2$-shifts 
for $D^*$ from Class HA1 and by $\bbH_2$-shifts or reflections for $D^*$ from Class HA2.  
Cf. Figure 6. 

\FigureF6

For instance, if $D^2=19$ then $(D^*)^2=21$, and the number of the EGMs for 
$D^2=19$ equals $28$. On the other hand, for $D^2=43$ the value $(D^*)^2$ is $48$. Therefore, for  
$D^2=43$ the number of EGMs equals $32$. Cf. Theorem 11 in Section \ref{SubSec3.5}.
 
If $D^*$ is a value of type HB then again the dominance analysis is needed to determine 
which PGSs generate EGMs.

Finally, consider $D^2$ from Classes HD or HE. From now on we refer to the values
of $D^2$ from these classes as {\it exceptional}. (In fact, these values are  
exceptions from Class HC.) It is convenient to divide Class HD into two sub-classes:
HD1: $D^2=$ 1, 13, 28, 49, 64, 97, 157; HD2: $D^2=$ 16, 256.
For $D=D^2=1$ we have a single PGS where all sites in 
$\bbH_2$ are occupied. There is just one EGM for all values of $u$ (not only for 
$u$ large), which is a Bernoulli random field over $\bbH_2$, with probability 
for a site being empty/vacant $1/(1+u)$ and occupied $u/(1+u)$. We will treat the case
$D=1$ as trivial and omit it from the forthcoming discussions. 

\FigureG7
 
%Suppose that $D^2=a^2+b^2+ab$ and $D^2=$ 13, 28, 49, 64, 97, 157 (sub-class HD2). 
Take $D^2=13, 28, 49, 64, 97, 157$ (sub-class HD2) and write $D^2=a^2+b^2+ab$ where $a, b$ 
are non-negative integers. Then we have a 
particular structure of a PGS related to a quadrilateral $OABC$ in $\bbH_2$ where (i) two 
adjacent sides $OA$ and $OC$ have length $D$ and form an angle $2\pi/3$, (ii) two other 
sides $AB$ and $BC$ (also adjacent to each other) have $|AB|^2=D^2+2a+b+1$ and 
$|BC|^2=D^2-a+b+1$, and (iii) the shorter diagonal $OB$ has $|OB|^2=D^2 + a+2b+1$. 
We place particles at the vertices 
of such a quadrilateral and then extend this pattern to the whole of $\bbH_2$, generating
a picture with intermittent stripes parallel to $OB$. Such configurations are called
$\beta$-configurations or configurations of type $\beta$ (in short: $\beta$-ACs);
these configurations are not constructed from sub-lattices.
These configurations yield PGSs in Class HD1; accordingly, we refer to them as $\beta$-PGSs. Cf. 
Figure 7 where $D^2=13$. %Figure 20 for $D^2=49$ . 
There are 
$66$ PGSs for $D^2=13$, $132$ for $D^2=28$, $222$ for $D^2=49$, $288$ for $D^2=64$,
$426$ for $D^2=97$,
and $678$ for $D^2=157$. For a large fugacity $u$, each PGS gives rise to a different $D$-EGM, 
and the number of the $D$-EGMs matches that of the $D$-PGSs. See Theorem 12(i)
from Section \ref{SubSec3.6}. 
\FigureH8

Another particular PGS-structure arises for $D^2=$ 16, 256 (sub-class HD2). For these
integers $D^2=a^2+b^2+ab$ where $a=0$, $b=D$. Then we take 
the value ${\ov D}=2D+1$ which is divisible by $3$.
Hence, an equilateral triangle, $DEO$, can be inscribed in $\bbH_2$, with side-length $\ov D$
and a horizontal base $OE$. Now, consider an inner equilateral triangle $ABC$ with side-length 
$\left(D^2+D+1\right)^{1/2}$ inscribed in $DEO$: the vertices $A,B,C$ lie in the 
sides $OD$, $DE$ and $EO$ and divide them at the ratio $D:(D+1)$. Let us put particles at the vertices 
of $A,B,C$ and $D,E,O$. The PGSs for $D^2=$16, 256 arise from triangles 
congruent to $DEO$, each carrying the above particles, via the extension to the whole of $\bbH_2$.
For this type of configurations we use the terms $\gamma$-configurations and $\gamma$-PGSs. Cf. 
Figure 8. There are $54$ $D$-PGSs for $D^2=16$ and 
726 for $D^2=256$. As above, each PGS gives rise to a different EGM, 
and the number of the EGMs matches that of the PGSs. See Theorem 12(ii) from Section \ref{SubSec3.6}. %in Section 7.
\FigureI9

For $D^2=$ 67 (Class HE) we have a competition between two types of PGSs: (a) 50 PGSs as 
in Class HA ($\alpha$-configurations) with the squared exclusion diameter $75$ and (b) 300 PGSs as in 
class HD1 ($\beta$-configurations). This is formally stated in Theorem 13 in Section \ref{SubSec3.7}.
See Figure 9. We conjecture that the $\alpha$-PGSs of type (a) are dominant.

We see that, for values $D$ from exceptional classes HD and HE on $\bbH_2$,
we have PGSs that are not generated from sub-lattices (apart from $D=1$), yet these cases do not lead to sliding. In
contrast, on $\bbZ^2$, if for a given attainable $D$ there exists a non-lattice PGS then this
value of $D$ exhibits sliding.

\section{Formal preliminaries and basic facts}\label{Sec2}

\subsection{The H-C model on $\bbA_2$}\label{SubSec2.1} We refer to a two-dimensional unit 
triangular lattice $\bbA_2$ as the set of points $\bx =(x_1;x_2)\in\bbR^2$ (sites of the lattice) 
with Euclidean co-ordinates
\be\label{triangularL} x_1=m-\diy\frac{1}{2}n\;\hbox{ and }\;x_2=\diy\frac{\sqrt 3}{2}n,\;\hbox{ 
where }\;m,n\in\bbZ.\ee 
Every site $\bx\in\bbA_2$ has six neighboring sites $\bx'$ such that the distance $\rho(\bx,\bx')$ 
equals $1$. In future we write $\bx\in\bbA_2$ for brevity. Alternatively to $\bx =(x_1;x_2)$, we 
also write $\bx \simeq (m,n)\in\bbA_2$. Geometrically, 
points $(1;0)\simeq (1,0)$ and $(-1/2;{\sqrt 3}/2)\simeq (0,1)$ lead to a natural basis for $\bbA_2$.

Here and below, $\rho (=\rho_2)$ stands for the 2D Euclidean metric: for
$\bx =(x_1;x_2),\by =(y_1;y_2)\in\bbR^2$, the distance $\rho (\bx ,\by)=\left[\rho_1(x_1,y_1)^2
+\rho_1(x_2,y_2)^2\right]^{1/2}$, where $\rho_1(x,y)=|y-x|$, $x,y\in\bbR$.
Alternatively, for $\bx \simeq (m,n),\by\simeq (u,v)\in\bbA_2$, 
$$\rho (\bx ,\by)^2=(m-u)^2+(n-v)^2-(m-u)(n-v).$$ 
Given $\bu\in\bbA_2$, we designate $\tT_\bu:\bbA_2\to\bbA_2$ 
to be an $\bbA_2$-shift by $\bu$, with $\tT_\bu\bx=\bx+\bu$,
$\bx\in\bbA_2$. Similarly, \ $\tR:\bbA_2\to\bbA_2$ stands for the reflection map about the 
horizontal axis.

Given a number $D\geq 1$, consider $D$-{\it admissible} configurations ($D$-ACs, or, in short,
ACs): 
$$\phi_{\bbA_2}:\bx\in\bbA_2\mapsto\phi_{\bbA_2}(\bx)\in\{0,1\}\;\hbox{ (or
shortly $\phi:\bbA_2\to\{0,1\}^{\bbA_2}$)}$$ 
such that for any two {\it occupied} sites $\bx$ and $\by$ with $\phi (\bx)=\phi (\by)=1$ the 
distance $\rho (\bx ,\by)\geq D$. (We can think that $\phi (\bx)=1$ means site $\bx$ is occupied 
in $\phi$ by a particle, and $\phi (\bx)=0$ that $\bx$ is vacant in $\phi$. Particles are 
treated as non-overlapped open disks of diameter $D$ with the centers placed at lattice sites.) 
We write 
$$\hbox{$\bx\in\phi$ if $\phi (\bx)=1$ and identify $\phi$ with the subset in $\bbA_2$ where 
$\phi (\bx)=1$.}$$ 
The value $D$ is called as an H-C exclusion diameter. The set of admissible configurations is denoted 
by $\cA=\cA(D,\bbA_2)$. As was said in Introduction, we can assume that $D^2$ is a L\"oschian
number:
\be\label{Loesch} \hbox{$D^2\in\bbN$ \ and \
$D^2= a^2 + b^2 + ab$ \ where \ $a,b \in {\bbZ}$;} \ee
it means that $D$ is attainable, i.e., there are sites $\bx,\by\in\bbA_2$ with
$\rho (\bx,\by )=D$.
Assumption \eqref{Loesch} does not restrict generality, as any other $D' >1$ can be replaced
by the smallest $D \ge D'$ satisfying \eqref{Loesch} without changing the set
$\cA$.

Set $\cA$ is a closed subset in the Cartesian product $\sX:=\{0,1\}^{\bbA_2}$ 
(the set of all $0,1$-configurations) in the Tykhonov topology. For $D=1$, $\cA =\sX$.

The notion of an AC can be defined for any $\bbV\subset\bbA_2$; accordingly, 
one can use the notation $\cA (\bbV)=\cA (D,\bbV)$. The restriction of configuration
$\phi\in\cA$ to $\bbV$ is denoted by $\phi\upharpoonright_\bbV$.

We are interested in some particular probability measures $\bmu$ on $(\sX,\fB(\sX))$ 
sitting on $\cA$ (i.e., such that $\bmu (\cA )=\bmu (\sX)=1$) where $\fB(\sX)$ is the Borel
$\sigma$-algebra in $\sX$. As we repeatedly stressed, the measures of interest are extreme
Gibbs/DLR, probability measures for high densities/large
fugacities, which are formally defined below. 

Let $\bbV\subset\bbA_2$ be a finite set and $\phi\in{\mathcal A}$.
We say that a finite configuration $\psi^\bbV\in\cA (\bbV)$
is $(\phi ,\bbV)$-compatible if the concatenated configuration
$\psi^\bbV\vee (\phi \upharpoonright_{\bbA_2\setminus\bbV})\in{\mathcal A}$.
The set of $(\phi ,\bbV)$-compatible configurations is denoted by $\cA(\bbV\|\phi )$.

Given $u>0$, consider a probability measure
$\mu_{\bbV}(\;\cdot\; \|\phi )$ on $\{0,1\}^\bbV$ given by
\be\label{GibbsV} \mu_{\bbV}(\psi^\bbV\,\|\phi )=\begin{cases}\diy\frac{u^{\sharp
(\psi^\bbV)-\sharp (\phi^{\,\bbV})}}{\bZ(\bbV\|\phi )},&\hbox{if $\psi^\bbV\in\cA(\bbV\|\phi )$,}\\
0,&\hbox{if $\psi^\bbV\in\{0,1\}^\bbV\setminus\cA(\bbV\|\phi )$.}\end{cases}\ee
Here $\sharp (\psi^\bbV)$ and $\sharp (\phi^\bbV)$ stand for the number of particles in 
$\psi^\bbV$ and $\phi^\bbV$:
$$\sharp (\psi^\bbV):=\#\big\{x\in\bbV:\;\psi (x)=1\big\},\;\;\sharp (\phi^\bbV):=\#\big\{x\in\bbV:\;\phi (x)=1\big\}.$$
Next, $\bZ(\bbV\;\|\phi )$
is the {\it partition function} in $\bbV$ with the boundary condition $\phi$\,:
\be\label{PartFnctnV}\bZ(\bbV\|\phi )=\sum\limits_{\psi^{\,\bbV}\in\cA (\bbV\|\phi )}u^{\sharp
(\psi^{\,\bbV})-\sharp (\phi^{\,\bbV})}.\ee
Measure $\mu_{\bbV}(\;\cdot\; \|\phi )$ sits on $\cA(\bbV\|\phi )$.
Parameter $u>0$ is called {\it fugacity} or {\it activity} (of an occupied site).

A probability measure $\bmu$ on $(\sX,\fB(\sX))$ is called a $D$-H-C
{\it Gibbs/DLR measure} (in short, $D$-H-C GM or GM when the reference 
to $D$ can be omitted)\ if (i) $\bmu ({\mathcal A})=1$, (ii)
$\forall$ \ finite $\bbV\subset\bbA_2$ and a function $f:\phi\in
\sX\mapsto f(\phi )\in{\mathbb C}$ depending only on the
restriction $\phi\upharpoonright_\bbV$, the integral 
$\bmu (f)=\int_\sX f(\phi )\rd\bmu (\phi )$ has the form
\be\label{GibbsInt}\begin{array}{c}\bmu (f)=\diy\int_{\sX}\int_{\{0,1\}^\bbV}
f(\psi^\bbV\vee\phi\upharpoonright_{\bbA_2\setminus\bbV})\rd\mu_\bbV(\psi^\bbV\,
\|\phi )\rd\bmu (\phi ).\end{array}\ee
One can say that under such measure $\bmu$, the probability of a configuration
$\psi^\bbV$ in a finite volume $\bbV\subset\bbA_2$, conditional on
a configuration
$\phi\upharpoonright_{\bbA_2\setminus\bbV}$, coincides with $\mu_{\bbV}(\psi^\bbV\,\|\phi)$,
for $\bmu$-a.a. $\phi\in\{0,1\}^{\bbH_2}$.

In the literature, equality \eqref{GibbsInt} is often referred to as the DLR equation for a 
measure $\bmu$ (in fact, it represents a system of equations labeled by $\bbV$ and $f$).
For the general theory of Gibbs measures, see the monograph \cite{Ge}, Chapters 3, 4, 5--8.

The $D$-H-C GMs form a {\it Choquet simplex} (in the weak-convergence
topology on the set of probability measures on $(\sX,\fB(\sX))$), which we denote by
$\sG=\sG(D,u,\bbA_2)$. An {\it extreme} $D$-H-C GM $\bmu$ is a $D$-H-C GM 
which does not admit a non-trivial decomposition $\bmu =a\bmu^{(1)}
+(1-a )\bmu^{(2)}$ in terms of other $D$-H-C GMs $\bmu^{(i)}$, $i=1,2$,
with $a\in (0,1)$. As was said, the extreme $D$-H-C GMs ($D$-EGMs or briefly EGMs) 
represent pure phases. The collection of $D$-EGMs 
is denoted by $\sE (D)=\sE(D,u,\bbA_2)$. (Argument $u$ will be systematically omitted.)
Any $D$-H-C Gibbs measure $\bmu$ is a barycenter/mixture for some unit mass distribution 
over $\sE (D)$.

\br%{\bf Remark.} 
{\rm The simplest version of the partition function is $\Xi (\bbV\|\varnothing )$,
with an empty boundary condition:
\be\label{(1)}
\bZ (\bbV\|\varnothing )=\sum\limits_{\psi^{\,\bbV}\in\cA (\bbV )}u^{\sharp
(\psi^{\,\bbV})}.\ee
Despite a straightforward (and appealing) form of $\Xi (\bbV\|\varnothing )$,
it is not always convenient (or at least not the most convenient) for the rigorous analysis in the
{\it thermodynamic limit}, for a sequence of volumes $\bbV_k\nearrow\bbA_2$
in the Van Hove sense.  %with $\bbV_{k+1}\supset\bbV_k$ and 
%$\operatornamewithlimits{\cup}\limits_k\bbV_k=\bbA_2$). 
The corresponding limit Gibbs
measure (if it exists) depends on the particular shape of volumes $\bbV_k$ which can be in a
`good'  or `bad' agreement with {\it symmetries} of the hard-core model on $\bbA_2$. In this paper
we concentrate on the partition function $\Xi (\bbV\|\vphi )$ with a PGS boundary condition $\vphi$. 
We also analyze a periodic version of \eqref{(1)}.} %Cf. \cite{CD}. 
\hfill $\blacktriangle$
\er

A {\it ground state} $\vphi\in\cA (D)$ in the H-C model with $u>1$ is defined by the property that 
one cannot remove finitely many particles from $\vphi$ and replace 
them by a larger number of particles without breaking $D$-admissibility. In other words, 
one cannot find a finite subset $\bbV\subset\bbA_2$ and a configuration 
$\psi^\bbV\in\cA (\bbV\|\vphi )$ 
such that $\sharp\psi^\bbV>\sharp\vphi\upharpoonright_\bbV$. 

A crucial role belongs to {\it periodic ground states} (PGSs). A $D$-AC $\phi\in\cA$ is 
said to be periodic if there exist two linearly independent vectors 
$\bfe^{(1)},\bfe^{(2)}\in\bbA_2$ such that $\phi (\bx)=\phi (\bx+\bfe^{(i)})$ $\forall$
$\bx\in\bbA_2$, $i=1,2$. In terms of $\bbA_2$-shifts $\tT_\bu$ it means that 
$\tT_{\bfe^{(i)}}\phi =\phi$ for $i=1,2$. The collection 
of PGSs for a given $D$ is denoted by $\sP(=\sP (D)=\sP (D,\bbA_2))$. 

The PGSs on $\bbA_2$
are relatively straightforward and obtained from $D$-sub-lattices. %The form of PGSs 
%on $\bbH_2$ depends on $D$ in a more subtle way and is discussed below.
%This simplifies %the analysis of he EGMs on $\bbA_2$ greatly.

%Given $\phi\in\cA (D)$ and $\vphi\in\sP (D)$,
%we say that a site $\bx\in\bbA_2$ is $\vphi$-{\it correct} in $\phi$ if 
%$\phi (\by )=\vphi (\by )$ \ $\forall$ $\by\in\bbA_2$ with $\rho (\bx,\by)\leq D$. This is 
%a step towards the specification of the definition in \cite{Za}, P. 561, for our model.

Now we turn to arithmetic properties of a given $D$. Any ordered pair of integers $(a, b)$ 
which is a solution to equation \eqref{Loesch} defines a $D$-sub-lattice of $\bbA_2$ containing 
the origin and the following 6 sites:
\be\label{(10)} (-a,b),\; (b, a+b),\; (a+b, a),\; (a,-b),\; (-b, a+b),\; (a+b, -a),\ee
which all are the solutions to \eqref{Loesch} as ordered pairs of integers. If  $ab=0$ or 
$a = b$  then the pair $(a, b)$ defines a single $D$-sub-lattice of $\bbA_2$ which is mapped 
into itself under the reflection $\tR$ (Class TA1). If $ab\neq 0$ and $a \neq b$
then the pair $(b, a)$ also defines a $D$-sub-lattice of $\bbA_2$ which 
is a reflection by $\tR$ of the sub-lattice defined by $(a, b)$ (Classes TA2 and TB).
For each $D$-sub-lattice of $\bbA_2$ 
generated by a solution to \eqref{Loesch} there are exactly $D^2$ distinct  
$\bbA_2$-shifts $\tT_\bu$ as there are exactly $D^2$ lattice sites inside the 
fundamental parallelogram of the $D$-sub-lattice. All shifted configurations
are PGSs. Moreover, all PGSs corresponding to a given $D$ are obtained as $\mathbb{A}_2$ shifts of $D$-sub-lattices generated by the solutions to \eqref{Loesch}.

\subsection{The H-C model on $\bbH_2$}\label{SubSec2.2}  Formally, $\bbH_2$ can be defined as the 
set-theoretical difference where we remove, from lattice $\bbA_2$, the sub-lattice $\bbA_2({\sqrt 3})$
with a fundamental parallelogram $\{(0,0),(1,2),(2,1),(1,-1)\}$:
\be\label{bbH2a}\beal \bbH_2=\bbA_2\setminus\bbA_2({\sqrt 3}),\\
\quad \hbox{ where }\;\diy\bbA_2({\sqrt 3})=\Big\{m\cdot (1,2)
+n\cdot (2,1):\;m,n\in\bbZ \Big\}.\ena\ee
Equivalently, 
\be\label{bbH2b}\bbH_2=\tT_{(1,0)}\bbA_2({\sqrt 3})\cup\tT_{(0,1)}\bbA_2({\sqrt 3}).\ee
Each site in $\bbH_2$ has three neighboring sites, at the Euclidean distance $1$. 
%Viz., sites $(\pm 1,0)$ have neighbors $(\pm 2;0)$ and $(\pm 1/2;\pm{\sqrt 3})$ while 
%site $(1/2;{\sqrt 3}/2)$ has neighbors  $(-1/2;{\sqrt 3}/2)$, $(1;0)$, $(1;{\sqrt 3})$.  
Lattice $\bbA_2$ is represented as the union of three disjoint congruent subsets:
$\bbA_2=\bbA_2({\sqrt 3})\cup \tT_{(1,0)}\bbA_2({\sqrt 3})\cup\tT_{(0,1)}\bbA_2({\sqrt 3})$.
If $D^2\equiv 0 \mod 3$ then all 3 vertices of an equilateral $\bbA_2$-triangle $\triangle$
with the side-length $D$ lie in the same subset. Otherwise $D^2\equiv 1 \mod 3$, and all  
vertices of $\triangle$ lie in different subset. Hence, for every L\"oschian number $D^2$
there are pairs of vertices $\bx,\bx'\in\bbH_2$ for which $\rho (\bx,\bx')=D$.

Lattice $\bbA_2$ is represented as a non-disjoint union
\be\label{bbA2union}\bbA_2=\bbH_2\cup\left(\tT_{(-1;0)}\bbH_2\right)
=\bbH_2\cup\left(\tT_{(1;0)}\bbH_2\right).\ee 
As above, 
we use the notation $\bx=(x;x')\in\bbH_2$ and $\bx \simeq (m,n)\in\bbH_2$.

%Clearly, the squared Euclidean distance $\rho (\bx ,\by)^2$ between two vertices 
%$\bx,\bx'\in\bbH_2$ is always a L\"oschian number.  Conversely, let $D^2$ be a
%a L\"oschian number and site $\bx=(x_1;x_2)\in\bbA_2$ have $\rho ({\mathbf 0},\bx )^2=D^2$. 
%Then, by virtue of \eqref{bbA2union}, either $\tT_{(1;0)}\bx=(x_1+1;x_2))\in\bbH_2$
%or $\tT_{(-1;0)}\bx=(x_1-1;x_2))\in\bbH_2$. In the first case we get 
%$\rho ((1;0),\tT_{(1,0)}\bx )^2=D^2$ whereas in the second  $\rho ((-1;0),\tT_{(-1,0)}\bx )^2=D^2$.
%Hence, for every L\"oschian number $D^2$
%there are pairs of vertices $\bx,\bx'\in\bbH_2$ for which $\rho (\bx,\bx')=D$.
% \vskip .5cm
  
 We use the term an $\bbH_2$-shift for any $\bbA_2$-shift $\tT_\bu$ where $\bu \simeq (m,n)$
 has both $u_1,u_2$ divisible by 3. Also, $\tR$ stands for the reflection about the horizontal
 axis: $\tR\bx =  (n-m,-m)$ for $\bx =(m,n)\in\bbH_2$. 
 
The definitions of admissible configurations, compatibility, partition functions, Gibbs measures and extreme Gibbs  
measures on $\bbH_2$ are similar to those on $\bbA_2$, and we do not repeat them. We also 
continue using a similar notation $\cA=\cA (D, \bbH_2)$, $\cA(\bbV)=\cA (D, \bbV)$, $\cA(\bbV\|\phi)$. 
The definition of a ground state and a periodic ground state on $\bbH_2$ are direct repetitions of
their counterparts on $\bbA_2$. 

%$\vphi$ is defined as above, by requiring that one 
%cannot find a finite subset $\bbV\subset\bbH_2$ and a configuration $\psi^\bbV\in\cA (\bbV\|\vphi )$ 
%compatible with $\vphi$ such that $\sharp\psi^\bbV>\sharp\vphi\upharpoonright_\bbV$. 
%Next, we say that a configuration $\phi$ on $\bbH_2$ is periodic if there exist two linearly independent
%vectors $\bfe^{(1)},\bfe^{(2)}\in\bbA_2$ with components $e^{(j)}_1$, $e^{(j)}_2$ divisible
%by  3 such that $\phi (\bx)=\phi (\bx +\bfe^{(i)})$ $\forall$ $\bx\in\bbH_2$, $i=1,2$. 
%In terms of $\bbH_2$-shifts $\tT_\bu$ it means that $\tT_{\bfe^{(i)}}\phi =\phi$ for $i=1,2$.

As in the case of $\bbA_2$, the crucial notion is a periodic ground state (PGS). However, 
on $\bbH_2$ a PGS is not necessarily obtained from a sub-lattice (although it is the
case for Classes HA, HB and HC). The set of PGSs for 
a given value $D$ is again denoted by $\sP (D)=\sP (D,\bbH_2)$ and that of EGMs by  
$\sE (D)=\sE (D,u,\bbH_2)$, respectively. (As above, argument 
$u$ will be systematically omitted.)

%Next, for a given $D$, we again can introduce the notion of a $D$-template: pictorially, it is 
%an $\bbH_2$-rhombus $F$ (with vertices in $\bbH_2$) such that $F$ contains the {\it basic cells} 
%of all PGSs; see below. For  $D$ from Classes HA, HB and HC
%the formal definition of a $D$-template given on $\bbA_2$ can be  repeated
%verbatim. For exceptional classes HD, HE, the definition will be modified below. 

\section{Main theorems}\label{Sec3}

\subsection{Templates. Contour definitions}\label{SubSec3.1} First, 
let us consider the case of $\bbA_2$. For a given $D$, 
%the fundamental parallelograms 
%($D$-rhombuses) for different $D$-sub-lattices are different but the 
{\it templates} $F_{k,l}=F^{\bbA_2}_{k,l}$ are defined by
%independent of the choice of a sub-lattice. Here
\be\label{3.1}F_{k,l}
:= \{ (m, n) \in \bbA_2:\;  kD^2 \le m < (k+1)D^2,\; lD^2 \le n < (l+1)D^2\},\; k,l \in {\bbZ}. \ee
Each template contains $D^4$ points.
Note that sites $(kD^2,lD^2)$ form a sub-lattice $\bbA_2(D^2)$,
and all $D$-PGSs are periodic relative to it.

%The concept of a template will allow us to establish 
%the Peierls bound which is a staple of the proof of our results. The meaning of a template 
%is that it is a shift of a common parallelogram for all $D$-PGSs with a given $D$. More precisely,
%each $D$-PGS on $\bbA_2/\bbH_2$ is periodic relative to a sub-lattice. These sub-lattices are 
%different for different equivalence classes. The intersection of these sub-lattices over all equivalence
%classes is also a non-empty sub-lattice; its fundamental parallelogram is the common parallelogram. 
%The lattice $\bbA_2$ can be partitioned into pair-wise disjoint shifted common parallelograms;
%each parallelogram in the partition represents a {\it template}. 

The family $\{F_{k,l}\}$ forms a partition of $\bbA_2$.
The template $F_{0,0}$, treated as a $D^2 \times D^2$-torus, is partitioned into $D^2$ 
rhombuses, one partition for each PGS-equivalence class. We frequently omit the indices $k, l$ in 
the notation for a template when their values are not important or are evident from the context. Figure 10 shows examples of
templates.

\FigureJ10 %Fig10 %\figurea

In what follows, we suppose that volume $\bbV\subset\bbA_2$ is a finite union of 
templates; such a set $\bbV$ is called a {\it basic lattice polygon} (briefly, a basic polygon).

Given a PGS $\vphi\in\sP$ and a basic polygon $\bbV\subset\bbA_2$, the
partition function $\bZ (\bbV\|\vphi)$ in \eqref{PartFnctnV} gives 
rise to a Gibbs probability distribution $\mu_{\bbV}(\;\cdot\; \|\phi )$  
on $\{0,1\}^\bbV$ concentrated on $\cA (\bbV\|\vphi)$. %In a standard way, 
%$\mu^\circ_{\bbV}(\;\cdot\; \|\phi )$ is treated as a probability distribution on 
%$\{0,1\}^{\bbA_2}$, and we can consider measures that are weak limit 
%points for $\mu_{\bbV_k}(\;\cdot\; \|\phi )$ where
%$\bbV_k$ is an increasing sequence of basic polygons. Such measures are examples of Gibbs 
%measures for the H-C model on $\bbA_2$ with diameter $D$; they sit on the set $\cA$. 
We say a PGS $\vphi\in\sP$ {\it generates} a GM $\bmu_\vphi$ if, \ $\forall$ \ sequence of basic polygons 
$\bbV_k\nearrow\bbA_2$ satisfying the Van Hove condition, 
\be\label{gener}\bmu_\vphi =\lim\limits_{k\to\infty}\mu_{\bbV_k}(\;\cdot\; \|\vphi ).\ee
Equivalently, we say that $\bmu_\vphi$ is {\it generated} by $\vphi$.

A specific construction of a GM exploits periodic boundary conditions, in toric 
volumes $\bbV_k=\bbT (k)$, $k = 1, 2,\ldots$ Here $\bbT (k)=\bbT^{\bbA_2} (k)$ is given by
\be\label{(9.2)}\beal\bbT (k) = \{ (m, n) \in \bbA_2:\;  -kD^2 \le m < kD^2,\; -kD^2 \le n < kD^2;\;\hbox{ with}\\
\qquad\quad\hbox{identification }\; (kD^2, n) \equiv (-kD^2 , n)\;\hbox{ and }\;(m,kD^2)\equiv(m,-kD^2)\}.
\ena\ee
To determine the admissible configurations in a torus we use the condition that 
$\rho^{(k)}(\bx,\by)\geq D$. Here the metric $\rho^{(k)}$ is the toric metric on $\bbT (k)$ defined by
$$\rho^{(k)}(\bx,\by)=\left[\rho^{(k)}_1(x_1,y_1)^2 + \rho^{(k)}_1(x_2,y_2)^2\right]^{1/2},$$
where $\bx=(x_1,x_2)$, $\by=(y_1,y_2)$. In turn, $\rho^{(k)}_1$ is a metric on the interval
$[-kD^2,kD^2)$, with $\rho^{(k)}_1(x,y)=\min\,\{y-x,x+2kD^2-y\}$ for $-kD^2 \leq x\le y<kD^2$. 
The set of admissible configurations in $\bbT (k)$
is denoted by $\cA_{\rm{per},k}=\cA_{\rm{per},k}(\bbT (k))$. In the same spirit as \eqref{(1)}, 
the partition function
in $\bbT (k)$ with periodic boundary condition is determined by
\be\label{(101)}\bZ_{\rm{per}}(\bbT (k))=\sum\limits_{\phi_{\bbT (k)}\in\cA_{\rm{per},k}}\;\prod_{\bx\in\bbT_k}
u^{\phi_{\bbT(k)}(\bx )}.\ee
This in turn defines the Gibbs distribution $\mu_{{\rm{per}},k}$. Next, we set 
\be\label{(110)}\bmu_{\rm{per}}=\lim\limits_{k\to\infty}\mu_{{\rm{per}},k}\ee
provided that the limit measure exists. 

\FigureK11 %Fig11% \figureb 

The concept of a template can be extended without changes to the case of $\bbH_2$ when
$D^2$ is from Classes HA, HB or HC (where the PGSs are $\alpha$-configurations). Here 
templates $F_{k,l}=F^{\bbH_2}_{k,l}$ are defined by \eqref{3.1} with the requirement 
$(m,n)\in\bbH_2$; in other words, $F^{\bbH_2}_{k,l}=F^{\bbA_2}_{k,l}\cap\bbH_2$. (Of course, 
in Class HC the number $D$ has to be replaced by $D^*$.) Then the torus $\bbT (k)=\bbT^{\bbH_2} (k)$
is introduced, as $\bbT^{\bbH_2} (k)=\bbT^{\bbA_2} (k)\cap\bbH_2$ where $\bbT^{\bbA_2} (k)$ is
defined in \eqref{(9.2)}. The set $\cA_{\rm{per},k}$ and measures $\mu_{{\rm{per}},k}$ and 
$\bmu_{\rm{per}}$ are defined as above: see \eqref{(101)} and \eqref{(110)}.

For the exceptional values $D^2$, the above construction needs the following modifications. 
For each PGS $\vphi$ we have a 
period parallelogram $\Pi (\vphi )$ with sides $\bfe^{(1)}(\vphi )
=\left(e^{(1)}_1(\vphi ),e^{(1)}_2(\vphi )\right)$ 
and $\bfe^{(2)}(\vphi )=\left(e^{(2)}_1(\vphi ),e^{(2)}_2(\vphi )\right)$ where 
$\tT_{\bfe^{(i)}(\vphi )}\vphi =\vphi$, $i=1,2$:
$$\beal\Pi (\vphi )=\Big\{(m,n)\in\bbH_2:\, (m,n)=\eps_1\bfe^{(1)}(\vphi )+\eps_2\bfe^{(2)}(\vphi )\;\hbox{ for 
$0\leq\eps_i<1$, $i=1,2$}\Big\}.\ena$$
Template $F_{0,0}=F^{\bbH_2}_{0,0}$ for an exceptional $D^2$ can be defined as a parallelogram
\be\label{ExceTempl}\beal F_{0,0}=\Big\{(m,n)\in\bbH_2:\, (m,n)=\eps_1{\ov\bfe}^{(1)}(D)+\eps_2{\ov\bfe}^{(2)}(D)\;\hbox{ for 
$0\leq\eps_i<1$, $i=1,2$}\Big\}.\ena\ee
Here ${\ov\bfe}^{(i)}(D)=\left({\ov e}^{(i)}_1(D),{\ov e}^{(i)}_2(D)\right)$, ${\ov e}^{(i)}_j(D)
={\rm{LCM}}\{e^{(i)}_j(\vphi ):\,\vphi\in\sP (D)\}$, $i,j=1,2$.
Finally, we set $F_{k,l}=\tT_{k{\ov\bfe}^{(1)}(D)+l{\ov\bfe}^{(2)}(D)}F_{0,0}$, $k,l\in\bbZ$, to form a partition of 
$\bbH_2$.

Pictorially, each PGS $\vphi$ is periodic relative to the sub-lattice defined by vectors $\bfe^{(i)}(\vphi )$,
$i=1,2$. Template $F_{0,0}=F^{\bbH_2}_{0,0}$ is the fundamental parallelogram for the
lattice, that is, the intersection of the above lattices for all $\vphi\in\sP (D)$. 

%To state and prove the Peierls bound for an arbitrary $D$, we have to use additional technical 
%constructions. 

Let $\vphi\in\sP(D)$ be a $D$-PGS and $\phi\in\cA(D)$ be an admissible 
configuration, on $\bbA_2$ or $\bbH_2$. Following the definition of correctness on P. 561 in \cite{Za} we say that a template
$F_{k, l}$ is {\it $\vphi$-correct} in $\phi$ if $\phi (\bx)=\vphi (\bx)$ for every site 
$\bx$ lying in 9 templates $F_{k + i, l + j}$, where $i, j = -1, 0, 1$. The 9 templates 
include the initial template $F_{k, l}$ and 8 neighboring templates considered as
connected to $F_{k, l}$. Cf. Figure 11.

A {\it contour support} in a configuration  $\phi\in\cA$ is defined as a connected component 
of the union of templates which are not $\vphi$-correct in $\phi$ for any $\vphi\in\sP$.  
A {\it contour} in $\phi$ is defined as a pair 
\be\label{Gamma}\Gam =\left(\rSp\,(\Gam ), \phi\upharpoonright_{\rSp\,(\Gam )}\right)\ee 
consisting of a contour support $\rSp\,(\Gam )$ and the restriction $\phi\upharpoonright_{\rSp\,(\Gam )}$. 
These definitions are specifications, for the H-C model, of general 
definitions on P. 561 in \cite{Za}. Accordingly, we define sets $\rInt\,(\Gam )$, 
$\qquad\qquad \rInt_\vphi (\Gam )$ and $\rExt\,(\Gam )$ by using Eqn (1.5) from \cite{Za}.
\FigureL12 %\figurec

Inside each of $\rExt\,(\Gam )$, $\rSp\,(\Gam )$ and $\rInt_\vphi (\Gam )$ we can specify a 
{\it boundary layer}: it is a connected set of templates where each parallelogram has a
neighboring template outside the corresponding $\rExt\Gam$, $\rSp\,(\Gam )$, or $\rInt_\vphi\Gam$. 
Each of $\rExt\,(\Gam )$ and $\rInt_\vphi (\Gam )$ has a single corresponding 
boundary layer while $\rSp\,(\Gam )$ has several of them. Every boundary layer in $\rSp\,(\Gam )$ 
has a corresponding (dual) boundary layer inside $\rExt\,(\Gam )$ or $\rInt_\vphi(\Gam )$. 
Moreover, in every boundary layer all occupied sites belong to the same $\vphi$ (which 
justifies the notation $\rInt_\vphi (\Gam )$). Finally, following Eqn (1.5) from \cite{Za}, a contour 
for which a boundary layer of $\rExt\, (\Gam )$ belongs to the
PGS $\vphi$ is called a $\vphi$-{\it contour}. See Figure 12.

Physically speaking, a $\vphi$-contour emerges when we add to $\vphi$ an amount of 
particles at some `inserted' sites and simultaneously remove the particles from $\vphi$
which are `repelled' by the inserted particles. The latter will be referred to as removed 
sites/particles. The whole procedure should of course maintain admissibility. In fact, for
any attainable $D^2$ on $\bbA_2$ and any attainable $D^2\neq 4, 7, 31, 133$ on $\bbH_2$
for $u$ large enough, every EGM $\bmu$ has the property that, with $\bmu$-probability 
$1$ the AC $\phi$ has no infinite contours. Cf. Theorem III(iv) in Section 3.2.

%A key property of the contour definitions is as follows. 
Let $\bbV$ be a basic polygon
and $\vphi\in\sP$ be a PGS. Then the partition function $\BZ (\bbV\|\vphi )$
can be written in the form
\be\label{PFthruCs}
\BZ (\bbV\|\vphi )=\sum_{\{\Gam_i\}\;{\rm{in}}\;\bbV}\prod\limits_iw (\Gam_i)
\ee 
Here and below, $w(\Gam )$ stands for the statistical {\it weight} of contour $\Gam$:
\be\label{SWoC}
w (\Gam )=u^{\sharp (\psi_\Gam )-\sharp (\vphi_\Gam )}.
\ee  
Further, the summation in Eqn \eqref{PFthruCs} is extended to collections of contours 
$\Gam_i$ in $\bbV$ compatible in the sense of the (general) PS theory; cf. \cite{Si}, \cite{Za}.
Here we say that $\Gam$ is a contour in $\bbV$ if the set $\rSp\,(\Gam )\setminus\bbV$ 
is empty or lies in the boundary layer of Ext $(\Gam )$.

\subsection{Pirogov--Sinai theory for the H-C model}\label{SubSec3.2} In Theorems I and II below we 
summarize our results on PGSs and EGMs on both lattices, $\bbA_2$ and $\bbH_2$.
These theorems form a prerequisite for the use of the PS theory. Applying the PS theory, we 
obtain Theorem III which holds true for all Classes of $D$ except for HS. In Theorems 1--13
we identify the PGSs and the EGMs for the respective Classes of values of $D$.

%{\bf Theorem I.} {\sl 
\bthmI\label{ThmI}
\begin{description}
\item[(i)] For any attainable $D>1$, the set $\sP (D,\bbA_2 )$ consists of 
$D$-sub-lattices and their shifts and reflections. In particular, set $\sP (D,\bbA_2 )$ is finite.

\item[(ii)] For any attainable non-sliding $D>1$, the set $\sP (D,\bbH_2 )$ is finite. If $D^2$ is
from Classes {\rm{HA}} or {\rm{HB}} then $\sP (D,\bbH_2 )$ consists of $(\alpha, D)$-configurations. 
If $D^2$ is from Class {\rm{HC}} then $\sP (D,\bbH_2 )$ consists of $(\alpha, D^*)$-configurations.
If $D^2$ is
from Class {\rm{HD1}} then $\sP (D,\bbH_2 )$ consists of $(\beta, D)$-configurations. 
If $D^2$ is
from Class {\rm{HD2}} then $\sP (D,\bbH_2 )$ consists of $(\gamma, D)$-configurations. 
For $D^2=67$ 
(Class {\rm{HE}}), set $\sP (D,\bbH_2 )$ consists of $(\beta, D)$- and $(\alpha, D^*)$-configurations
where $(D^*)^2=75$.
\end{description}
\ethmI

The proof of Theorem I involves the material of Section \ref{Sec4} and completed in Section \ref{SubSec4.6}.

%The Peierls bound for the H-C model on $\bbA_2$ and $\bbH_2$ is given in Theorem II. 
Let us define 
\be\label{SoD}\beal S=S(D)=D^2{\sqrt 3}/2=2\times\Big(\hbox{the area of a $D$-triangle}\Big).\ena\ee
Next, for exceptional non-sliding
values $D^2 = 13, 16, 28, 49, 64, 67, 97,157, 256$ on $\bbH_2$ we set:
\be\beac S^\RD (\sqrt{13})=16.5 {\sqrt 3}/2,\;S^\RD (\sqrt{16})=20.25{\sqrt 3}/2,\;S^\RD
(\sqrt{28})=33{\sqrt 3}/2,\\ 
S^\RD (\sqrt{49}) =55.5{\sqrt 3}/2,\;\; S^\RD (\sqrt{64}) =72{\sqrt 3}/2, \;\;S^\RD (\sqrt{67}) 
= 75{\sqrt 3}/2, \\
S^\RD (\sqrt{97}) = 106.5{\sqrt 3}/2,S^\RD (\sqrt{157}) = 169.5{\sqrt 3}/2,S^\RD (\sqrt{256}) 
= 272.25 {\sqrt 3}/2.\ena\ee
Here the notation $S^\RD(D)$ refers to a minimal re-distributed triangle area for a given 
value $D$ from the above list. The general concept of a re-distributed area will be introduced in 
Section \ref{SubSec4.2}. %Cf. Section 4.5, particularly, Lemma 4.12.

The Peierls bound is established in Theorem II below. It refers to the quantity
\be\label{NoTs}\|\rSp (\Gam) \|:=\hbox{the number of (incorrect) templates in}\;\rSp (\Gam).\ee

%{\bf Theorem II.}
\bthmII \label{ThmII} (The Peierls bound for contours) The weight $w(\Gam )$ of contour \newline
$\Gam =\left(\rSp\,(\Gam ), \phi\upharpoonright_{\rSp\,(\Gam )}\right)$ obeys the bound
\be\label{PeierlsB}
\qquad\qquad\qquad w(\Gam) \le u^{- p\|\rSp\;(\Gam )\|}.\ee
Here $p >0$ (the Peierls constant) satisfies 

$\bullet$ on $\bbA_2$:
\be\label{pon1} p=p(D,\bbA_2)\geq {\sqrt 3}/(288S(D)),\ee

$\bullet$ on $\bbH_2$:
\be\label{pon2} p=p(D,\bbH_2)\geq {\sqrt 3}/(288S(D)),\;\hbox{ if $D^2$ falls in Class 
\rH\rA \ or \rH\rB,}\ee
\be\label{pon3}\beal p=p(D,\bbH_2)\geq {\sqrt 3}/(288S(D^*)),\;\hbox{if $D^2$ falls in Class \rH\rC, 
where $(D^*)^2$}\\
\qquad\hbox{is the closest L\"oschian number with $(D^*)^2>D^2$ such that  $3\big|(D^*)^2$,}\ena
\qquad\ee
\be\label{pon4} \beal p=p(D,\bbH_2)\geq {\sqrt 3}/(288S^\RD(D))\;\hbox{ for $D^2 = 13, 16, 28, 49, 64, 67 ,97,$}\\ 
\qquad \hbox{$157, 256$} \quad (\hbox{Classes \rH\rD, \rH\rE}).\ena
\ee
\ethmII

\bra%{\bf Remark.} 
{\rm The bounds upon $p$ in Eqns \eqref{pon1}--\eqref{pon4} are far from optimal and 
have been selected in a universal form for simplicity. The value $p$ can be improved at an 
expense of additional technicalities.} \hfill $\blacktriangle$
\era

The proof of Theorem II starts in Section \ref{Sec4} and is completed in Section \ref{SubSec5.1}.

\bthmIII \label{ThmIII}
For all $D$ exists a value
$u_0(D,\bbA_2)\in (0,\infty )$, and for all $D\neq 4, 7, 31, 133$ exists a value
$u_0(D,\bbH_2)\in (0,\infty )$ such that, for $u\geq u_0 (D,\,\cdot\,)$, on  
$\bbA_2$ or $\bbH_2$, respectively, the following assertions hold true.
\begin{description}
\item[{\rm{(i)}}] Each {\rm{EGM}} $\bmu\in\sE (D)$ is generated by a {\rm{PGS}}. That is, each {\rm{EGM}} is of the
form $\bmu_\vphi$ for some $\vphi\in\sP (D)$. If {\rm{PGS}}s $\vphi_i$ generate 
{\rm{EGM}}s $\bmu_{\vphi_i}$, $i=1,2$, and $\vphi_1\neq\vphi_2$ then $\bmu_{\vphi_1}\perp
\bmu_{\vphi_2}$. The {\rm{EGM}}s inherit symmetries between the {\rm{PGS}}s: if {\rm{PGS}}s 
$\vphi_i$ generate 
{\rm{EGM}}s $\bmu_{\vphi_i}$, $i=1,2$, and $\vphi_1=\tT_{\bu}\vphi_2$ or $\vphi_1=\tR\vphi_2$
then $\bmu_{\vphi_1}=\tT_{\bu}\bmu_{\vphi_2}$ or $\bmu_{\vphi_1}=\tR\bmu_{\vphi_2}$, 
respectively.

\item[{\rm{(ii)}}] {\rm{EGM}}-generation is a class property: if a $\vphi\in\sP (D)$ 
generates an EGM $\bmu_\vphi$ then every {\rm{PGS}} $\wt\vphi$ from the same 
 {\rm{PGS}}-equivalence 
class generates an EGM $\bmu_{\wt\vphi}$. Such a class is referred to as dominant. If an 
equivalence class is unique, it is dominant.

\item[{\rm{(iii)}}] Measure $\bmu_{\rm{per}}$ exists and is a uniform mixture of the measures $\bmu_\vphi$
where $\vphi$ runs through all dominant {\rm{PGS}}-equivalence classes. If there is a single 
equivalence class then $\bmu_{\rm{per}}$ is a uniform mixture of all measures $\bmu_\vphi$ where
$\vphi\in\sP$.

\item[{\rm{(iv)}}] Each {\rm{EGM}} $\bmu_\vphi$ exhibit the following properties. For $\bmu_\vphi$-almost all $\phi\in\cA$: 
\begin{description}
\item[{\rm{(A)}}] All contours $\Gam$ in $\phi$ are finite. 
\item[{\rm{(B)}}] For 
any site $\bx$ there exist only finitely many contours $\Gam$ (possibly none) such that 
$\bx\in\rInt (\Gam )\cup\rSp\,(\Gam )$. 
\item[{\rm{(C)}}] There are countably many disjoint connected sets 
of $\vphi$-correct templates one of which is infinite and all remaining ones are finite. 
\item[{\rm{(D)}}] For any $\vphi'\in\sP\setminus\{\vphi\}$, there are countably many disjoint connected sets 
of $\vphi'$-correct templates, and they all are finite.
\end{description}

\item[{\rm{(v)}}]  Measure $\bmu_\vphi$ admits a polymer expansion and consequently has an exponential 
decay of correlations.

\item[{\rm{(vi)}}] As $u\to\infty$, measure $\bmu_\vphi$ converges weakly to a measure sitting
on a single {\rm{AC}} $\vphi$.
\end{description}
\ethmIII

%{\bf Proof of Theorem III.} \; 
\bprfIII(i, ii) According to Corollary on P. 565 in \cite{Za}, there exists at least one
PGS $\vphi\in\sP (D)$ which generates an EGM $\bmu_\vphi\in\sE (D)$. It is obvious that every PGS 
$\wt\vphi$ from the same equivalence class generates an EGM $\bmu_{\wt\vphi}$ and that  EGMs
$\bmu_\vphi$, $\bmu_{\wt\vphi}$ are related with the same symmetry as PGSs $\vphi$,
$\wt\vphi$. The fact that each EGM $\bmu$ is generated by a PGS follows from Corollary on P. 578 in \cite{Za}, 
completed with Theorem 1  from \cite{DoS}. The mutual singularity of measures $\bmu_{\vphi_1}$
and $\bmu_{\vphi_2}$ for $\vphi_1\neq\vphi_2$ can be deduced from assertion (iv) which is deduced
from \cite{Za} below.

Passing to assertion (iii), the main concern is the existence of contours {\it winding} around the torus
$\bbT_k$. However, the $\mu_{{\rm{per}},k}$-probability of the event
$\cW_k\subset\cA_{{\rm{per}},k}$ that such contour is present in an admissible configuration
$\phi_{\bbT^{(1)}_k}$ becomes negligible as $k\to\infty$, as winding contours
are too large. On the remaining event, ${\ov\cW}_k =\cA_{\rm{per},k}\setminus\cW_k$,
the statistics of the random configuration is described in terms of the ensemble of {\it external contours}.
Furthermore, event ${\ov\cW}_k$ can be partitioned into $\sharp \sE$ parts,
${\ov\cW}_{k,\bmu_\vphi}$, $\bmu_\vphi\in\sE$, so that on ${\ov\cW}_{k,\bmu_\vphi}$ all external 
contours are $\vphi$-contours, and each ${\ov\cW}_{k,\bmu_\vphi}$ will have the same limit 
probability: $\lim\limits_{k\to\infty}
\mu_{{\rm{per}},k}\left({\ov\cW}_{k,\bmu_\vphi}\right)=
\diy\frac{1}{\sharp \sE}$.  $\bigg($Here we use the property that the ratio
$\diy\frac{\mu_{{\rm{per}},k}\left({\ov\cW}_{k,\bmu_\vphi}\right)}{
\mu_{{\rm{per}},k}\left({\ov\cW}_{k,\bmu_{\vphi'}}\right)}$ tends to $1$ as $k\to\infty$ for any
choice of EGMs $\bmu_\vphi,\bmu_{\vphi'}\in\sE$. This follows from the contour representation
for the sum $\sum\limits_{\phi_{\bbT (k)}\in{\ov\cW}_{k,\mu_\vphi}}\;\prod\limits_{\bx\in\bbT_k}
u^{\phi_{\bbT(k)}(\bx )}=\bZ_{\rm{per}}(\bbT (k))\mu_{{\rm{per}},k}\left({\ov\cW}_{k,\bmu_\vphi}\right)$.$\bigg)$
This argument leads to the formula
$\lim\limits_{k\to\infty}\mu_{{\rm{per}},k}={\diy\frac{1}{\sharp \sE}}\sum\limits_{\bmu_\vphi\in\sE}\bmu_\vphi$. 

(iv) Statements (A, B) follow from the fact that in an EGM $\bmu_\vphi$, the probability of a contour 
$\Gam$ is $\leq u^{\|\rSp (\Gam )\| p(D)/3}$. Cf. \cite{Za}, Theorem on P. 564. In turn, (C, D) follow from
(A, B). 

Statements (v, vi) follow from \cite{Za}, Theorem on P. 564. %$\qquad\blacksquare$
\eprfIII

In  assertions Theorems 1--13 below we work under the condition $u >u_0(D,\,\cdot\,)$ assumed in Theorem III.

\subsection{PGSs and EGMs for Class TA}\label{SubSec3.3} Here and below, we say that a sub-lattice 
in $\bbA_2$ has type $(a,b)$, or is  
an $(a,b)$-sub-lattice when it is a $D$-sub-lattice generated by the sites $(0,0)$ 
and $(a+b,-b)$, where $a,b\in\bbZ$, $D^2=a^2+b^2+ab$. If $ab=0$, the corresponding 
$D$-sub-lattice is horizontal; if $a=b$, it is vertical. In cases where $ab\neq 0$
and $a\neq b$, we have inclined $D$-sub-lattices. Correspondingly, a PGS-equivalence 
class is called an $(a,b)$-class or class $(a,b)$ if it contains a sub-lattice of type
$(a,b)$. We denote this class by $\sP(a,b)$, and its PGSs are called $(a,b)$-PGSs. 
We also use the term an $\alpha$-configuration and its specifications: 
a $(D,\alpha )$-configuration or $((a,b),\alpha )$-configuration, or $\alpha$-configurations of 
type $(a,b)$, intermittently. %Our results in this section are the following Theorems 1 and 2. 

\bthma%{\bf Theorem 1.} 
(Class {\rm{TA1}}) \begin{description}
\item[{\rm{(i)}}] 
Let $D$ be an integer not divisible by primes of the form $3v+1$. Then in $\bbA_2$ there 
is a unique $D$-sub-lattice which is horizontal and has type $(D,0)$. Thus, on $\bbA_2$ 
there is a single {\rm{PGS}}-equivalence class, which contains $D^2$ different {\rm{PGS}}s. 
The {\rm{PGS}}s are horizontal $\big((D,0),\alpha\big)$-configurations, hence
reflection-invariant. Different {\rm{PGS}}s are obtained from each other by $\bbA_2$-shifts. 
Consequently, for such $D$ the number of {\rm{EGM}}s on $\bbA_2$ equals $D^2$.

\item[{\rm{(ii)}}] Let $D /\sqrt{3}$ be an integer not divisible by primes of the form $3v+1$. Then 
in $\bbA_2$ there is a unique $D$-sub-lattice which is vertical and has type 
$(\frac{D}{\sqrt 3},\frac{D}{\sqrt 3})$. Thus, on $\bbA_2$  there is a single 
{\rm{PGS}}-equivalence class, which contains $D^2$ different {\rm{PGS}}s.
The {\rm{PGS}}s are vertical $\Big((\frac{D}{\sqrt 3},\frac{D}{\sqrt 3}),\alpha\Big)$-configurations, 
hence reflection-invariant. Different {\rm{PGS}}s are obtained from each other by $\bbA_2$-shifts. 
Consequently, for such $D$ the number of {\rm{EGM}}s on $\bbA_2$ equals $D^2$.
\end{description}
\ethma

\bthmb%{\bf Theorem 2.} 
(Class {\rm{TA2}}) Let $D^2$ be 
an integer whose prime decomposition contains {\rm{(i)}} a factor 
$3$ in any power, {\rm{(ii)}} primes of the form $3v+2$, in even powers, possibly zero, and 
{\rm{(iii)}} a single prime of the form $3v+1$. Then in $\bbA_2$ there are exactly two $D$-sub-lattices, 
which are inclined and taken to each other by reflections. Hence, on $\bbA_2$ there is a single  
{\rm{PGS}}-equivalence class, which contains $2D^2$ {\rm{PGS}}s. The {\rm{PGS}}s are 
inclined $(D,\alpha)$-configurations, hence not reflection-invariant. Different {\rm{PGS}}s are 
obtained from each other by $\bbA_2$-shifts and reflections. Consequently, for such $D$ the number of 
{\rm{EGM}}s on $\bbA_2$ equals $2D^2$.
\ethmb

\subsection{PGSs and EGMs for Class TB}\label{SubSec3.4} For a generic $D$ from Class TB, there are at least 
three $D$-sub-lattices among PGSs on $\bbA_2$. 

\bthmc%{\bf Theorem 3.} 
(Class {\rm{TB}}) Suppose that the prime decomposition of $D^2$ contains 
{\rm{(i)}} a factor $3$ in any power, {\rm{(ii)}} primes of the form $3v+2$, in even powers, possibly 
zero, and {\rm{(iii)}}
$M\geq 2$ primes of the form $3v+1$, some of which may coincide. Then the 
following assertions hold true. 
\begin{description}
\item[{\rm{(i)}}] The number of {\rm{PGS}}-equivalence classes 
on $\bbA_2$ is $\geq 2$, and it increases when $\lceil M/2\rceil$ increases.    
\item[{\rm{(ii)}}] At most one class contains $D^2$ {\rm{PGS}}s. It consists of horizontal $(D,0)$-{\rm{PGS}}s 
if $D$ is integer, or of vertical $(\frac{D}{\sqrt 3},\frac{D}{\sqrt 3})$-{\rm{PGS}}s if $D/{\sqrt 3}$ is integer.
All other equivalence classes contain two inclined $D$-sub-lattices and $2D^2$ {\rm{PGS}}s each; these 
sub-lattices are taken to each other by reflections. 
\item[{\rm{(iii)}}] Furthermore, a measure 
$\bmu_\vphi$ is reflection-invariant iff the {\rm{PGS}} $\vphi$ comes from a dominant equivalence 
class of cardinality $D^2$. 
\item[{\rm{(iv)}}] Let $J=J(D,\bbA_2)$ denote the number of 
dominant equivalence classes labeled by $1,\ldots ,J$ in an arbitrary order. Let $m_jD^2$ 
stand for the number of {\rm{PGS}}s in the dominant class $j$, where $m_j=1, 2$ and $1\leq j\leq J$. 
Then the total number of {\rm{EGM}}s equals $D^2\sum\limits_{j=1}^Jm_j$.
\end{description}
\ethmc

As we said earlier, we conjecture that the number of dominant classes $J(D,\bbA_2)=1$.
This is confirmed in several examples considered in Theorems 4-6 on $\bbA_2$
and Theorem 10 on $\bbH_2$.

The analysis of dominance for selected examples of $D^2$ from Class TB is given in Theorems 4--6.
For these examples, we reach the level of Theorems 1, 2 in the description of the structure of
EGMs. The examples have been selected to demonstrate different outcomes of the 
competition between inclined and horizontal or vertical equivalence classes.

\bthmd%{\bf Theorem 4.} 
For $D^2=49$ on $\bbA_2$, there are $147$ {\rm{PGS}}s divided in two 
equivalence classes: horizontal $(7,0)$ and inclined $(5,3)$. The $(7,0)$-class $\sP (7,0)$
consists of $49$ {\rm{PGS}}s and is the only one dominant. 
The $(7,0)$-{\rm{PGS}}s are reflection-invariant and obtained from each other by $\bbA_2$-shifts.
Consequently, we have in total $49$ {\rm{EGM}}s, and they all are of the form $\bmu_\vphi$ 
where $\vphi\in\sP (7,0)$.
\ethmd

In the next result the choice of the dominant PGS class between the inclined and horizontal 
ones is inverted.

\bthme%{\bf Theorem 5.} 
For $D^2=169$ on $\bbA_2$, there are $507$ {\rm{PGS}}s divided in two 
equivalence classes: inclined $(8,7)$ and horizontal $(13,0)$. The $(8,7)$-class $\sP (8,7)$ 
consists of $338$ {\rm{PGS}}s and is the only one dominant. The $(8,7)$-{\rm{PGS}}s 
are not reflection-invariant; they are obtained from each other by $\bbA_2$-shifts and 
reflections. Consequently, we have in total $338$ {\rm{EGM}}s $\bmu_\vphi$, and they all 
are of the form $\bmu_\vphi$ where $\vphi\in\sP (8,7)$.
\ethme

Finally, we discuss a case where we have one vertical class, $(7,7)$, and one inclined, $(11,2)$.

\bthmf%{\bf Theorem 6.} 
For $D^2=147$ on $\bbA_2$, there are $441$ {\rm{PGS}}s divided in two equivalence
classes, the vertical $(7,7)$ and the inclined $(11,2)$. The $(7,7)$-class $\sP (7,7)$ consists 
of $147$ {\rm{PGS}}s and is the only one that is dominant. The $(7,7)$-PGSs are reflection-invariant 
and obtained from each other by $\bbA_2$-shifts. Consequently, we have in total $147$ {\rm{EGM}}s,
and they all are of the form $\bmu_\vphi$ where $\vphi\in\sP (7,7)$.
%are -periodic but not reflection-invariant.}
\ethmf

\subsection{PGSs and EGMs for Classes HA, HB and HC}\label{SubSec3.5}
The results for Classes HA and HB go in parallel to those for Classes TA and TB. Recall, for 
a value $D^2$ from Classes HA and HB, the PGSs on $\bbH_2$ are $(D,\alpha )$-PGSs
restricted to $\bbH_2$; see Theorem I(ii). We will use on $\bbH_2$ the same 
terminology as on $\bbA_2$.

\bthmg%{\bf Theorem 7.} 
(Class {\rm{HA1}}) 
\begin{description}
\item[{\rm{(i)}}] Let $D/3$ be an integer not divisible by primes of the 
form $3v+1$. Then on $\bbH_2$ there is a single {\rm{PGS}}-equivalence class, which contains
$2D^2/3$ {\rm{PGS}}s. The {\rm{PGS}}s are horizontal $((D,0),\alpha)$-configurations, hence
reflection-invariant. Different {\rm{PGS}}s are obtained from each other by $\bbH_2$-shifts.
Consequently, the number of {\rm{EGM}}s on $\bbH_2$ equals $2D^2/3$.

\item[{\rm{(ii)}}] Let $D/\sqrt{3}$ be an integer not divisible by primes of the form $3v+1$. Then 
on $\bbH_2$  there is a single {\rm{PGS}}-equivalence class, which contains $2D^2/3$ {\rm{PGS}}s.
The {\rm{PGS}}s are vertical $\left((\frac{D}{\sqrt 3},\frac{D}{\sqrt 3}),\alpha\right)$-configurations, 
hence reflection-invariant. Different {\rm{PGS}}s are obtained from each other by $\bbH_2$-shifts. 
Consequently, the number of {\rm{EGM}}s on $\bbH_2$ equals $2D^2/3$.
\end{description}
\ethmg

\bthmh%{\bf Theorem 8.} 
(Class {\rm{HA2}}) Let $D^2$ be an integer whose prime decomposition contains 
{\rm{(i)}} at least one factor 3, {\rm{(ii)}} primes of the form $3v+2$, in even powers, possibly zero, and 
{\rm{(iii)}} a single prime of the form $3v+1$. Then on $\bbH_2$  there is a 
single {\rm{PGS}}-equivalence class, which contains $4D^2/3$ {\rm{PGS}}s. The {\rm{PGS}}s are
inclined $(D,\alpha)$-configurations, hence not 
reflection-invariant. Different {\rm{PGS}}s are obtained from each other by shifts 
and reflections. Consequently,  the number of {\rm{EGM}}s on $\bbH_2$ equals $4D^2/3$.
\ethmh

In Theorem 9 we use the same terminology of dominant classes as in Theorem 3

\bthmi%{\bf Theorem 9.} 
(Class {\rm{HB}}) Suppose that the prime decomposition of $D^2$ contains 
{\rm{(i)}} at least one factor 3, {\rm{(ii)}} primes of the form $3v+2$, in even powers, possibly zero, and 
{\rm{(iii)}} at least two prime factors of the form $3v+1$, some of which may coincide. 
Then all {\rm{PGS}}s are 
$(\bbH_2,D,\alpha )$-configurations obtained as the restrictions to
$\bbH_2$ of their $(\bbA_2,D,\alpha )$-counterparts. Thus, the number of
 {\rm{PGS}}-equivalence classes on $\bbH_2$ is the same as on $\bbA_2$.
Furthermore, the assertions {\rm{(ii)}}--{\rm{(iv)}}
of Theorem $3$ are transferred from $\bbA_2$  to $\bbH_2$, with the proviso that the
number of {\rm{PGS}}s in the equivalence classes is $2D^2/3$ in place of $D^2$ and
$4D^2/3$ in place of $2D^2$. Hence, in assertion {\rm{(iv)}}, the total number 
of $D$-{\rm{PGS}}s on $\bbH_2$ should be equal to $D^2\sum\limits_{j=1}^Jm_j$ 
where $m_j=2/3,4/3$, $1\leq j\leq J$, and $J=J(D,\bbH_2)$.
\ethmi

We conjecture that, for any $D^2$ from Class HB, $J(D,\bbH_2)=1$. 
%As before, a corollary of Theorem 9 is that, in the above notation, the cardinality 
%$\sharp\sE (D)$ equals $\big(2D^2/3\big)\sum\limits_{j=1}^Jm_j$.

\smallskip
An analog of Theorem 6 is

\bthmj%{\bf Theorem 10.}  
For $D^2=147$ (Class {\rm HB}), there are $294$ {\rm{PGS}}s 
divided in two equivalence classes, the vertical $(7,7)$ and the inclined 
$(11,2)$. The $(7,7)$-class $\sP(7,7)$ consists of $98$ {\rm{PGS}}s
and is the only one dominant. The {\rm{PGS}}s $\vphi\in\sP (7,7)$ are reflection-invariant and
obtained from each other by $\bbH_2$-shifts. Consequently, we have in total $98$ {\rm{EGM}}s 
$\bmu_\vphi$, where $\vphi\in\sP (7,7)$.
\ethmj

Results for Class HC are given in Theorem 11 below.

\bthmk%{\bf Theorem 11.} 
Assume $D^2$ is not divisible by 3 and not from Classes {\rm HD}, {\rm HE} or {\rm HS}. 
Consider the number $D^*$ such that $D^*>D$ and $(D^*)^2$ is the 
closest L\"oschian number to $D^2$ divisible by $3$. Then the {\rm PGS}s on $\bbH_2$ 
are the $(D^*,\alpha )$-configurations.
\begin{description}
\item[{\rm (A1)}] Suppose the value $D^*$ belongs to Class {\rm HA1}. Then the assertions of 
Theorem $7$ can be repeated with $D$ replaced by $D^*$. In particular, the number of 
{\rm EGM}s on $\bbH_2$ equals  $2(D^*)^2/3$.

\item[{\rm (A2)}] Suppose the value $D^*$ belongs to Class {\rm HA2}. Then 
the assertions of Theorem $8$  can be repeated 
with $D$ replaced by $D^*$. In particular, the number of 
{\rm EGM}s on $\bbH_2$ equals  $4(D^*)^2/3$.

\item[{\rm (B)}] Suppose that the above value $(D^*)^2>D^2$ belongs to Class {\rm HB}.  Then 
the assertions of Theorem $9$ can be repeated with $D$ replaced by $D^*$.
In particular, the total number of {\rm EGM}s on $\bbH_2$ equals $(D^*)^2\sum\limits_{j=1}^Jm_j$ 
where $m_j=2/3,4/3$, $1\leq j\leq J$, and $J=J(D^*,\bbH_2)$.
\end{description}
\ethmk

%An example of dominance analysis in Class HC is presented in Theorem 12. It is
%related to Theorem 10 but the outcomes of the two theorems are opposite.
%\vskip .5cm

%{\bf Theorem 12.} {\sl For $D^2=139$ on $\bbH_2$ (Class {\rm HC}), the value $(D^*)^2=147$.
%There are $294$ $D$-{\rm{PGS}}s on  $\bbH_2$, and they are $(D^*,\alpha)$-configurations 
%referred to in Theorem $10$. The inclined $(11,2)$-class $\sP (11,2)$ is dominant, with 
%$196$ {\rm{PGS}}s. The $(11,2)$-PGSs are not reflection-invariant; they are obtained from 
%each other by $\bbH_2$-shifts and reflections. Consequently, there are $196$ $D$-{\rm{EGM}}s 
%$\bmu_\vphi$ where $\vphi\in\sP (11,2)$.}
%\vskip .5cm

\subsection{PGSs and EGMs for Class HD}\label{SubSec3.6} To conclude our results, it remains to consider 
exceptional non-sliding values $D^2$. For Class HD we have the following

\bthml%{\bf Theorem 12.} 
Assume $D^2>1$ is from Class  {\rm HD}. 
\begin{description}
\item[{\rm{(i)}}] For  $D^2=13,\, 28,\, 49,\, 64,\, 97,\, 157$ (sub-class {\rm{HD1}}): the number of 
{\rm{PGS}}s on $\bbH_2$ equals $66$, $132$, $222$, $288$, $426$ and $678$, 
respectively, and they are $(D,\beta )$-configurations. The {\rm{PGS}}s are not reflection-invariant
and are obtained from each other by shifts and reflections. The number of the {\rm{EGM}}s 
matches that of the {\rm{PGS}}s. 

\item[{\rm{(ii)}}] For $D^2=16, 256$  (sub-class {\rm{HD2}}): the number of {\rm{PGS}}s on $\bbH_2$ equals  
$54$ and $726$, respectively, and they are $(D,\gamma )$-configurations. The {\rm{PGS}}s 
are not reflection-invariant and are obtained from each other by shifts and reflections.  
The number of the $D$-{\rm{EGM}}s matches that of the {\rm{PGS}}s.
\end{description}
\ethml

\subsection{PGSs and EGMs for Class HE}\label{SubSec3.7} Class HE ($D^2=67$) is the one where the description of 
EGMs requires new techniques and is not given in this paper. However, the analysis of the 
PGSs can be done.

\bthmm%{\bf Theorem 13.} {\sl 
For the value $D^2=67$ (Class  {\rm HE}): there are $300$ {\rm{PGS}}s of 
type $(D,\beta )$ and $50$ {\rm{PGS}}s of type $(D^*,\alpha )$, with  $(D^*)^2=75$. 
\ethmm

As was said earlier, we conjecture that $(D^*,\alpha )$-PGSs form the only dominant equivalence 
class, and so the number of $D$-{\rm{EGM}}s equals $50$. %However, it remains an open question.

\section{The PGSs on $\bbA_2$  and $\bbH_2$ via MRA-triangles}\label{Sec4}

%The aim of this section is to prove Theorem I.
%the following Lemma IV:
%\vskip .5cm

%{\bf Lemma IV.} (i) {\sl For any attainable $D$ on $\bbA_2$ and any attainable  
%$D$ from Classes \rH\rA or \rH\rB  on $\bbH_2$, every {\rm{PGS}} is a 
%a $D$-sub-lattice or its $\bbA_2/\bbH_2$-.  For any attainable $D$ 
%from Class \rH\rC \ on  $\bbH_2$,  every {\rm{PGS}} is a a $D$-sub-lattice or its 
%$\bbH_2$-. Here and below, $D^*$ stands for
%a nearest attainable value with $(D^*)^2>D^2$ and $3|D^2$.}

%(ii) {\sl For values $D^2$ from Classes \rH\rD \ and
%\rH\rE on $\bbH_2$, the PGSs are as listed in Theorems $13$, $14$.}

To prove Theorem I, we develop a united approach to the analysis of PGSs
covering the whole variety of cases in Theorems 1--13. It is based on the notion of a re-distributed 
area of a triangle in the Delaunay triangulation (DT) for a $D$-AC $\phi\in\cA$ and the concept
of a MRA-triangle minimizing a re-distributed area. In our approach, we have been inspired 
by (i) an idea of a local energy minimizer serving as an indicator of a PGS (see \cite{HS}) and 
(ii) a specific choice of a minimizer as a triangle area in a DT and a related notion of a saturated 
configuration (see \cite{ChW}). Elements of such an approach have been used for lattice $\bbZ^2$ 
in \cite{MSS1}, Sect. \ref{Sec3}.

\subsection{V-cells, C-triangles and saturated configurations}\label{SubSec4.1} 
The key point of our construction is that maximizing the number of particles in an AC $\phi$ 
can be done through minimizing triangle areas in the Delaunay triangulation of $\phi$; see 
below. One caveat here is that minimization should exclude `sliver' obtuse triangles (as their 
area can be arbitrarily small). The other caveat is that minimization is applied not to the
`standard' triangle area but to its modification which we call a re-distributed (RD) area 
$s^\RD (\triangle )$. And finally, one has to verify that the triangles minimizing the RD-area 
(MRA-triangles) form a tessellation of the whole $\bbA_2/\bbH_2$.
 
Let us pass to a formal argument.  
Consider an arbitrary set ${\bbE} \subset {\bbR}^2$, with at least two points,
such that $\rho(\bx,\by) \ge D$ for any two distinct $\bx, \by \in {\bbE}$. For each $\bx \in {\bbE}$
define the {\it Voronoi cell} $\cV (\bx, \bbE)$ as the set of points $\bz \in {\bbR}^2$ satisfying $\rho(\bx,\bz)
\le \rho(\by,\bz)$, $\forall$ $\by \in {\bbE} \setminus \{\bx\}$. The Voronoi cells (V-cells, for short)
are always convex polygons.

\FigureM13

We apply the above definition to a given $D$-AC  $\phi\in\cA (D)$ with at least two particles; this 
yields a collection of {\it Voronoi cells} $\cV (\bx,\phi )$ constructed for the occupied sites 
$\bx\in\phi$. Here $\cA (D)$ may stand for $\cA (D,\bbA_2/\bbH_2)$ or
$\cA (D,\bbR^2)$.  If $\phi$ has no unbounded
V-cells then to each cell $\cV(\bx,\phi )$ there is assigned a
finite set of circles centered at the vertices of $\cV(\bx,\phi )$
and passing through $\bx$. We call them {\it V-circles} in $\phi$. Each
$\bx\in\phi$ lies in at least one of V-circles but no
$\bx\in\phi$ falls inside a circle. The sites $\by\in\phi$ lying in a given
V-circle form the vertices of a {\it constituting polygon}. These polygons
form a tessellation of $\bbR^2$: they have disjoint interiors, and the union of
their closures gives the entire plane. If a constituting polygon has
$\geq 3$ vertices, it can be divided (non-uniquely) into {\it constituting triangles}
(in short: \rC-triangles); this produces the {\it Delaunay triangulation} (DT) of
$\phi$ (and of $\bbR^2$). See Figure 13 (a).

\bl \label{Lem4.1}%{\bf Lemma 4.1.} 
Let $\triangle$ be a \rC-triangle in a $D$-\rA\rC $\phi$ and consider $3$ pair-wise 
disjoint disks of diameter $D$ centered at the vertices of $\triangle$. Consider $3$ sectors in these 
disks which are intersections of the circles with the angles of $\triangle$ and let $\bbS (\triangle )$
denote the union of these sectors. Then the area of $\bbS (\triangle)$, i.e., the sum of the areas of 
these $3$ sectors, equals $\pi D^2/8$.
\el

\bp Let us stress that $\bbS (\triangle )$ not necessarily lies  completely inside triangle 
$\triangle$. Nevertheless, the sets  $\bbS (\triangle )$ where $\triangle$ runs over \rC-triangles 
of $\phi$ form a partition of the union of the disks
$\operatornamewithlimits{\cup}\limits_{\bx\in\phi}\bbD(\bx,D/2)$ (modulo a set of measure $0$).
Here $\bbD(u,r)$ stands for the disk of radius $r>0$ centered at $\bu\in\bbR^2$: $\bbD(\bu, r)
=\{\by\in\bbR^2:\rho (\bu,\by)\leq r\}$.

For each angle of size $\alpha$ in $\triangle$ the intersection with the corresponding disk is 
a full sector with the angular measure $\alpha$ and area $\alpha D^2/8$. The sum of the 
triangle angles equals $\pi$. Cf. Figure 13 (b). 
\ep

Lemma \ref{Lem4.1} establishes a principal fact that the number of particles in an AC $\phi$ equals 
the doubled number of C-triangles in the DT. Hence, to maximize the number of particles
one would like to minimize the area of C-triangles. However, since the triangular areas in a DT
can be arbitrarily small (when a C-triangle is obtuse and has a large circumradius, i.e, the radius 
of the corresponding V-circle), we use the idea of saturation allowing us to discard C-triangles 
that have the area close to 0; see Lemma \ref{Lem4.2}. 
 
%\vskip .5cm

%The following text establishes some related definitions and known facts. In our own 
%opinion, it is quite elementary although seems rather tedious. It should be an easier
%reading when reader is familiar with the concept of the Delaunay triangulation.  The 
%problems are two-fold: (a) we have to separate obtuse
%triangles and (b) M-triangles may not generate a PGS on the whole $\bbA_2$ or $\bbH_2$.

%\vskip .5cm
\def\a{{\alpha}}

A $D$-AC $\phi$ is called {\it saturated} if no occupied site can be added to it without breaking
admissibility. A {\it saturation} of a given $D$-AC $\phi$ is a completion of $\phi$ (in some 
uniquely defined way) with the maximal possible amount of added occupied sites.

Clearly, every {\rm{PGS}} configuration is saturated (this is also true for non-periodic GSs). 
Saturated configurations are convenient as a natural initial step in a procedure of identifying 
PGSs within the set ${\mathcal A}(D)$ of admissible configurations.
The use of saturated configurations also makes more transparent the
derivation of the Peierls bound in Section \ref{SubSec5.1}.

The idea of a saturated configuration worked well in the study of dense-packed circle 
configurations in $\bbR^2$; cf. \cite{ChW}. We attempt to emulate a similar approach on 
$\bbA_2/\bbH_2$. It generates some technical complications that are addressed in Lemmas \ref{Lem4.2} - \ref{Lem4.10}.

\bl\label{Lem4.2}%{\bf Lemma 4.2.}\label{LemmaXX} 
A saturated configuration on $\bbA_2/\bbH_2$
does not have \rV-circles of radius $\geq D+1$.
\el

\bp Suppose there exists a V-circle of radius $\geq D+1$. The center of the 
V-circle may not lie in $\bbA_2/\bbH_2$ but is at distance $\leq 1$ from one of the 
$\bbA_2/\bbH_2$-sites. Then an additional particle can be added at this site without 
breaking admissibility. This contradicts the saturation assumption. 
\ep

We would like to note a difference between Lemma \ref{Lem4.2} and Lemma 2 from \cite{ChW}.
We have a lower bound $D+1$ whereas in \cite{ChW}, Lemma 2, one has $D$. This
creates a particular technical complication arising on $\bbA_2/\bbH_2$ compared with 
$\bbR^2$.

Lemma \ref{Lem4.2} enables us to discard C-triangles which have a circumradius  $> D+1$ and 
focus on those with a circumradius $\leq D+1$ in our analysis of PGSs. The remaining 
obtuse C-triangles are tackled via a routine of the area re-distribution. More precisely, 
C-triangles with circumradius $\leq D-1$ are tackled in Lemmas \ref{Lem4.4}, 4.8 and 4.10,
depending upon the class of the value $D^2$. In Lemmas 4.5.1 - 4.5.3 we 
treat C-triangles on $\bbA_2/\bbH_2$ with a circumradius  between $D-1$ and $D+1$. 
Such a C-triangle, let us denote it by $\triangle$, can have area $<S(D)/2$ (when $\triangle$ is  
obtuse).  However, it turns out that in this case there will be an adjacent \rC-triangle $\triangle'$ 
(sharing a side with $\triangle$) with a rather large area, so that the area of the union 
$\triangle\cup\triangle'$ is $\geq S(D)+1$. It may also happen that two or three \rC-triangles
$\triangle_j$, of area $<S(D)/2$ each, share a common adjacent  triangle $\triangle'$; in this 
case there will again be a lower bound upon the area of their union. Such an 
observation allows us to circumspect obtuse C-triangles via Lemmas 4.5.1 - 4.5.3. 
For formal definitions, see Section \ref{SubSec4.2}.  
%justify the constraint $\alpha_i\leq\pi/2$ in the minimization  problem (4.1).

\subsection{Redistributed areas for triangles}\label{SubSec4.2} In this section we introduce re-distributed 
areas assigned to a \rC-triangle $\triangle$, which can be conveniently lower-bounded. 
One, $s^\RD(\triangle )$, characterizes the triangle {\it per se},  the other, 
$\varSigma (\triangle ,\phi)$, considers it in a $D$-AC $\phi$. The bounds involve the 
quantities $S(D)$ and $S^\RD (D)$ determined 
in \eqref{SoD} and \eqref{SRD(D)}, respectively. This will enable us  to analyze the PGSs 
for the whole array of the situations on $\bbA_2/\bbH_2$, including the exceptional non-sliding 
values $D^2=13, 16, 28, 49, 64, 67, 97, 157, 256$ (Classes HD and HE). 

An $\bbA_2/\bbH_2$-triangle $ABC$ is called a {\it qualifying triangle} if all its side-lengths 
are $\geq D$ while the circumradius is $\leq D+1$. All triangles we consider 
from now are supposed to be qualifying. 

A collection of two triangles $ABC$ and $ACE$ with mutually disjoint interiors is called 
a 2-{\it triangle group} if all sides and diagonals of quadrilateral $ABCE$ are not shorter than 
$D$, vertex $E$ does not lie inside the circumcircle of $ABC$, and vertex $B$ does not lie 
inside the circumcircle of $ACE$.

A collection of three triangles $ABC$, $CDE$ and $ACE$ with mutually 
disjoint interiors is called a 3-{\it triangle group} if all sides and diagonals of 
pentagon $ABCDE$ are not shorter than $D$, vertices $D, E$  do not lie inside the 
circumcircle of triangle $ABC$, vertices $A, B$ do not lie inside the circumcircle 
of $CDE$, and vertices $B, D$ do not lie inside the circumcircle of $ACE$.

A collection of four triangles $ABC$, $CDE$, $EFA$ and $ACE$ with 
mutually disjoint interiors is called a 4-{\it triangle group} if all sides and diagonals 
of hexagon $ABCDEF$ are not shorter than $D$, vertices $D, E, F$  do not lie
inside the circumcircle of $ABC$, vertices $F, A, B$ do not lie inside the circumcircle 
of $CDE$, vertices $B, C, D$ do not lie inside the circumcircle of $EFA$, and vertices 
$B, D, F$ do not lie inside the circumcircle of triangle $ACE$.

For each triangle group one can calculate the corresponding average triangle 
area which we call the {\it re-distributed group area}. 

For any triangle $ABC$ one can consider all triangle groups containing this triangle 
such that side $AB$ is shared with another triangle in the group but sides $BC$ and $CA$ 
are not shared. The minimal redistributed group area among 
all such groups is called the {\it $AB$-re-distributed area} of $ABC$ and 
denoted by $s^{\RD}_{AB}(ABC)$. The $BC$-re-distributed area 
$s^\RD_{BC}(ABC)$ and $CA$-re-distributed
area $s^{\RD}_{CA}(ABC)$ are defined in a similar way.

The quantity
\be\label{RDarea-s}s^\RD(ABC) = \max\big(s(ABC),s^{\RD}_{AB}(ABC),s^{\RD}_{BC}(ABC),
s^{\RD}_{CA}(ABC)\big)\ee
is called the {\it re-distributed area of triangle} $ABC$. Here and below $s(ABC)$ stands
for the area of $ABC$; a similar meaning will have the notation $s(\triangle)$, $s(\triangle\cup\triangle')$ 
and so on. If the maximum in \eqref{RDarea-s} is achieved
at $s^{ \RD}_\bullet (ABC)$ then the corresponding triangle side is called a  
{\it re-distributing side} (of $ABC$) and denoted by $\sigma (ABC)$.

\FigureN14

The doubled minimal redistributed area is denoted by $S^\RD (D)=S^\RD (D,\bbA_2/\bbH_2)$: 
\be\label{SRD(D)}\bear S^\RD (D)=2\times\min\,\Big[s^\RD(\triangle ):\;\triangle\;\hbox{ runs over
 the triangles on $\bbA_2/\bbH_2$}\Big].\ena\ee
A triangle $\triangle$ with minimal re-distributed area, i.e., with 
$s^\RD (\triangle )=S^\RD (D)/2$, is called an {\it MRA-triangle}. 

Note that if $\triangle$ has an area $s(\triangle )<
s^\RD (\triangle )$ then $\triangle$ has a re-distributing side.

An important observation is that, by virtue of Lemma \ref{Lem4.1} and the definition of an MRA-triangle,
any PGS consists of MRA-triangles. Cf. \cite{HS}, Criterium on P179.  %More precisely, we offer the
%following Lemma 4.3.
More precisely, we say that a $D$-AC is MRA-perfect if its C-triangles are all 
MRA-triangles.

\bl \label{Lem4.3}%{\bf Lemma 4.3.} {\sl 
Given an attainable $D^2$, suppose that there exists a perfect $D$-{\rm{AC}}
$\phi$. Then any $D$-{\rm{PGS}} is a perfect configuration.
\el

\bp Owing to Lemma \ref{Lem4.1}, the particle density in $\phi$ is $1/S^\RD (D)$.
Then any periodic $D$-AC has the particle density $\leq 1/S^\RD (D)$. Next, let $\psi$ be any periodic 
$D$-AC containing a non-MRA-triangle. Then, in a large basic quadrilateral
polygon $\bbV$, one can
construct a perturbation of $\psi$ having more particles in $\bbV$ than $\psi$ has.
In fact, such a perturbation will have the same pattern as $\phi$ in $\bbV$. We will
have to remove some particles from $\psi$ along the boundary $\partial\bbV$ 
but will gain an amount of particle proportional to the number of sites in $\bbV$.

Consequently, any PGS should consist of MRA-triangles. 
\ep

\brb {\rm In the course of this section, we will check that for any attainable $D^2$
on both $\bbA_2$ and $\bbH_2$ there exists at least one perfect $D$-AC.
Moreover, as we show further in this section, for any non-sliding $D$ the number of
perfect configurations (and hence that of the PGSs) is finite. Furthermore, it will be shown 
that every perfect configuration is periodic, hence a PGS. 
 
Consequently, any non-periodic ground state contains at least one infinite connected 
component of non-MRA triangles and no finite ones. Moreover, the number of non-MRA 
triangles in a $\bbA_2/\bbH_2$-hexagon $\bbV (L)$ of side-length $L$ can only grow linearly 
with $L$; this means that non-MRA triangles form, effectively, a one-dimensional array. Let us
repeat once more that, according to \cite{DoS}, non-periodic ground states do not generate EGMs
on $\bbA_2/\bbH_2$.} \hfill $\blacktriangle$
\erb

\bl\label{Lem4.4}%{\bf Lemma 4.4.}  {\sl 
For any $D^2$ on $\bbA_2$ and for any $D^2$ divisible by $3$ 
on $\bbH_2$: 
\begin{description}
\item[{\rm{(i)}}] for a $D$-triangle $\triangle^\circ$,
we have $s^\RD (\triangle^\circ )={\sqrt 3}D^2/4=S(D)/2$ (cf. Eqn \eqref{SoD}), 
\item[{\rm{(ii)}}] $\forall$ $D$-admissible $\bbA_2/\bbH_2$-triangle $\triangle$ non-congruent to $\triangle^\circ$
such that the circumradius of $\triangle$ is $\leq D-1$, we have 
\be\label{sRDtrianlgebound} \qquad\qquad\qquad\qquad\qquad s^\RD (\triangle)\geq s^\RD (\triangle^\circ )+\frac{\sqrt 3}{8}.\ee
\end{description}
\el

\bp (i) A $D$-triangle $ABC$ can be complemented by
its reflection about a given side, say $AB$, to form a 2-triangle group. Hence,
$s^\RD_{AB}(ABC)\leq S(D)/2$. By definition \eqref{RDarea-s}, it follows that
$s^\RD (ABC)= S(D)/2$.

(ii) The triangles under consideration in assertion (ii) have the maximum angle strictly
between $\pi/3$ and $2\pi/3$. The sinus of such an angle is $>{\sqrt 3}/2$.
Hence, $s (\triangle )> \diy\frac{D^2{\sqrt 3}}{4}=\frac{S(D)}{2}$. This implies 
\eqref{sRDtrianlgebound} for the area of an $\bbA_2$-triangle multiplied by $8/{\sqrt 3}$ is 
integer.
\ep

%For very obtuse triangles, with  For triangles with 

Consider an arbitrary saturated $D$-AC $\phi$ on $\bbA_2/\bbH_2$ and identify all triangles 
in $\phi$ with an area $<S^\RD (D)/2$. For each such triangle $\triangle$ consider a 2-triangle 
group $\triangle\cup\triangle'$
formed by a triangle $\triangle'$, called a {\it donor}, adjacent to $\triangle$ along the 
redistributing side $\sigma (\triangle )$. If several such 2-triangle groups 
have a common donor $\triangle'$ then we unite them into a single 3-triangle 
group or 4-triangle group. (By construction, a donor $\triangle'$
has area $>S^\RD (D)/2$.) In case of a 3-triangle group we have 
a donor $\triangle'$ with area $\geq S^\RD (D)/2$ grouped 
with two adjacent triangles of area $<S^\RD (D)/2$. In case of a 4-triangle group 
we have a donor $\triangle'$ of area 
$s (\triangle')\geq S^\RD (D)/2$ grouped with three adjacent triangles of area 
$<S^\RD (D)/2$. By construction, each triangle $\triangle$ in $\phi$ belongs to at 
most one group. Furthermore, the grouping uniquely assigns the {\it redistributed group 
area} $\varSigma (\triangle ,\phi)$ to each triangle $\triangle$ in the AC $\phi$. 
Namely, $\varSigma (\triangle ,\phi)$ is the total area of the triangles 
in the group divided by the number of the triangles in the group
containing $\triangle$ in $\phi$.

Next, if $\triangle$ is not a donor then $\varSigma (\triangle ,\phi)\geq s^\RD (\triangle )\geq
S^\RD (D)/2$. If $\triangle$ is a donor we have that $s^\RD (\triangle )\geq s(\triangle )
>\varSigma (\triangle ,\phi)\geq S^\RD (D)/2$; the last inequality holds since 
$\varSigma (\,\bullet\, ,\phi)$ is the same for all members in the group. Finally, 
if $\triangle$ does not belong to any group in $\phi$ then $s^\RD (\triangle )\geq s (\triangle )
=\varSigma (\triangle ,\phi) \geq S^\RD (D)/2$; the equality $s (\triangle )
=\varSigma (\triangle ,\phi)$ and inequality $\varSigma (\triangle ,\phi) \geq S^\RD (D)/2$ 
follow directly from the way in which $\triangle$ is identified in $\phi$.

%; the last equality and inequality express the way how this $\triangle$ is defined.}

%As soon as there exist at least one $D$-admissible configuration consisting solely 
%of MRA-triangles (an RD-perfect configuration), all theory constructed in Lemmas 
%4.2--4.10 for M-triangles is applicable verbatim to MRA-triangles. 

\subsection{Equality $S^\RD (D)=S(D)$ on $\bbA_2$ and -- for $3|D^2$ -- on $\bbH_2$} \label{SubSec4.3} In this section we 
give three lemmas, 4.5.1--4.5.3, treating C-triangles on $\bbA_2/\bbH_2$ with 
a circumradius $r=D+\delta$, $\delta\in [-1,1]$. Then we proceed with Lemma 4.6 which,
together with Lemma \ref{Lem4.3}, 
establishes the equality $S^\RD (D)=S(D)$ under some conditions upon $D^2$.

 \def\rC{{\rm C}}
 
\blema \label{Lem4.5.1}%{\bf Lemma 4.5.1.} {\sl 
Suppose that a \rC-triangle $\triangle$ has the circumradius
$r=D+\delta$ where $-1\leq\delta\leq 1$. Then 
\be\label{eq:Tsep'}
\qquad\qquad\qquad s(\triangle )\geq \diy\frac{D^3}{2r}\sqrt{1-\frac{D^2}{4r^2}}
>\diy\frac{{\sqrt 3}D^2}{4}-\frac{D\delta}{2\sqrt 3}.\ee 
Here $\diy\frac{D^3}{2r}\sqrt{1-\frac{D^2}{4r^2}}$ is the area of an isosceles triangle 
with circumradius $r$ and two side-lengths $D$.  
The longest side in this triangle has length $\diy <D{\sqrt 3}+\frac{\delta}{\sqrt 3}$.
\elema

\bp Suppose a \rC-triangle $\triangle$ with vertices $A, B, C$ satisfies the
assumptions of the lemma. Let the side-lengths be $AB=l_0$, $BC=l_1$, $CA=l_2$,
with $D\leq l_0\leq l_1\leq l_2\leq 2r$.
If two side-lengths are $>D$, say $l_1 ,l_2>D$, then the area of the $\triangle$ can be
made smaller by moving vertex $C$ along the circumcircle towards $B$, until
the length of side $BC$ becomes $D$. Indeed,
in the process of motion $l_0$ remains fixed but the height from $C$ to $AB$ shortens.
Thus, the area of $\triangle$
is lower-bounded by the area of an isosceles triangle with two side-lengths $D$ and the
remaining side-length $\diy 2D\sqrt{1-\frac{D^2}{4r^2}}$. (On  $\bbH_2$, it 
is not necessarily an $\bbH_2$-triangle.) A direct calculation shows that for $D\geq 1$
and $\delta\in (-1,1)$ the bound $\diy 2D\sqrt{1-\frac{D^2}{4r^2}}
< D{\sqrt 3}+\frac{\delta}{\sqrt 3}$ holds true.
(The right-hand side is simply the Taylor expansion in $\delta$ up to order $1$.)
The area of such a triangle equals 
$\diy\frac{D^3}{2r}\sqrt{1-\frac{D^2}{4r^2}}$. Finally, 
$\diy\frac{D^3}{2r}\sqrt{1-\frac{D^2}{4r^2}}= \frac{D^3}{2(D+\delta )}\sqrt{1-\frac{D^2}{4(D+\delta)^2}}
>\frac{{\sqrt 3}D^2}{4}-\frac{D\delta}{2\sqrt 3}$.
\ep

\blemb\label{Lem4.5.2}%{\bf Lemma 4.5.2.} {\sl 
Suppose that a \rC-triangle $\triangle$ with side-lengths 
$l_0, l_1, l_2$ has the circumradius $r=D+\delta$ where 
$-1\leq\delta\leq 1$. Consider an adjacent \rC-triangle $\triangle'$ that shares
with $\triangle$ the longest side (of length $l_2$). Then the area 
$s(\triangle\cup\triangle' )$ is lower-bounded by the area of a trapeze 
inscribed in a circle of radius $r$, with three sides being of length $D$. 
Furthermore, for $D^2\geq 400$ we have $s(\triangle\cup\triangle')\geq 
\diy\frac{3{\sqrt 3}D^2}{4}-2\delta^2$.
\elemb

\bp Again, we assume
$D\leq l_0\leq l_1\leq l_2\leq 2r$. Two vertices of triangle $\triangle'$ are the
end-points of the side of length $l_2$ and lie in the V-circle of radius $r$ circumscribing
$\triangle$. The third vertex of $\triangle'$ cannot lie inside this V-circle
but can be placed on the circle. It also should lie outside the circles of radius $D$
centered at the end-points of the side of length $l_2$. Under these restrictions, the minimal area
of $\triangle'$ is not less than the area of a triangle inscribed in the V-circle which shares the side
 of length $l_2$ with $\triangle$ and has the other side of length $D$. (Cf. the proof of Lemma 
4.5.1.) If we now minimize the area of $\triangle$, we obtain a pair $\triangle$, $\triangle'$ 
 forming a trapeze, as specified in the assertion of Lemma 4.5.2. (Again, on $\bbH_2$
 the resulting triangle is not necessarily an $\bbH_2$-triangle.)

The area of the trapeze in question equals the sum of the areas of 4 triangles, 3
 of which are identical. The area of each of these identical triangles is $\diy\frac{r^2}{2}\sin (2\alpha )$
 where $\sin (\alpha )=\diy\frac{D}{2r}$. The area of the fourth triangle 
 is $\diy\frac{r^2}{2}\sin (2\pi -6\alpha)$.
 All-in-all, the area of the trapeze is $\diy\frac{r^2}{2}4\sin^3(\alpha )$, which equals
$$\frac{2D^3}{r}\left(\sqrt{1-\frac{D^2}{4r^2}}\right)^3=\frac{3{\sqrt 3}D^2}{4}-{\sqrt 3}\delta^2
+\frac{19\delta^3}{3{\sqrt 3}D}-\frac{113\delta^4}{9{\sqrt 3}D^2}+\ldots $$
A  straightforward calculation asserts that for $D^2\geq 400$ and
$-1\leq\delta\leq 1$ this expression is
$\diy\geq\frac{3{\sqrt 3}D^2}{4}-2\delta^2$, as claimed in the lemma. 
\ep

\blemc \label{Lem4.5.3} %{\bf Lemma 4.5.3.} {\sl  
Suppose that a $\rC$-triangle $\triangle$ has the circumradius 
$r=D+\delta$ where $-1\leq\delta\leq 1$. Let $\triangle'$  be the adjacent \rC-triangle 
sharing the longest side with $\triangle$ (cf. Lemma {\rm{4.5.2}}).
\begin{description}
\item[{\rm{(i)}}] Suppose that
$\triangle'$ is adjacent to another \rC-triangle, $\triangle_1$, with circumradius
$r_1=D+\delta_1$ where $-1\leq\delta_1\leq 1$. Then we have $s(\triangle')
\geq  3D^2/4$.

\item[{\rm{(ii)}}] Further, suppose $\triangle'$ is adjacent to other two \rC-triangles,
$\triangle_1$ and $\triangle_2$, with
circumradii $r_1=D+\delta_1$ and $r_2=D+\delta_2$ where $-1\leq\delta_1,\delta_2\leq 1$.
Then $s(\triangle' )\geq D^2$.
\end{description}
\elemc

\bp  (i) Here the triangle $\triangle'$ has one side-length $\geq D$ and two others
$\geq D\sqrt 3$ by construction. On the other hand, the side-lengths 
are $\leq 2D+2$ since the circumradius is \ $\leq D+1$. Therefore, the area of $\triangle'$
is greater than or equal to the area of a triangle with side-lengths $D$, $D{\sqrt 3}$, $D{\sqrt 3}$.
The area of such a triangle is, clearly,  $\geq 3D^2/4$.

(ii) In this case all side-lengths of $\triangle'$ are $\geq D\sqrt 3$. Hence, 
the area of $\triangle'$ is $\geq D^2$. 
\ep

\blemd\label{Lem4.6}%{\bf Lemma 4.6.} {\sl 
For any $D^2$ on $\bbA_2$ and  for any $D^2$ divisible by $3$ on 
$\bbH_2$, we have that 
\be\label{SRD=S} \qquad\qquad\qquad\qquad S^\RD (D)=S(D) ,\ee
and the equality $S^\RD (D)=2s^\RD (\triangle )$ is attained only when $\triangle$
is a $D$-triangle. For each of these values of $D$, the corresponding
{\rm{MRA}}-perfect configuration exists and has type $(D,\alpha )$.
\elemd

\bp %For $D^2\geq 400$, the area $S^\RD(D)$  satisfies
%\be\label{MinAreaS(D)}S(D)\leq S^\RD(D)<S(D)+9\sqrt{3D}.\ee
%Consequently, the following holds true.  
In the situation of Lemma 4.4 we have the bound 
\be\label{sRDbound}2s^\RD (\triangle)\geq S(D)+\frac{\sqrt 3}{4}.\ee

Next, in the situation of Lemma {\rm{4.5.2}} (in particular, for $D^2\geq 400$) we have:
\be\label{MinAreaS(D)1}2s^\RD (\triangle )\geq \frac{3}{2}S(D)-2\delta^2\geq S(D)+\frac{\sqrt 3}{2}.\ee
Next, in case {\rm{(i)}} of Lemma {\rm{4.5.3}},
\be\label{MinAreaS(D)2}\diy 3s^\RD (\triangle )\geq 2\left (2S(D)
-\frac{D\delta}{2\sqrt 3}\right)+\frac{\sqrt 3}{2}S(D)\geq
\frac{3}{2}\left(S(D)+\frac{\sqrt 3}{2}\right).\ee
Finally, in case {\rm{(ii)}} of Lemma {\rm{4.5.3}},
\be\label{MinAreaS(D)3}\diy 4s^\RD (\triangle )\geq 3\left(2S(D)-\frac{D\delta}{2\sqrt 3}\right)
+\frac{2}{\sqrt 3}S(D)\geq 2\left(S(D)+\frac{\sqrt 3}{2}\right).\ee
Together, \eqref{sRDbound}, \eqref{MinAreaS(D)1}, \eqref{MinAreaS(D)2}, \eqref{MinAreaS(D)3} 
imply that for $D^2\geq 400$ and every C-triangle $\triangle$ different from a $D$-triangle: 
$2s^\RD (\triangle )>S(D)$. This implies the assertion of Lemma 4.6 for $D^\geq 400$.

For $1\leq D^2< 400$ the proof is done by a computer enumeration. 
\ep

\subsection{MRA-triangles for Class HC on $\bbH_2$}\label{SubSec4.4}
Next, we analyze the situation on $\bbH_2$, for values $D^2$ not divisible by $3$.

\bleme \label{Lem4.7}%{\bf Lemma 4.7.} {\sl 
For any L\"oschian $D^2\geq 300$ there exists a L\"oschian number that 
is $\geq D^2$, is divisible by $3$ and is at distance at most $18{\sqrt D}$ from $D^2$. 
\eleme

\bp Consider L\"oschian numbers divisible by 3 of the form
$$(l-3k)^2 + (l+3k)^2 + (l-3k)(l+3k) = 3l^2+9k^2.$$
(It is simply the set of all  L\"oschian numbers scaled 3 times.) Now, take an arbitrary  L\"oschian number
$D^2$ and find $l$ such that $3l^2 \le D^2 \le 3(l+1)^2$. Then find $k$ such that
$$3l^2 + 9k^2 \le D^2 \le \min\big(3(l+1)^2,\; 3l^2 + 9(k+1)^2\big)$$
Then $9k^2 \le 6l+3$, i.e. $k \le \sqrt{(2l+1)/3}$.
The distance from $D^2$ to $\min\big(3(l+1)^2,\; 3l^2 + 9(k+1)^2\big)$ is at most
\be\beal 9(k+1)^2 -9k^2 = 18k+9  \\
\qquad\qquad <  18 \sqrt{(2l+1)/ 3} + 9 < 18 \sqrt{l} \le 18 \sqrt{\sqrt{D^2 /3}}  
\le 18{\sqrt D},\ena\ee 
where  inequality involving $l$ in the middle is true for $l > 9$. 
\ep
 
In what follows, we refer to $D^*(=D^*(D))$ as the nearest L\"oschian number not less than $D$
 such that $3|(D^*)^2$.  
 
\blemf\label{Lem4.8}%{\bf Lemma 4.8.} {\sl  
Any non-equilateral $D$-admissible $\bbH_2$-triangle 
$\triangle$ with circumradius \ $\leq D-1$ and the shortest side-length $<D^*$
has at least one side with squared length 
 $\geq D^2 + D+1$. %and $\diy (2D)^2\left(1-\frac{D^2}{4(D-1)^2}\right)$. 
Consequently, for the double area $2s(\triangle )$ we have: 
\be\label{eq:L4.8_1} \beal 2s(\triangle )\geq h(D)\;\hbox{ where}\\
\diy\qquad h(D)=\min\,\bigg[\frac{D^3}{D-1}
\sqrt{1-\frac{D^2}{4(D-1)^2}},\;\frac{1}{2}\sqrt{(3D^2-D-1)(D^2+D+1)}\bigg].
\ena\ee 
Furthermore, for $D^2\geq 12$: 
\be\label{eq:L4.8_2}
\diy \qquad\qquad h(D)>\frac{\sqrt 3}{2}D^2 + \frac{D}{2{\sqrt 3}}.\ee
\elemf

\FigureO15

\bp Referring to Figure 15, suppose that triangle $\triangle$ is $OCB''$, and its shortest
side is $OC$, with $D\leq |OC|<D^*$. Hence, $|OC|^2$ is not divisible by 3. 

Let $B'$ be the vertex of an equilateral triangle, with $|OB'|=|CB'|=|OC|$. Then $B'$ will be at
the center of a unit hexagon, as shown in Figure 15. Consequently, $|B'B''| \ge 1$, and 
at least one of the triangles $OB'B''$ or $CB'B''$ is obtuse with the corresponding obtuse angle $\geq 2\pi/ 3$. 
Therefore, by the cosine theorem the squared length of the longest side of the obtuse triangle is at least $D^2 + D + 1$. 
Hence, as long as triangle $OCB''$ is acute, we have that 
$s(OCB'')\geq\diy\frac{1}{2} \sqrt{(3D^2-D-1)(D^2+D+1)}$. The latter value is the area of an 
isosceles triangle with side-lengths $D$, $D$ and $\sqrt{D^2+D+1}$ (and the circumradius $D-1$). 

On the other hand, if $OCB''$ is obtuse then, according to Lemma 4.5.1, 
$s(OCB'' )\geq \diy\frac{D^3}{2(D-1)}\sqrt{1-\frac{D^2}{4(D-1)^2}}$. 

The last assertion of the lemma is straightforward for $D^2\geq 12$. 
\ep

%\vskip .5cm
%\newpage

%{\bf Lemma 4.9.} {\sl  Suppose $D^2 > 54^4$ and $(D^*)^2\geq D^2$ is the L\"oschian number 
%closest to $D^2$
%which is divisible by 3. Then for the doubled area of the equilateral $D^*$-triangle $\triangle^*$ we have
%\be 2s(\triangle^*)= \frac{{\sqrt 3}(D^*)^2}{2} \le 
%\frac{{\sqrt 3}(D^2 + 18{\sqrt D})}{2} = \frac{\sqrt 3}{2}\big(D^2 + 18{\sqrt D}\big).\ee}
%\vskip .5cm

% {\bf Proof.} If  $(D^*)^2\geq D^2$ is the L\"oschian number closest to $D^2$
% which is divisible by 3 then the doubled area of the equilateral $D^*$-triangle is
% \be 2 (D^*)^2 \le 2(D^2 + 18{\sqrt D}) = 2D^2 + 36{\sqrt D} < 2D^2 +\frac{2\sqrt{D^2}}{3}\ee
% as soon as $D^2 > 54^4$.  \qquad $\blacksquare$
%\vskip .5cm

%{\bf Lemma 4.10.} {\sl For any value $D^2$ of Class \rH\rC \ on $\bbH_2$:} (i) 
%{\sl for an equilateral $D^*$-triangle $\triangle^*$
%we have $s^\RD (\triangle^* )={\sqrt 3}(D^*)^2/4=S(D^*)/2$ (cf. Eqn \eqref{SoD}),} (ii)
%{\sl $\forall$ $D$-admissible $\bbH_2$-triangle $\triangle$ not congruent to $\triangle^*$
%and such that the circumradius of $\triangle$ is not longer than $ D+1$, we have $s^\RD (\triangle)
%\geq s^\RD (\triangle^*)+1/2$.} \vskip .5cm

%{\bf Proof.} The argument used in the proof of Lemma 4.4(i) is also valid for the proof of 
%assertion (i), and we omit it to avoid a repetition. 

%(ii)

%\vskip .5cm
 
\blemg \label{Lem4.9} %{\bf Lemma 4.9.} {\sl 
For any value $D^2$ of Class \rH\rC \ on $\bbH_2$
we have that 
\be\label{(4.11)} \qquad\qquad\qquad S^\RD (D)=S(D^*),\ee
and the equality $S^\RD (D)=2s^\RD (\triangle )$ is attained only when $\triangle$
is congruent to a $D^*$-triangle $\triangle^*$. Moreover,
for each value of $D$ of Class \rH\rC, the corresponding {\rm{MRA}}-perfect configuration
exists and has type $\alpha (D^*)$.
\elemg

\bp By construction, $s^\RD (\triangle )\geq s(\triangle )$ for any $D$-admissible 
$\bbA_2/\bbH_2$-triangle. First, we consider triangles satisfying the conditions of 
Lemma 4.8. For any such $\triangle$ we have
\be\label{(4.12)}2s^\RD (\triangle )\geq 2s (\triangle )\geq h(D).\ee
If $D^2\geq (54)^4$ then, by Lemma 4.7,  for any such $\triangle$, 
\be h(D)\geq \frac{\sqrt 3}{2}\big(D^2 + 18{\sqrt D}\big)\geq\frac{{\sqrt 3}(D^*)^2}{2}=2s(\triangle^*),\ee
and therefore $s^\RD (\triangle )>s(\triangle^*)$.
Hence, such a triangle cannot be an MRA-triangle if $D^2\geq (54)^4$. When $D^2< (54)^4$ then, 
instead of utilizing Lemma 4.7 we verify numerically that, apart from 184 values,
every $D^2< (54)^4$ satisfies the bound $\big(D^2 + D/3\big) >(D^*)^2$, which again implies 
that $s^\RD (\triangle )> s(\triangle^*)$, with the help of \eqref{eq:L4.8_1}.  Cf. 
Section \ref{Sec9} and Program 1
{\tt NearestLoschianNumber} in the ancillary file.
The non-exceptional $D$ among 184 remaining values are tackled by a separate computer 
program which calculates  $S^\RD (D)$ and verifies \eqref{(4.11)}. Cf. Section \ref{Sec9} and Program 
2 {\tt SpecialD} in the ancillary file. 

Next, we discuss the case where we have a non-equilateral triangle $\triangle$ with 
circumradius $\leq D-1$ and the shortest side-length $\geq D^*$. Here, the inequality $2s(\triangle )
\geq S^\RD (D^*)$ is straightforward if $\triangle$ is acute and follows from Lemma 4.5.1 if
$\triangle$ is obtuse.

Finally, a $D$-admissible triangle with circumradius between $D-1$ 
and $D+1$ cannot be an MRA-triangle by virtue of an argument similar to the one in the 
proof of Lemma 4.6. A lower bound $s^\RD (\triangle )> s(\triangle^*)$ for such a
triangle $\triangle$ is obtained by 
repeating the proof of Lemma 4.6 where we use analogs of inequalities 
\eqref{MinAreaS(D)1}, \eqref{MinAreaS(D)2} and \eqref{MinAreaS(D)3} with the value 
$S(D)$ in the RHS replaced by  $S(D^*)$. Namely, in the situation of  
Lemma {\rm{4.5.2}}, 
\be\label{(4.15)} 2s^\RD (\triangle )\geq \frac{3}{2}S(D)-2\delta^2\geq S(D^*)+\frac{\sqrt 3}{2},\ee
in case {\rm{(i)}} of Lemma {\rm{4.5.3}},
\be\label{(4.16)} 3s^\RD (\triangle )\geq 2\left (2S(D)
-\frac{D\delta}{2\sqrt 3}\right)+\frac{\sqrt 3}{2}S(D)\geq
\frac{3}{2}\left(S(D^*)+\frac{\sqrt 3}{2}\right),\ee
and in case {\rm{(ii)}} of Lemma {\rm{4.5.3}},
\be\label{(4.17)} 4s^\RD (\triangle )\geq 3\left(2S(D)-\frac{D\delta}{2\sqrt 3}\right)
+\frac{2}{\sqrt 3}S(D)\geq 2\left(S(D^*)+\frac{\sqrt 3}{2}\right).\ee

For $D^2\geq (54)^4$ bounds \eqref{(4.15)}-\eqref{(4.17)} are a consequence of 
Lemma 4.7. For non-exceptional $D^2< (54)^4$, \eqref{(4.15)} and \eqref{(4.16)} are 
verified numerically -- cf. Section \ref{Sec9} and Program 1 {\tt NearestLoschianNumber} in the ancillary file --
while second bound in \eqref{(4.17)} follows from the second bound in \eqref{(4.16)}. 
This implies the desired estimate $s^\RD (\triangle )> s(\triangle^*)$ for a 
$D$-admissible triangle $\triangle$ with circumradius between $D-1$ 
and $D+1$.

Thus, it is established that for any $\triangle$ with circumradius $\leq D+1$ not
congruent to $\triangle^*$ we have the bound 
\be\label{sRD>sRD*} \qquad\qquad\qquad s^\RD (\triangle )>s^\RD (\triangle^*).\ee 
This leads to the assertions of Lemma 4.9. 
\ep

\subsection{MRA-triangles for Classes HD and HE on $\bbH_2$}\label{SubSec4.5} In this section we establish 
the values $S^\RD (D)$ when $D$ is exceptional and non-sliding. We use the notation $[l^2_0|l^2_1|l^2_2]$,
referred to as a triangle type, to indicate a triangle with side-lengths $l_0\leq l_1\leq l_2$. For example, in Figure 7, triangles 
$AOB$, $AOH$, $HFO$, $CED$  have type $[13|19|21]$ 
whereas triangles $OBC$, $OCE$, $OFE$ have type $[13|16|21]$.

\blemh\label{Lem4.10} %{\bf Lemma 4.10.} {\sl  
The \RD-perfect configurations exist for all exceptional $D^2=$ 
$13$, $16$, $28$, $49$, $64$, $67$, $97$, $157$, $256$ and are periodic. The corresponding 
values of $S^\RD (D)$
and triangle groups on $\bbH_2$ at which these values are achieved are as follows:
\be\label{SRDexceptional}\begin{array}{llc}
S^\RD (\sqrt{13}) = 16.5 {\sqrt 3}/2  &\{[13|16|21],[13|19|21]\} &\beta\\
S^\RD (\sqrt{16}) = 20.25 {\sqrt 3}/2 &\{[21|21|21],[16|21|25], &\\
&\;\;[16|21|25],[16|21|25]\}&\gamma\\
S^\RD (\sqrt{28}) = 33{\sqrt 3}/2 &\{[28|31|39],[28|37|39]\}&\beta\\ 
S^\RD (\sqrt{49}) =55.5{\sqrt 3}/2 &\{[49|52|63],[49|61|63]\}&\beta\\
S^\RD (\sqrt{64}) =72{\sqrt 3}/2 &\{[64|73|81]\}&\beta\\
S^\RD (\sqrt{67}) = 75{\sqrt 3}/2 &\{[75|75|75]\}&\alpha (\sqrt{75})\\
S^\RD (\sqrt{67}) = 75{\sqrt 3}/2 &\{[67|73|84],[67|79|84]\}&\beta \\
S^\RD (\sqrt{97}) = 106.5{\sqrt 3}/2 &\{[97|103|117],[97|112|117]\}&\beta \\ 
S^\RD (\sqrt{157}) = 169.5{\sqrt 3}/2 &\{[157|169|183],[157|172|183]\}&\beta\\ 
S^\RD (\sqrt{256}) = 272.25 {\sqrt 3}/2 &\{[273|273|273],[256|273|289], &\\  
&\;\;[256|273|289],[256|273|289]\}&\gamma . \ena\ee
respectively. For each of these values of $D$, the corresponding {\rm{MRA}}-perfect configuration
exists, and its type is listed in the right column.
\elemh

\bp The calculation of $S^\RD (D)$ involves a finite number of 
qualified triangles is performed by Program 2 {\tt SpecialD}. Cf. Section \ref{Sec9}
and ancillary file.
\ep

\brc {\rm Lemma 4.10 indicates that the exceptional values of $D$ emerge when 
the MRA-triangles are non-unique or non-equilateral. As a result, we have 
PGSs which are not obtained from max-dense sub-lattices. E.g., for $D^2=67$ we have
(i) an MRA-triangle that is a $D$-triangle with $D^2=75$, and (ii) a 2-triangle group formed by
non-equilateral MRA-triangles.

Another notable case is $D^2=64$ where an MRA-triangle is unique but not equilateral
(and forms a group on its own). Here all occupied sites in the $\beta$-PGSs have V-cells 
of area $72{\sqrt 3}/2$; these V-cells are congruent hexagons. However, there exist ACs 
where some V-cells (still hexagons) have area $71.5{\sqrt 3}/2$.} \hfill $\blacktriangle$
\erc

\subsection{Proof of Theorem I}\label{SubSec4.6} 

\bprfI Owing to Lemma 4.3,  
Theorem I follows from Lemma 4.6 and 4.9 and 4.10 establishing the existence of RD-perfect
configurations for all non-sliding values of $D$. 
\eprfI

\brdd {\rm A corollary of Theorem I is that the particle density in a PGS (per a unit Euclidean
area) equals $1/S^\RD (D)$.} \hfill $\blacktriangle$
\erdd
 
\section{The Peierls bound}\label{Sec5}

\subsection{The Peierls bound via MRA-triangles}\label{SubSec5.1} As was said, an application of the PS theory needs a
Peierls bound. Here we establish the Peierls bound by using the machinery of MRA-triangles.
We again begin with some auxiliary notions and statements. Throughout Section \ref{SubSec5.1} we assume 
that on $\bbA_2$, $D$ takes any attainable value while on $\bbH_2$ the value $D$ is non-sliding
(i.e., not from Class HS).

Let $\phi^*$ be a saturation of a given $D$-AC $\phi$. If an added occupied site
$x\in\phi^*\setminus\phi$ lies in a template then, clearly, this template is incorrect (more
precisely, non-$\vphi$-correct in $\vphi$ for each $\vphi\in\sP$). We say that such a template
is an {\it s-defect} (in $\phi$). Another possibility for a defect is where, in the saturation
$\phi^*$, a template has a non-empty intersection with one of C-triangles
that is not an MRA-triangle. We call it a {\it t-defect} (again in $\phi$). Finally, an incorrect 
template can be simply a neighbor of an s- or a t-defect. 
We call it an {\it n-defect} (still in $\phi$). Observe that 
any triangle intersecting the support of the n-defect template is an MRA-triangle.

We would like to note that C-triangles considered in Lemmas 4.5.1--4.5.3
lead to t-defects by definition.

\bl \label{Lem5.1}%{\bf Lemma 5.1.} 
{\rm{(A Peierls bound in terms of defects)}} Let $D$ be not from Class {\rm{HS}}. Consider 
a $\vphi$-contour $\Gam =\big({\rm{Supp}}\,(\Gam ),\phi\upharpoonright_{{\rm{Supp}}\,(\Gam )}\big)$ 
containing $m=\|({\rm{Supp}}\,(\Gam )\|$ incorrect templates. 
%enclosed by $n$ adjacent $\wt\vphi$-correct templates belonging to
%the boundary layers of $\Gam$ and  (for the same
%or different values of $\vphi$). 
Additionally, assume that $m = i + j + k$ where $i, j, k$
give the amount of \rs-, \rt- and \rn-defects in $\phi$, respectively. Then for the weight $w(\Gam )$
we have that
\be \qquad\qquad\qquad w(\Gam  )\leq u^{-i-j{\sqrt 3}/(32 S^\RD (D))}.\ee

%\beal\hbox{the amount of occupied 
%sites in $\phi\upharpoonright_{{\rm{Supp}}\,\Gam}$  $\Gam$ is}\\
%\quad\hbox{$\leq$ $\Big($the amount of occupied 
%sites in $\vphi\upharpoonright_{{\rm{Supp}}\,\Gam}\Big)$ $-$ 
%$\diy \left(i+\frac{j}{8S^\RD(D)}\right)$.}\ena\ee}
\el

\bp  The integer value $i$ contributed by s-defects is straightforward, so we
consider the saturation $\phi^*$ and its t-defects only. 

Observe that, according to Lemmas 4.4, 4.6, 4.9, 4.10, the re-distributed area of any C-triangle
that is not an MRA-triangle is at least $\diy\frac{1}{2}\left(S^\RD (D)+\frac{\sqrt 3}{2}\del (D)\right)$, 
where $\del (D) \ge 1/2$. (The overall minimal value $1/2$ for $\del (D)$ is attained, e.g., for 
$D^2=16$ in Lemma 4.10.) 

Further, a C-triangle that is not an MRA-triangle can be shared by at most 4 templates. 
Therefore, $j$ templates with t-defects contain (in the $\bbR^2$-sense) at least $j/4$ 
C-triangles that are not MRA-triangles. Consider a torus $\bbT$ formed by an integer number 
of templates and containing ${\rm{Supp}}\,(\Gam )$. Then $\bbT$ contains at most $s(\bbT )/S^\RD (D)$ C-triangles 
where $s(\bbT )$ is the area of $\bbT$.  
  
On the other hand, the maximal possible amount of C-triangles in 
$\phi^*\upharpoonright_{\bbT}$ is $\leq\big(s(\bbT )- j{\sqrt 3}/16\big)\big/S^\RD (D)$.
Next, owing to Lemma \ref{Lem4.1}, the number of particles in $\phi^*\upharpoonright_{\bbT}$
and $\vphi\upharpoonright_{\bbT}$ is obtained by dividing the amount of C-triangles by a 
factor 2. Finally, we can pass from $\bbT$ to ${\rm{Supp}}\,(\Gam )$ as the number
of particles in $\phi^*\upharpoonright_{\bbT\setminus{\rm{Supp}}\,(\Gam )}$
and $\vphi\upharpoonright_{\bbT\setminus{\rm{Supp}}\,(\Gam )}$ is the same by construction.
\ep

Informally, Lemma \ref{Lem5.1} states that the increment of `energy' (i.e., decrease in the number
of particles) caused by a deviation from a PGS
is lower-bounded proportionally to the `size' of the deviation. This is the gist of Peierls
bounds used in the Pirogov--Sinai theory and its applications.

\bl \label{Lem5.2}%{\bf Lemma 5.2.} {\sl 
Let $\vphi', \vphi''\in{\sP}(D)$ be two distinct
{\rm{PGS}}s. Consider a $D$-AC $\phi$ containing a connected component $\Lam$ of 
$\vphi'$-correct templates enclosed by a connected component of $\vphi''$-correct 
templates. Then $\phi$ contains a closed chain of adjacent non-{\rm{MRA}} \rC-triangles enclosing \ $\Lam$.
\el

\bp On $\bbA_2$ and on $\bbH_2$ for $D$ from Classes HA, HB, HC,
the MRA-triangles are equilateral. Such triangles from two distinct PGSs
cannot share a side in a $D$-AC. For $D$ from Classes HD and HE on $\bbH_2$
the assertion is verified case-by-case. 
\ep
%(A Peierls estimate for contours)

%{\bf Lemma V.} {\rm{}} {\it Consider a finite $\vphi$-contour  $\Gam$
%containing $m$ $\vphi$-incorrect templates. Then %the amount of occupied sites inside this 
%contour is at most
%\be w(\Gam )\leq u^{-\|{\rm{Supp}}\,(\Gam )\|{\sqrt 3}/288}.\ee}
%\vskip .5cm

%{\bf Proof of Theorem II.}
\bprfII
The theorem is a direct consequence of Lemmas \ref{Lem5.1} and \ref{Lem5.2} with an
additional factor
$1/9$ accounting for the possibility for each s- or t-defect to be surrounded by 8
n-defects. 
\eprfII

%Lemma V completes the verification of assumptions of the PS theory for the hard-core models
%on $\bbA_2/\bbH_2$. We are now in position to prove Theorem II.

%\vskip .5cm

%{\bf Proof of Theorem II.} Theorem II follows directly from Lemma V.

%\vskip .5cm

\subsection{A Peierls bound via Voronoi cells}\label{SubSec5.2} An alternative method of establishing the Peierls
bound is to use V-cells: it works on $\bbA_2$ and -- when $3|D^2$ -- on $\bbH_2$ (Classes HA
and HB). Thus, from now on until the end of Section \ref{SubSec5.2} we assume that the attainable
value $D^2$ is arbitrary on $\bbA_2$ and is divisible by 3 on $\bbH_2$. Consequently, the
PGSs are configurations of type $(D,\alpha )$ (obtained from $D$-sub-lattices). 
The V-cell method is considerably shorter than the MRD-triangle method but it has a drawback 
that the obtained Peierls constant $\ovp =\ovp (D)$ is not explicit. 

It is known \cite{F, Hs} that for any given $D$, a V-cell with the minimal possible area 
among $D$-ACs $\phi\in\cA (D,\bbR^2)$ 
is a perfect hexagon with the side length $D/\sqrt{3}$ and area $S=S(D)$ defined in Eqn \eqref{SoD}.
We call it a {\it perfect} V-cell.  A $D$-AC $\phi\in\cA (D,\bbR^2)$ 
is called V-{\it perfect} if it contains only perfect V-cells. The only V-perfect 
$D$-AC $\phi\in\cA (D,\bbR^2)$ are triangular lattices in $\bbR^2$ with 
the distance $D$ between neighboring lattice sites (see again \cite{Hs}).

\bl \label{Lem5.3}%{\bf Lemma 5.3.} {\sl 
For each $D$ there exists a number $\odel =\odel (D, \bbA_2/\bbH_2)>0$ such 
that the area of a non-perfect \rV-cell in any {\rm{AC}} $\phi\in\cA$ is $>S + \odel (D)$. 
\el

\bp As follows from \cite{Hs}, to analyze optimal and next-to-optimal V-cells for $\bx$
in $\phi\in\cA$, it suffices to consider sites at distance at most $4D$ from $\bx$ which
yields finitely many possibilities of drawing V-cells on $\bbA_2$ or  $\bbH_2$.
There are always an optimal and a next-to-optimal cells as not all possibilities are the same.
\ep

\bree {\rm It is precisely the fact that  $\odel $ is not determined explicitly that leads to a 
non-explicit Peierls constant $\ovp$ in Lemma \ref{Lem5.4}.} \hfill $\blacktriangle$
\eree

Given a basic polygon $\bbV$, a PGS $\vphi\in\sP$ and an AC $\phi_\bbV\in\cA (\bbV\|\vphi )$, 
we have the set-theoretical identity
\be\label{(4)}
\qquad\qquad\qquad \bigcup_{\bx \in\phi_V}\cV (\bx, \phi_\bbV)= \bigcup_{\bx \in\vphi\upharpoonright_\bbV}
\cV (\bx, \vphi\upharpoonright_\bbV).\ee
Therefore, since PGS $\vphi$ is an $(D,\alpha )$-configuration, for the partition function \eqref{PartFnctnV} we 
have that
\be\label{(5)}
\qquad\qquad\qquad \bZ (\bbV\|\vphi ) = \sum_{\phi_\bbV\in\cA (\bbV\|\vphi )}\; \prod_{\bx \in\phi_\bbV} u^{-S^{-1} 
\left( \left | \cV (\bx,\phi_\bbV)\right| - S \right) }.   \ee
Here, and in Lemma \ref{Lem5.4} below, we use the notation $\left| \cV (\bx,\phi_\bbV) \right|$ and
$\left| \cV (\bx,\phi\upharpoonright_\Gam) \right|$ for the 
area of $\cV (\bx,\phi_\bbV)$ and $\cV (\bx,\phi\upharpoonright_\Gam )$ where, in turn,
$\phi\upharpoonright_\Gam :=\phi\upharpoonright_{{\rm{Supp}}\,(\Gam )}$. We also write 
$ |\rSp\,(\Gam )|$ for the area of $\rSp\,(\Gam )$.

Recall, the quantities $\|\rSp\,(\Gam )\|$ and $w(\Gam )$ are defined in Eqns \eqref{NoTs} 
and \eqref{SWoC}, respectively. 

\bl\label{Lem5.4}%{\bf Lemma 5.4.} 
{\rm{(A Peierls bound via V-cells)}} There exists a constant $\ovp=\ovp(D) > 0$ 
such that for any contour \ $\Gam =({\rm{Supp}}\,(\Gam ),\phi\upharpoonright_{\Gam})$ we have 
\be\label{(6)}
\qquad\qquad\qquad w(\Gam ) = \prod_{\bx \in\phi\upharpoonright_\Gam} u^{-S^{-1}\left( \left|
\cV (\bx,\phi\upharpoonright_\Gam )\right| - S \right)}\leq  u^{-\ovp(D)\|\rSp\,(\Gam )\|}.\ee
\el

\bp The equality in Eqn \eqref{(6)} is simply a re-writing of \eqref{SWoC}. Further, we need to
consider sites $\bx$ where $\left| \cV (\bx,\phi\upharpoonright_\Gam) \right| >S$; otherwise
(i.e., when $\left| \cV (\bx,\phi\upharpoonright_\Gam) \right| =S$) site 
$\bx$ does not contribute into \eqref{(6)}. Observe that
$$\hbox{if }\;\left|\cV (\bx,\phi\upharpoonright_\Gam) \right| - S \ge S\;\hbox{ then }\;
\left| \cV (\bx,\phi\upharpoonright_\Gam) \right| - S \ge\frac{1}{2} \left|\cV (\bx,\phi\upharpoonright_\Gam)\right|\,.$$
On the other hand, by Lemma \ref{Lem5.3},
$$\hbox{if }\;\left |\cV (\bx,\phi\upharpoonright_\Gam)\right| - S < S\;\hbox{ then }\;
\left |\cV (\bx,\phi\upharpoonright_\Gam)\right| - S \ge \odel \ge\frac{\odel}{2S} 
\left|\cV (\bx,\phi\upharpoonright_\Gam)\right|\,. $$
According to the definition of a $\vphi$-correct template, we have an inequality
$$\sum\limits_{\bx \in\phi\upharpoonright_\Gam} 
\left|\cV (\bx,\phi\upharpoonright_\Gam)\right|{\mathbf 1}\Big(
\left|\cV (\bx,\phi\upharpoonright_\Gam)\right| >S \Big) \ge\diy\frac{1}{9 D^2} |\rSp\,(\Gam )|\,.$$ 
Also,
$||\rSp\,(\Gam )|| =\kappa |\rSp\,(\Gam )|$ where $\kappa =1/S^2=4/(3D^4)$. Thus, we can take
\be\label{(6A)}\diy\ovp (D)=\diy \frac{\kappa}{9 D^2} \min \left(\frac{1}{2}, \frac{\odel}{\sqrt 3} 
D^{-2}\right).\ee
\ep

\subsection{Proof of Theorems 1--3 (Classes TA1, TA2, TB), Theorems 7--9 (Classes HA1, HA2, HB), 
Theorem 11
(Classe HC), Theorems 12 (Classes HD1 and HD2) and 13 (Class HE)}\label{SubSec5.3} Owing to Theorems I and II, the 
proof of the listed theorems is reduced to
an explicit description of the PGS-equivalence classes. The structure of these classes
follows from the arithmetic properties of the value $D^2$ to which the conditions
of each theorem explicitly refer.  Once again, in Theorem 13 we restrict ourselves to the
analysis of PGSs only. %\label{Sec5}

\section{\bf Proof of Theorems 4--6 and 10}\label{Sec6}

\subsection{Dominance for the H-C model}\label{SubSec6.1} Our study of dominance follows an approach 
developed in \cite{Sl1}, \cite{Za} and \cite{BS}. In particular, we use an appropriate family of {\it small}
\ contours (see Definition~1 on page 566 in \cite{Za}) and then compare the {\it free energies} 
of the corresponding {\it truncated models} to decide which PGS-equivalence class is dominant.

In the examples of dominance presented in Theorems 4--6 (Class TB), and 10 (Class HB) 
we have that all PGSs are $\alpha$-configurations; this fact generates a number of similarities in the analysis of these
examples. Let us first give a common summary of our construction. 

The `smallest' contour in a PGS is generated by the removal of a single particle. The statistical
weight of such a contour is $u^{-1}$; we can say that it represents a $u^{-1}$-excitation. The 
{\it density} of such $u^{-1}$-contours is the same in each of the PGSs $\vphi\in\sP$. Similarly, 
the removal of two particles at distance $D$ from each other generates a contour of statistical weight $u^{-2}$. 
Again, the density of such $u^{-2}$-contours is the same in every PGS $\vphi\in\sP$.

The next category of a small contour is generated when three particles are removed at the 
vertices of a $D$-triangle $\triangle$, and one particle is inserted at a site inside $\triangle$. Here the new 
occupied site should lie at a distance $\ge D$ from any other sub-lattice site. As before, the corresponding 
contour has statistical weight $u^{-2}$. We can speak of a {\it single insertion} repelling 
3 particles from a PGS. Next, we will have to deal with double, triple, and quadruple admissible insertions 
maintaining the weight $u^{-2}$ for the emerging contour. To stress the latter property, we will often 
speak of $u^{-2}$-insertions.  

Figures 16--22 show the structure of $u^{-2}$-insertions for the cases considered in Theorems 
4--6 and 10. Double admissible insertions occur when 4 particles are removed from
the vertices of a $D$-rhombus formed by two adjacent $D$-triangles $\triangle_1$, 
$\triangle_2$, and 2 particles are inserted inside $\triangle_1\cup\triangle_2$. (For $D^2=49$ on
lattice $\bbA_2$ (Theorem 4), the single and double $u^{-2}$ suffice.) Next, triple
admissible insertions occur when 5 particles are removed from the boundary of a trapeze
formed by three pair-wise adjacent $D$-triangles $\triangle_1$, $\triangle_2$, $\triangle_3$,  
and 3 particles are inserted inside $\triangle_1\cup\triangle_2\cup\triangle_3$. Finally, quadruple 
admissible insertions occur when 6 particles are removed from the boundary of an $2D$-triangle 
by four pair-wise adjacent $D$-triangles $\triangle_1$, $\triangle_2$, 
$\triangle_3$, $\triangle_4$,  
and 4 particles are inserted inside $\triangle_1\cup\triangle_2\cup\triangle_3\cup\triangle_4$.

Any other contour in the truncated model for the considered examples has statistical weight at 
most $u^{-3}$. This statement requires a certain effort to verify (including a substantial computer 
assistance); it is done in Section \ref{Sec7} in the form of Lemmas \ref{Lem7.1}--\ref{Lem7.4}. 
With Lemmas \ref{Lem7.1}--\ref{Lem7.4} at hand, we can use in a standard way the polymer expansions for
the free energies of the truncated models (see, e.g., Sections 1.7, 2.1 in \cite{Za} or Section~3.a 
in \cite{Se}). This allows us to upper-bound the contribution to these free energies from 
contours with weight  $\leq u^{-3}$ by $cu^{-3}$ where $c > 0$ is an absolute constant.
Thus, for $u$ large enough the determination of a dominant class is reduced to the count
of densities of single, double, triple and quadruple insertions.

\subsection{Proof of Theorem 4}\label{SubSec6.2} 
For $D^2=49$ on $\bbA_2$, we have two PGS-equivalence 
classes (inclined and horizontal); they are 
determined by inclined  $D$-sub-lattices containing sites $(3, 5)$ or $(5 ,3)$ and a horizontal one 
containing site $(7, 0)$. We will use pairs $(5, 3)$ and $(7, 0)$ for referring to these sub-lattices and 
their associated PGSs. We want to check that the horizontal $(7, 0)$-PGSs are dominant and the
inclined  $(5,3)$-PGSs are not.
\FigureP16  %\figuree

Any inclined $D$-triangle for $D^2=49$ covers 12 sites where we have a single $u^{-2}$-insertion 
repelling precisely 3 particles at the vertices of this triangle. In Figure 16 (a) 
these sites are marked by orange balls, while the 
repelled sites from a PGS are marked by black balls. The orange balls are placed at sites covered 
by closed concave circular triangles. On the other hand, 
open bi-convex lenses indicate positions where an inserted particle repels 4 black balls at
the vertices of a $D$-rhombus which yields a $u^{-3}$-insertion. 
%For a future use, inserted sites 
%(gray balls) in Figure 16 (a) are divided into 
%4 categories. Category one is marked by a quadruple $(4,27,31)$ and includes three inserted
%sites.

Similarly, any horizontal triangle also covers 12 sites where a single insertion repels 3 particles 
at the vertices of a $D$-triangle. We use the same legend to mark these possibilities 
in Figure 16 (c): closed concave circular triangles cover single $u^{-2}$-insertions, 
open bi-convex lenses indicate positions where an inserted particle repels 4 black balls at
the vertices of a $D$-rhombus, yielding an $u^{-3}$-insertion.   

Thus, both inclined and horizontal PGSs have the same density of single $u^{-2}$-insertions.

The small contour which detects a difference is constructed when 4 particles at the vertices of
a $D$-rhombus are removed and 2 particles inside the rhombus are inserted, maintaining
admissibility. The statistical weight of this contour also equals $u^{-2}$. Figures 16 (b) and 16 (d)
show examples of double $u^{-2}$-insertions marked by red. 
For any inclined $(5,3)$-rhombus %, e.g. with vertices
%$$(0, 0), (5, 3), (-3, 8), {\rm and\,} (2, 11),$$
there are 6 such pairs of sites. For any $(7,0)$-rhombus
%e.g. with vertices
%$$(0, 0), (7, 0), (0, 7), {\rm and\,} (7, 7),$$
there are 7 such pairs.

Any other contour in the truncated model for $D^2=49$ has statistical weight at most $u^{-3}$; this
is proven in Lemma \ref{Lem7.1} in Section \ref{Sec7}. Therefore,
only the horizontal PGS-equivalence class contains the dominant PGSs.  \hfill$\rule{1ex}{1ex}$\par\medskip %\hfill $\footnotesize{\blacksquare}$ 

\subsection{Proof of Theorem 5}\label{SubSec6.3} The argument in the proof of Theorem~5 is similar to that of 
Theorem~4, except for specific numbers of small contours. Here we distinguish between 
inclined (8, 7)- and horizontal (13, 0)-PGSs. The first difference with Theorem~4 is in the
categories of contours having the statistical weight $u^{-2}$. As before, we have single and 
double admissible $u^{-2}$-insertions; see Figure 17. 
\FigureQ17 %\figuref
\FigureR18 %\fgri
In addition, we can place 3 particles in a trapeze and also 4 particles in a $2D$-triangle: see Figure 18. 
The number of single insertions equals 39 per a triangle or 78 per a $D$-rhombus in 
both PGS types. However, in the remaining three categories of $u^{-2}$-insertions, the (8, 7)-PGSs 
dominate distinctively, with
113 vs 78 doubles in a $D$-rhombus, 61 vs 20 triples in a trapeze and 39 vs 3 quadruples in a
$2D$-triangle. The enumeration of cases above can be performed manually but we also present
a Java routine which automates this task. Cf. Program 3 {\tt CountExcitations} in the ancillary file
and Section \ref{Sec9}.

The verification that no other contour with statistical weight smaller than $u^{-3}$ exists
for the horizontal $D$-sub-lattice is done in Lemma \ref{Lem7.2} in Section \ref{Sec7}; 
it requires a more massive enumeration than in Lemma \ref{Lem7.1} and therefore
relies on a computer-assisted argument. \hfill$\rule{1ex}{1ex}$\par\medskip %\hfill $\footnotesize{\blacksquare}$  

\subsection{Proof of Theorem 6}\label{SubSec6.4} The argument for $D^2=147$ on $\bbA_2$ repeats that for $D^2=49,169$ 
and is again based on an exact count of $u^{-2}$-insertions. Here we distinguish 
between two PGS types referred to as vertical and (7, 7)- and inclined (11, 2)-PGSs. 
\FigureS19
Single $u^{-2}$-insertions do not favor any specific PGS: their number equals 34 per triangle or 68 per rhombus
in all PGSs. However, the double, triple and and quadruple $u^{-2}$-insertions favor the (7, 7)-PGSs. 
Viz., the number of double insertions is 86 in a (7, 7)- and 51 in an (11, 2)-rhombus, the number 
of triple insertions is 39 in a (7, 7)- and 1 in an (11, 2)-triangle, and the number of quadruple 
insertions is 1 in a (7, 7)- and 0 in an (11, 2)-triangle. Again, this enumeration can be performed manually
or by executing the Java routine  from Program 3 {\tt CountExcitations}.
Cf. Figures 19 and 20. The proof is completed 
by applying Lemma \ref{Lem7.3} guaranteeing that the $u^{-2}$-insertions are the only ones listed 
above. \hfill$\rule{1ex}{1ex}$\par\medskip %\hfill $\footnotesize{\blacksquare}$ 
%\vskip .5 truecm

%We would like to note that, for a general value of $D$ on $\bbA_2$, it is possible to verify 
%that the double insertions 
%used in Theorems 4--6 are, in a sense, the smallest possible contours which 
%can distinguish between PGSs. It is because the number of single $u^{-2}$-insertions
%removing 3 vertices of a covering triangle is the same for all PGSs occurring for a given $D$.

%For  the case $D^2 =49,169, 147$ it was checked by a direct enumeration in the
%proofs of Theorems 4--6. % \newpage %\vskip .5cm

\FigureT20

%In fact, for an arbitrary attainable $D$ the amount of single $u^{-2}$-insertions is exactly one half of 
%the amount of lattice sites inside the disk of radius $D$, again centered at the origin. This 
%difference is calculated under the following conditions. (i) The disk is open (in $\bbR^2$). (ii) The 
%hexagon is half-open (also in $\bbR^2$), i.e., 3 vertices and 3 open sides of
%the hexagon boundary are included while the other 3 sides and 3 vertices are not. Additionally, (iii)
%the vertices of the hexagon are placed at the sites of the sub-lattice that generates a given class
%of PGSs. It is easy to see that under the conditions (i)--(iii) both quantities participating in the
%difference are independent on a choice of the $D$-sub-lattice.

We do not know if the contours of weight $u^{-2}$ always suffice to determine
dominant PGSs and if a dominant PGS class is always unique. A numerical calculation covering 
$D^2 \leq 100000$ confirms that there is only one $D$-sub-lattice which dominates in the amount
of $u^{-2}$-contours. We conjecture that it is the sub-lattice which has an orange ball 
in a site at the shortest distance from the triangle vertex among all $D$-sub-lattices.

\subsection{Proof of Theorem 10}\label{SubSec6.5} Once more, we follow the established scheme of counting the 
$u^{-2}$-insertions. As in Theorem 6, we distinguish between the inclined (11, 2)- 
and vertical (7, 7)-PGSs types, now on $\bbH_2$. The number of vertical PGSs equals $98$ 
while the number of inclined PGSs is $196$. 

As before, the  $u^{-1}$-contours do not make a distinction. The analysis of dominance focuses
on admissible $u^{-2}$-insertions. 

A vertical (7, 7)-PGS is shown in Figure 21.  As earlier, single 
$u^{-2}$-insertions remove 3 particles at the vertices of a $D$-triangle and add 
one inside the same triangle. They are again marked by orange balls in frame (a). The number 
of such insertions is 21 in triangles
$OAB$, $OCD$  (and also $BFC$ in frame (b)), and  $25$ in triangles $OBC$, $ODE$ (and 
also $CDG$ in frame (b)). In total, we have $46$ single insertions in each of five rhombuses 
$OABC$, $OBFC$, $OBDC$, $OCGD$, $OCDE$ featured  in frame (b).

\FigureU21

Double $u^{-2}$-insertions remove 4 particles at the vertices of a $D$-rhombus and add 2 particles 
inside the same rhombus.  In Figure 21 (b), a double $u^{-2}$-insertion is marked by a red bar. 
The number of admissible double insertions inside every $D$-rhombus equals $108$. 

Triple and quadruple admissible $u^{-2}$-insertions for vertical PGSs are also shown 
in Figure 21 (b). In a triple insertion 5 particles are removed and 3 added (blue balls joined 
by blue bars), whereas in a quadruple insertion 6 particles are removed and 4 added (green 
balls joined by tripods of green bars), following the same geometric pattern as before (a trapeze 
or a $2D$-triangle). %As earlier, a convenient way of counting triple insertions is to refer to the middle triangle in a trapeze. 
%For trapeze $OABCD$, with triangle $OBC$ 
%in the middle, the number  of admissible triple insertions equals 
%$17$. The same number emerges for trapezes $OBFCD$ and $OBCGD$. Thus, 
%the number of triple insertions involving $OBC$ as a middle $D$-triangle 
%is $51$. For trapezes $OBCDE$, $OCGDE$ and $OBCGD$, with triangle $OCD$ 
%in the middle, the number  of admissible triple insertions equals $4$. Thus, the number  of triple 
%insertions involving $OCD$ as a middle $D$-triangle is $12$. 
In total, we have $63$ triple
$u^{-2}$-insertions per a $D$-rhombus.

Quadruple admissible $u^{-2}$-insertions cannot occur inside the $2D$-triangle $EBG$. 
However, for triangle $AFD$ they can occur, and their number equals $9$. Hence, the number of
quadruple insertions with the middle point of a tripod inside triangle $OBC$ equals $9$.
Thus, the total number of admissible quadruple $u^{-2}$-insertions in a $D$-rhombus is $9$.
According to Lemma \ref{Lem7.4}, the list of all admissible $u^{-2}$-insertions is exhausted 
by the aforementioned possibilities.  All-in-all, the above 
count yields $226$ admissible $u^{-2}$-insertions per a $D$-rhombus in a 
vertical PGS.
\FigureV22

The situation with an inclined (11, 2)-PGSs for $D^2=147$ on $\bbH_2$ is shown in 
Figure 22. In frame (a)
we again put orange balls in positions where a single insertion repels three occupied sites
at the vertices of the covering $D$-triangle. The number of such insertions is $24$ in 
triangles $OAB$ and $OCD$ and $22$ in triangles $OBC$ and $ODE$, with $46$ insertion per a $D$-rhombus. 
Next, in frame (b) we mark by a red bar an admissible double insertion removing 4 vertices of the covering 
$D$-rhombus. The number of double $u^{-2}$-insertions is $23$ in all rhombuses in the inclined PGS.
Thus, the total amount of double insertions equals $69$ per a $D$-rhombus.

As before, triple insertions repelling 5 vertices in an inclined PGS occur when 3 particles 
are put in a trapeze, one insertion for each involved triangle. The total number of 
admissible triple $u^{-2}$-insertions per a $D$-rhombus is $3$.

%Again, 
%a convenient way of counting triple insertions is to refer to the middle triangle in 
%a trapeze. For each of trapezes $OABCD$, $OEFAB$ and $OCDEF$, with the middle 
%triangles $OBC$, $OFA$ and $ODE$, respectively, the number of triple insertions 
%equals 1 whereas for $OFABC$, $OBCDE$ and $ODEFA$ there are no such insertions.

Lastly, quadruple insertions repelling 6 vertices could have occurred when 4 particles are 
put in a $2D$-triangle. However, in an inclined PGS such insertions do not exist. 
According to Lemma \ref{Lem7.4}, the list of admissible $u^{-2}$-insertions in an inclined PGS is 
exhausted by the above types. All-in-all, the number of admissible $u^{-2}$-insertions 
per a $D$-rhombus in an inclined PGS equals $118$. Hence, for $D^2=147$ on $\bbH_2$, the 
vertical PGS class is dominant, and for $u$ large enough we have $98$ {\rm{EGM}}s generated 
by the (7, 7)-PGSs. \hfill$\rule{1ex}{1ex}$\par\medskip   %\hfill $\tiny{\blacksquare}$ 

\section{\bf Proof of technical assertions from Theorems 4--6, 10}\label{Sec7}

In this section we verify that the small contours with a weight $\geq u^{-2}$ which were 
determined in Section \ref{Sec6} are the only ones possible for the selected values of $D^2$, and 
any other contour has the statistical weight $\le u^{-3}$. The corresponding statements are 
Lemmas \ref{Lem7.1}--\ref{Lem7.3} on $\bbA_2$ and Lemma \ref{Lem7.4} on $\bbH_2$; the latter, essentially, implied
by Lemma \ref{Lem7.3}.

The argument is based on the following construction. Without loss of generality 
we can assume that the underlying PGS $\vphi\in\sP (D)$
has $\vphi ({\mathbf 0})=1$ (i.e., $\vphi $ is a $D$-sub-lattice in $\bbA_2$). 
Recall, a $\vphi$-contour $\Gam$ can 
be obtained by adding finitely many particles at some inserted sites, and then removing the 
particles from $\vphi$ which are repelled by the inserted ones (removed sites/particles). The 
resulting admissible configuration is denoted by $\phi$. One can also remove from $\phi$ any 
additional particles but such an unforced removal can only decrease the weight $w(\Gam )$ and 
therefore will be disregarded. %That is, we 

%If $\sharp =\sharp (\phi )$ is the difference between the numbers of removed and inserted sites then
%$w(\Gam )=u^{-\sharp}$. 
As we saw earlier, every inserted site repels
from $\vphi $ either 3 or 4 removed sites. The inserted sites which repel 3 removed sites
are located inside closed concave circular triangles identified in the proof of Theorems 4--6 (orange balls in Figures 16, 17, 19). 
%(the light-gray areas in Figure 16--20). 
The complement (in $\mathbb{R}^2$) to these triangles consists of
mutually disjoint open circular bi-convex lenses (gray areas in Figures 16, 17, 19). A particle
inserted in a lattice site belonging to a lens repels 4 removed sites.

Consider a $D$-connected component $\Delta$ of the set of removed sites (together with the 
corresponding inserted sites). Let $\sharp(\Delta)$ denote the difference between the
numbers of removed and inserted sites in $\Delta$. Our goal is to
verify that the weight $u^{-\sharp(\Delta)}$ of any such component is at most $u^{-3}$,
i.e., $\sharp(\Delta) \ge 3$. Note that a contour support has been defined in Section \ref{SubSec3.1} by 
using the notion of a template; hence it can include more than one $\Delta$. In that case 
the statistical weight of the contour is the product of the statistical weights of constituting $D$-connected
components, and for our purposes it is enough to estimate the weight of a single
component $\Delta$.

To evaluate $\sharp (\Delta )$, it is convenient to introduce a {\it total repelling force} $\rF (\bx)=\rF (\bx ,\phi )$ 
acting (in the resulting 
AC $\phi\in\cA (D)$) upon a removed site $\bx\in\vphi$. Such a force is accumulated
from all inserted sites $\by_i\in\phi$ that repel site $\bx$: $\rF (\bx) = \sum\limits_i \rF (\bx, \by_i, \phi)$.
We require that every summand $\rF (\bx, \by_i, \phi)$ is non-negative and depends only
on the Euclidean distance $\rho(\bx, \by_i)$ between $\bx$ and $\by_i$. The square of this
distance $\rho(\bx, \by_i)^2$ is always a positive integer, and we use a shorthand notation
$f_r$ for $\rF (\bx, \by, \phi)$, with $r = \rho(\bx, \by)^2\in\bbN$, $\bx\in\vphi$, $\by\in\phi$. 
With this notation at hand, $\forall$ \ $\bx\in\vphi$,
$$\rF (\bx):=\sum\limits_{r<D^2,\;\by\in\phi}f_r{\mathbf 1}\Big(\hbox{$\by$ removes $\bx$,
and $r=\rho (\by ,\bx)^2$}\Big),\;
\hbox{ where }\;f_r\geq 0. \eqno (7.1.\rA)$$
The coefficient $f_r$ is referred to as a {\it local repelling force} at distance ${\sqrt r}$. A dual quantity 
$\rG (\by) =\rG (\by,\phi )$ represents the{\it  total repelling force} generated by an inserted site $\by$:
$$\rG (\by):=\sum\limits_{r<D^2,\;\bx\in\vphi}f_r{\mathbf 1}\Big(\rho (\bx,\by)^2=r,\;\bx\;
\hbox{is removed by}\;\by\Big),\;\;\;\by\in\phi .\eqno (7.1.\rB)$$

Our aim is to find $f_r$ such that, for any site $\bx\in\vphi$ and any site $\by\in\phi$ removing 
3 or 4 sites from $\vphi$, 
$${\rm{(a)}}\quad\rF (\bx ,\phi )\leq 1,\qquad{\rm{(b)}}\quad\rG (\by,\phi )=1. \eqno (7.2)$$

Owing to (7.2), if the
{\it deficit} $\delta (\bx )=\delta (\bx ,\phi )$ of the removed site $\bx$ is calculated as $1-\rF (\bx)$
then $\delta (\bx )\geq 0$, and 
$$\sharp(\Delta )=\sum_{\bx\in\vphi}\delta (\bx, \phi )
{\mathbf 1}\Big(\hbox{site $\bx$ is removed when passing from $\vphi$ to 
$\phi$}\Big).\eqno (7.3)$$

From now on we assume that the configuration $\phi$ has a single $D$-connected component 
$\Delta$, and the rest of the argument deals with this $\Delta$. Figure 23 shows a fragment 
of a set $\Delta$ with a collection of inserted and removed sites. 
\WFigure25
%Accordingly, we write that $\sharp (\phi )=\sharp (\Delta )$.

%The family $\{ f_r \}$ of local repelling forces is called {\it proper} if for any $\Delta$ and each
%removed site $\bx \in \Delta$ the deficit $\delta (\bx )\ge 0$. A proper family of repelling forces
%s a convenient tool for estimating the value of $\sharp(\Delta)$ as the sum of deficits $\delta (\bx )$
%collected over some removed sites in $\Delta$ can only increase when the remaining removed
%sites in $\Delta$ are also taken into account.

The next observation is that set $\Delta$ consists of {\it internal} sites for which all 6 sub-lattice
neighbors also belong to $\Delta$ and {\it boundary} sites which have at least one occupied
$D$-sub-lattice neighbor (obviously, not belonging to $\Delta$). Each $D$-connected component
of the boundary sites in $\Delta$ defines a closed broken line in $\mathbb{R}^2$, and the set
$\Delta$ can be understood as $\mathbb{R}^2$-polygon with the boundary $\partial\Delta$
formed by these broken lines.

In general, the boundary $\partial\Delta$ can have
several connected components: one external and zero or more internal ones.
An ambiguous situation arises when 4 $D$-segments from $\partial\Delta$ meet at the same
boundary site (i.e., this site has 2 opposite $D$-neighbors that are occupied). In that
case we fictitiously cut this site along the short line segment (of length less than 1) which passes
trough this site and has both ends inside $\Delta$ (viewed as an open polygon in $\mathbb{R}^2$).
This removes the ambiguity, and the exterior and the interior of $\Delta$ become uniquely defined.

It is clear that, as an $\mathbb{R}^2$-polygon, $\Delta$ can only have vertices with angles $\pi/3$,
$2\pi/3$ and $4\pi/3$. We say that the corresponding removed sites from $\partial\Delta$ are of
type $\pi/3$, $2\pi/3$ and $4\pi/3$ respectively. The remaining sub-lattice sites from
$\partial\Delta$ correspond to the angle $\pi$; we say that such a site has type $\pi$.

If a vertex $\bx\in\partial\Delta$ is repelled only by a single inserted site $\by$ then imagine the
particle at $\by$ being deleted. Then vertex $\bx$ also disappears from $\Delta$ (as nothing
repels it anymore), and the value $\sharp (\Delta )$ does not increase. (Actually, $\sharp (\Delta)$ 
remains intact if $\by$ repels a single vertex $\bx$ in $\Delta$.)  In Figures 24, 25 we refer to 
such a site $\bx$ as {\it deletable}. 
\FigureY26

In view of the above definition, every polygon $\Delta$ that can be reduced, by the process of deletion, 
to an irreducible polygon $\Delta^0$, for which $\sharp (\Delta^0)\leq \sharp (\Delta)$. By definition, a
polygon $\Delta$ with a single inserted site is irreducible. 
The simplest form of $\Delta^0$ is a $D$-triangle with a single inserted site, where
$\sharp (\Delta^0)=2$. We would like to: (i) list all $\Delta$s that are reduced to a $D$-triangle
(possibly, with the help of a computer), and (ii) demonstrate that for all other irreducible polygons 
$\Delta^0$, we have $\sharp (\Delta^0)\geq 3$.

In fact, the next irreducible case is where $\Delta^0$ is a $D$-rhombus with a single 
inserted site: it has $\sharp (\Delta^0) \ge 3$, in agreement with property (ii).
%(Since an irreducible $D$-rhombus has one inserted site removing all 4 vertices of the rhombus.)
%verifi (ii) in this case.
%%Suppose that for a given $D^2$ we are able to verify that adding a repelling
%site to a $D$-rhombus always produces a polygon with $\sharp \ge 3$. Then, 
 %\fgraa
\FigureZ27
For any other (larger) irreducible polygon $\Delta^0$, the boundary $\partial\Delta^0$ must have 
(i) no vertex 
of type $\pi/3$ and (ii) at least 6 vertices of type $2\pi/3$. Each of the latter 6 vertices is repelled 
by exactly
two inserted sites. Our goal in the lemmas below is to find, for the corresponding value of $D^2$, 
a collection of 
repelling forces $\{f_r\}$ such that 
$$\delta (\bx) > 1/3\;\hbox{ for any $\bx \in \partial\Delta^0$
of type $2\pi/3$.}\eqno (7.4)$$ 
This would imply the desired assertions, as $6 \delta (\bx) > 2$. 

Let us now pass to specific cases. The proofs of Lemmas \ref{Lem7.1}--\ref{Lem7.4} require a finite 
enumeration which was done by computer. Cf. Programs 4 {\tt VerifyRepellingForces} and 
5 {\tt  CountMinDelta} in Section \ref{Sec9} and in the ancillary file.

The first case is $D^2=49$ on $\bbA_2$. Define the following family $\{f_r\}$:
$$\bealllll
 f_{1}  = 44/56,&  f_{3} = 40/56,&  f_{4} = 40/56,&  f_{7} = 31/56,&  f_{9} = 31/56,\\
 f_{12} = 22/56,&   f_{13} = 22/56,& f_{16} = 17/56,&     f_{19} = 17/56,& f_{21} = 17/56,\\
 f_{25} = 8/56,& f_{27} = 8/56,&   f_{28} = 8/56,& f_{31} = 8/56,& f_{36} = 4/56,\\
 f_{37} = 4/56,&   f_{39} = 4/56,& f_{43} = 4/56,&   f_{48} = 4/56.
\ena\eqno (7.5)$$
The values $r=$ 1, 3, 4, 7,  9, 12, 13, 16, 19, 21, 25, 27, 28, 31, 36, 37, 39, 43, 48 in (7.5)
represent all squared Euclidean distances from ${\mathbf 0}$ to the $\bbA_2$-sites within 
an open $\mathbb{R}^2$-disk of radius~$7$.

\bl\label{Lem7.1}%{\bf Lemma 7.1.} {\sl 
The family $(7.5)$ gives a collection of local repelling forces  for 
$D^2=49$ on $\bbA_2$ satisfying $(7.2)$ and $(7.4)$, for both 
horizontal and inclined {\rm{PGS}}s. 
More precisely, for this collection, $\forall$ irreducible polygon $\Delta^0$ and vertex 
$\bx \in \partial\Delta^0$ of type $2\pi/3$, 
$$\delta (\bx )\geq 1 - f_{19}- f_{31} = 1 - f_{21}- f_{27} = 31 / 56. \eqno (7.6)$$
\el

\brff {\rm For the proof of Theorem 4, it suffices to find a collection $\{f_r\}$ only for 
inclined PGSs. Cf. Lemma \ref{Lem7.2}. However, it turns out that the family $\{f_r\}$ in (7.5) serves
both types of PGSs.} \hfill $\blacktriangle$
\erff

%{\bf Proof of Lemma 7.1.} The shortest proof is given with the help of a computer.
%See ... ... ... .. . \qquad $\blacksquare$

%\vskip 1truecm

Next, we deal with $D^2=$169. Here we consider the values $r$ 
representing the squared Euclidean distances from ${\mathbf 0}$ to all $\bbA_2$-sites within 
an open $\mathbb{R}^2$-disk of radius~$13$. For such values $r$ we set:%\newpage

$$\bealllll
 f_{1}  = 131 / 135,&  f_{3} = 127 / 135,&  f_{4} = 251 / 270,&  f_{7} = 241 / 270,&  f_{9} = 116 / 135,\\
 f_{12} = 37 / 45,&   f_{13} = 221 / 270,& f_{16} = 7 / 9,&     f_{19} = 133 / 180,& f_{21} = 383 / 540,\\
 f_{25} = 179 / 270,& f_{27} = 19 / 30,&   f_{28} = 169 / 270,& f_{31} = 317 / 540,& f_{36} = 281 / 540,\\
 f_{37} = 14 / 27,&   f_{39} = 131 / 270,& f_{43} = 41 / 90,&   f_{48} = 37 / 90,&   f_{49} = 43 / 108,\\
 f_{52} = 103 / 270,& f_{57} = 35 / 108,&  f_{61} = 53 / 180,&  f_{63} = 151 / 540,& f_{64} = 5 / 18,\\
 f_{67} = 7 / 27,&    f_{73} = 119 / 540,& f_{75} = 11 / 54,&   f_{76} = 109 / 540,& f_{79} = 11 / 60,\\
 f_{81} = 22 / 135,&  f_{84} = 83 / 540,&  f_{91} = 2 / 15,&    f_{93} = 31 / 270,&  f_{97} = 29 / 270,\\
f_{100} = 1 / 9,&    f_{103} = 4 / 45,&   f_{108} = 2 / 27,&   f_{109} = 2 / 27,&   f_{111} = 7 / 108,\\
f_{112} = 1 / 15,&   f_{117} = 2 / 45,&   f_{121} = 11 / 270,& f_{124} = 23 / 540,& f_{127} = 11 / 270,\\
f_{129} = 4 / 135,&  f_{133} = 4 / 135,&  f_{139} = 2 / 135,&  f_{144} = 1 / 54,&   f_{147} = 2 / 135,\\
f_{148} = 2 / 135,&  f_{151} = 1 / 270,&  f_{156} = 1 / 270,&  f_{157} = 1 / 180&   f_{163} = 0.
\ena\eqno (7.7)$$

\bl\label{Lem7.2}%{\bf Lemma 7.2.} {\sl 
The family $(7.7)$ gives a collection of local repelling forces  for 
$D^2=169$ on $\bbA_2$ satisfying  $(7.2)$ and $(7.4)$, for horizontal  $(13,0)$-{\rm{PGS}}s. 
More precisely, for  this collection, $\forall$ irreducible polygon $\Delta^0$ in a horizontal 
$(13,0)$-{\rm{PGS}} and vertex $\bx \in \partial\Delta^0$ of type $2\pi/3$, 
$$\delta (\bx )\geq 1 - f_{28} - f_{133} = 93 / 270. \eqno (7.8)$$
\el

%{\bf Proof of Lemma 7.2.} 

Finally, we consider the example of $D^2=$147 on $\bbA_2$. Set:
$$\bealllll
f_{1}  = 24 / 24,& f_{3}  = 24 / 24,& f_{4}  = 24 / 24,& f_{7}  = 23 / 24,& f_{9}  = 22 / 24,\\
f_{12} = 21 / 24,& f_{13} = 21 / 24,& f_{16} = 20 / 24,& f_{19} = 19 / 24,& f_{21} = 18 / 24,\\
f_{25} = 16 / 24,& f_{27} = 15 / 24,& f_{28} = 15 / 24,& f_{31} = 14 / 24,& f_{36} = 12 / 24,\\
f_{37} = 12 / 24,& f_{39} = 11 / 24,& f_{43} = 10 / 24,& f_{48} =  9 / 24,& f_{49} =  8 / 24,\\
f_{52} =  7 / 24,& f_{57} =  6 / 24,& f_{61} =  5 / 24,& f_{63} =  4 / 24,& f_{64} =  4 / 24,\\
f_{67} =  4 / 24,& f_{73} =  3 / 24,& f_{75} =  3 / 24,& f_{76} =  2 / 24,& f_{79} =  2 / 24,\\
f_{81} =  2 / 24,& f_{84} =  2 / 24,& f_{91} =  1 / 24,& f_{93} =  1 / 24,& f_{97} =  1 / 24,
\ena\eqno (7.9\rA)$$
with 
$$f_r = 0\;\hbox{ for }\;r > 97.\eqno (7.9\rB)$$
This yields a family of values $f_r$ where $r$ represents the squared Euclidean distance 
from ${\mathbf 0}$ to all $\bbA_2$-sites within an open $\mathbb{R}^2$-disk of radius~$\sqrt{147}$.

\bl\label{Lem7.3}%{\bf Lemma 7.3.} {\sl 
The family $(7.9\rA,\rB)$ gives a collection of local repelling forces $\{f_r\}$ for 
$D^2=147$ on $\bbA_2$ satisfying $(7.2)$ and $(7.4)$, for inclined $(11,2)$-{\rm{PGS}}s. 
More precisely, for this collection, $\forall$ irreducible polygon $\Delta^0$ in an inclined $(11,2)$-{\rm{PGS}} 
and vertex $\bx \in \partial\Delta^0$ of type $2\pi/3$, 
$$\delta (\bx ) \geq 1 - f_{37} - f_{100} = 1/2. \eqno (7.10)$$
\el

Next, we extend our analysis to $\bbH_2$. Consider the values $f_r$ given by Eqns  (7.9A,B)
for $r$ representing the squared Euclidean distance $\leq 147$ between two 
$\bbH_2$-sites. We call it the $\bbH_2$-projected family (7.9A,B).

\bl\label{Lem7.4}%{\bf Lemma 7.4.}  {\sl 
The $\bbH_2$-projected family $(7.9\rA,\rB)$ gives a collection of local 
repelling forces $\{f_r\}$ for 
$D^2=147$ on $\bbH_2$ satisfying $(7.2)$ and $(7.4)$, for inclined $(11,2)$-{\rm{PGS}}s. 
More precisely, for this collection, $\forall$ irreducible polygon $\Delta^0$ in an inclined 
$(11,2)$-{\rm{PGS}} and vertex $\bx \in \partial\Delta^0$ of type $2\pi/3$, the bound
$(7.10)$ holds true.
\el

%{\bf Lemma 7.5.}  {\sl The $\bbH_2$-projected family $(7.9\rA,\rB)$ gives a collection of local 
%repelling forces $\{f_r\}$ for 
%$D^2=139$ on $\bbH_2$ and $(\alpha ,D^*)$-PGSs $\vphi\in\sP (D,\bbH_2)$ where $(D^*)^2=147$. 
%This family satisfies $(7.4)$, for vertical $(7,7)$-{\rm{PGS}}s. 
%More precisely, for  this collection, $\forall$ irreducible polygon $\Delta^0$ in a vertical {\rm{PGS}} 
%$\vphi\in\sP (7,7)$ and vertex $\bx \in \partial\Delta^0$ of type $2\pi/3$, the bound
%$(7.10)$ holds true.} \vskip .5 truecm

\brgg {\rm In essence, the local repelling forces $f_r$ are related to an attempt to improve a Peierls 
constant for the listed values $D^2=$ 49, 147, 169.  In our opinion, this method in its present form 
can work only for moderate values of $D^2$.} \hfill $\blacktriangle$
\ergg
 
\section{A brief note on sliding on $\bbH_2$}\label{Sec8}

As was noted, the values $D^2=$ 4, 7, 31, 133 from Class HS exhibit sliding on $\bbH_2$.
This is characterized by a cost-free passage from one type of PGSs to another.
\Figurea1
For  $D^2=$ 4 we have two types of PGSs: (a) one formed by hexagons with side-length 2,
and (b) the other formed by $\beta$-configurations. These patterns can be intermittent in a 
stripe-like fashion, which generates countably many PGSs with no loss in the 
weight in the course of transition. See Figure 26 (a).

Similarly, for  $D^2=$ 7 we have the following types of PGSs: (a) an $(\alpha ,D)$-configuration for
$D^2=9$, (b) an assortment of $\beta$-configurations of various shapes and orientations.
Again, it is possible to combine these patterns and generate countably many
PGSs with no loss in the weight in the course of transition. See Figure 26 (b).

For $D^2=$ 31, we have a competition between strips formed by $\wt D$-triangles with 
${\wt D}^2=36$
and triangles with squared side-lengths 31, 36, 43. Such triangles share common sides and 
have equal areas. Again, the separation line does not incur a loss in weight.

Similarly, for $D^2=$ 133, we have a competition between strips formed by $\wt D$-triangles with 
${\wt D}^2=144$
and triangles with squared side-lengths 133, 144, 157. Such triangles again share a common side and 
have equal areas. As earlier, the separation line does not incur any loss in weight. 

A pattern resembling that for $D^2=$ 31, 133
is typical for sliding on $\bbZ^2$; cf. \cite{MSS1}.
\Figureb2

\section{Comments on computer programs used in the proofs}\label{Sec9}

The ancillary file to this paper contains five Java programs and their outputs.
These programs are used to (a) assist the proof of Lemmas  
4.9, 4.10, (b) identify distinguishing small contours (admissible $u^{-2}$-insertions) for the values
$D^2=49$, $D=169$ and $D^2=147$ and (c) assist the proof of Lemmas \ref{Lem7.1}--\ref{Lem7.4}. Lemmas 4.9 and 4.10 are parts of the proof of Theorem I(ii) (in the part
concerning Classes HD and HE) and Theorems 12 and
13. Lemmas \ref{Lem7.1}--\ref{Lem7.4} are parts of the proof of Theorems 4--6 and 10,
for $D^2=49$, $D=169$ and $D^2=147$. The routines can be executed on any computer hosting Java 
Development Kit version 1.4 and later. The routine from Program 5 requires 3GB of RAM 
while other routines are not resource hungry. Executions of all routines take from few 
seconds to 60 minutes. 

Program 1 {\tt NearestLoschianNumber} is a routine working on $\bbH_2$.
It checks that, apart from 184 values, for each $D^2< (54)^4$ such that $D^2$ is not divisible 
by 3, (i) the inequality involving the RHS of Eqn \eqref{eq:L4.8_1}
holds true: $\diy\frac{\sqrt 3}{2}D^2 + \frac{D}{2{\sqrt 3}}>\frac{{\sqrt 3}(D^*)^2}{2}=2s (\triangle^* )$, 
(ii) inequalities  \eqref{(4.15)}, \eqref{(4.16)} are satisfied. It leads to the conclusion that for each 
$D^2< (54)^4$, apart from the above 184 values, if $D^2$ is not divisible by 3 then this value $D^2$ belongs 
to Class HC: the PGSs are 
$(\alpha ,D^*)$-configurations, where $D^*>D$ is the nearest L\"oschian number such that 
$3|(D^*)^2$. This assists the proof of Lemmas 4.9 and 4.10.
  
 % the PGSs are triangular 
 % $(\alpha ,D^*)$-configurations where $D^*>D$ is the attainable minimal value such that
 % $3|\,(D^*)^2$. 

%(identifies the pairs $D$ and $D^*$ with $D^2\leq (54)^2$ which violates inequalities $(D^2+D)/3\geq (D^*)^2$ and (4.5), (4.6) with replacing $D$ by $D^*$ in the right most side. )

Program 2 {\tt SpecialD} is a routine analyzing the 184 values $D^2$ detected by Program 1.
It specifies the values forming Classes HD, HE and HS on $\bbH_2$. This routine (i) extracts 
the exceptional values $D^2=$ 4, 7, 13, 16, 28, 31, 49, 64, 67, 97, 133, 157, 256 (Classes HD, HE, HS) 
and identifies the PGSs for the exceptional non-sliding $D^2\neq$ 4, 7, 31, 133,
(ii) checks that each $D^2$ among the above 184 values which is not from Classes HD, HE, HS
belongs to Class HC. Again, this assists the proof of Lemmas 4.9 and 4.10.

Program 3 {\tt CountExcitations} is a routine calculating, for a given sub-lattice in 
$\bbA_2$ or $\bbH_2$, the 
amount of distinguishing small contours (admissible $u^{-2}$-insertions) of the types 
used in the proofs of Theorems 4--6, 10. The execution results are presented only for
$D^2=49$, $D=169$ and $D^2=147$ and the sub-lattices used in the proof of Lemmas \ref{Lem7.1}--\ref{Lem7.4}.

Program 4 {\tt VerifyRepellingForces} is a routine that verifies, for $D^2=$ 49, 169, 147, and a
family of local repelling forces 
$\{f_r\}$, if Eqn (7.2) is satisfied. 
%As above, the execution results are presented for $D=7$, $D=13$  and $\sqrt{147}$ only.

Program 5 {\tt CountMinDelta} is a routine that checks inequality (7.4), for a given sub-lattice and a 
family of local repelling forces, for $D^2=49$, $D=169$ and $D^2=147$. 

Programs 4, 5 assist the proof of Lemmas \ref{Lem7.1}--\ref{Lem7.4}.

\vskip 0.5cm
\noindent
{\bf Acknowledgement} IS and YS thank the Math Department, Penn State University, for hospitality
and support. YS thanks St John's College, Cambridge, for long-term support.

%\WFigure25 \XFigure24A \YFigure25 \ZFigure26


\begin{thebibliography}{10}

\bibitem{BKZZ}
Barbier, J., Krzakala, F., Zdeborova, L., Zhang, P.
The hard-core model on random graphs revisited. {\it J.  Phys.}: Conf. Ser. {\bf 473},  012021
(2013).

\bibitem{Ba} Baxter, R. Hard hexagons: exact solutions. {\it J. Phys.} A: Math. Gen., {\bf 13}
(1980), L61-L70.

\bibitem{BS}
Bricmont, J., Slawny, J. Phase transitions in systems with a finite number of dominant ground states.
{\it J. Stat. Physics}, {\bf 54}, (1989), 89-161.

\bibitem{ChW} Chang, H.-C., Wang, L.-C. A simple proof of Thue's theorem on circle packing.
arXiv:1407.2199, 2010.

\bibitem{CKPUZ} Charbonneau, P., Kurchan, J., Parisi, G., Urbani, P.,  Zamponi, F. Glass and
jamming transitions: from exact results to finite-dimensional descriptions. {\it Annual Rev. Cond.
Mat. Phys.}, {\bf 8} (2016), 1-26.

%\bibitem{CoK}
%Cohn, H.,  Kumar, A.D. Universally optimal distributions of points on spheres.
%{\it J. Amer. Math. Soc.}, {\bf 20} (2007), 99-148.

%\bibitem{CKRV1}
%Cohn, H.,  Kumar, A.D., Miller, S.D., Radchenko, D., Viazovska, M. The sphere packing problem
%in dimension 24. {\it Annals of Math.}, {\bf 185} (2017), 1017-1033.

%\bibitem{CKRV2}
%Cohn, H.,  Kumar, A.D., Miller, S.D., Radchenko, D., Viazovska, M.
%Universal optimality of the E8 and Leech lattices and interpolation formulas.
%arXiv:1902.05438.

\bibitem{CD} Connelly, R., W. Dickinson. Periodic planar disc packings. {\it Phil. Trans. Roy. Soc.,}
{\bf A372}, 20120039. http://dx.doi.org/10.1098/rsta.2012.0039

\bibitem{CS} Conway, J., Sloane, N. {\it Sphere packings, lattices and groups}, 3rd Ed (With
 contributions by E. Bannai, R.E. Borcherds, J. Leech, S.P. Norton, A.M. Odlyzko, R.A. Parker,
 L. Queen and B.B. Venkov.) New York: Springer, 1999.

%\bibitem{DLRS} De Loera, J.A., Rambau, J., Santos, F. {\it Triangulations. Structures for
%algorithms and applications.} Heidelberg: Springer, 2010.

\bibitem{Dob}
Dobrushin R.L. The problem of uniqueness of a Gibbsian random field and
the problem of phase transitions. {\it Funct. Anal. Appl.} {\bf 4}:4 (1968), 302--312.

\bibitem{DoS} Dobrushin R, Shlosman S., The problem of translation invariance of Gibbs states at low
temperatures. {\it Mathematical Physics Reviews} {\bf 5}, 53-195. Soviet Sci. Rev. Sect. C Math. Phys.
Rev., 5, Harwood Academic Publ., Chur, 1985.

\bibitem{F}
Fejes T\'{o}th, L. Some packing and covering theorems. {\it Acta Sci. Math.}, {\bf 12A}
(1950), 62-67.

\bibitem{Ge}
Georgii, H.O. {\it Gibbs measures and phase transitions.} De Gryuter, 2011.

%\bibitem{Ha}
%Hales, T.C. A proof of the Kepler conjecture. {\it Annals of Math.}, {\bf 162} (2005), 1065-1185.

\bibitem{HeP}
Heilmann, O.J., Praestgaard, E. Phase transition of hard hexagons on a triangular lattice.
{\it J. Stat. Physics}, {\bf 9} (1973), 1-22.

\bibitem{HS} Holsztynski, W., Slawny, J. Peierls condition and number of ground states.
{\it Comm. Math Phys.} {\bf 61} (1978), 177-190.

\bibitem{Hs}
Hsiang, W-Y. Simple proof of a theorem of Thue on the maximal density of circle
packings in $E^2$, {\it L'Einseignement Math\'ematique}, {\bf 38} (1992), 125-131.

\bibitem{JaL}
Jauslin, I., Lebowitz, J.L. High-fugacity expansion, Lee-Yang zeros and
order-disorder transitions in hard-core lattice systems. {\it Comm. in Math. Phys.},
{\bf 364}:2 (2018), 655-682.

\bibitem{KMRTSZ}  Krzakała, F., Montanari, F., Ricci-Tersenghi, F., Semerjian, G., Zdeborova, L.
Gibbs states and the set of solutions of random constraint satisfaction problems. {\it Proc. Nat.
Acad. Sci. USA}, {\bf 104}:25 (2007), 10318-10323.

%\bibitem{MSS1} Mazel, A., Stuhl, I., Suhov, Y. Hard-core configurations on a triangular lattice 
%and Eisenstein primes. arXiv:1803.04041v1, 2018. 

\bibitem{MSS1} Mazel, A., Stuhl, I., Suhov, Y.  High-density hard-core model on $\bbZ^2$
and norm equations in ring $\bbZ [{\sqrt[6] -1}]$. arXiv:1909.11648v2, 2019.

\bibitem{M}
Misaghian, M. Factor rings and their decompositions in the Eisenstein integers ring
$Z [\omega ]$. {\it Armenian J. of Mathematics,} {\bf 5}:1 (2013), 58--68.

\bibitem{N} Nair, U.P.  Elementary results on the binary quadratic form $a^2 + ab + b^2$.
arXiv: math.0408107v1, 2004

\bibitem{PaZ}
Parisi, G., Zamponi, F. Mean-field theory of hard-sphere glasses and jamming.
{\it Rev. Mod. Phys.}, {\bf 82} (2010), 789--845.

\bibitem{PeS}
Peled, R., Samotij, W. Odd cutsets and the hard-core model on $\bbZ^2$.
{\it  Annales de l'IHP}, Probabilités et Statistiques, {\bf 50} (2014), 975-998.

\bibitem{PiS}
Pirogov, S.A., Sinai, Ya.G. Phase diagrams of classical lattice systems. {\it Teor. Mat. Fiz.}
{\bf 25} (1975), 1185-1192; {\bf 26} (1976), 61-76.

\bibitem{Se} Seiler, E. {\it Gauge theories as a problem of Constructive Quantum field theory
and Statistical mechanics.} Lecture Notes in Physics, {\bf 159}. Berlin: Springer, 1982.

\bibitem{Si} Sinai, Ya. G. {\it Theory of phase transitions: rigorous results.}
Oxford {\it et al.}: Pergamon Press, 1982.

\bibitem{Sl1} Slawny, J. Low-temperature expansion for lattice systems with many ground states. 
{\it J. of Stat. Physics}, {\bf 20} (1979), 711-717.

%\bibitem{V}
%Vyazovskaya, M. The sphere packing problem in dimension 8. {\it Annals of Math.}, {\bf 185}
%(2017), 991-1015.

\bibitem{Za}
Zahradnik, M. An alternate version of Pirogov-Sinai theory. {\it Comm. Math Phys.} {\bf 93}
(1984), 559-581.

\end{thebibliography}
\end{document}